\documentclass{amsart}
\usepackage{amsmath}
\usepackage{todonotes} 
  
\input xy
\xyoption{pdf} 
\xyoption{all}
\CompileMatrices

\usepackage{xcolor,sseq,amssymb, caption}

\usepackage{wasysym}
\usepackage{array}
\usepackage{pdflscape}

\usepackage{pdfpages}    
\usepackage{longtable}  

\usepackage{tikz}
\newcommand{\ssdiff}[1]{\ssarrow{-1}{#1}\ssmove{1}{-#1}}
\newcommand{\ssreddiff}[1]{\ssarrow[color=red]{-1}{#1}\ssmove{1}{-#1}}

\newcommand{\ssblueline}[1]{\ssline[color=blue]{0}{#1}}
\newcommand{\ssgreenline}[1]{\ssline[color=green]{0}{#1}}

\newcommand{\sscurvedblueline}[1]{\ssline[color=blue, curve=-.1]{0}{#1}}

\def\circbox{
     \begin{tikzpicture}[x=1.2ex,y=1.2ex]
        \draw (0,0) rectangle +(1.,1);
        \draw (.5,.5)  circle (.4mm); 
     \end{tikzpicture}
            }
\def\JJbox{
     \begin{tikzpicture}[x=1.2ex,y=1.2ex]
        \draw (0,0) rectangle +(1.,1);
        \fill (0,1) -- (.5,0) -- (1,1) -- cycle;
     \end{tikzpicture}
            }
\def\JJJbox{
     \begin{tikzpicture}[x=1.2ex,y=1.2ex]
        \draw (0,0) rectangle +(1.,1);
        \fill (0,0) -- (.5,1) -- (1,0) -- cycle;
     \end{tikzpicture}
            }

\def\diagbox{
     \begin{tikzpicture}[x=1.2ex,y=1.2ex]
        \draw (0,0) rectangle +(1,1);
        \fill (0,0) -- (0,1) -- (1,0) -- cycle;
     \end{tikzpicture}
            }
\def\twobox{\diagbox}
\def\fourbox{
     \begin{tikzpicture}[x=1.2ex,y=1.2ex]
        \draw (0,0) rectangle +(1,1);
        \fill (0,0) -- (1,1) -- (1,0) -- cycle;
     \end{tikzpicture}
            }
\def\circoldiagbox{
     \begin{tikzpicture}[x=1.2ex,y=1.2ex]
        \draw (.5,1.4)  circle (.4mm); 
        \draw (0,0) rectangle +(1,1);
        \fill (0,0) -- (0,1) -- (1,0) -- cycle;
     \end{tikzpicture}
            }

\def\JJoldiagbox{
     \begin{tikzpicture}[x=1.2ex,y=1.2ex]
        \fill (.1, 1.8) -- (.8, 1.8) -- (.5, 1.2) -- cycle;
        \draw (0,0) rectangle +(1,1);
        \fill (0,0) -- (0,1) -- (1,0) -- cycle;
     \end{tikzpicture}
            }

\newcount\xval
\newcount\yval
\newcount\difflim

\newcount\ssdoublexval
\newcount\ssdoubleyval

\def\ssray#1#2#3#4#5#6#7{
\ifnum#3>#1
\else 
\ifnum#4>#2
\else
   \ssmoveto{#3}{#4}\ssdrop{#5}
   \xval=#3
   \yval=#4
   \advance \xval by #6
   \advance \yval by #7
   \ssray{#1}{#2}{\xval}{\yval}{#5}{#6}{#7}
\fi
\fi
}

\def\sslineray#1#2#3#4#5#6#7#8{
\ifnum#3>#1
\else 
\ifnum#4>#2
\else
   \ssmoveto{#3}{#4}\ssline{#5}{#6}\ssmove{-#5}{-#6}
   \xval=#3
   \yval=#4
   \advance \xval by #7
   \advance \yval by #8
   \sslineray{#1}{#2}{\xval}{\yval}{#5}{#6}{#7}{#8}
\fi
\fi
}

\def\sseightfanray#1#2#3#4#5#6{
\ifnum#3>#1\else \ifnum#4>#2\else
   \xval=#3
   \yval=#4
   \ssmoveto{#3}{#4}
   \advance \xval by 1
   \advance \yval by 7
   \ifnum\xval>#1\else \ifnum\yval>#2\else
         \ssline{1}{7}\ssmove{-1}{-7}
   \fi\fi      
   \advance \xval by 2
   \advance \yval by -2
   \ifnum\xval>#1\else \ifnum\yval>#2\else
         \ssline{3}{5}\ssmove{-3}{-5}
   \fi\fi      
   \advance \xval by 2
   \advance \yval by -2
   \ifnum\xval>#1\else \ifnum\yval>#2\else
         \ssline{5}{3}\ssmove{-5}{-3}
   \fi\fi      
   \advance \xval by 2
   \advance \yval by -2
   \ifnum\xval>#1\else \ifnum\yval>#2\else
         \ssline{7}{1}\ssmove{-7}{-1}
   \fi\fi      
   \advance \xval by -7
   \advance \yval by -1
   \advance \xval by #5
   \advance \yval by #6
   \sseightfanray{#1}{#2}{\xval}{\yval}{#5}{#6}
\fi
\fi
}

\def\sscoloreightfanray#1#2#3#4#5#6{
\ifnum#3>#1\else \ifnum#4>#2\else
   \xval=#3
   \yval=#4
   \ssmoveto{#3}{#4}
   \advance \xval by 1
   \advance \yval by 7
   \ifnum\xval>#1\else \ifnum\yval>#2\else
         \ssline{1}{7}\ssmove{-1}{-7}
   \fi\fi      
   \advance \xval by 2
   \advance \yval by -2
   \ifnum\xval>#1\else \ifnum\yval>#2\else
         \ssline[color=blue]{3}{5}\ssmove{-3}{-5}
   \fi\fi      
   \advance \xval by 2
   \advance \yval by -2
   \ifnum\xval>#1\else \ifnum\yval>#2\else
         \ssline[color=green]{5}{3}\ssmove{-5}{-3}
   \fi\fi      
   \advance \xval by 2
   \advance \yval by -2
   \ifnum\xval>#1\else \ifnum\yval>#2\else
         \ssline[color=cyan]{7}{1}\ssmove{-7}{-1}
   \fi\fi      
   \advance \xval by -7
   \advance \yval by -1
   \advance \xval by #5
   \advance \yval by #6
   \sscoloreightfanray{#1}{#2}{\xval}{\yval}{#5}{#6}
\fi
\fi
}

\def\sseightfandoubleray#1#2#3#4#5#6#7#8{
\ifnum#3>#1\else \ifnum#4>#2\else
   \sseightfanray{#1}{#2}{#3}{#4}{#5}{#6}
   \ssdoublexval=#3
   \ssdoubleyval=#4
   \advance \ssdoublexval by #7
   \advance \ssdoubleyval by #8
   \sseightfandoubleray{#1}{#2}{\ssdoublexval}{\ssdoubleyval}{#5}{#6}{#7}{#8}
\fi\fi
}

\def\sscoloreightfandoubleray#1#2#3#4#5#6#7#8{
\ifnum#3>#1\else \ifnum#4>#2\else
   \sscoloreightfanray{#1}{#2}{#3}{#4}{#5}{#6}
   \ssdoublexval=#3
   \ssdoubleyval=#4
   \advance \ssdoublexval by #7
   \advance \ssdoubleyval by #8
   \sscoloreightfandoubleray{#1}{#2}{\ssdoublexval}{\ssdoubleyval}
                            {#5}{#6}{#7}{#8}
\fi\fi
}

\def\ssantiray#1#2#3#4#5#6#7{
\ifnum#3<#1
\else 
\ifnum#4<#2
\else
   \ssmoveto{#3}{#4}\ssdrop{#5}
   \xval=#3
   \yval=#4
   \advance \xval by #6
   \advance \yval by #7
   \ssantiray{#1}{#2}{\xval}{\yval}{#5}{#6}{#7}
\fi
\fi
}

\def\ssraydiff#1#2#3#4#5#6#7#8{
\difflim=#2
\advance \difflim by -#8
\ifnum#3>#1
\else 
\ifnum#4>#2
\else
   \ssmoveto{#3}{#4}\ssdrop{#5}
   \ifnum#4>\difflim 
   \else
       \ssdiff{#8}
   \fi
   \xval=#3
   \yval=#4
   \advance \xval by #6
   \advance \yval by #7
   \ssraydiff{#1}{#2}{\xval}{\yval}{#5}{#6}{#7}{#8}
\fi
\fi
}

\def\ssrayreddiff#1#2#3#4#5#6#7#8{
\difflim=#2
\advance \difflim by -#8
\ifnum#3>#1
\else 
\ifnum#4>#2
\else
   \ssmoveto{#3}{#4}\ssdrop{#5}
   \ifnum#4>\difflim 
   \else
       \ssreddiff{#8}
   \fi
   \xval=#3
   \yval=#4
   \advance \xval by #6
   \advance \yval by #7
   \ssrayreddiff{#1}{#2}{\xval}{\yval}{#5}{#6}{#7}{#8}
\fi
\fi
}

\def\ssdoubleray#1#2#3#4#5#6#7#8#9{
\ifnum#3>#1
\else 
\ifnum#4>#2
\else
   \ssray{#1}{#2}{#3}{#4}{#5}{#6}{#7}
   \ssdoublexval=#3
   \ssdoubleyval=#4
   \advance \ssdoublexval by #8
   \advance \ssdoubleyval by #9
   \ssdoubleray{#1}{#2}{\ssdoublexval}{\ssdoubleyval}
               {#5}{#6}{#7}{#8}{#9}
\fi
\fi
}

\def\ssantidoubleray#1#2#3#4#5#6#7#8#9{
\ifnum#3<#1
\else 
\ifnum#4<#2
\else
   \ssantiray{#1}{#2}{#3}{#4}{#5}{#6}{#7}
   \ssdoublexval=#3
   \ssdoubleyval=#4
   \advance \ssdoublexval by #8
   \advance \ssdoubleyval by #9
   \ssantidoubleray{#1}{#2}{\ssdoublexval}{\ssdoubleyval}{#5}{#6}{#7}{#8}{#9}
\fi
\fi
}

\def\ssrayblueline#1#2#3#4#5#6#7{
\ifnum#3>#1
\else 
\ifnum#4>#2
\else
   \ssmoveto{#3}{#4}\ssdrop{#5}\ssblueline{-2}
   \xval=#3
   \yval=#4
   \advance \xval by #6
   \advance \yval by #7
   \ssrayblueline{#1}{#2}{\xval}{\yval}{#5}{#6}{#7}
\fi
\fi
}

\def\ssraygreenline#1#2#3#4#5#6#7{
\ifnum#3>#1
\else 
\ifnum#4>#2
\else
   \ssmoveto{#3}{#4}\ssdrop{#5}\ssgreenline{-2}
   \xval=#3
   \yval=#4
   \advance \xval by #6
   \advance \yval by #7
   \ssraygreenline{#1}{#2}{\xval}{\yval}{#5}{#6}{#7}
\fi
\fi
}

\def\ssdoublerayblueline#1#2#3#4#5#6#7#8#9{
\ifnum#3>#1
\else 
\ifnum#4>#2
\else
   \ssrayblueline{#1}{#2}{#3}{#4}{#5}{#6}{#7}
   \ssdoublexval=#3
   \ssdoubleyval=#4
   \advance \ssdoublexval by #8
   \advance \ssdoubleyval by #9
   \ssdoublerayblueline{#1}{#2}{\ssdoublexval}{\ssdoubleyval}
                       {#5}{#6}{#7}{#8}{#9}
\fi
\fi
}

\def\ssdoubleraygreenline#1#2#3#4#5#6#7#8#9{
\ifnum#3>#1
\else 
\ifnum#4>#2
\else
   \ssraygreenline{#1}{#2}{#3}{#4}{#5}{#6}{#7}
   \ssdoublexval=#3
   \ssdoubleyval=#4
   \advance \ssdoublexval by #8
   \advance \ssdoubleyval by #9
   \ssdoubleraygreenline{#1}{#2}{\ssdoublexval}{\ssdoubleyval}
                       {#5}{#6}{#7}{#8}{#9}
\fi
\fi
}

\def\ssrayblueblueline#1#2#3#4#5#6#7{
\ifnum#3>#1
\else 
\ifnum#4>#2
\else
   \ssmoveto{#3}{#4}\ssdrop{#5}\ssblueline{-4}
   \xval=#3
   \yval=#4
   \advance \xval by #6
   \advance \yval by #7
   \ssrayblueblueline{#1}{#2}{\xval}{\yval}{#5}{#6}{#7}
\fi
\fi
}

\def\ssraygreengreenline#1#2#3#4#5#6#7{
\ifnum#3>#1
\else 
\ifnum#4>#2
\else
   \ssmoveto{#3}{#4}\ssdrop{#5}\ssgreenline{-4}
   \xval=#3
   \yval=#4
   \advance \xval by #6
   \advance \yval by #7
   \ssraygreengreenline{#1}{#2}{\xval}{\yval}{#5}{#6}{#7}
\fi
\fi
}

\def\ssdoublerayblueblueline#1#2#3#4#5#6#7#8#9{
\ifnum#3>#1
\else 
\ifnum#4>#2
\else
   \ssrayblueblueline{#1}{#2}{#3}{#4}{#5}{#6}{#7}
   \ssdoublexval=#3
   \ssdoubleyval=#4
   \advance \ssdoublexval by #8
   \advance \ssdoubleyval by #9
   \ssdoublerayblueblueline{#1}{#2}{\ssdoublexval}{\ssdoubleyval}
                           {#5}{#6}{#7}{#8}{#9}
\fi
\fi
}

\def\ssdoubleraygreengreenline#1#2#3#4#5#6#7#8#9{
\ifnum#3>#1
\else 
\ifnum#4>#2
\else
   \ssraygreengreenline{#1}{#2}{#3}{#4}{#5}{#6}{#7}
   \ssdoublexval=#3
   \ssdoubleyval=#4
   \advance \ssdoublexval by #8
   \advance \ssdoubleyval by #9
   \ssdoubleraygreengreenline{#1}{#2}{\ssdoublexval}{\ssdoubleyval}
                           {#5}{#6}{#7}{#8}{#9}
\fi
\fi
}

\def\ssreddiffray#1#2#3#4#5#6#7{
\difflim=#2
\advance \difflim by -#7
\ifnum#3>#1  
\else
\ifnum#4>\difflim 
\else
       \ssmoveto{#3}{#4}\ssreddiff{#7}
       \xval=#3
       \yval=#4
       \advance \xval by #5
       \advance \yval by #6
       \ssreddiffray{#1}{#2}{\xval}{\yval}{#5}{#6}{#7}
\fi
\fi
}

\def\reddiffdoubleray#1#2#3#4#5#6#7#8#9{
\difflim=#2
\advance \difflim by -#7
\ifnum#3 > #1  
\else
\ifnum#4 > \difflim 
\else
       \ssreddiffray{#1}{#2}{#3}{#4}{#5}{#6}{#7}
       \ssdoublexval=#3
       \ssdoubleyval=#4
       \advance \ssdoublexval by #8
       \advance \ssdoubleyval by #9
       \reddiffdoubleray{#1}{#2}{\ssdoublexval}{\ssdoubleyval}{#5}{#6}{#7}{#8}{#9}
\fi
\fi
}

\def\ssreddiffantiray#1#2#3#4#5#6#7{
\difflim=#1
\advance \difflim by 1
\ifnum#3< \difflim  
\else
\ifnum#4<#2
\else
       \ssmoveto{#3}{#4}\ssreddiff{#7}
       \xval=#3
       \yval=#4
       \advance \xval by #5
       \advance \yval by #6 
       \ssreddiffantiray{#1}{#2}{\xval}{\yval}{#5}{#6}{#7}
\fi
\fi
}

\def\reddiffantidoubleray#1#2#3#4#5#6#7#8#9{
\difflim=#2
\ifnum#3<#1  
\else
\ifnum#4<\difflim 
\else
       \ssreddiffantiray{#1}{#2}{#3}{#4}{#5}{#6}{#7}
       \ssdoublexval=#3
       \ssdoubleyval=#4
       \advance \ssdoublexval by #8
       \advance \ssdoubleyval by #9
       \reddiffantidoubleray{#1}{#2}{\ssdoublexval}{\ssdoubleyval}{#5}{#6}{#7}{#8}{#9}
\fi
\fi
}

\def\ssbluelineray#1#2#3#4#5#6#7{
\difflim=#2
\advance \difflim by -#7
\ifnum#3>#1  
\else
\ifnum#4>\difflim 
\else
       \ssmoveto{#3}{#4}\ssblueline{#7}
       \xval=#3
       \yval=#4
       \advance \xval by #5
       \advance \yval by #6
       \ssbluelineray{#1}{#2}{\xval}{\yval}{#5}{#6}{#7}
\fi
\fi
}

\def\sscurvedbluelineray#1#2#3#4#5#6#7{
\difflim=#2
\advance \difflim by -#7
\ifnum#3>#1  
\else
\ifnum#4>\difflim 
\else
       \ssmoveto{#3}{#4}\sscurvedblueline{#7}
       \xval=#3
       \yval=#4
       \advance \xval by #5
       \advance \yval by #6
       \sscurvedbluelineray{#1}{#2}{\xval}{\yval}{#5}{#6}{#7}
\fi
\fi
}

\def\ssbluelineantiray#1#2#3#4#5#6#7{
\ifnum#3<#1
\else
\ifnum#4<#2
\else
       \ssmoveto{#3}{#4}\ssblueline{#7}
       \xval=#3
       \yval=#4
       \advance \xval by #5
       \advance \yval by #6
       \ssbluelineantiray{#1}{#2}{\xval}{\yval}{#5}{#6}{#7}
\fi
\fi
}

\def\ssgreenlineantiray#1#2#3#4#5#6#7{
\ifnum#3<#1
\else
\ifnum#4<#2
\else
       \ssmoveto{#3}{#4}\ssgreenline{#7}
       \xval=#3
       \yval=#4
       \advance \xval by #5
       \advance \yval by #6
       \ssgreenlineantiray{#1}{#2}{\xval}{\yval}{#5}{#6}{#7}
\fi
\fi
}

\def\bluelinedoubleray#1#2#3#4#5#6#7#8#9{
\difflim=#2
\advance \difflim by -#7
\ifnum#3>#1  
\else
\ifnum#4>\difflim 
\else
       \ssbluelineray{#1}{#2}{#3}{#4}{#5}{#6}{#7}
       \ssdoublexval=#3
       \ssdoubleyval=#4
       \advance \ssdoublexval by #8
       \advance \ssdoubleyval by #9
       \bluelinedoubleray{#1}{#2}{\ssdoublexval}{\ssdoubleyval}{#5}{#6}{#7}{#8}{#9}
\fi
\fi
}

\def\bluelineantidoubleray#1#2#3#4#5#6#7#8#9{
\difflim=#2
\advance \difflim by -#7
\ifnum#3<#1  
\else
\ifnum#4<\difflim 
\else
       \ssbluelineantiray{#1}{#2}{#3}{#4}{#5}{#6}{#7}
       \ssdoublexval=#3
       \ssdoubleyval=#4
       \advance \ssdoublexval by #8
       \advance \ssdoubleyval by #9
       \bluelineantidoubleray{#1}{#2}{\ssdoublexval}{\ssdoubleyval}{#5}{#6}{#7}{#8}{#9}
\fi
\fi
}

\def\greenlineantidoubleray#1#2#3#4#5#6#7#8#9{
\difflim=#2
\advance \difflim by -#7
\ifnum#3<#1  
\else
\ifnum#4<\difflim 
\else
       \ssgreenlineantiray{#1}{#2}{#3}{#4}{#5}{#6}{#7}
       \ssdoublexval=#3
       \ssdoubleyval=#4
       \advance \ssdoublexval by #8
       \advance \ssdoubleyval by #9
       \greenlineantidoubleray{#1}{#2}{\ssdoublexval}{\ssdoubleyval}{#5}{#6}{#7}{#8}{#9}
\fi
\fi
}

\def\bluessray#1#2#3#4#5#6#7{
\ifnum#3>#1
\else 
\ifnum#4>#2
\else
   \ssmoveto{#3}{#4}\ssdrop[color=blue]{#5}
   \xval=#3
   \yval=#4
   \advance \xval by #6
   \advance \yval by #7
   \bluessray{#1}{#2}{\xval}{\yval}{#5}{#6}{#7}
\fi
\fi
}

\def\violetssray#1#2#3#4#5#6#7{
\ifnum#3>#1
\else 
\ifnum#4>#2
\else
   \ssmoveto{#3}{#4}\ssdrop[color=violet]{#5}
   \xval=#3
   \yval=#4
   \advance \xval by #6
   \advance \yval by #7
   \violetssray{#1}{#2}{\xval}{\yval}{#5}{#6}{#7}
\fi
\fi
}

\def\violetssantiray#1#2#3#4#5#6#7{
\ifnum#3<#1
\else 
\ifnum#4<#2
\else
   \ssmoveto{#3}{#4}\ssdrop[color=violet]{#5}
   \xval=#3
   \yval=#4
   \advance \xval by #6
   \advance \yval by #7
   \violetssantiray{#1}{#2}{\xval}{\yval}{#5}{#6}{#7}
\fi
\fi
}

\def\bluessantiray#1#2#3#4#5#6#7{
\ifnum#3<#1
\else 
\ifnum#4<#2
\else
   \ssmoveto{#3}{#4}\ssdrop[color=blue]{#5}
   \xval=#3
   \yval=#4
   \advance \xval by #6
   \advance \yval by #7
   \bluessantiray{#1}{#2}{\xval}{\yval}{#5}{#6}{#7}
\fi
\fi
}

\def\ssbluelinedoubleray#1#2#3#4#5#6#7#8#9{
\ifnum#3>#1
\else 
\ifnum#4>#2
\else
   \ssbluelineray{#1}{#2}{#3}{#4}{#5}{#6}{#7}
   \ssdoublexval=#3
   \ssdoubleyval=#4
   \advance \ssdoublexval by #8
   \advance \ssdoubleyval by #9
   \ssbluelinedoubleray{#1}{#2}{\ssdoublexval}{\ssdoubleyval}
                       {#5}{#6}{#7}{#8}{#9}
\fi
\fi
}

\def\ssgreenlineray#1#2#3#4#5#6#7{
\difflim=#2
\advance \difflim by -#7
\ifnum#3>#1  
\else
\ifnum#4>\difflim 
\else
       \ssmoveto{#3}{#4}\ssgreenline{#7}
       \xval=#3
       \yval=#4
       \advance \xval by #5
       \advance \yval by #6
       \ssgreenlineray{#1}{#2}{\xval}{\yval}{#5}{#6}{#7}
\fi
\fi
}

\def\ssgreenlinedoubleray#1#2#3#4#5#6#7#8#9{
\ifnum#3>#1
\else 
\ifnum#4>#2
\else
   \ssgreenlineray{#1}{#2}{#3}{#4}{#5}{#6}{#7}
   \ssdoublexval=#3
   \ssdoubleyval=#4
   \advance \ssdoublexval by #8
   \advance \ssdoubleyval by #9
   \ssgreenlinedoubleray{#1}{#2}{\ssdoublexval}{\ssdoubleyval}
                       {#5}{#6}{#7}{#8}{#9}
\fi
\fi
}

\newcommand{\alphax }{\overline{\nu } }
\newcommand{\betax }{\xi}
\usepackage{tmf5}   

\begin{document}

\title[The slice spectral sequence for the $C_{4}$ analog of real
$K$-theory]{The slice spectral sequence for the
$C_{4}$ analog of real $K$-theory}

\author{Michael~A.~Hill}
\address{Department of Mathematics \\ University of California Los Angeles\\
 Los Angeles, CA 90095}
\email{mikehill@math.ucla.edu}

\author{Michael~J.~Hopkins}
\address{Department of Mathematics \\ Harvard University 
\\Cambridge, MA 02138}
\email{mjh@math.harvard.edu}

\author{Douglas~C.~Ravenel}
\address{Department of Mathematics \\ University of Rochester 
\\Rochester, NY 14627}
\email{doug@math.rochester.edu}

\thanks{The authors were supported by DARPA Grant 
FA9550-07-1-0555 and NSF Grants DMS-0905160 , DMS-1307896\dots }

\date{\today}

\begin{abstract}
We describe the slice spectral sequence
of a 32-periodic $C_{4}$-spectrum $\KH$ related to the
$C_{4}$ norm 
\begin{displaymath}
{MU^{((C_{4}))}=N_{C_{2}}^{C_{4}}MU_{\reals}}
\end{displaymath}
 
\noindent of the real cobordism spectrum $MU_{\reals}$.  We will give
it as a {\SS} of Mackey functors converging to the graded Mackey
functor $\upi_{*}\KH$, complete with differentials and exotic
extensions in the Mackey functor structure.  

The slice spectral sequence for the 8-periodic real $K$-theory
spectrum $K_{\reals}$ was first analyzed by Dugger.  The $C_{8}$
analog of $\KH$ is 256-periodic and detects the Kervaire invariant
classes $\theta_{j}$.  A partial analysis of its slice spectral
sequence led to the solution to the Kervaire invariant problem, namely
the theorem that $\theta_{j}$ does not exist for $j\geq 7$.
\end{abstract}



\maketitle

\tableofcontents
\listoftables

\section{Introduction}\label{sec-intro}

In \cite{HHR} we derived the main theorem about the Kervaire invariant
elements from some properties of a $C_{8}$-{\eqvr} spectrum we called
$\Omega $ constructed as follows.  We started with the
$C_{2}$-spectrum $MU_{\reals}$, meaning the usual complex cobordism
spectrum $MU$ equipped with a $C_{2}$ action defined in terms of
complex conjugation.  

Then we defined a functor $N_{C_{2}}^{C_{8}}$, the norm of
\cite[\S2.2.3]{HHR} which we abbreviate here by $N_{2}^{8}$, from the
category of $C_{2}$-spectra to that of $C_{8}$-spectra.  Roughly
speaking, given a $C_{2}$-spectrum $X$, $N_{2}^{8}X$ is underlain by
the fourfold smash power $X^{\wedge 4}$ where a generator $\gamma $ of
$C_{8}$ acts by cyclically permuting the four factors, each of which
is invariant under the given action of the subgroup $C_{2}$.  In a similar
way one can define a functor $N_{H}^{G}$ from $H$-spectra to
$G$-spectra for any finite groups $H\subseteq G$.

A $C_{8}$-spectrum such as $N_{2}^{8}MU_{\reals}$, which is a
commutative ring spectrum, has {\eqvr} homotopy groups indexed by $RO
(C_{8})$, the orthogonal {\rep} ring for the group $C_{8}$.  One
element of the latter is $\rho_{8}$, the regular {\rep}.  In
\cite[\S9]{HHR} we defined a certain element $D\in
\pi_{19\rho_{8}}N_{2}^{8}MU_{\reals}$ and then formed the associated
mapping telescope, which we denoted by $\Omega_{\mathbb{O}}$.  The
symbol $\mathbb{O}$ was chosen to suggest a connection with the
octonions, but there really is none apart from the fact that the
octonions are 8-dimensional like $\rho _{8}$.  

$\Omega_{\mathbb{O}}$ is also a
$C_{8}$-{\eqvr} commutative ring spectrum.  We then proved that it is
{\eqvr}ly equivalent to $\Sigma^{256}\Omega_{\mathbb{O}}$; we call
this result the Periodicity Theorem.  Then our spectrum $\Omega $ is
$\Omega_{\mathbb{O}}^{C_{8}}$, the fixed point spectrum of
$\Omega_{\mathbb{O}}$.

It is possible to do this with $C_{8}$ replaced by $C_{2^{n}}$ for any
$n$. The dimension of the periodicity is then $2^{1+n+2^{n-1}}$.  For
example it is 32 for the group $C_{4}$ and $2^{13}$ for $C_{16}$.  We
chose the group $C_{8}$ because it is the smallest that suits our
purposes, namely it is the smallest one yielding a fixed point
spectrum that detects the Kervaire invariant elements $\theta_{j}$.

We know almost nothing about $\pi_{*}\Omega $, only that it is
periodic with periodic 256, that $\pi_{-2}=0$ (the Gap Theorem of
\cite[\S8]{HHR}), and that when $\theta_{j}$ exists its image in
$\pi_{*}\Omega $ is nontrivial (the Detection Theorem of
\cite[\S11]{HHR}).

We also know, although we did not say so in \cite{HHR}, that more
explicit computations would be much easier if we cut
$N_{2}^{8}MU_{\reals}$ down to size in the following way.  Its
underlying homotopy, meaning that of the spectrum $MU^{\wedge 4}$, is
known classically to be a polynomial algbera over the integers with
four generators (cyclically permuted up to sign by the group action)
in every positive even dimension.  This can be proved with methods
described by Adams in \cite{Ad:SHGH}. For the cyclic group $C_{2^{n}}$
one has $2^{n-1}$ generators in each positive even degree.
Specific generators $r_{i,j}\in \pi_{2i}MU^{\wedge 2^{n-1}}$ for $i>0$
and $0\leq j<^{n-1}$ are defined in \cite[\S5.4.2]{HHR}.

{\em There is a way to kill all the generators above dimension 2k}
that was described in \cite[\S2.4]{HHR}.  Roughly speaking, let $A$ be
a wedge of suspensions of the sphere spectrum, one for each monomial in
the generators one wants to kill.  One can define a multiplication and
group action on $A$ corresponding to the ones in $\pi_{*}MU^{\wedge
4}$.  Then one has a map $A\to MU^{\wedge 4}$ whose restriction to
each summand represents the corresponding monomial, and a map $A\to
S^0$ (where the target is the sphere spectrum, not the space $S^{0}$)
sending each positive dimensional summand to a point.  This leads to
two maps
\begin{displaymath}
S^{0}\wedge A\wedge MU^{\wedge 4}\rightrightarrows S^{0}\wedge MU^{\wedge 4}
\end{displaymath}

\noindent whose coequalizer we denote by $S^{0}\smashove{A}MU^{\wedge
4}$.  Its homotopy is the quotient of $\pi_{*}MU^{\wedge 4}$ obtained
by killing the polynomial generators above diimension $2k$.  The
construction is {\eqvr}, meaning that $S^{0}\smashove{A}MU^{\wedge 4}$
underlies a $C_{8}$-spectrum.  

In \cite[\S7]{HHR} we showed that for $k=0$ the spectrum we get is the
integer {\SESM} spectrum $H\Z$; we called this result the Reduction
Theorem. In the non{\eqvr} case this is obvious. We are in effect
attaching cells to $MU^{\wedge 4}$ to kill all of its homotopy groups
in positive dimensions, which amounts to constructing the 0th
Postnikov section.  In the {\eqvr} case the proof is more
delicate.

Now consider the case $k=1$, meaning that we are killing the
polynomial generators above dimension 2.  Classically we know that
doing this to $MU$ (without the $C_{2}$-action) produces the
connective complex $K$-theory spectrum, some times denoted by $k$,
$bu$ or (2-locally) $BP\langle 1\rangle$.  Inverting the Bott element
via a mapping telescope gives us $K$ itself, which is of course
2-periodic.  In the $C_{2}$-{\eqvr} case one gets the ``real $K$-theory''
spectrum $K_{\reals}$ first studied by Atiyah in \cite{Atiyah:KR}.  It
turns out to be 8-periodic and its fixed point spectrum is $KO$, which
is also referred to in other contexts as real $K$-theory.

The spectrum we get by killing the generators above dimension 2 in
the\linebreak $C_{8}$-spectrum $N_{2}^{8}MU_{\reals}$ will be denoted
analogously by $k_{[3]}$.  We can invert the image of $D$ by forming a
mapping telescope, which we will denote by $K_{[3]}$.  More
generally we denote by $k_{[n]}$ the spectrum obtained from
$N_{C_{2}}^{C_{2^{n}}}MU_{\reals}$ by killing all generators above
dimension 2. In particular $k_{[1]}=k_{\reals}$.  Then we denote the
mapping telescope (after defining a suitable $D$) by $K_{[n]}$ and its
fixed point set by $KO_{[n]}$.  

For $n\geq 3$, $KO_{[n]}$ also has a Periodicity, Gap and Detection
Theorem, so it could be used to prove the Kervaire invariant theorem.

{\em Thus $K_{[3]}$ is a substitute for $\Omega _{\mathbb{O}}$
with much smaller and therefore more tractable homotopy groups.}  A
detailed study of them might shed some light on the fate of
$\theta_{6}$ in the 126-stem, the one hypothetical Kervaire invariant
element whose status is still open.  {\em If we could show that
$\pi_{126}KO_{[3]}=0$, that would mean that
$\theta_{6}$ does not exist.}

The computation of the {\eqvr} homotopy $\upi_{*}K_{[3]}$ at this time
is daunting.  {\em The purpose of this paper is to do a similar
computation for the group $C_{4}$ as a warmup exercise.}  In the
process of describing it we will develope some techniques that are
likley to be needed in the $C_{8}$ case.  We start with
$N_{2}^{4}MU_{\reals}$, kill its polynomial generators (of which there
are two in every positive even dimension) above dimension 2 as
described previously, and then invert a certain element in
$\pi_{4\rho_{4}}$.  We denote the resulting spectrum by $\KH$, see
\ref{def-kH} below.  This spectrum is known to be 32-periodic.  In an
earlier draft of this paper it was denoted by $K_{\mathbf{H}}$.

The computational tool for finding these homotopy groups is the slice
spectral sequence introduced in \cite[\S4]{HHR}.  Indeed we do not
know of any other way to do it.  For $K_{\reals}$ it was first
analyzed by Dugger \cite{Dugger} and his work is described below in
\S\ref{sec-Dugger}.  In this paper we will study the slice spectral
sequence of Mackey functors associated with $\KH$.  We will rely
extensively on the results, methods and terminology of \cite{HHR}.

{\em We warn the reader that the computation for $\KH$ is more
intricate than the one for $\KR$.}  For example, the slice spectral
sequence for $\KR$, which is shown in Figure \ref{fig-KR}, involves
five different Mackey functors for the group $C_{2}$.  We abbreviate
them with certain symbols indicated in Table \ref{tab-C2Mackey}.  The
one for $\KH$, partly shown in Figure \ref{fig-KH}, involves over
twenty Mackey functors for the group $C_{4}$, with symbols indicated
in Table \ref{tab-C4Mackey}.

\begin{figure}
\begin{center}
\includegraphics[width=11cm]{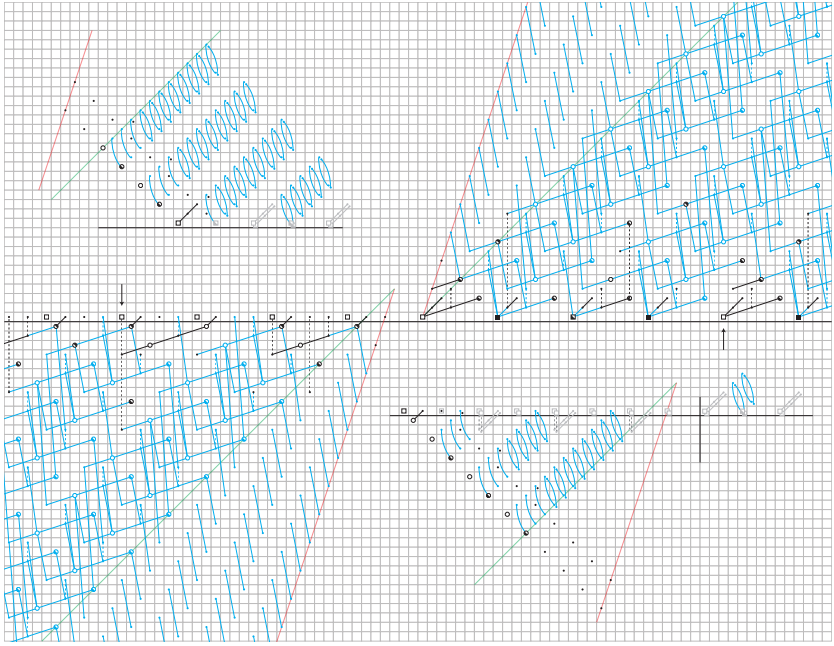} \caption[The 2008
poster.]  {The 2008 poster. The first and third quadrants show
$\EE_{4} (G/G)$ of the slice {\SS} for $\KH$ with the elements of
Prop. \ref{prop-perm} excluded.  The second quadrant indicates
$d_{3}$s as in Figures \ref{fig-sseq-7a} and \ref{fig-sseq-7b}.  The
fourth quadrant indicates comparable $d_{3}$s in the third quadrant of
the slice {\SS} as in Figures \ref{fig-sseq-8a} and
\ref{fig-sseq-8b}.}  \label{fig-poster}
\end{center}
\end{figure}

Part of this {\SS} is also illustrated in an unpublished poster produced in
late 2008 and shown in Figure \ref{fig-poster}.  It shows the {\SS}
converging to the homotopy of the fixed point spectrum $\KH^{C_{4}}$.
The corresponding {\SS} of Mackey functors converges to the graded
Mackey functor $\upi_{*}\KH$.  

In both illustrations some patterns of $d_{3}$s and families of
elements in low filtration are excluded to avoid clutter. In the
poster, representative examples of these are shown in the second and
fourth quadrants, the {\SS} itself being concentrated in the first and
third quadrants.  In this paper those patterns are spelled out in
\S\ref{sec-C2diffs} and \S\ref{sec-C4diffs}.

\bigskip

We now outline the rest of the paper.  Briefly, the next five sections
introduce various tools we need. Our objects of study, the spectra
$\kH$ and $\KH$, are formally introduced in \S\ref{sec-kH}.  Dugger's
computation for $\KR$ is recalled in \S\ref{sec-Dugger}.  The final
six sections describe the computation for $\kH$ and $\KH$.

In more detail, \S\ref{sec-nonsense} collects some notions from
{\eqvr} stable homotopy theory with an emphasis on Mackey functors.
Definition \ref{def-graded} introduces new notation that we will
ocasionally need.

\S\ref{sec-HZZ} concerns the {\eqvr} analog of the homology of a point
namely, the $RO (G)$-graded homotoy of the integer {\SESM} spectrum
$H\Z$. In particular Lemma \ref{lem-aeu} describes some relations
among certain elements in it including the ``gold relation'' between
$a_{V}$ and $u_{V}$.
   
\S\ref{sec-gendiffs} describes some general properties of spectral
sequences of Mackey functors.  These include Theorem \ref{thm-exotic}
about the relation between differential and exotic extensions in the
Mackey functor structure and Theorem \ref{thm-normdiff} on the norm of
a differential.

\S\ref{sec-C4}  lists some concise symbols for various specific Mackey functors
for the groups $C_{2}$ and $C_{4}$ that we will need.  Such functors
can be spelled out explicitly by means of Lewis diagrams
(\ref{eq-Lewis}), which we usually abbreviate by symbols shown in
Tables \ref{tab-C2Mackey} and \ref{tab-C4Mackey}.

In \S\ref{sec-chain} we study some chain complexes of Mackey functors
that arise as cellular chain complexes for $G$-CW complexes of the
form $S^{V}$.

In \S\ref{sec-kH} we formally define (in \ref{def-kH}) the
$C_{4}$-spectra of interest in this paper, $\kH$ and $\KH$.

In \S\ref{sec-Dugger} we describe the slice {\SS} for an easier case,
the $C_{2}$-spectrum for real $K$-theory, $K_{\reals}$.  This is due to
Dugger \cite{Dugger} and serves as a warmup exercise for us. It turns
out that everything in the {\SS } is formally determined by the
structure of its $\EE_{2}$-term and Bott periodicity.

In \S\ref{sec-more} we introduce various  elements
in the homotopy groups of $\kH$ and $\KH$.  They are collected in Table
\ref{tab-pi*}, which spans several pages.
In \S\ref{sec-slices} we determine the $\EE_{2}$-term of the
slice {\SS } for $\kH$ and $\KH$.

In \S\ref{sec-diffs} we use the Slice Differentials Theorem of
\cite{HHR} to determine some differentials in our {\SS}.

In \S\ref{sec-C2diffs} we examine the $C_{4}$-spectrum $\kH$ as a
$C_{2}$-spectrum.  This leads to a calculation only slightly more
complicated than Dugger's.  It gives a way to remove a lot of clutter
from the $C_{4}$ calculation.

In \S\ref{sec-C4diffs} we determine the $\EE_{4}$-term of our {\SS}.
It is far smaller than $\EE_{2}$ and the results of
\S\ref{sec-C2diffs} enable us to ignore most of it.  What is left is
small enough to be shown legibly in the {\SS } charts of Figures
\ref{fig-E4} and \ref{fig-KH}.  They illustrate integrally graded (as
opposed to $RO (C_{4})$-graded) spectral sequences of Mackey functors,
which are discussed in \S\ref{sec-C4}. In order to read these charts
one needs to refer to Table \ref{tab-C4Mackey} which defines the
``heiroglyphic'' symbols we use for the specific Mackey functors that
we neede.

We finish the calculation in \S\ref{sec-higher} by dealing with the
remaining differentials and exotic Mackey functor extensions.  It
turns out that they are all formal consequences of $C_{2}$
differentials of the previous section along with the results of
\S\ref{sec-gendiffs}.

The result is a complete description of the {\em integrally graded}
portion of $\upi_{\star}\kH$.  It is best seen in the {\SS } charts of
Figures \ref{fig-E4} and \ref{fig-KH}.  Unfortunately we do not have a
clean description, much less an effective way to display the full $RO
(C_{4})$-graded homotopy groups.  

For $G=C_{2}$, the two irreducible orthogonal {\rep}s are the trivial
one of degree 1, denoted by the symbol 1, and the sign {\rep} denoted
by $\sigma $.  Thus $RO (G)$ is additvley a free abelain group of rank
2, and the {\SS } of interest is trigraded. In the $RO (C_{2})$-graded
homotopy of $\KR$, a certain element of degree $1+\sigma $ (the degree
of the regular {\rep} $\rho_{2}$) is invertible.  This means that each
component of $\upi_{\star}\KR$ is canonically isomorphic to a Mackey
functor indexed by an ordinary integer.  See Theorem
\ref{thm-ROG-graded} for a more precise statement.  Thus the full
(trigraded)$RO (C_{2})$-graded slice {\SS } is determined by bigraded
one shown in Figure \ref{fig-KR}.

For $G=C_{4}$, the {\rep} ring $RO (G)$ is additively a free abelian group of
rank 3, so it leads to a quadrigraded {\SS }.  The three irreducible
{\rep}s are the trivial and sign {\rep}s 1 and $\sigma $ (each having
degree one) and a degree two {\rep} $\lambda $ given by a rotation of
the plane $\reals^{2}$ of order 4.  The regular {\rep} $\rho_{4}$ is
isomorphic to $1+\sigma +\lambda $.  As in the case of $\KR$, there is
an invertible element $\normrbar_{1}$ (see Table \ref{tab-pi*}) in
$\upi_{\star}\KH$ of degree $\rho_{4}$.  This means we can reduce the
quadigraded slice {\SS } to a trigraded one, but finding a full
description of it is a problem for the future.

\section{Recollections about {\eqvr} stable homotopy theory}
\label{sec-nonsense}

We first discuss some structure on the {\eqvr} homotopy groups of a\linebreak 
$G$-spectrum $X$.  {\em We will assume throughout that $G$ is a finite
cyclic $p$-group.}  This means that its subgroups are well ordered by
inclusion and each is uniquely determined by its order.  The results
of this section hold for any prime $p$, but the rest of the paper
concerns only the case $p=2$.  We will define several maps indexed by
pairs of subgroups of $G$.  {\em We will often replace these indices
by the orders of the subgroups, sometimes denoting $|H|$ by $h$.}

 The homotopy groups  can be defined in terms of finite $G$-sets $T$
Let
\begin{displaymath}
\upi_{0}^{G}X (T) = [T_{+}, X]^{G},
\end{displaymath}

\noindent be the set of homotopy classes of {\eqvr} maps from $T_{+}$,
the suspension spectrum of the union of $T$ with a disjoint base
point, to the spectrum $X$.  We will often omit $G$ from the notation
when it is clear from the context.  For an orthogonal {\rep} $V$ of
$G$, we define
\begin{displaymath}
\upi_{V}X (T) = [S^{V}\wedge T_{+}, X]^{G}.
\end{displaymath}

\noindent As an $RO (G)$-graded contravariant abelian group valued
functor of $T$, this converts disjoint unions to direct sums. This
means it is determined by its values on the sets $G/H$ for subgroups
$H\subseteq G$.  

Since $G$ is abelian, $H$ is normal and $\upi_{V}X (G/H)$ is a $Z[G/H]$-module.

Given subgroups $K\subseteq H\subseteq G$, one has pinch and fold maps
between the $H$-spectra $H/H_{+}$ and $H/K_{+}$.  This leads to a diagram
\begin{numequation}\label{eq-pinch-fold}
\begin{split}
\xymatrix
@R=1mm
@C=4mm
{
    &H/H_{+}\ar@<.5ex>[rr]^(.5){\mbox{pinch} }
        &{}
            &H/K_{+}\ar@<.5ex>[ll]^(.5){\mbox{fold} }\\
    &   &\ar@{=>}[ddd]^(.6){G_{+}\smashove{H} (\cdot)}\\
    &   &\\
    &   &\\
    &   &\\
G/H_{+}\ar@{=}[r]^(.5){}
    &G_{+}\smashove{H}H/H_{+}\ar@<.5ex>[rr]^(.5){\mbox{pinch} }
        &{} &G_{+}\smashove{H}H/K_{+}\ar@{=}[r]^(.5){}
                                     \ar@<.5ex>[ll]^(.5){\mbox{fold} } 
                &G_{+}\smashove{K}K/K_{+}\ar@{=}[r]^(.5){} 
                   &G/K_{+}.
}
\end{split}
\end{numequation}%

\noindent Note that while the fold map is induced by a map of
$H$-sets, the pinch map is not.  It only exists in the stable
category.  

\begin{defin}\label{def-fixed-point-maps}
{\bf The Mackey functor structure maps in $\upi_{V}^{G}X$.}  
The fixed point transfer
and restriction maps
\begin{displaymath}
\xymatrix
@R=5mm
@C=15mm
{
{\upi_{V}X (G/H)}\ar@<-.5ex>[r]_(.5){\Res_{K}^{H} }
    &{\upi_{V}X (G/K)}\ar@<-.5ex>[l]_(.5){\Tr_{K}^{H} }
}
\end{displaymath}

\noindent are the ones induced by the composite maps in the bottom row
of (\ref{eq-pinch-fold}).
\end{defin}

These satisfy the formal properties needed to make $\upi_{V}X$ into a
Mackey functor; see \cite[Def. 3.1]{HHR}.  They are usually referred
to simply as the transfer and restriction maps.  We use the words
``fixed point'' to distinguish them from another similar pair of maps
specified below in Definition \ref{def-gp-action}.

We remind the reader that a Mackey functor $\underline{M}$ for a
finite group $G$ assigns an abelian group $\underline{M} (T)$ to every
finite $G$-set $T$.  It converts disjoint unions to direct sums. It  is
therefore determined by its values on orbits, meaning $G$-sets for the
form $G/H$ for various subgroups $H$ of $G$.  For subgroups
$K\subseteq H\subseteq G$, one has a map of $G$-sets $G/K\to G/H$.  In
categorical language $\underline{M}$ is actually a pair of functors,
one covariant and one contravariant, both behaving the same way on
objects.  Hence we get maps both ways between $\underline{M} (G/K)$
and $\underline{M} (G/H)$.  For the Mackey functor $\upi_{V}X$, these
are the two maps of \ref{def-fixed-point-maps}.

One can generalize the definition of a Mackey functor by replacing the
target category of abelian groups by one's favorite abelian category,
such as that of $R$-modules over graded abelian groups.

\begin{defin}\label{def-green}
A {\bf graded Green functor $\underline{R}_{*}$ for a group $G$} is a
Mackey functor for $G$ with values in the category of graded abelian
groups such that $\underline{R}_{*} (G/H)$ is a graded commutative
ring for each subgroup $H$ and for each pair of subgroups $K\subseteq
H\subseteq G$, the restriction map $\Res_{K}^{H}$ is a ring
homomorphism and the transfer map $\Tr_{K}^{H}$ satisfies the {\bf
Frobenius relation}
\begin{displaymath}
\Tr_{K}^{H}(\Res_{K}^{H}(a)b)=a (\Tr_{K}^{H}(b))\qquad
\mbox{for $a\in \underline{R}_{*} (G/H)$ and $b\in
\underline{R}_{*} (G/K)$}.
\end{displaymath}
\end{defin}

When $X$ is a ring spectrum, we have the {\em fixed point Frobenius
relation}
\begin{numequation}
\label{eq-Frob} \Tr_{K}^{H}(\Res_{K}^{H}(a)b)=a (\Tr_{K}^{H}(b))\qquad
\mbox{for $a\in \upi_{\star}X (G/H)$ and $b\in
\upi_{\star}X (G/K)$}.
\end{numequation}%
In particular this means that 
\begin{numequation}
\label{eq-Frobcor}
a(\Tr_{K}^{H}(b)) = 0\qquad \mbox{when } \Res_{K}^{H}(a)=0.
\end{numequation}%

For a {\rep} $V$ of $G$, the group 
\begin{displaymath}
\upi_{V}^{G}X (G/H)=\pi_{V}^{H}X=[S^{V},X]^{H}
\end{displaymath}

\noindent is isomorphic to 
\begin{displaymath}
[S^{0},S^{-V}\wedge X]^{H}=\pi_{0} (S^{-V}\wedge X)^{H}.
\end{displaymath}

\noindent However fixed points do not respect smash products, so we
cannot equate this group with
\begin{displaymath}
\pi_{0} (S^{-V^{H}}\wedge X^{H}) = [S^{V^{H}},X^{H}]=\pi_{|V^{H}|}X^{H}
=\upi^{G}_{|V^{H}|}X (G/H).
\end{displaymath}

Conversely a $G$-{\eqvr} map $S^{V}\to X$ represents an element in
\begin{displaymath}
[S^{V},X]^{G}=\pi_{V}^{G}X=\upi^{G}_{V}X (G/G).
\end{displaymath}

The following notion is useful.

\begin{defin}\label{def-MF-induction}
{\bf Mackey functor induction and restriction.}  For  s subgroup $H$ of $G$
and an $H$-Mackey functor $\underline{M}$, the
induced $G$-Mackey functor $\uparrow_{H}^{G}\underline{M}$ is given by
\begin{displaymath}
\uparrow_{H}^{G}\underline{M} (T)=\underline{M} (i_{H}^{*}T)
\end{displaymath}

\noindent for each finite $G$-set $T$, where $i_{H}^{*}$ denotes the
forgetful functor from $G$-sets (or spaces or spectra) to $H$-sets.  

For a $G$-Mackey functor $\underline{N}$, the restricted  $H$-Mackey functor
$\downarrow_{H}^{G}\underline{N}$ is given by
\begin{displaymath}
\downarrow_{H}^{G}\underline{N} (S)=\underline{N} (G\times_{H}S)
\end{displaymath}

\noindent for each finite $H$-set $S$.
\end{defin}

This notation is due to Th\'evenaz-Webb \cite{Thevenaz-Webb}. They put
the decorated arrow on the right and denote $G\times_{H}S$ by
$S\uparrow_{H}^{G}$ and $i_{H}^{*}T$ by $T\downarrow_{H}^{G}$.

We also need notation for $X$ as an $H$-spectrum
for subgroups $H\subseteq G$.  For this purpose we will enlarge the
orthogonal {\rep} ring of $G$, $RO (G)$, to the {\rep} ring Mackey
functor $\underline{RO} (G)$ defined by $\underline{RO} (G) (G/H)=RO
(H)$.  This was the motivating example for the definition of a Mackey
functor in the first place.  In it the transfer map on a {\rep} $V$ of
$H$ is the induced {\rep} of a supergroup $K\supseteq H$, and its
restriction to a subgroup is defined in the obvious way. In particular
the restriction of the transfer of $V$ is $|K/H|V$.

More generally for a finite $G$-set $T$, $\underline{RO} (G) (T)$ is
the ring (under pointwise direct sum and tensor product) of functors
to the category of finite dimensional orthogonal real vector spaces
from $B_{G}T$, the split groupoid (see \cite[A1.1.22]{Rav:MU}) whose
objects are the elements of $T$ with morphisms defined by the action
of $G$.

\begin{defin}\label{def-graded}
{\bf $\underline{RO} (G)$-graded homotopy groups.}  For each
$G$-spectrum $X$ and each pair $(H,V)$ consisting of a subgroup
$H\subseteq G$ and a virtual orthogonal {\rep} $V$ of $H$, let the
$G$-Mackey functor $\upi_{{H,V}} (X)$ be defined by
\begin{displaymath}
\upi_{{H,V}} (X) (T)
:= \left[(G_{+}\smashove{H}S^{V})\wedge T_{+}, X \right]^{G}
\cong  \left[S^{V}\wedge i_{H}^{*}T_{+}, i_{H}^{*}X \right]^{H}
 = \upi_{V}^{H} (i_{H}^{*}X) (i_{H}^{*}T),
\end{displaymath}

\noindent for each finite $G$-set $T$.  Equivalently, $\upi_{{H,V}}
(X)=\uparrow_{H}^{G}\upi_{V}^{H} (i_{H}^{*}X)$ (see
\ref{def-MF-induction}) as Mackey functors.  We will often denote
$\upi_{G,V}$ by $\upi_{V}^{G}$ or $\upi_{V}$.
\end{defin}

We will be studying the $\underline{RO} (G)$-graded slice {\SS }
$\left\{\EE_{r}^{s,\star} \right\}$ of Mackey functors with
$r,s\in \Z$ and $\star\in \underline{RO} (G)$.  We will use the
notation $\EE_{r}^{s,(H,V)}$ for such Mackey functors,
abbreviating to $\EE_{r}^{s,V}$ when the subgroup is
$G$. Most of our spectral sequence charts will display the values of
$\EE_{2}^{s,t}$ for integral values of $t$ only.

The following definition should be compared with \cite[(2.3)]{Ad:Prereq}.

\bigskip
\begin{defin}\label{def-homeo}
{\bf An {\eqvr} homeomorphism.} Let $X$ be a $G$-space and $Y$ an
$H$-space for a subgroup $H\subseteq G$.  We define the {\eqvr}
homeomorphism
\begin{displaymath}
\tilde{u}_{H}^{G} (Y,X):G\times_{H} (Y\times i_{H}^{*}X) \to 
                        (G\times_{H}Y)\times X
\end{displaymath}

\noindent by $(g,y,x)\mapsto (g,y,g (x))$ for $g\in G$, $y\in Y$ and
$x\in X$. We will use the same notation for a
similarly defined homeomorphism
\begin{displaymath}
\tilde{u}_{H}^{G} (Y,X):G_{+}\smashove{H} (Y\wedge  i_{H}^{*}X) \to 
                        (G_{+}\smashove{H}Y)\wedge  X
\end{displaymath}

\noindent for a $G$-spectrum $X$ and $H$-spectrum $Y$.  We will abbreviate 
\begin{displaymath}
\tilde{u}_{H}^{G} (S^{0},X):G_{+}\smashove{H}  i_{H}^{*}X \to 
                        G/H_{+}\wedge  X
\end{displaymath}

\noindent by $\tilde{u}_{H}^{G} (X)$.

For {\rep}s $V$ and $V'$ of $G$ both restricting to $W$ on $H$, but
having distinct restrictions to all larger subgroups, we define
$\tilde{u}_{V-V'}=\tilde{u}_{H}^{G} (S^{V})\tilde{u}_{H}^{G} (S^{V'})^{-1}$,
so the following diagram of {\eqvr} homeomorphisms commutes:
\begin{numequation}\label{eq-u{V-V'}}
\begin{split}
\xymatrix
@R=2mm
@C=10mm
{   &   &G/H\wedge S^{V}\\
G_{+}\smashove{H}S^{W}\ar[rru]^(.5){\tilde{u}_{H}^{G} (S^{V})}
                      \ar[rrd]_(.5){\tilde{u}_{H}^{G} (S^{V'})}\\
    &   &G/H\wedge S^{V'}\ar[uu]_(.5){\tilde{u}_{V-V'}}.
}
\end{split}
\end{numequation}%

\noindent When $V'=|V|$ (meaning that $H=G_{V}$ acts trivially on
$W$), then we abbreviate $\tilde{u}_{V-V'}$ by $\tilde{u}_{V}$.
\end{defin}

If $V$ is a {\rep} of $H$ restricting to $W$ on $K$, we can smash the
diagram (\ref{eq-pinch-fold}) with $S^{V}$ and get
\begin{numequation}\label{eq-pinch-fold-VW}
\begin{split}
\xymatrix
@R=1mm
@C=5mm
{
    S^{V}\ar@<.5ex>[rr]^(.4){\mbox{pinch} }
        &{}
            & H/K_{+}\wedge S^{V}\ar@<.5ex>[ll]^(.6){\mbox{fold} }\\
       &\ar@{=>}[ddd]^(.6){G_{+}\smashove{H} (\cdot)}\\
       &\\
       &\\
       &\\
    G_{+}\smashove{H}S^{V}\ar@<.5ex>[rr]^(.4){\mbox{pinch} }
       &{} &G_{+}\smashove{H} (H/K_{+}\wedge S^{V} )\ar[r]^(.55){\approx }
                                     \ar@<.5ex>[ll]^(.6){\mbox{fold} } 
                &G_{+}\smashove{H} (H_{+}\smashove{K} S^{W} )\ar@{=}[r]^(.5){}
                    &G_{+}\smashove{K}S^{W},
}
\end{split}
\end{numequation}%

\noindent where the homeomorphism is induced by that of Definition
\ref{def-homeo}.

\begin{defin}\label{def-gp-action}
{\bf The group action restriction and transfer maps.} For subgroups
$K\subseteq H\subseteq G$,  let $V\in RO (H)$ be a virtual {\rep} of $H$ restricting to $W\in RO (K)$. The group action  transfer
and restriction maps
\begin{displaymath}
\xymatrix
@R=5mm
@C=15mm
{ 
{\uparrow_{H}^{G}\upi_{V}^{H} (i_{H}^{*}X) } \ar@{=}[r]^(.5){}
    &{\upi_{H,V}X }\ar@<-.5ex>[r]_(.5){\underline{r}_{K}^{H} }
        &{\upi_{K,W}X }\ar@<-.5ex>[l]_(.5){\underline{t}_{K}^{H,V} }
                       \ar@{=}[r]^(.5){}
            &{\uparrow_{K}^{G}\upi_{W}^{K} (i_{K}^{*}X) } 
}
\end{displaymath}

\noindent (see \ref{def-MF-induction}) are the ones induced by the
composite maps in the bottom row of (\ref{eq-pinch-fold-VW}). The
symbols $t$ and $r$ here are underlined because they are maps {\em of}
Mackey functors rather than maps within Mackey functors.
\end{defin}

We include $V$ as an index for the group action transfer
$\underline{t}_{K}^{H,V}$ because its target is not determined by its
source.

Thus we have abelian groups $\upi_{{H',V}} (X) (G/H'')$ for all
subgroups $H', H''\subseteq G$ and {\rep}s $V$ of $H'$.  Most of them
are redundant in view of Theorem \ref{thm-module} below. In what
follows, we will use the notation $H_{\cup}=H'\cup H''$ and
$H_{\cap}=H'\cap H''$. 

\begin{lem}\label{lem-module}
{\bf An {\eqvr} module structure.} For a $G$-spectrum $X$ and
$H'$-spectrum $Y$,
\begin{displaymath}
[G_{+}\smashove{H'}Y, X]^{H''}
    =\Z[G/H_{\cup }]\otimes [H_{\cup +}\smashove{H'}Y,X]^{H''}
\end{displaymath}

\noindent as $\Z[G/H'']$-modules.
\end{lem}

\proof
As abelian groups,
\begin{align*}
[G_{+}\smashove{H'}Y, X]^{H''}
 & = [i_{H''}^{*} (G_{+}\smashove{H'}Y), X]^{H''}  \\
 & = \left[\bigvee_{|G/H_{\cup }|}H_{\cup +}\smashove{H'}Y, 
            X \right]^{H''}\\
 & = \bigoplus_{|G/H_{\cup }|}[H_{\cup +}\smashove{H'}Y, X]^{H''}
\end{align*}

\noindent and $G/H''$ permutes the wedge summands of
$\bigvee_{|G/H_{\cup }|}H_{\cup +}\smashove{H'}Y$ as it
permutes the elements of $G/H_{\cup }$.  \qed\bigskip

\begin{thm}\label{thm-module}
{\bf The module structure for $\underline{RO (G)}$-graded
homotopy\linebreak groups.} For subgroups $H', H''\subseteq G$ with
$H_{\cup}=H'\cup H''$ and $H_{\cap}=H'\cap H''$, and a virtual {\rep}
$V$ of $H'$ restricting to $W$ on $H_{\cap}$,
\begin{displaymath}
\upi_{H',V}X (G/H'')
 \cong  \Z[G/H_{\cup }]\otimes \upi_{H_{\cap},W}X (G/G)  
 \cong  \Z[G/H_{\cup }]\otimes \upi_{W}^{H_{\cap}}
                       i_{H_{\cap}}^{*}X (H_{\cap}/H_{\cap})
\end{displaymath}

\noindent as $\Z[G/H'']$-modules.

Suppose that $H''$ is a proper  subgroup of $H'$ and $\gamma \in H'$
is a generator.  Then as an element in $\Z[G/H'']$, $\gamma $ induces
multiplication by $-1$ in $\upi_{H',V}X (G/H'')$ iff $V$ is
nonorientable.
\end{thm}

\proof We start with the definition and use the homeomorphism of
Definition \ref{def-homeo} and the module structure of Lemma \ref{lem-module}.
\begin{align*}
\upi_{H',V}X (G/H'')
 & =  [(G_{+}\smashove{H'} S^{V})\wedge G/H''_{+},X]^{G} \\
 & =  [G_{+}\smashove{H''}(G_{+}\smashove{H'}S^{V}),X]^{G} \\ 
 & =  [G_{+}\smashove{H'}S^{V},X]^{H''} 
   =  \Z[G/H_{\cup }] \otimes [H_{\cup +}\smashove{H'}S^{V}, X]^{H''}\\
\aand 
[H_{\cup +}\smashove{H'}S^{V}, X]^{H''}
 & =  [S^{W}, X]^{H_{\cap}}
   =  [G_{+}\smashove{H_{\cap}}S^{W}, X]^{G}\\
 & =  \upi_{W}^{H_{\cap}} (i_{H_{\cap}}^{*}X) (H_{\cap}/H_{\cap})
   =  \upi_{H_{\cap},W}X (G/G).
\end{align*}

For the statement about nonoriented $V$, we have 
\begin{displaymath}
\upi_{H',V}X (G/H'') = \Z[G/H']\otimes \upi_{W}^{H''}i_{H''}^{*}X (H''/H'')
= \Z[G/H']\otimes  [S^{W},X]^{H''}.
\end{displaymath}

\noindent 
Then $\gamma $ induces a map of degree $\pm 1$ on the sphere depending
on the orientability of $V$.  \qed
\bigskip

Theorem \ref{thm-module} means that we need only consider the groups 
\begin{displaymath}
\upi_{H,V}X (G/G) \cong  \upi_{V}^{H}i^{*}_{H}X (H/H).
\end{displaymath}

When $H\subset G$ and $V$ is a virtual {\rep} of $G$, we have
\begin{numequation}\label{eq-easy-iso}
\begin{split}
\upi_{V}X (G/H)
 \cong  \upi_{H,i^{*}_{H}V}X (G/G)
 \cong  \upi^{H}_{i^{*}_{H}V}i^{*}_{H}X (H/H).
\end{split}
\end{numequation}%

\noindent This isomorphism makes the following diagram commute for
$K\subseteq H$.
\begin{displaymath}
\xymatrix
@R=5mm
@C=10mm
{
{\upi_{V}X (G/H)}\ar[r]^(.45){\cong }
              \ar@<-.5ex>[d]_(.5){\Res_{K}^{H} }
    &{\upi_{H,i^{*}_{H}V}X (G/G)}\ar[r]^(.5){\cong }
                  \ar@<-.5ex>[d]_(.5){\underline{r}_{K}^{H} }
        &{\upi^{H}_{i^{*}_{H}V}i^{*}_{H}X (H/H)}
                    \\
{\upi_{V}X (G/K)}\ar[r]^(.45){\cong }
              \ar@<-.5ex>[u]_(.5){\Tr_{K}^{H} }
    &{\upi_{K,i_{K}^{*}V}X (G/G)}\ar[r]^(.5){\cong }
                  \ar@<-.5ex>[u]_(.5){\underline{t}_{K}^{H,i^{*}_{H}V} }
        &{\upi^{K}_{i_{K}^{*}V}i^{*}_{K}X (K/K)}
}
\end{displaymath}

\noindent {\em We will use the three groups of (\ref{eq-easy-iso})
interchangeably as convenient and use the same notation for elements
in each related by this canonical isomorphism.} Note that the group on
the left is indexed by $RO (G)$ while the two on the right are indexed
by $RO (H)$.  This means that if $V$ and $V'$ are {\rep}s of $G$ each
restricting to $W$ on $H$, then $\upi_{V}X (G/H)$ and $\upi_{V'}X
(G/H)$ are canonically isomorphic.  The first of these is
\begin{displaymath}
[G/H_{+}\wedge S^{V}, X]^{G}\cong 
[G_{+}\smashove{H}S^{W}, X]^{G} \cong 
[S^{W}, i_{H}^{*}X]^{H}
\end{displaymath}

\noindent where the first isomorphism is induced by the homeomorphism
$\tilde{u}_{H}^{G} (X)$ of Definition \ref{def-homeo} and the second
is the fact that $G_{+}\smashove{H} (\cdot )$ is the left adjoint of
the forgetful functor $i_{H}^{*}$.

\begin{rem}\label{rem-via}
{\bf Factorization via restriction}.  For a ring spectrum $X$, such as
the one we are studying in this paper, an indecomposable element in
$\upi_{\star}X (G/H)$ may map to a product $xy\in \upi_{H,\star}X
(G/G)$ of elements in groups indexed by {\rep}s of $H$ that are not
restrictions of {\rep}s of $G$.  When this happens we may denote the
indecomposable element in $\upi_{\star}X (G/H)$ by $[xy]$.  This
factorization can make some computations easier.
\end{rem}

\section{The $RO (G)$-graded homotopy of $\HZ$}\label{sec-HZZ}

We describe part of the $RO (G)$-graded Green functor
$\upi_{\star}(\HZ)$, where $\HZ $ is the integer {\SESM} spectrum
$\HZ$ in the $G$-{\eqvr} category, for some finite cyclic 2-group $G$.  For
each actual (as opposed to virtual) $G$-{\rep} $V$ we have an {\eqvr}
reduced cellular chain complex $C^{V}_{*}$ for the space $S^{V}$.  It
is a complex of $\ints [G]$-modules with $H_{*} (C^{V})=H_{*}
(S^{|V|})$.

One can convert such a chain complex $C_{*}^{V}$ of $\Z[G]$-modules to
one of Mackey functors as follows.  Given a $\Z[G]$-module $M$, we get
a Mackey functor $\uM$ defined by
\begin{numequation}
\label{eq-fpm}
\uM (G/H)= M^{H}\qquad \mbox{for each subgroup $H\subseteq G$.} 
\end{numequation}%

\noindent We call this a {\em fixed point Mackey functor}. In it each
restriction map $\Res_{K}^{H}$ (for $K\subseteq H\subseteq G$) is one
to one.  When $M$ is a permutation module, meaning the free abelian
group on a $G$-set $B$, we call $\uM$ a {\em permutation Mackey
functor} \cite[\S3.2]{HHR}.

In particular the $Z[G]$-module $\Z$ with trivial group action (the
free abelian group on the $G$-set $G/G$) leads to a Mackey functor
$\underline{\Z}$ in which each restriction map is an isomorphism and
the transfer map $\Tr_{K}^{H}$ is multiplication by $|H/K|$.  For each
Mackey functor $\underline{M}$ there is an {\SESM} spectrum
$H\underline{M}$ \cite[\S5]{Greenlees-May}, and $\HZZ$ is the same as
$\HZ$ with trivial group action.

Given a finite $G$-CW spectrum $X$, meaning one
built out of cells of the form $G_{+}\smashove{H} e^{n}$, we get a
reduced cellular chain complex of $\Z[G]$-modules $C_{*}X$, leading to
a chain complex of fixed point Mackey functors $\uC_{*}X$.
Its homology is a graded Mackey functor $\uH_{*}X$ with
\begin{displaymath}
\uH_{*}X (G/H)
 = \upi_{*} (X\wedge \HZZ) (G/H)
 = \pi_{*} (X\wedge \HZZ)^{H}.
\end{displaymath}

\noindent In particular $\uH_{*}X (G/\ee) = H_{*}X$, the
underlying homology of $X$.  In general $\uH_{*}X (G/H)$
is not the same as $H_{*} (X^{H})$ because fixed points do not commute
with smash products.  

For a finite cyclic 2-group $G=C_{2^{k}}$, the irreducible {\rep}s are
the 2-dimensional ones $\lambda (m)$ corresponding to rotation through
an angle of $2\pi m/2^{k}$ for $0<m<2^{k-1}$, the sign {\rep} $\sigma
$ and the trivial one of degree one, which we denote by 1.  The
2-local {\eqvr} homotopy type of $S^{\lambda (m)}$ depends only on the
2-adic valuation of $m$, so we will only consider $\lambda (2^{j})$
for $0\leq j\leq k-2$ and denote it by $\lambda_{j}$.  The planar
rotation $\lambda_{k-1}$ though angle $\pi$ is the same {\rep} as
$2\sigma $. {\em We will denote $\lambda (1)=\lambda_{0}$ simply by
$\lambda $.}

We will describe the chain complex $C^{V}$ for
 \begin{displaymath}
V=a+b\sigma +\sum_{2\leq j\leq k}c_{j}\lambda_{k-j}.
\end{displaymath}

\noindent  for nonnegative integers $a$, $b$ and $c_{j}$.
The isotropy group of $V$ (the largest subgroup fixing all of $V$) is 
\begin{displaymath}
G_{V}=\mycases{
C_{2^{k}}=G 
      &\mbox{for }b=c_{2}=\dotsb =c_{k}=0\\
C_{2^{k-1}}=:G' 
      &\mbox{for $b>0$ and  $c_{2}=\dotsb =c_{k}=0$}\\
C_{2^{k-\ell }}
      &\mbox{for  $c_{\ell }>0$  and $c_{1+\ell }=\dotsb =c_{k}=0$}
}
\end{displaymath}

The sphere $S^{V}$ has a $G$-CW structure with reduced cellular chain complex 
$C^{V}$ of the form
\begin{numequation}\label{eq-CVn}
\begin{split}
C^{V}_{n}=\mycases{
\ints 
       &\mbox{for }n=d_{0}\\
\ints[G/G'] 
       &\mbox{for }d_{0}<n\leq d_{1}\\
\ints[G/C_{2^{k-j}}] 
       &\mbox{for $d_{j-1}<n\leq d_{j}$ and $2\leq j\leq \ell $}\\
0      &\mbox{otherwise.}
}
\end{split}
\end{numequation}%

\noindent where 
\begin{displaymath}
d_{j}=\mycases{
a       &\mbox{for $j=0$}\\
a+b
       &\mbox{for $j=1$}\\
a+b+2c_{2}+\dotsb +2c_{j}
       &\mbox{for $2\leq j\leq \ell $,}
}
\end{displaymath}

\noindent so $d_{\ell }=|V|$. 

The boundary map $\partial_{n}:C_{n}^{V}\to C_{n-1}^{V}$ is determined
by the fact that\linebreak ${H_{*}(C^{V})=H_{*}(S^{|V|})}$.  More
explicitly, let $\gamma $ be a generator of $G$ and
\begin{displaymath}
\zeta_{j}=\sum_{0\leq t<2^{j }}\gamma^{t}
\qquad \mbox{for }1\leq j \leq k.
\end{displaymath}

\noindent   Then we have 
\begin{displaymath}
\partial_{n}=\mycases{
\nabla
       &\mbox{for }n=1+d_{0}\\
(1-\gamma )x_{n}
       &\mbox{for $n-d_{0}$ even and $2+d_{0}\leq n\leq d_{n}$}\\
x_{n}
       &\mbox{for $n-d_{0}$ odd and $2+d_{0}\leq n\leq d_{n}$}\\
0      &\mbox{otherwise,}
}
\end{displaymath}

\noindent where $\nabla$ is the fold map sending $\gamma \mapsto 1$,
and $x_{n}$ denotes multiplication by an element in $\Z[G]$ to be
named below. We will use the same symbol below for the quotient map
$\Z[G/H]\to\Z[G/K]$ for $H\subseteq K\subseteq G$.  The elements
$x_{n}\in \Z[G]$ for $2+d_{0}\leq n\leq |V|$ are determined
recursively by $x_{2+d_{0}}=1$ and
\begin{displaymath}
x_{n}x_{n-1}= \zeta_{j } \qquad \mbox{for }2+d_{j-1}<n\leq 2+d_{j}. 
\end{displaymath}

\noindent It follows that $H_{|V|}C^{V}=\Z$ generated by either
$x_{1+|V|}$ or its product with $1-\gamma $, depending on the parity
of $b$.

This complex is 
\begin{displaymath}
C^{V} = \Sigma^{|V_{0}|} C^{V/V_{0}}
\end{displaymath}

\noindent where $V_{0}=V^{G}$.  This means we can assume without loss
of generality that $V_{0}=0$.

An element 
\begin{displaymath}
x\in H_{n} \uC^{V} (G/H)
 = \uH_{n} S^{V} (G/H)
\end{displaymath}

\noindent
corresponds to an element $x\in \upi_{n-V} \HZ (G/H) $.

We will denote the dual complex $\Hom_{\ints } (C^{V},\ints )$ by
$C^{-V}$.  Its chains lie in dimensions $-n$ for $0\leq n\leq |V|$.
An element $x\in \uH_{-n} (S^{-V}) (G/H)$ corresponds to an
element $x\in \upi_{V-n} \HZ (G/H)$.

The method we have just described determines only a portion of the\linebreak 
$RO(G)$-graded Mackey functor $\upi_{(G, \star)}\HZZ$, namely
the groups in which the index differs by an integer from an actual
{\rep} $V$ or its negative.  For example it does not give us
$\upi_{\sigma -\lambda }\HZZ$ for $|G|\geq 4$.

We leave the proof of the following as an exercise for the reader.

\begin{prop}\label{prop-top}
{\bf The top (bottom) homology groups for $S^{V}$ ($S^{-V}$).} Let $G$
be a finite cyclic 2-group and $V$ a nontrivial {\rep} of $G$ of
degree $d$ with $V^{G}=0$ and isotropy group $G_{V}$.  Then
$C_{d}^{V}=C_{-d}^{-V}=\Z[G/G_{V}]$ and
\begin{enumerate}
\item [(i)] If $V$ is oriented then
$\uH_{d}S^{V}=\underline{\Z}$, the constant $\Z$-valued
Mackey functor in which each restriction map is an isomorphism and
each transfer $\Tr_{H}^{K}$ is multiplication by $|K/H|$.

\item [(ii)] 
$\uH_{-d}S^{-V}=\underline{\Z} (G,G_{V})$, the constant
$\Z$-valued Mackey functor in which
\begin{displaymath}
\Res_{H}^{K}=\mycases{
1       &\mbox{for }K\subseteq G_{V}\\
{|K/H|}   &\mbox{for }G_{V}\subseteq H
}
\end{displaymath}

\noindent and 
\begin{displaymath}
\Tr_{H}^{K}=\mycases{
{|K/H| }  &\mbox{for }K\subseteq G_{V}\\
1       &\mbox{for }G_{V}\subseteq H.
}
\end{displaymath}

\noindent (The above completely describes the cases where
$|K/H|=2$, and they determine all other restrictions and transfers.)  The
functor $\underline{\Z}(G,e)$ is also known as the dual
$\underline{\Z}^{*}$.  These isomorphisms are induced by the maps
\begin{displaymath}
\xymatrix
@R=0mm
@C=15mm
{
{\uH_{d}S^{V}}\ar@{=}[dd]^(.5){}
    &   &{\uH_{-d}S^{-V}}\ar@{=}[dd]^(.5){}\\
    &    &\\
{\underline{\Z}}\ar[r]^(.4){\Delta}
    &{\underline{\Z[G/G_{V}]}}\ar[r]^(.5){\nabla}
        &{\underline{\Z} (G,G_{V}).}
}
\end{displaymath}


\item [(iii)] If $V$ is not oriented then
$\uH_{d}S^{V}=\underline{\Z}_{-}$, where
\begin{displaymath}
\underline{\Z}_{-} (G/H)=\mycases{
0   &\mbox{for }H=G\\
\Z_{-}:=\Z[G]/ (1+\gamma )
    &\mbox{otherwise}
}
\end{displaymath}

\noindent where each restriction map $\Res_{H}^{K}$ is an isomorphism
and each transfer $\Tr_{H}^{K}$ is multiplication by $|K/H|$ for each
proper subgroup $K$.  

\item [(iv)] We also have
$\uH_{-d}S^{-V}=\underline{\Z} (G,G_{V})_{-}$, where
\begin{displaymath}
\underline{\Z} (G,G_{V})_{-} (G/H)=\mycases{
0       &\mbox{for $H=G$ and $V=\sigma $}\\
\Z/2    &\mbox{for $H=G$ and $V\neq \sigma $}\\
\Z_{-}  &\mbox{otherwise}
}
\end{displaymath}

\noindent with the same restrictions and transfers as $\underline{\Z}
(G,G_{V})$. These isomorphisms are induced by the maps 
\begin{displaymath}
\xymatrix
@R=0mm
@C=15mm
{
{\uH_{d}S^{V}}\ar@{=}[dd]^(.5){}
   &    &{\uH_{-d}S^{-V}}\ar@{=}[dd]^(.5){}\\
   &    &\\
{\underline{\Z}_{-}}\ar[r]^(.4){\Delta_{-} }
   &{\underline{\Z[G/G_{V}]}}\ar[r]^(.5){\nabla_{-}}
        &{\underline{\Z} (G,G_{V})_{-}.}
}
\end{displaymath}


\end{enumerate}
\end{prop}

The Mackey functor $\underline{\Z(G,G_{V})}$ is one of those defined
(with different notation) in \cite[Def. 2.1]{HHR:RO(G)}.

\begin{defin}\label{def-aeu}
{\bf Three elements in $\upi_{\star}^{G}(\HZ)$.}  Let $V$
be an actual (as opposed to virtual) {\rep} of the finite cyclic
2-group $G$ with $V^{G}=0$ and isotropy group $G_{V}$.
\begin{enumerate}
\item [(i)] The {\eqvr} inclusion ${S^{0}\to S^{V}}$ defines an
element in $\upi_{-V}S^{0} (G/G)$ via the isomorphisms
\begin{displaymath}
\upi_{-V}S^{0} (G/G)=
\upi_{0}S^{V} (G/G)= \pi_{0}S^{V^{G}}=\pi_{0}S^{0}=\Z,
\end{displaymath}

\noindent and we will use the symbol $a_{V}$ to denote its image in
$\upi_{-V}\HZZ (G/G)$.

\item [(ii)] The underlying equivalence $S^{V}\to S^{|V|}$ defines an
element in
\begin{displaymath}
\upi_{V}S^{|V|} (G/G_{V}) = \upi_{V-|V|}S^{0} (G/G_{V})
\end{displaymath}

\noindent and we will use the symbol $e_{V}$ to denote its image in
$\upi_{V-|V|}\HZZ (G/G_{V})$.

\item [(iii)] If $W$ is an oriented {\rep} of $G$ (we do not require
that $W^{G}=0$), there is a map
\begin{displaymath}
\Delta :\ints  \to C^{W}_{|W|}= \Z[G/G_{W}]
\end{displaymath}

\noindent as in Proposition \ref{prop-top} giving an element 
\begin{displaymath}
u_{W}\in \uH_{|W|}S^{W} (G/G) = \upi_{|W|-W} \HZ(G/G).
\end{displaymath}

For nonoriented $W$, Proposition \ref{prop-top} gives a map
\begin{displaymath}
\Delta_{-} :\ints_{-}  \to C^{W}_{|W|}
\end{displaymath}

\noindent and an element 
\begin{displaymath}
u_{W}\in \uH_{|W|}S^{W} (G/G') = \upi_{|W|-W} \HZ(G/G').
\end{displaymath}

\noindent   
\end{enumerate}
\end{defin}

The element $u_{W}$ above is related to the element $\tilde{u}_{V}$ of
(\ref{eq-u{V-V'}}) as follows. 

\begin{lem}\label{lem-uV}
{\bf The restriction of $u_{W}$ to a unit and permanent cycle.}  Let $W$ be a
nontrivial {\rep} of $G$ with $H=G_{W}$.  Then the homeomorphism
\begin{displaymath}
\Sigma^{-W}\tilde{u}_{W}:G/H_{+}\wedge S^{|W|-W}\to G/H_{+}
\end{displaymath}

\noindent of
(\ref{eq-u{V-V'}}) induces an isomorphism $\upi_{0}\HZZ (G/H)\to
\upi_{|W|-W}\HZZ (G/H)$ sending the unit to $\Res_{H}^{K}(u_{W})$
for $u_{W}$ as defined in (iii) above and $K=G$ or $G'$ depending on
the orientability of $W$.

The product
\begin{displaymath}
\Res_{H}^{K}(u_{W})e_{W}\in \upi_{0}\HZZ (G/H)=\Z
\end{displaymath}

\noindent is a generator, so $e_{W}$ and $\Res_{H}^{K}(u_{W})$ are
units in the ring $\upi_{\star}\HZZ (G/H)$, and $\Res_{H}^{K}(u_{W})$
is in the Hurewicz image of $\upi_{\star}S^{0} (G/H)$.

\end{lem}

\proof
The diagram 
\begin{displaymath}
\xymatrix
@R=5mm
@C=10mm
{G/K_{+}\wedge S^{|W|-W}
    &G/H_{+}\wedge S^{|W|-W}\ar[r]^(.6){\tilde{u}_{W}}
                            \ar[l]_(.5){\mbox{fold} }
        &G/H_{+}
}
\end{displaymath}

\noindent induces (via the functor $[\cdot , \HZZ]^{G}$)
\begin{displaymath}
\xymatrix
@R=5mm
@C=10mm
{{\upi_{|W|-W}\HZZ (G/K)}\ar[r]^(.5){\Res_{H}^{K}}\ar@{=}[d]^(.5){}
    &{\upi_{|W|-W}\HZZ (G/H)}\ar@{=}[d]^(.5){}
        &{\upi_{0}\HZZ (G/H)}\ar[l]_(.4){\cong }\ar@{=}[d]^(.5){} \\
{\underline{H}_{|W|}S^{W} (G/K)}           
    &{\underline{H}_{|W|}S^{W} (G/H)}           
        &\Z
}
\end{displaymath}

\noindent The restriction map is an isomorphism by Proposition
\ref{prop-top} and the group on the left is generated by $u_{W}$.

The product is the composite of $H$-maps
\begin{displaymath}
\xymatrix
@R=5mm
@C=10mm
{
{S^{W}}\ar[r]^(.5){e_{W}}
    &{S^{|W|}}\ar[r]^(.5){\Res_{H}^{K}(u_{W})}
        &{\Sigma^{W}\HZZ,}
}
\end{displaymath}

\noindent which is the standard inclusion.
\qed\bigskip

Note that $a_{V}$ and $e_{V}$ are induced by maps to {\eqvr} spheres
while $u_{W}$ is not.  This means that in any {\SS} based on a
filtration where the subquotients are {\eqvr} $\HZ $-modules, elements
defined in terms of $a_{V}$ and $e_{V}$ will be permanent cycles,
while multiples and powers of $u_{W}$ can support nontrivial
differentials.  Lemma \ref{lem-uV} says a certain restriction of
$u_{W}$ is a permanent cycle.

Each nonoriented $V$ has the form $W+\sigma $ where $\sigma $ is the
sign {\rep} and $W$ is oriented.  It follows that
\begin{displaymath}
u_{V}=u_{\sigma}\Res_{G'}^{G}(u_{W})\in \upi_{|V|-V}H\Z (G/G').
\end{displaymath}







Note also that $a_{0}=e_{0}=u_{0}=1$.  The trivial
representations contribute nothing to $\pi_{\star}(H\Z)$.  We can limit
our attention to {\rep}s $V$ with $V^{G}=0$.  Among such {\rep}s
of cyclic 2-groups, the oriented ones are precisely the
ones of even degree.

\begin{lem}\label{lem-aeu}
{\bf Properties of $a_{V}$, $e_{V}$ and $u_{W}$.} The elements
$a_{V}\in \upi_{-V}\HZ(G/G)$, $e_{V}\in \underline{\pi
}_{V-|V|}\HZ(G/G_{V})$ and $u_{W}\in \upi_{|W|-W}\HZ(G/G)$ for
$W$ oriented of Definition \ref{def-aeu} satisfy the following.

\begin{enumerate}

\item \label{aeu:i} 
$a_{V+W}=a_{V}a_{W}$ and $u_{V+W}=u_{V}u_{W}$.

\item \label{aeu:ii} 
$|G/G_{V}|a_{V}=0$ where $G_{V}$ is the isotropy group of $V$.

\item \label{aeu:iii} For oriented $V$, $\Tr_{G_{V}}^{G}(e_{V})$
and $\Tr_{G_{V}}^{G'}(e_{V+\sigma})$ have infinite order, while\linebreak 
$\Tr_{G_{V}}^{G}(e_{V+\sigma})$ has order 2 if
$|V|>0$ and $\Tr_{G_{V}}^{G}(e_{\sigma})=\Tr_{G'}^{G}(e_{\sigma})=0$. 

\item \label{aeu:iv} For oriented $V$ and $G_{V}\subseteq H\subseteq G$
%
\begin{align*}
\Tr_{G_{V}}^{G}(e_{V})u_{V}
 & = |G/G_{V}|\in \upi_{0}\HZZ (G/G)=\Z  \\
\aand 
\Tr_{G_{V}}^{G'}(e_{V+\sigma })u_{V+\sigma }
 & =  |G'/G_{V}|\in \upi_{0}\HZZ (G/G')=\Z
 \qquad \mbox{for }|V|>0. 
\end{align*}

\item \label{aeu:v}
$a_{V+W}\Tr_{G_{V}}^{G}(e_{V+U})=0$ if $|V|>0$.

\item \label{aeu:vi}
For $V$ and $W$ oriented, $u_{W}\Tr_{G_{V}}^{G}(e_{V+W})
=|G_{V}/G_{V+W}|\Tr_{G_{V}}^{G}(e_{V})$.

\item \label{aeu:vii}
{\bf The gold (or $au$) relation. } For $V$ and $W$ oriented {\rep}s
of degree 2 with $G_{V}\subseteq G_{W}$, ${a_{W}u_{V} =
|G_{W}/G_{V}|a_{V}u_{W}}$.

\end{enumerate}

For nonoriented $W$ similar statements hold in $\upi_{\star}\HZ
(G/G')$. $2W$ is oriented and $u_{2W}$ is defined in
$\upi_{2|W|-2W}\HZZ (G/G)$ with ${\Res_{G'}^{G}(u_{2W})=u_{W}^{2}}$.

\end{lem}

\proof (\ref{aeu:i}) 
This follows from the existence of the pairing
$C^{V}\otimes C^{W}\to C^{V+W}$.  It induces an isomorphism in $H_{0}$
and (when both $V$ and $W$ are oriented) in $H_{|V+W|}$.

(\ref{aeu:ii})  This holds because $H_{0} (V)$ is killed by $|G/G_{V}|$.

(\ref{aeu:iii}) This follows from Proposition \ref{prop-top}.

(\ref{aeu:iv}) 
Using the Frobenius relation we have 
%
%
\begin{align*}
\Tr_{G_{V}}^{G}(e_{V})u_{V}
 & = \Tr_{G_{V}}^{G}(e_{V}\Res_{G_{V}}^{G}(u_{V}))
   = \Tr_{G_{V}}^{G}(1) \qquad \mbox{by \ref{lem-uV} } \\
 &  = |G/G_{V}|  \\
\Tr_{G_{V}}^{G'}(e_{V+\sigma })u_{V+\sigma }
 & = \Tr_{G_{V}}^{G'}(e_{V+\sigma }\Res_{G_{V}}^{G'}(u_{V+\sigma }))
   = \Tr_{G_{V}}^{G'}(1) = |G'/G_{V}|.
\end{align*}

(\ref{aeu:v}) We have 
\begin{displaymath}
a_{V+W}\Tr_{G_{V}}^{G}(e_{V+U}):S^{-|V|-|U|} \to S^{W-U}.
\end{displaymath}

\noindent It is null because the bottom cell of ${S^{W-U}}$ is in
dimension ${-|U|}$.

(\ref{aeu:vi}) Since $V$ is oriented, then we are computing in a
torsion free group so we can tensor with the rationals.  It follows
from (\ref{aeu:iv}) that
\begin{align*}
\Tr_{G_{V+W}}^{G}(e_{V+W})
 & =   \frac{|G/G_{V+W}|}{u_{V}u_{W}}\\                        
\aand 
\Tr_{G_{V}}^{G}(e_{V})
   & =   \frac{|G/G_{V}|}{u_{V}}\\                        
\sso
u_{W}\Tr_{G_{V+W}}^{G}(e_{V+W} )
   & =   \frac{|G/G_{V+W}|}{u_{V}} 
    =   |G_{V}/G_{V+W}|\Tr_{G_{V}}^{G}(e_{V}).
\end{align*}

(\ref{aeu:vii}) 
For $G=C_{2^{n}}$, each oriented {\rep} of degree 2 is $2$-locally
equivalent to a $\lambda_{j}$ for $0\leq j<n$.  The isotropy group is
$G_{\lambda_{j}} = C_{2^{j}}$.  Hence the assumption that
$G_{V}\subset G_{W}$ is can be replaced with $V=\lambda_{j}$ and
$W=\lambda_{k}$ with $0\leq j<k<n$.  the statment we wish to prove is 
\begin{displaymath}
a_{\lambda_{k}}u_{\lambda_{j}}=2^{k-j}a_{\lambda_{j}}u_{\lambda_{k}}.
\end{displaymath}

One has a map $S^{\lambda_{j}}\to S^{\lambda_{k}}$ which is the
suspension of the $2^{k-j}$th power map on the equatorial circle.
Hence its underlying degree is $2^{k-j}$.  We will denote it by
$a_{\lambda_{k}}/a_{\lambda_{j}}$ since there is a diagram
\begin{displaymath}
\xymatrix
@R=3mm
@C=20mm
{
  &S^{\lambda_{j}}\ar[dd]^(.5){a_{\lambda_{k}}/a_{\lambda_{j}}}\\
S^{0}\ar[ru]^(.5){a_{\lambda_{j}}}\ar[rd]_(.5){a_{\lambda_{k}}}\\
  &S^{\lambda_{k}}.
}
\end{displaymath}

We claim there is a similar diagram 
\begin{numequation}\label{eq-au-claim}
\begin{split}
\xymatrix
@R=3mm
@C=20mm
{
  &S^{\lambda_{k}}\wedge \HZZ\ar[dd]^(.5){u_{\lambda_{j}}/u_{\lambda_{k}}}\\
S^{2}\ar[ru]^(.5){u_{\lambda_{k}}}\ar[rd]_(.5){u_{\lambda_{j}}}\\
  &S^{\lambda_{j}}\wedge \HZZ.
}
\end{split}
\end{numequation}%

\noindent in which the underlying degree of the vertical map is one.

Smashing $a_{\lambda_{k}}/a_{\lambda_{j}}$ with $\HZZ$ and composing
with $u_{\lambda_{j}}/u_{\lambda_{k}}$ gives a factorization of the
degree $2^{k-j}$ map on $S^{\lambda_{j}}\wedge \HZZ$.  Thus we have 
\begin{align*}
\frac{u_{\lambda_{j}}}{u_{\lambda_{k}}}
\frac{a_{\lambda_{k}}}{a_{\lambda_{j}}}
 & = 2^{k-j}  \\
u_{\lambda_{j}}a_{\lambda_{k}}
 & = 2^{k-j} u_{\lambda_{k}}a_{\lambda_{j}} 
\end{align*}

\noindent as desired.

The vertical map in (\ref{eq-au-claim}) would follow from a map 
\begin{displaymath}
S^{\lambda_{k}-\lambda_{j}}\to \HZZ
\end{displaymath}

\noindent with underlying degree one. Let $G=C_{2^{n}}$ and $G \supset
H=C_{2^{j}}$.  Then $S^{-\lambda_{j}}$ has a cellular structure of the
form
\begin{displaymath}
G/H_{+}\wedge S^{-2} \cup G/H_{+} \wedge e^{-1} \cup  e^{0}.
\end{displaymath}

\noindent We need to smash this with $S^{\lambda_{k}}$.  Since
$\lambda_{k}$ restricts trivially to $H$, 
\begin{displaymath}
G/H_{+}\wedge
S^{\lambda_{k}}=G/H_{+}\wedge S^{2}.
\end{displaymath}

\noindent  This means 
\begin{displaymath}
S^{\lambda_{k}-\lambda_{j}}
 = S^{\lambda _{k}}\wedge S^{-\lambda_{j}} 
 = G/H_{+}\wedge S^{0} \cup
     G/H_{+} \wedge e^{1} \cup e^{0}\wedge S^{\lambda_{k}}.
\end{displaymath}

\noindent  Thus its cellular chain complex has the form
\begin{displaymath}
\xymatrix
@R=5mm
@C=15mm
{
2   &\Z[G/K] \ar[d]^(.5){1-\gamma }\ar[drr]^(.5){\Delta }\\
1   &\Z[G/K] \ar[d]^(.5){\nabla }\ar[drr]^(.5){-\Delta }
        &   &\Z[G/H]\ar[d]^(.5){1-\gamma }\\
0   &\Z &   &\Z[G/H]
}
\end{displaymath}

\noindent where $K=G/C_{p^{k}}$ and the left column is the chain
complex for $S^{\lambda_{k}}$.  

There is a corresponding chain complex of fixed point Mackey functors.
Its value on the $G$-set $G/L$ for an arbitrary subgroup $L$ is
\begin{displaymath}
\xymatrix
@R=5mm
@C=15mm
{
2   &\Z[G/\max(K,L)] \ar[d]^(.5){1-\gamma }\ar[drr]^(.5){\Delta }\\
1   &\Z[G/\max(K,L)] \ar[d]^(.5){\nabla }\ar[drr]^(.5){-\Delta }
        &   &\Z[G/\max(H,L)]\ar[d]^(.5){1-\gamma }\\
0   &\Z &   &\Z[G/\max(H,L)]
}
\end{displaymath}

\noindent For each $L$ the map $\Delta $ is injective and maps the
kernel of the first $1-\gamma $ isomorphically to the kernel of the
second one.  This means we can replace the above by a diagram of the
form
\begin{displaymath}
\xymatrix
@R=5mm
@C=15mm
{
1   &\coker (1-\gamma ) \ar[d]^(.5){\nabla }\ar[drr]^(.5){-\Delta }\\
0   &\Z &   &\coker (1-\gamma )
}
\end{displaymath}

\noindent where each cokernel is isomorphic to $\Z$ and each
map is injective.

This means that $\underline{H}_{*}S^{\lambda_{k}-\lambda_{j}}$ is
concentrated in degree 0 where it is the pushout of the diagram above,
meaning a Mackey functor whose value on each subgroup is $\Z$.  Any
such Mackey functor admits a map to $\underline{\Z}$ with underlying
degree one.  This proves the claim of (\ref{eq-au-claim}). \qed\bigskip

The $\Z$-valued Mackey functor
$\underline{H}_{0}S^{\lambda_{k}-\lambda_{j}}$ is discussed in more
detail in \cite{HHR:RO(G)}, where it is denoted by $\underline{\Z} (k,j)$.

\section{Generalities on differentials and Mackey functor extensions}
\label{sec-gendiffs}

Before proceeding with a discussion about sepctral sequences, we need
the following.

\begin{rem}\label{rem-abuse}
{\bf Abusive {\SS } notation}.  When $d_{r} (x)$ is a nontrivial
element of order 2, the elements $2x$ and $x^{2}$ both survive to
$\underline{E}_{r+1}$, but in that group they are not the products
indicated by these symbols since $x$ itself is no longer present.
More geneally if $d_{r} (x)=y$ and $\alpha y=0$ for some $\alpha $,
then $\alpha x$ is present in $\underline{E}_{r+1}$.  This abuse of
notation is customary because it would be cumbersome to rename these
elements when passing from $\underline{E}_{r}$ to
$\underline{E}_{r+1}$.  We will sometimes denote them by $[2x]$,
$[x^{2}]$ and $[\alpha x]$ respectively to emphasize their
indecomposability.
\end{rem}

Now we make some observations about the relation between exotic
transfers and restriction with certain differentials in the slice
{\SS}.  By ``exotic'' we mean in a higher filtration.  In a {\SS} of
Mackey functors converging to $\upi_{\star}X$, it can happen that an
element $x\in \upi_{V}X (G/H)$ has filtration $s$, but its restriction
or transfer has a higher filtration.  {\em In the {\SS} charts in this
paper, exotic transfers and restrictions will be indicated by
{\color{blue}blue} and {\color{green} dashed green} lines respectively.}

\begin{lem}\label{lem-hate}
{\bf Restriction kills $a_{\sigma}$ and $a_{\sigma }$ kills
transfers.}  Let $G$ be a finite cyclic 2-group with sign {\rep}
$\sigma $ and index 2 subgroup $G'$, and let $X$ be a $G$-spectrum.
Then in $\upi_{*}X (G/G)$ the image of $\Tr_{G'}^{G}$ is the kernel of
multiplication by $a_{\sigma}$, and the kernel of $\Res_{G'}^{G}$ is
the image of multiplication by $a_{\sigma}$.

Suppose further that 4 divides the order of $G$ and let $\lambda $ be
the degree 2 representation sending a generator $\gamma \in G$ to a
rotation of order 4.  Then restriction kills $2a_{\lambda }$ and
$2a_{\lambda }$ kills transfers.
\end{lem}
 
\proof
Consider the cofiber sequence obtained by smashing $X$ with
\begin{numequation}\label{eq-hate}
\begin{split}
\xymatrix
@R=5mm
@C=10mm
{
S^{-1}\ar[r]^(.5){a_{\sigma}}
    &S^{\sigma -1}\ar[r]^(.5){}
        &G_{+}\smashove{G'}S^{0}\ar[r]^(.5){}
            &S^{0}\ar[r]^(.5){a_{\sigma}}
                &S^{\sigma}
}
\end{split}
\end{numequation}%

\noindent Since $(G_{+}\smashove{G'}X)^{G}$ is equivalent to $X^{G'}$,
passage to fixed point spectra gives
\begin{displaymath}
\xymatrix
@R=5mm
@C=10mm
{
\Sigma^{-1}X^{G}\ar[r]^(.5){}
    &\left(\Sigma^{\sigma -1}X \right)^{G}\ar[r]^(.5){}
        &X^{G'}\ar[r]^(.5){}
            &X^{G}\ar[r]^(.5){}
                &\left(\Sigma^{\sigma}X \right)^{G}
}
\end{displaymath}

\noindent so the exact sequence of homotopy groups is
\begin{center}
\includegraphics[width=12cm]{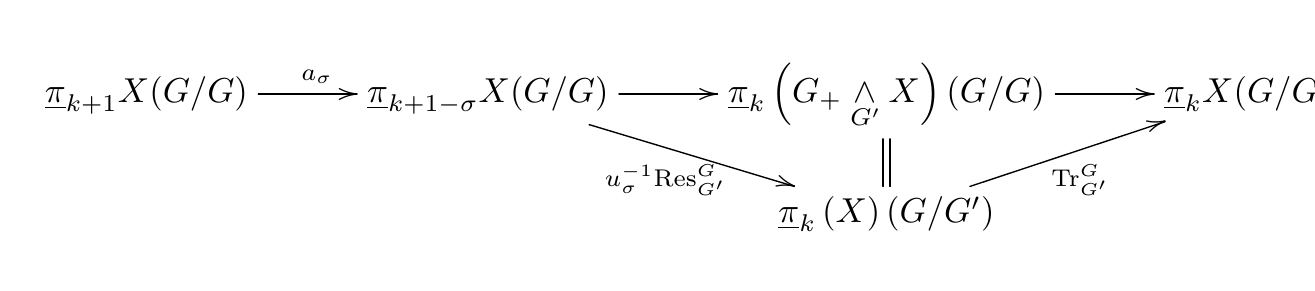}
\end{center}

\noindent Note that the isomorphism $u_{\sigma}$ is invertible. This
gives the exactness required by both statements.

For the statements about $a_{\lambda }$, note that $\lambda $ restricts
to $2\sigma_{G'}$, where $\sigma_{G'}$ is the sign {\rep} for the
index 2 subgroup $G'$.  It follows that $\Res_{G'}^{G}(a_{\lambda
})=a_{\sigma_{G'}}^{2}$, which has order 2.  Using the Frobenius
relation, we have for $x\in \upi_{*}X (G/G')$,
\begin{displaymath}
2a_{\lambda }\Tr_{G'}^{G}(x)
 = \Tr_{G'}^{G}(\Res_{G'}^{G}(2a_{\lambda })x)
 = \Tr_{G'}^{G}(2a_{\sigma_{G'} }^{2}x)
 = 0.
\end{displaymath}
\qed\bigskip

This implies that when $a_{\sigma}x$ is killed by a differential but
$x\in \EE_{r} (G/G)$ is not, then $x$ represents an element
that is $\Tr_{G'}^{G}(y)$ for some $y$ in lower filtration.  Similarly
if $x$ supports a nontrivial differential but $a_{\sigma}x$ is a
nontrivial permanent cycle, then the latter represents an element with
a nontrivial restriction to $G'$ of higher filtration.  In both cases
the converse also holds.

\begin{thm}\label{thm-exotic}
{\bf Exotic transfers and restrictions in the $RO (G)$-graded slice
{\SS}.} Let $G$ be a finite cyclic 2-group with index 2 subgroup $G'$
and sign {\rep} $\sigma $, and let $X$ be a $G$-{\eqvr} spectrum
with $x\in \EE_{r}^{s,V}X (G/G)$ (for $V\in RO (G)$) in the
slice {\SS} for $X$.  Then
\begin{enumerate}
\item [(i)] Suppose there is a permanent cycle $y'\in
\EE_{r}^{s+r,V+r-1}X (G/G')$.  Then there is a nontrivial
differential $d_{r} (x)=\Tr_{G'}^{G}(y')$ iff $[a_{\sigma}x]$ is a
permanent cycle with $\Res_{G'}^{G}(a_{\sigma}x)=u_{\sigma}y'$. In
this case $[a_{\sigma }x]$ represents the Toda bracket $\langle
a_{\sigma },\,\Tr_{G'}^{G} ,\,y' \rangle$.

\item [(ii)] Suppose there is a permanent cycle $y\in
\EE_{r}^{s+r-1 ,V+r+\sigma -2}X (G/G)$. Then there is a
nontrivial differential $d_{r}(x)=a_{\sigma}y$ iff $\Res_{G'}^{G}(x)$
is a permanent cycle with
$\Tr_{G'}^{G}(u_{\sigma}^{-1}\Res_{G'}^{G}(x))=y$.  In this case
$\Res_{G'}^{G}(x)$ represents the Toda bracket $\langle
\Res_{G'}^{G},\,a_{\sigma } ,\,y \rangle$.
\end{enumerate}

\end{thm}

In each case a nontrivial $d_{r}$ is equivalent to a Mackey functor
extension raising filtration by $r-1$.  In (i) the permanent cycle
$a_{\sigma}x$ is not divisible in $\upi_{\star}X$ by $a_{\sigma}$
and therefore could have a nontrivial restriction in a higher
filtration.  Similarly in (ii) the element denoted by
$\Res_{G'}^{G}(x)$ is not a restriction in $\upi_{\star}X$, so we
cannot use the Frobenius relation to equate
$\Tr_{G'}^{G}(u_{\sigma}^{-1}\Res_{G'}^{G}(x))$ with
$\Tr_{G'}^{G}(u_{\sigma}^{-1})x$.

We remark that the proof below makes no use of any properties
specific to the slice filtration.  The result holds for any {\eqvr}
filtration with suitable formal properties.

Before giving the proof we need the following.

\begin{lem}\label{lem-formal}
{\bf A formal observation}. Suppose we have a
commutative diagram up to sign
\begin{displaymath}
\xymatrix
@R=5mm
@C=10mm
{
A_{0,0}\ar[r]^(.5){a_{0,0}}\ar[d]^(.5){b_{0,0}}
    &A_{0,1}\ar[r]^(.5){a_{0,1}}\ar[d]^(.5){b_{0,1}}
        &A_{0,2}\ar[r]^(.5){a_{0,2}}\ar[d]^(.5){b_{0,2}}
            &\Sigma A_{0,0}\ar[d]^(.5){b_{0,0}}\\
A_{1,0}\ar[r]^(.5){a_{1,0}}\ar[d]^(.5){b_{1,0}}
    &A_{1,1}\ar[r]^(.5){a_{1,1}}\ar[d]^(.5){b_{1,1}}
        &A_{1,2}\ar[r]^(.5){a_{1,2}}\ar[d]^(.5){b_{1,2}}
            &\Sigma A_{1,0}\ar[d]^(.5){b_{1,0}}\\
A_{2,0}\ar[r]^(.5){a_{2,0}}\ar[d]^(.5){b_{2,0}}
    &A_{2,1}\ar[r]^(.5){a_{2,1}}\ar[d]^(.5){b_{2,1}}
        &A_{2,2}\ar[r]^(.5){a_{2,2}}\ar[d]^(.5){b_{2,2}}
            &\Sigma A_{2,0}\ar[d]^(.5){b_{2,0}}\\
\Sigma A_{0,0}\ar[r]^(.5){a_{0,0}}
    &\Sigma A_{0,1}\ar[r]^(.5){a_{0,1}}
        &\Sigma A_{0,2}\ar[r]^(.5){a_{0,2}}
            &\Sigma^{2} A_{0,0}
}
\end{displaymath}

\noindent in which each row and column is a cofiber sequence.  Then
suppose that from some spectrum $W$ we have a map $f_{3}$ and
hypothetical maps $f_{1}$ and $f_{2}$ making the following diagram
commute up to sign, where $c_{i,j}=b_{i,j+1}a_{i,j}=a_{i+1,j}b_{i,j}$.
\begin{numequation}\label{eq-formal}
\begin{split}
\xymatrix
@R=5mm
@C=10mm
{
W\ar[rrr]^(.5){f_{3}}\ar@{-->}[ddr]^(.5){f_{1}}
 \ar@{-->}[drr]^(.5){f_{2}}\ar[ddd]^(.5){f_{3}}
    &   &   &\Sigma A_{0,0}\ar[d]^(.5){b_{0,0}}\ar[dr]^(.5){c_{0,0}}\\
    &   &A_{1,2}\ar[r]^(.5){a_{1,2}}\ar[d]^(.5){b_{1,2}}\ar[dr]^(.5){c_{1,2}}
            &\Sigma A_{1,0}\ar[d]^(.5){b_{1,0}}\ar[r]^(.5){a_{1,0}}
                &\Sigma A_{1,1}\ar[d]^(.5){b_{1,1}}\\
    &A_{2,1}\ar[r]^(.5){a_{2,1}}\ar[d]^(.5){b_{2,1}}\ar[dr]^(.5){c_{2,1}}
        &A_{2,2}\ar[r]^(.5){a_{2,2}}\ar[d]^(.5){b_{2,2}}
            &\Sigma A_{2,0}\ar[r]^(.5){a_{2,0}}
                &\Sigma A_{2,1}\\
\Sigma A_{0,0}\ar[r]^(.5){a_{0,0}}\ar[rd]_(.5){c_{0,0}}
    &\Sigma A_{0,1}\ar[r]^(.5){a_{0,1}}\ar[d]^(.5){b_{0,1}}
        &\Sigma A_{0,2}\ar[d]^(.5){b_{0,2}}
            &\\
    &\Sigma A_{1,1}\ar[r]^(.5){a_{1,1}}
        &\Sigma A_{1,2}
}
\end{split}
\end{numequation}
 
\noindent Then $f_{1}$ exists iff $f_{2}$ does. When this happens,
$c_{0,0}f_{3}$ is null  and we have Toda brackets
\begin{displaymath}
\langle a_{1,1},\, c_{0,0} ,\,f_{3} \rangle \ni f_{2}
\qquad \aand 
\langle b_{1,1},\, c_{0,0} ,\,f_{3} \rangle \ni f_{1}.
\end{displaymath}

\end{lem}

\proof Let $R$ be the pullback of $a_{2,1}$ and $b_{1,2}$, so we have a diagram 
\begin{displaymath}
\xymatrix
@R=5mm
@C=10mm
{
{}
    &{A_{0,2}}\ar@{=}[r]^(.5){}\ar[d]^(.5){}
        &{A_{0,2}}\ar[d]^(.5){b_{0,2}}
            &{}\\
{A_{2,0}}\ar@{=}[d]^(.5){}\ar[r]^(.5){}
    &{R}\ar[d]^(.5){}\ar[r]^(.5){}
        &{A_{1,2}}\ar[d]^(.5){b_{1,2}}\ar[r]^(.5){c_{1,2}}
            &{\Sigma A_{2,0}}\ar@{=}[d]^(.5){}\\
{A_{2,0}}\ar[r]^(.5){a_{2,0}}
    &{A_{2,1}}\ar[d]^(.5){c_{2,1}}\ar[r]^(.5){a_{2,1}}
        &{A_{2,2}}\ar[d]^(.5){b_{2,2}}\ar[r]^(.5){a_{2,2}}
            &{\Sigma A_{2,0}}\\
{}
    &{\Sigma A_{0,2}}\ar@{=}[r]^(.5){}
        &{\Sigma A_{0,2}}
            &{}
}
\end{displaymath}

\noindent in which each row and column is a cofiber sequence.  Thus we
see that $R$ is the fiber of both $c_{1,2}$ and $c_{2,1}$.  If $f_{1}$
exists, then
\begin{displaymath}
c_{2,1}f_{1}=a_{0,1}b_{2,1}f_{1}=a_{0,1}a_{0,0}f_{3}
\end{displaymath}

\noindent which is null homotopic, so $f_{1}$ lifts to $R$, which
comes equipped with a map to $A_{1,2}$, giving us $f_{2}$.  Conversely
if $f_{2}$ exists, it lifts to $R$, which comes equipped with a map to
$A_{2,1}$, giving us $f_{1}$.

The statement about Toda brackets follows from the way they are defined.
\qed\bigskip 

\noindent {\em Proof of Theorem \ref{thm-exotic}}.
For a $G$-spectrum $X$ and integers $a<b<c\leq
\infty $ there is a cofiber sequence
\begin{displaymath}
\xymatrix
@R=5mm
@C=10mm
{
P_{b+1}^{c}X\ar[r]^(.5){i}
    &P_{a}^{c}X\ar[r]^(.5){j}
        &P_{a}^{b}X\ar[r]^(.5){k}
            &\Sigma P_{b+1}^{c}X.
}
\end{displaymath}

\noindent When $c=\infty $, we omit it from the notation.  We will
combine this and the one of (\ref{eq-hate}) to get a diagram similar
to (\ref{eq-formal}) with $W=S^{V}$ to prove our two statements.

For (i) note that $x\in \EE_{1}^{s,t}X (G/G)$ is by definition an
element in $\upi_{t-s}P^{s}_{s}X (G/G)$.  We will assume for
simplicity that $s=0$, so $x$ is represented by a map from some
$S^{V}$ to $(P_{0}^{0}X)^{G}$.  Its survival to $\EE_{r}$ and
supporting a nontrivial differential means that it lifts to
$(P_{0}^{r-2}X)^{G}$ but not to $(P_{0}^{r-1}X)^{G}$.  The value of
$d_{r} (x)$ is represented by the composite $kx$ in the diagram below,
where we can use Lemma \ref{lem-formal}.
\begin{displaymath}
\xymatrix
@R=4mm
@C=12mm
{ 
S^{V-1}\ar[ddd]^(.5){y'}\ar@{-->}[ddr]^(.5){x}
     \ar@{-->}[drr]^(.5){w}\ar[rrr]^(.5){y'}
    &   &   &(P_{r-1}X)^{G'} \ar[d]^(.5){i}\\
    &   &(\Sigma^{\sigma -1}P_{0}X)^{G}\ar[d]^(.5){j}
                 \ar[r]^(.55){u_{\sigma}^{-1}\Res_{G'}^{G}}
            &(P_{0}X )^{G'}\ar[d]^(.5){j}\\
    &(\Sigma^{-1}P_{0}^{r-2}X)^{G}\ar[d]^(.5){k}\ar[r]^(.5){a_{\sigma}}
        &(\Sigma^{\sigma-1 }P_{0}^{r-2}X)^{G}\ar[d]^(.5){k}
                 \ar[r]^(.55){u_{\sigma}^{-1}\Res_{G'}^{G}}
            &(P_{0}^{r-2}X)^{G'} \\
(P_{r-1}X)^{G'}\ar[r]^(.5){\Tr_{G'}^{G}}
    &(P_{r-1}X)^{G}\ar[r]^(.5){a_{\sigma}}
        &(\Sigma^{\sigma} P_{r-1}X)^{G}
}
\end{displaymath}

\noindent The commutativity of the lower left trapezoid is the
differential of (i),\linebreak ${d_{r} (x)=\Tr_{G'}^{G}(y')}$.  The
existence of the map $w$ making the diagram commute follows from that
of $x$ and $y'$. It is the representative of $a_{\sigma}x$ as a
permanent cycle, which represents the indicated Toda bracket.  The
commutativity of the upper right trapezoid identifies $y'$ as
$u_{\sigma}^{-1}\Res_{G'}^{G}(x)$ as claimed.  For the converse we
have the existence of $y'$ and $w$ and hence that of $x$.

The second statement follows by a similar argument based on the diagram
\begin{displaymath}
\xymatrix
@R=4mm
@C=13mm
{ 
S^{V+\sigma -1}\ar[ddd]^(.5){y}\ar@{-->}[ddr]^(.5){x}
     \ar@{-->}[drr]^(.5){w}\ar[rrr]^(.5){y}
    &   &   &(P_{r-1}X)^{G} \ar[d]^(.5){i}\\
    &   &(P_{0}X)^{G'}\ar[d]^(.5){j}
                 \ar[r]^(.5){\Tr_{G'}^{G}}
            &( P_{0}X )^{G}\ar[d]^(.5){j}\\
    &(\Sigma^{\sigma -1}P_{0}^{r-2}X)^{G}\ar[d]^(.5){k}
                       \ar[r]^(.55){u_{\sigma}^{-1}\Res_{G'}^{G}}
        &(P_{0}^{r-2}X)^{G'}\ar[d]^(.5){k}
                       \ar[r]^(.5){\Tr_{G'}^{G}}
            &( P_{0}^{r-2}X)^{G} \\
( P_{r-1}X)^{G}\ar[r]^(.5){a_{\sigma}}
    &(\Sigma^{\sigma} P_{r-1}X)^{G}
                \ar[r]^(.55){u_{\sigma}^{-1}\Res_{G'}^{G}}
        &(\Sigma P_{r-1}X)^{G'}.
            &
}
\end{displaymath}

\noindent Here $w$ represents $u_{\sigma }^{-1}\Res_{G'}^{G}(x)$ as a
permement cycle, so we get a Toda bracket containing
$\Res_{G'}^{G}(x)$ as indicated. \qed

\bigskip

Next we study the way differentials interact with the norm.
Suppose we have a subgroup $H\subset G$ and an $H$-{\eqvr} ring
spectrum $X$ with $Y=N_{H}^{G}X$.  Suppose we have {\SS}s converging
to $\upi_{\star}X$ and $\upi_{\star}Y$ based on towers
\begin{displaymath}
...\to P_{n}^{H}X \to P_{n-1}^{H}X \to \dotsb 
\qquad \aand 
...\to P_{n}^{G}Y \to P_{n-1}^{G}Y \to \dotsb 
\end{displaymath}
 
\noindent for functors $P_{n}^{H}$ and $P_{n}^{G}$ equipped with
suitable maps
\begin{displaymath}
P^{H}_{m}X \wedge P^{H}_{n}X \to  P^{H}_{m+n}X ,\quad 
P^{G}_{m}Y \wedge P^{G}_{n}Y \to  P^{G}_{m+n}Y 
\quad \mbox{and}\quad  
N_{H}^{G}P^{H}_{m}X \to P^{G}_{m|G/H|}Y.
\end{displaymath}

Our slice {\SS } for each of the spectra studied in this paper fits
this description.

\begin{thm}\label{thm-normdiff}
{\bf The norm of a differential.} Suppose we have spectral sequences
as described above and a differential $d_{r} (x)=y$ for $x\in
\EE_{r}^{s,\star}X (H/H)$.  Let $\rho=\Ind_{H}^{G}1 $ and suppose that
$a_{\rho }$ has filtration $|G/H|-1$.  Then in the {\SS } for
$Y=N_{H}^{G}X$,
\begin{displaymath}
d_{|G/H| (r-1) +1} (a_{\rho } N_{H}^{G}x)= N_{H}^{G}y
\in \EE_{|G/H| (r-1) +1}^{|G/H| (s+r),\star}Y (G/G).
\end{displaymath}
\end{thm}

\proof The differential can be represented by a diagram
\begin{displaymath}
\xymatrix
@R=5mm
@C=10mm
{
S^{V}\ar@{=}[r]^(.5){}
&S (1+V)\ar[r]^(.5){}\ar[d]_(.5){y}
    &D (1+V)\ar[r]^(.5){}\ar[d]^(.5){}
        &S^{1+V}\ar[d]^(.5){x}\\
&P_{s+r}^{H}X\ar[r]^(.5){}
    &P_{s}^{H}X\ar[r]^(.5){}
        &P^{H}_{s}X/P^{H}_{s+r}X
}
\end{displaymath}

\noindent for some orthogonal {\rep} $V$ of $H$, where each row is a
cofiber sequence.  We want to apply the norm functor $N_{H}^{G}$ to
it.  Let $W=\Ind_{H}^{G}V$.  Then we get
\begin{displaymath}
\xymatrix
@R=5mm
@C=10mm
{
S^{W}\ar@{=}[r]^(.5){}
&N_{H}^{G}S (1+V)\ar[r]^(.5){}\ar[d]_(.5){N_{H}^{G}y}
    &D (\rho +W)\ar[r]^(.5){}\ar[d]^(.5){}
        &S^{\rho +W}\ar[d]^(.5){N_{H}^{G}x}\\
&N_{H}^{G}P_{s+r}^{H}X\ar[r]^(.5){}
    &N_{H}^{G}P_{s}^{H}X\ar[r]^(.5){}
        &N_{H}^{G} (P_{s}^{H}X/P_{s+r}^{H}X).                   
}
\end{displaymath}

\noindent Neither row of this diagram is a cofiber sequence, but we
can enlarge it to one where the top and bottom rows are, namely
\begin{displaymath}
\xymatrix
@R=5mm
@C=10mm
{
S^{W}\ar[r]^(.5){}\ar@{=}[d]
    &D (1+W)\ar[r]^(.5){}\ar[d]^(.5){a_{\rho }}
        &S^{1 +W}\ar[d]^(.5){a_{\rho }}\\
S^{W}\ar[r]^(.5){}\ar[d]_(.5){N_{H}^{G}y}
    &D (\rho +W)\ar[r]^(.5){}\ar[d]^(.5){}
        &S^{\rho +W}\ar[d]^(.5){N_{H}^{G}x}\\
N_{H}^{G}P_{s+r}^{H}X\ar[r]^(.5){}\ar[d]^(.5){}
    &N_{H}^{G}P_{s}^{H}X\ar[r]^(.5){}\ar[d]^(.5){}
        &N_{H}^{G} (P_{s}^{H}X/P_{s+r}^{H}X)\ar[d]^(.5){}\\
P_{(s+r)|G/H|}^{G}Y\ar[r]^(.5){}                  
    &P_{s|G/H|}^{G}Y\ar[r]^(.5){}
        &P_{s|G/H|}^{G}Y/P_{(s+r) |G/H|}^{G}Y                
}
\end{displaymath}

\noindent Here the first two bottom vertical maps are part of the
multiplicative structure the {\SS } is assumed to have.  Composing the
maps in the three columns gives us the diagram for the desired
differential.  \qed\bigskip

Given a $G$-{\eqvr} ring spectrum $X$, let $X'=i^{*}_{H}X$ denote its
restriction as an $H$-spectrum.  Then $N_{H}^{G}X'=X^{(|G/H|)}$ and
the multiplication on $X$ gives us a map from this smash product to
$X$.  This gives us a map $\pi_{\star}X' \to \pi_{\star}X$ called the
{\em internal norm}, which we denote abusively by $N_{H}^{G}$.  The
argument above yields the following.

\begin{cor}\label{cor-normdiff}
{\bf The internal norm of a differential.}  With notation as above,
suppose we have a differential $d_{r} (x)=y$ for $x\in
\EE_{r}^{s,\star}X' (H/H)$.  Then
\begin{displaymath}
d_{|G/H| (r-1) +1} (a_{\rho } N_{H}^{G}x)= N_{H}^{G}y
\in \EE_{|G/H| (r-1) +1}^{|G/H| (s+r),\star}X (G/G).
\end{displaymath}
\end{cor}

The following is useful in making such calculations. It is very
similar to \cite[Lemma 3.13]{HHR}.

\begin{lem}\label{lem-norm-au}
{\bf The norm of $a_{V}$ and $u_{V}$}.  With notation as above, let
$V$ be a {\rep} of $H$ with $V^{H}=0$ and let $W=\Ind_{H}^{G}V$.  Then
$N_{H}^{G} (a_{V})=a_{W}$.  If $V$ is oriented (and hence even
dimensional, making $|V|\rho $ oriented), then 
\begin{displaymath}
u_{|V|\rho }N
(u_{V})=u_{W}.
\end{displaymath}
\end{lem}

\proof The element $a_{V}$ is represented by the map $S^{0}\to
S^{V}$, the inclusion of the fixed point set. Applying the norm
functor to this map gives
\begin{displaymath}
S^{0}=N_{H}^{G}S^{0}\to N_{H}^{G}S^{V}=S^{W},
\end{displaymath}

\noindent which is $a_{W}$.

When $V$ is oriented, $u_{V}$ is represented by a map $S^{|V|}\to
S^{V}\wedge \HZZ$.  Applying the norm functor and using the
multiplication in $\HZZ$ leads to a map
\begin{displaymath}
\xymatrix
@R=8mm
@C=5mm
{
{S^{|V|\rho }}\ar@{=}[r]^(.5){}
    &{N_{H}^{G}S^{|V|}}\ar[rr]^(.5){N_{H}^{G}u_{V}}
        &   &{S^{W}\wedge \HZZ}
}
\end{displaymath}

\noindent Now smash both sides with $\HZZ$, precompose with
$u_{|V|\rho }$ and follow with the multiplication on $\HZZ$, giving
\begin{displaymath}
\xymatrix
@R=8mm
@C=8mm
{
{S^{|V||\rho |}}\ar[r]^(.4){u_{|V|\rho }}
    &{S^{|V|\rho }\wedge \HZZ}\ar[rr]^(.45){N_{H}^{G}u_{V}\wedge \HZZ}
        &   &{S^{W}\wedge \HZZ\wedge \HZZ}\ar[r]^(.5){}
                &{S^{W}\wedge \HZZ,}
}
\end{displaymath}

\noindent which is $u_{W}$ since $|W|=|V||\rho |$.
\qed\bigskip

\section{Some Mackey functors for $C_{4}$ and $C_{2}$}\label{sec-C4}
 
We need some notation for Mackey functors to be used in {\SS} charts.
{\em In this paper, when a cyclic group or subgroup appears as an
index, we will often replace it by its order.}  We can specify Mackey
functors $\underline{M}$ for the group $C_{2}$ and $\underline{N}$ for
$C_{4}$ by means of Lewis diagrams (first introduced in
\cite{Lewis:ROG}),
\begin{numequation}\label{eq-Lewis}
\begin{split}
\xymatrix
@R=8mm
@C=10mm
{
{\underline{M} (C_{2}/C_{2})}\ar@/_/[d]_(.5){\Res_{1}^{2}}\\
{\underline{M} (C_{2}/e)}\ar@/_/[u]_(.5){\Tr_{1}^{2}} 
}
\qquad \aand 
\xymatrix
@R=4mm
@C=10mm
{
{\underline{N} (C_{4}/C_{4})}\ar@/_/[d]_(.5){\Res_{2}^{4}}\\
{\underline{N} (C_{4}/C_{2})}\ar@/_/[d]_(.5){\Res_{1}^{2}}
                           \ar@/_/[u]_(.5){\Tr_{2}^{4}}\\
{\underline{N} (C_{4}/e).}\ar@/_/[u]_(.5){\Tr_{1}^{2}}
}
\end{split}
\end{numequation}%


\noindent We omit Lewis' looped arrow indicating the Weyl group action
on $\underline{M} (G/H)$ for proper subgroups $H$.  This notation is
prohibitively cumbersome in {\SS } charts, so we will abbreviate
specific examples by more concise symbols.  These are shown in Tables
\ref{tab-C2Mackey} and \ref{tab-C4Mackey}.  {\em Admittedly some of
these symbols are arbitrary and take some getting used to, but we have
to start somewhere.}  Lewis denotes the fixed point Mackey functor
for a $\Z G$-module $M$ by $R (M)$.  He abbreviates $R (\Z)$and $R
(\Z_{-})$ by $R$ and $R_{-}$. He also defines (with similar
abbreviations) the orbit group Mackey functor $L (M)$ by
\begin{displaymath}
L (M) (G/H)=M/H.
\end{displaymath}

\noindent In this case each transfer map is the surjection of the
orbit space for a smaller subgroup onto that of a larger one.  The
functors $R$ and $L$ are the left and right adjoints of the forgetful
functor $\underline{M}\mapsto \underline{M} (G/e)$ from Mackey
functors to $\Z G$-modules.

Over $C_{2}$  we have short exact sequences
\begin{displaymath}
\xymatrix
@R=5mm
@C=10mm
{
0 \ar[r]^(.5){}
   &{\twobox} \ar[r]^(.5){}
      &\Box \ar[r]^(.5){}
        &\bullet \ar[r]^(.5){}
           &0\\
0\ar[r]^(.5){}
   &\bullet \ar[r]^(.5){}
      &{\dot{\Box}} \ar[r]^(.5){}
        &{\oBox} \ar[r]^(.5){}
           &0\\
0\ar[r]^(.5){}
   &\Box \ar[r]^(.5){}
      &{\widehat{\Box}} \ar[r]^(.5){}
        &{\oBox} \ar[r]^(.5){}
           &0
}
\end{displaymath}

\noindent We can apply the induction functor to each them to get a
short exact sequence of Mackey functors over $C_{4}$.

Five of the Mackey functors in Table \ref{tab-C4Mackey} are fixed
point Mackey functors (\ref{eq-fpm}), meaning they are fixed points of
an underlying $Z[G]$-module $M$, such as $\Z[G]$ or
\begin{displaymath}
\begin{array}[]{rcllcl}
\Z     & = & \Z[G]/ (\gamma -1)\qquad \qquad &
\Z[G/G']& = &\Z[G]/ (\gamma^{2} -1)\\
\Z_{-} & = &\Z[G]/ (\gamma +1)&
\Z[G/G']_{-}
       & = & \Z[G]/ (\gamma^{2}+1)
\end{array}
\end{displaymath}

\begin{table}
\caption[Some {$C_{2}$}-Mackey functors]{Some $C_{2}$-Mackey functors}
\label{tab-C2Mackey}
\begin{center}
\begin{tabular}{|p{2.1cm}|c|c|c|c|c|
c|}
\hline 
Symbol
&$\Box$
   &$\oBox $
      &$\bullet$
         &$\twobox$
            &$\dot{\Box} $
                     &$\widehat{\Box}$\\
\hline 
Lewis\newline   diagram
&$
\xymatrix
@R=5mm
@C=10mm
{
\Z\ar@/_/[d]_(.5){1}\\
\Z\ar@/_/[u]_(.5){2} 
}$ 
  &
$
\xymatrix
@R=4mm
@C=10mm
{
0\ar@/_/[d]_(.5){}\\
\Z_{-}\ar@/_/[u]_(.5){}
}$
         &
$
\xymatrix
@R=4mm
@C=10mm
{
\Z/2\ar@/_/[d]_(.5){}\\
0\ar@/_/[u]_(.5){}
}
$         &
$\xymatrix
@R=5mm
@C=10mm
{
\Z\ar@/_/[d]_(.5){2}\\
\Z\ar@/_/[u]_(.5){1}
}$           &
$\xymatrix
@R=5mm
@C=10mm
{
\Z/2\ar@/_/[d]_(.5){0}\\
\Z_{-}\ar@/_/[u]_(.5){1}
}$            
               &
$\xymatrix
@R=5mm
@C=10mm
{
\Z\ar@/_/[d]_(.5){\Delta}\\
\Z[C_{2}]\ar@/_/[u]_(.5){\nabla}
}$\\
\hline 
Lewis symbol
&$R$&$R_{-}$
        &$\langle \Z/2 \rangle$
            &$L$&$L_{-}$
                    &$R (\Z^{2})$                \\
\hline 
\end{tabular}
\end{center}
\end{table}

\begin{table}
\caption[Some {$C_{4}$}-Mackey functors]{Some $C_{4}$-Mackey functors,
where $G=C_{4}$ and $G'$ is its index 2 subgroup.  The notation
$\underline{\Z}(G,H)$ is defined in \ref{prop-top}(i).}
\label{tab-C4Mackey}
\begin{center} 
\includegraphics[width=12cm]{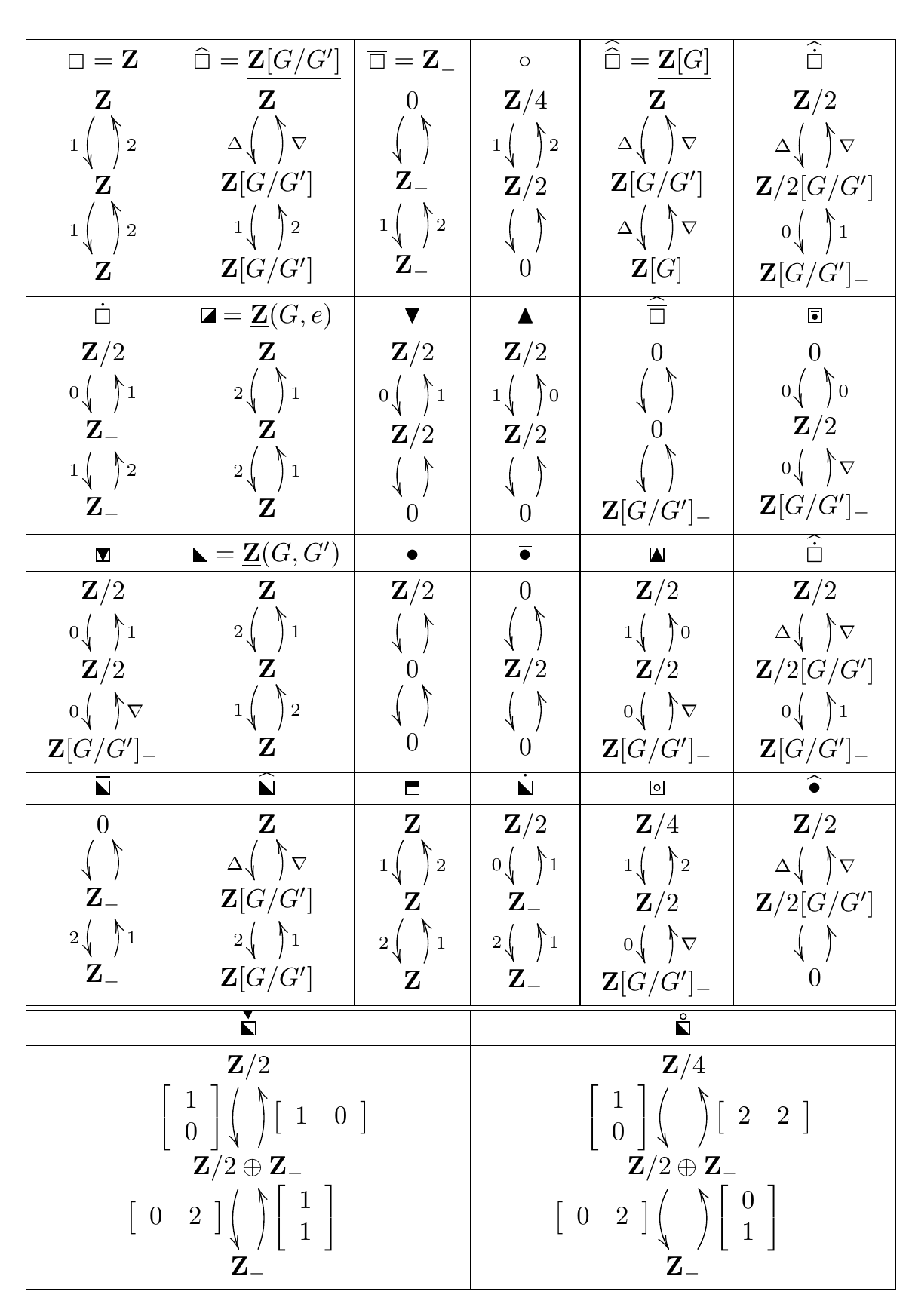}
\end{center}
\end{table}

We will use the following notational conventions for $C_{4}$-Mackey functors.
\begin{enumerate}
\item [$\bullet$] Given a $C_{2}$-Mackey functor $\underline{M}$ with Lewis diagram 
\begin{displaymath}
\xymatrix
@R=8mm
@C=10mm
{
{A}\ar@/_/[d]_(.5){\alpha }\\
{B}\ar@/_/[u]_(.5){\beta }
}
\end{displaymath}

\noindent with $A$ and $B$ cyclic, we will use the symbols
$\underline{M}$, $\underline{\overline{M} }$ and $\dot{\underline{M}}$
for the $C_{4}$-Mackey functors with Lewis diagrams
\begin{displaymath}
\xymatrix
@R=8mm
@C=10mm
{
{A}\ar@/_/[d]_(.5){\alpha }\\
{B}\ar@/_/[u]_(.5){\beta }\ar@/_/[d]_(.5){1 }\\
{B,}\ar@/_/[u]_(.5){2 }
}
\qquad 
\xymatrix
@R=8mm
@C=10mm
{
{0}\ar@/_/[d]_(.5){ }\\
{A_{-}}\ar@/_/[d]_(.5){\alpha }\ar@/_/[u]_(.5){ }\\
{B_{-}}\ar@/_/[u]_(.5){\beta }
}
\qquad \aand 
\xymatrix
@R=8mm
@C=10mm
{
{\Z/2}\ar@/_/[d]_(.5){0 }\\
{A_{-}}\ar@/_/[d]_(.5){\alpha }\ar@/_/[u]_(.5){\tau  }\\
{B_{-}}\ar@/_/[u]_(.5){\beta }
}\end{displaymath}

\noindent where a generator $\gamma \in C_{4}$ acts via multiplication
by $-1$ on $A$ and $B$ in the second two, and the transfer $\tau $ is
nontrivial. 

\item [$\bullet$] For a $C_{2}$-Mackey functor $\underline{M}$ we
will denote $\uparrow_{2}^{4}\underline{M}$ (see \ref{def-MF-induction}) by
$\widehat{\underline{M}}$.  For a Mackey functor $\underline{M}$
defined over the trivial group, we will denote
$\uparrow_{1}^{2}\underline{M}$ and $\uparrow_{1}^{4}\underline{M}$ by
$\widehat{\underline{M}}$ and $\widehat{\widehat{\underline{M}}}$.

\end{enumerate}

Over $C_{4}$, in addition to the short exact sequences induced up from
$C_{2}$, we have  
\begin{numequation}\label{eq-C4-SES}
\begin{split}
\xymatrix
@R=2mm
@C=10mm
{
0\ar[r]
   &{\bullet}\ar[r]
      &{\dot{\Box}}\ar[r]
         &{\oBox}\ar[r]
            &0\\
0\ar[r]
   &{\JJ}\ar[r]
      &{\circ }\ar[r]
         &\bullet\ar[r]
            &0\\
0\ar[r]
   &{\JJ}\ar[r]
      &{\JJbox}\ar[r]
         &{\widehat{\oBox}}\ar[r]
            &0\\
0\ar[r]
   &{\bullet}\ar[r]
      &{\circ }\ar[r]
         &\JJJ\ar[r]
            &0\\
0\ar[r]
   &{\fourbox}\ar[r]
      &\Box\ar[r]
         &\circ\ar[r]
            &0\\
0\ar[r]
   &{\twobox}\ar[r]
      &\Box\ar[r]
         &\bullet\ar[r]
            &0
}
\end{split}
\end{numequation}%

\begin{defin}\label{def-enriched}
{\bf A $C_{4}$-enriched  $C_{2}$-Mackey functor.} 
For a $C_{2}$-Mackey functor $\uM$ as above,
$\widetilde{\uM}$ will denote the $C_{2}$-Mackey functor
enriched over $\Z[C_{4}]$ defined by
\begin{displaymath}
\widetilde{\uM} (S)
    =\Z[C_{4}]\otimes_{\Z[C_{2}]}\uM (S)
\end{displaymath}

\noindent for a finite $C_{2}$-set $S$.  Equivalently, in the notation
of Definition \ref{def-MF-induction}, $\widetilde{\uM}=
\downarrow_{2}^{4}\uparrow_{2}^{4}\uM$.
\end{defin}

\section{Some chain complexes of Mackey functors}\label{sec-chain}

As noted above, a $G$-CW complex $X$, meaning one built out of cells
of the form $G_{+}\smashove{H} e^{n}$, has a reduced cellular chain
complex of $\Z[G]$-modules $C_{*}X$, leading to a chain complex of
fixed point Mackey functors (see (\ref{eq-fpm})) $\uC_{*}X$.  When
$X=S^{V}$ for a {\rep} $V$, we will denote this complex by
$\uC_{*}^{V}$; see (\ref{eq-CVn}).  Its homology is the graded Mackey
functor $\uH_{*}X$.  Here we will apply the methods of \S\ref{sec-HZZ}
to three examples.

\bigskip (i) Let $G=C_{2}$ with generator $\gamma $, and $X=S^{n\rho
}$ for $n>0$, where $\rho $ denotes the regular {\rep}.  We have seen
before \cite[Ex. 3.7]{HHR} that it has a reduced cellular chain complex $C$
with
\begin{numequation}
\label{eq-C2chain}
C_{i} ^{n\rho_{2}} =\mycases{
\Z[G]/ (\gamma -1)
       &\mbox{for }i=n\\
\Z[G]
       &\mbox{for }n<i\leq 2n\\
0      &\mbox{otherwise}.
}
\end{numequation}%

\noindent Let $c_{i}^{(n)}$ denote a generator of $C_{i}^{n\rho_{2}}$.
The boundary operator $d$ is given by
\begin{numequation}
\label{eq-C2boundary}
d (c_{i+1}^{(n)}) =\mycases{
c_{i}^{(n)}
       &\mbox{for }i=n\\
\gamma_{i+1-n} (c_{i}^{(n)})
       &\mbox{for }n<i\leq 2n\\
0      &\mbox{otherwise}
}
\end{numequation}%

\noindent where $\gamma_{i}=1- (-1)^{i}\gamma $.  For future reference, let 
\begin{displaymath}
\epsilon_{i}=1- (-1)^{i}=\mycases{
0       &\mbox{for $i$ even}\\
2       &\mbox{for $i$ odd.}
}
\end{displaymath}

\noindent  This chain complex has the form
\begin{displaymath}
\xymatrix
@R=3mm
@C=10mm
{
n  &n+1 
      &n+2 
         &n+3 
            &  &2n\\
\Box
   &{\widehat{\Box}}\ar[l]_(.5){\nabla}
      &{\widehat{\Box}}\ar[l]_(.5){\gamma_{2} }
         &{\widehat{\Box}}\ar[l]_(.5){\gamma_{3} }
            &\dotsb \ar[l]_(.5){}
               &{\widehat{\Box}}\ar[l]_(.5){\gamma_{n}}\\
\Z\ar@/_/[d]_(.5){1}
   &\Z\ar@/_/[d]_(.5){\Delta}\ar[l]_(.5){2}
      &\Z\ar@/_/[d]_(.5){\Delta}\ar[l]_(.5){0}
         &\Z\ar@/_/[d]_(.5){\Delta}\ar[l]_(.5){2}
            &\dotsb \ar[l]_(.5){}
               &\Z\ar@/_/[d]_(.5){\Delta}\ar[l]_(.5){\epsilon_{n}}\\
\Z\ar@/_/[u]_(.5){2}
   &\Z[G]\ar@/_/[u]_(.5){\nabla}\ar[l]_(.5){\nabla}
      &\Z[G]\ar@/_/[u]_(.5){\nabla}\ar[l]_(.5){\gamma_{2}}
         &\Z[G]\ar@/_/[u]_(.5){\nabla}\ar[l]_(.5){\gamma_{3} }
            &\dotsb \ar[l]_(.5){}
               &\Z[G]\ar@/_/[u]_(.5){\nabla}\ar[l]_(.5){\gamma_{n}}
}
\end{displaymath}

\noindent  Passing to homology we get 
\begin{displaymath}
\xymatrix
@R=3mm
@C=7mm
{
n  &n+1 
      &n+2 
         &n+3 
            &  &2n\\
\bullet
   &0 &\bullet
         &0  &\dotsb 
               &{\uH_{2n}} \\
\Z/2\ar@/_/[d]_(.5){}
   &0\ar@/_/[d]_(.5){ }
      &\Z/2\ar@/_/[d]_(.5){ }
         &0\ar@/_/[d]_(.5){ }
            &\dotsb 
               &{\uH_{2n} (G/G)}
                         \ar@/_/[d]_(.5){\Delta  }\\
0\ar@/_/[u]_(.5){}
   &0\ar@/_/[u]_(.5){}
      &0\ar@/_/[u]_(.5){}
         &0\ar@/_/[u]_(.5){}
            &\dotsb 
               &\Z[G]/ (\gamma_{n+1} )\ar@/_/[u]_(.5){\nabla}
}
\end{displaymath}

\noindent where 
\begin{displaymath}
\uH_{2n} (G/G)
=\mycases{
\Z
       &\mbox{for $n$ even}\\
0
       &\mbox{for $n$ odd}
}
\qquad \aand 
\uH_{2n}
=\mycases{
\Box
       &\mbox{for $n$ even}\\
\overline{\Box}
       &\mbox{for $n$ odd}
}
\end{displaymath}

\noindent Here $\Box$ and $\overline{\Box}$ are fixed point Mackey
functors but $\bullet$ is not.

\begin{figure} 
\begin{center}
\includegraphics[width=12cm]{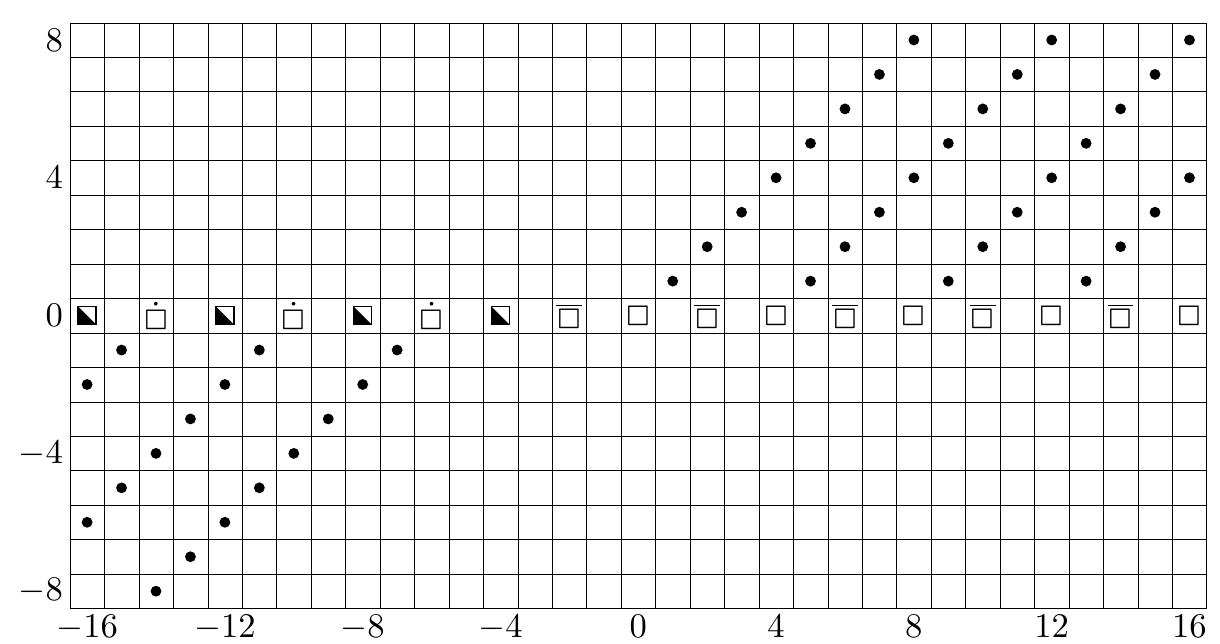} \caption[ The slice
spectral sequence for {${\bigvee_{n \in \Z}\Sigma^{n\rho_{2}}\HZZ }$}.]
{The (collapsing) Mackey functor slice spectral sequence for
${\bigvee_{n \in \Z}\Sigma^{n\rho_{2}}\HZZ }$.  The symbols are
defined in Table \ref{tab-C2Mackey}.  When the Mackey functor
$\upi_{(2-\rho_{2})n-s}\HZZ=\underline{H}_{2n-s}S^{n\rho_{2}}$ is
nontrivial, it is shown at ${(2n-s,s)}$ in the chart. Compare with
Figure \ref{fig-KR}.  } \label{fig-sseq-1}
\end{center}
\end{figure}

Similar calculations can be made for $S^{n\rho_{2}}$ for $n<0$.  The
results are indicated in Figure \ref{fig-sseq-1}.  This is originally
due to unpublished work of Stong and is reported in \cite[Theorem 2.1
and Table 2.2]{Lewis:ROG}.  This information will be used in
\S\ref{sec-Dugger}.

In other words the $RO (G)$-graded Mackey functor valued homotopy of
$\HZZ$ is as follows.  For $n\geq -1$ we have
\begin{displaymath}
\upi_{i}\Sigma^{n\rho_{2}}\HZZ
 =\upi_{i-n\rho_{2} }\HZZ = 
\mycases{
\Box
       &\mbox{for  $n$ even and $i=2n$}\\
\overline{\Box} 
       &\mbox{for $n$ odd and $i=2n$}\\
\bullet
       &\mbox{for $n\leq i<2n$ and $i+n$ even} \\
0      &\mbox{otherwise}
}
\end{displaymath}

\noindent For $n\leq -2$ we have
\begin{displaymath}
\upi_{i}\Sigma^{n\rho_{2}}\HZZ
 =\upi_{i-n\rho_{2} }\HZZ = 
\mycases{
\twobox 
       &\mbox{for $n$ even and $i=2n$}\\
\dot{\Box} 
       &\mbox{for $n$ odd and $i=2n$}\\
\bullet
       &\mbox{for $2n<i\leq n-3$ and $i+n$ odd}\\
0      &\mbox{otherwise}
}
\end{displaymath}

We can use Definition \ref{def-aeu} to name some elements of these groups.

Note that $\HZZ$ is a commutative ring spectrum, so there is a
commutative multiplication in $\upi_{\star}\HZZ$, making it a
commutative $RO (G)$-graded Green functor.  For such a functor $\uM$
on a general group $G$, the restriction maps are a ring homomorphisms
while the transfer maps satisfy the Frobenius relations
(\ref{eq-Frob}).
 
Then the generators of various groups in $\upi_{\star}\HZZ$ are
\begin{align*}
\mbox{\sc $(4m-2)$-slices for $m>0$ :}\hspace{-1cm} \\
a^{2m-1-2i}u^{i}
 & =      a_{(2m-1-2i)\sigma}u_{2i\sigma}\\
 & \in    \upi_{2m-1+2i }\Sigma^{(2m-1)\rho_{2}}\HZZ (G/G)\\
 & =      \upi_{2i- (2m-1)\sigma}\HZZ (G/G)\\                        
        &\qquad \mbox{for $0\leq i<m$}\\
x^{2m-1}=u_{(2m-1)\sigma}
 & \in    \upi_{4m-2}\Sigma^{(2m-1)\rho_{2}}\HZZ (G/\ee)\\
 &  =    \upi_{(2m-1) (1-\sigma ) }\HZZ (G/\ee)\\
        &\qquad \mbox{with  $\gamma (x)=-x$}\\
\mbox{\sc $4m$-slices for $m>0$ :} \\
a^{2m-2i}u^{i}
 & =      a_{(2m-2i)\sigma}u_{2i\sigma}\\
 & \in    \upi_{2m-1+2i }\Sigma^{(2m-1)\rho_{2}}\HZZ (G/G)\\
 & =      \upi_{2i- (2m-1)\sigma}\HZZ (G/G)\\                        
        &\qquad \mbox{for $0\leq i\leq m$}\\
 &\qquad \mbox{ and with }\Res (u)=x^{2}\\  
\mbox{\sc negative slices:}  \\ 
z_{n}=e_{2n\rho_{2} }
 & \in    \upi_{-4n}\Sigma^{-2n\rho_{2}}\HZZ (G/\ee)\\
 &     = \upi_{2n (\sigma -1)}\HZZ (G/\ee)
                \qquad \mbox{for $n>0$}\\
a^{-i}\Tr (x^{-2n-1})
  & \in  \upi_{-4n-2-i}\Sigma^{-(2n+1+i)\rho_{2}}\HZZ (G/G)\\
  &   = \upi_{(2n+1) (\sigma -1)+ i\sigma}\HZZ (G/G)\\
&\qquad 
\qquad \mbox{for $n>0$ and $i\geq 0$} .
\end{align*}

\noindent We have
relations
\begin{displaymath}
\begin{array}[]{rclrcl}
2a  &\hspace{-2.5mm} = & 0\qquad &                       
\Res (a) 
    &\hspace{-2.5mm} = & 0\\
z_{n} &\hspace{-2.5mm} = & x^{-2n}\qquad &
\Tr (x^{n}) &\hspace{-2.5mm} = &\mycases{
2u^{n/2}
        &\mbox{for $n$ even and $n\geq 0$}\\
\Tr_{}^{}(z_{-n/2})
        &\mbox{for $n$ even and $n< 0$}\\
0       &\mbox{for $n$ odd and $n>-3$.}
}
\end{array}
\end{displaymath}

\bigskip

(ii) Let $G=C_{4}$ with generator $\gamma $, $G'=C_{2}\subseteq G$, the
subgroup generated by $\gamma^{2}$, and $\widehat{S} (n,G')
=G_{+}\smashove{G'}S^{n\rho_{2}}$.  Thus we have
\begin{displaymath}
C_{*} (\widehat{S} (n,G')) = \Z[G]\otimes_{\Z[G']} C_{*}^{n\rho_{2}}
\end{displaymath}

\noindent with $C_{*}^{n\rho_{2}}$ as in (\ref{eq-C2chain}). The
calculations of the previous example carry over verbatim by the
exactness of Mackey functor induction of Definition \ref{def-MF-induction}.

\bigskip
\bigskip(iii) Let $G=C_{4}$ and $X=S^{n\rho_{4}}$.  Then the reduced
cellular chain complex (\ref{eq-CVn}) has the form
\begin{displaymath}
C_{i}^{n\rho_{4}} = \mycases{
\Z       &\mbox{for }i=n\\
\Z[G/G']
         &\mbox{for }n<i\leq 2n\\
\Z[G]
         &\mbox{for }2n<i\leq 4n\\
0        &\mbox{otherwise}
}
\end{displaymath}

\noindent in which generators $c_{i}^{(n)}\in C_{i}^{n\rho_{4}}$ satisfy
\begin{displaymath}
d (c_{i+1}^{(n)})= \mycases{
c_{i}^{(n)}
       &\mbox{for }i=n\\
\gamma_{i+1-n} c_{i}^{(n)}
       &\mbox{for $n<i\leq 2n$}\\
\mu_{i+1-n} c_{i}^{(n)}
       &\mbox{for $2n< i<4n$ and $i$ even}\\
\gamma_{i+1-n} c_{i}^{(n)}
       &\mbox{for $2n<i<4n$ and $i$ odd}\\
0       &\mbox{otherwise,}
}
\end{displaymath}

\noindent where
\begin{displaymath}
\mu_{i}=\gamma_{i} (1+\gamma^{2})= (1- (-1)^{i}\gamma ) (1+\gamma^{2}).
\end{displaymath}



\begin{figure}
\begin{center}
\includegraphics[width=12cm]{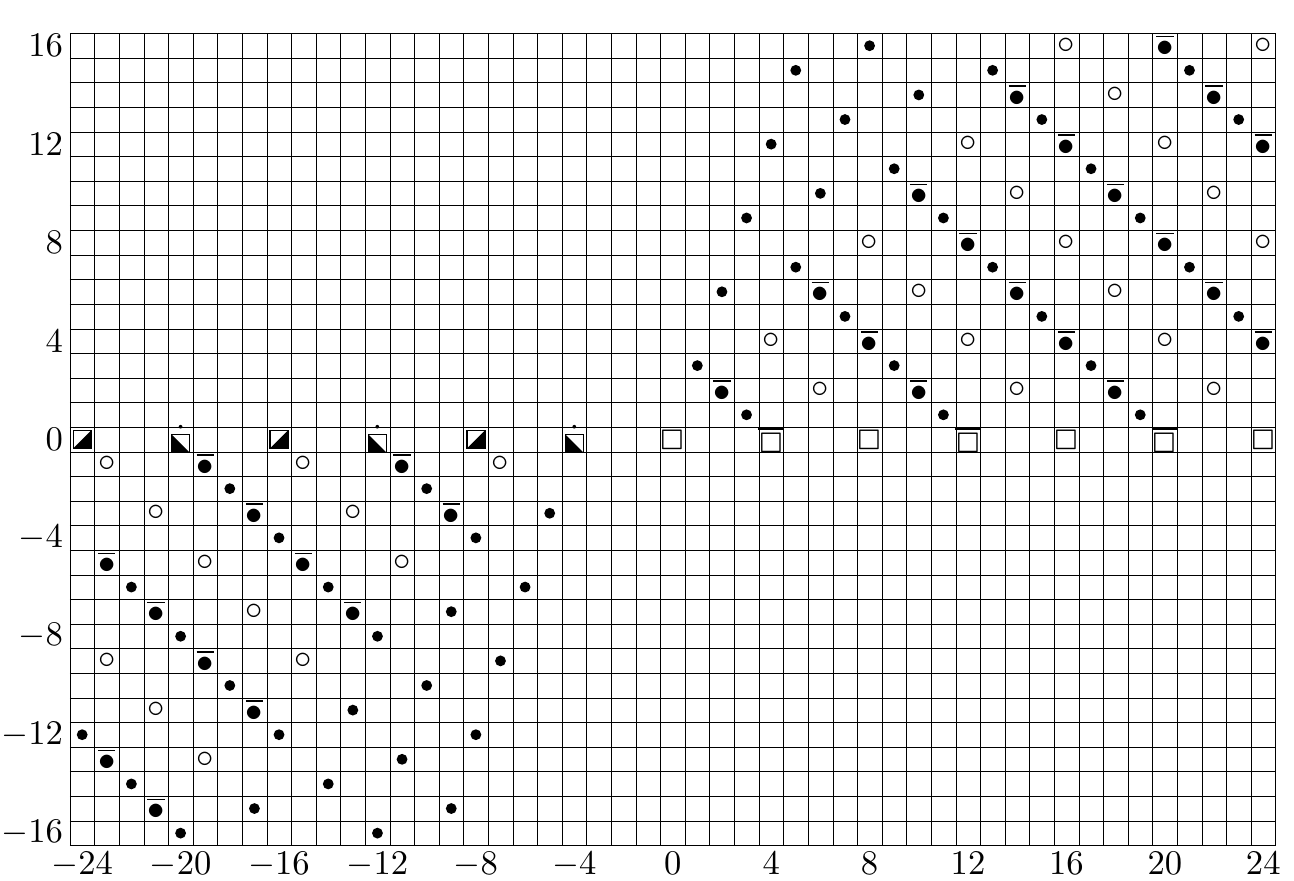} 
\caption[The slice spectral sequence for {${\bigvee_{n
\in \Z} \Sigma ^{n\rho_{4}}\HZZ }$}.]  {The Mackey functor slice
spectral sequence for ${\bigvee_{n \in \Z} \Sigma ^{n\rho_{4}}\HZZ }$.
The symbols are defined in Table \ref{tab-C4Mackey}.  The Mackey
functor at position $(4n-s,s)$ is $\upi_{n (4-\rho_{4})
-s}\HZ=\underline{H}_{4n-s}S^{n\rho_{4}}$.} 
\label{fig-sseq-3}
\end{center}
\end{figure}

The values of $\underline{H}_{*}S^{n\rho_{4}}$  are illustrated in Figure \ref{fig-sseq-3}.  The Mackey
functors in filtration 0 (the horizontal axis) are the ones described
in Proposition \ref{prop-top}.
\bigskip

As in (i), we name some of these elements.  Let $G=C_{4}$ and
$G'=C_{2}\subseteq G$.  Recall that the regular {\rep} $\rho_{4}$ is
$1+\sigma +\lambda$ where $\sigma $ is the sign {\rep} and $\lambda$
is the 2-dimensional {\rep} given by a rotation of order 4.

Note that while Figure \ref{fig-sseq-1} shows all of $\underline{\pi
}_{\star}\HZZ$ for $G=C_{2}$, Figure \ref{fig-sseq-3} shows only a
bigraded portion of this trigraded Mackey functor for $G=C_{4}$,
namely the groups for which the index differs by an integer from a
multiple of $\rho_{4}$.  We will need to refer to some elements not
shown in the latter chart, namely
\begin{numequation}\label{eq-au}
\begin{split}\left\{
\begin{array}[]{rlrlrl}
a_{\sigma } 
    &\hspace{-3mm} \in \underline{H}_{0}S^{\sigma } (G/G) 
        &a_{\lambda }  
            &\in \underline{H}_{0}S^{\lambda  } (G/G) 
                &\overline{a}_{\lambda }
                    &=\Res_{2}^{4}(a_{\lambda })
\\
u_{2\sigma } 
    &\hspace{-3mm} \in \underline{H}_{2}S^{2\sigma   } (G/G) ) 
        &u_{\sigma } 
            &\in \underline{H}_{1}S^{\sigma } (G/G') 
                &\overline{u} _{\sigma }  
                    &=\Res_{1}^{2}(u_{\sigma })
\\
u_{\lambda  } 
    &\hspace{-3mm} \in \underline{H}_{2}S^{\lambda  } (G/G) 
        &\overline{u}_{\lambda }
            &=\Res_{2}^{4}(u_{\lambda })
                &\overline{\overline{u}} _{\lambda }
                    &=\Res_{1}^{4}(u_{\lambda })
\end{array} \right.
\end{split}
\end{numequation}%

\noindent subject to the relations
\begin{numequation}\label{eq-au-rel}
\begin{split}
\left\{\begin{array}[]{rlrlrl}
2a_{\sigma } 
    &\hspace{-2.5mm} = 0
        &\Res_{2}^{4}(a_{\sigma })
            &=0  \\
4a_{\lambda }
    &\hspace{-2.5mm} =0
        &2\overline{a}_{\lambda }  
            &=0 
                &\Res_{1}^{4}(a_{\lambda }) 
                    &=0  \\
\Res_{2}^{4}(u_{2\sigma })
    &\hspace{-2.5mm} = u_{\sigma }^{2}
       &a_{\sigma }^{2}u_{\lambda }
            &=2a_{\lambda }u_{2\sigma } 
                &\mbox{ (gold relation)};\hspace{-1cm}
\end{array} \right.
\end{split}
\end{numequation}%

\noindent see Definition \ref{def-aeu} and Lemma \ref{lem-aeu}.  

We will denote the generator of $\EE_{2}^{s,t} (G/H)$ (when
it is nontrivial) by $x_{t-s,s}$, $y_{t-s,s}$ and $z_{t-s,s}$ for
$H=G$, $G'$ and $\ee$ respectively. Then the generators
for the groups in the 4-slice are
\begin{align*}
y_{4,0} = u_{\rho_{4}}
   =u_{\sigma}\Res^{4}_{2}(u_{\lambda})
  & \in   \upi_{4}\Sigma^{\rho_{4}}\HZZ (G/G')
        = \upi_{3-\sigma-\lambda}\HZZ (G/G')\\
  &   \qquad \mbox{with }\gamma (x_{4,0})=-x_{4,0}\\
x_{3,1} = a_{\sigma}u_{\lambda}
  & \in   \upi_{3}\Sigma^{\rho_{4}}\HZZ (G/G)
        = \upi_{2-\sigma-\lambda}\HZZ (G/G)\\
y_{2,2} = \Res^{4}_{2} (a_{\lambda})u_{\sigma} 
  & \in   \upi_{2}\Sigma^{\rho_{4}}\HZZ (G/G')
        = \upi_{1-\sigma -\lambda}\HZZ (G/G')\\
x_{1,3} = a_{\rho_{4}}
   = a_{\sigma} a_{\lambda}
  & \in   \upi_{1}\Sigma^{\rho_{4}}\HZZ (G/G)
        = \upi_{-\sigma -\lambda}\HZZ (G/G)
\end{align*}

\noindent and the ones for the 8-slice are
\begin{align*}
x_{8,0} = u _{2\lambda+2\sigma}=u_{2\rho_{4}}
  & \in   \upi_{8}\Sigma^{2\rho_{4}}\HZZ (G/G)
        = \upi_{6-2\sigma -2\lambda}\HZZ (G/G)\\
  &   \qquad \mbox{with }y_{4,0}^{2}=y_{8,0}=\Res_{2}^{4}(x_{8,0})\\
x_{6,2} = a_{\lambda} u _{\lambda+2\sigma}
  & \in   \upi_{6}\Sigma^{2\rho_{4}}\HZZ (G/G)
        = \upi_{4-2\sigma -2\lambda}\HZZ (G/G)\\
  &   \qquad \mbox{with }x_{3,1}^{2}=2x_{6,2}\\
  &   \qquad \mbox{and }y_{4,0}y_{2,2}=y_{6,2}=\Res_{2}^{4}(x_{6,2})\\
x_{4,4} = a_{2\lambda}u_{2\sigma}
  & \in   \upi_{4}\Sigma^{2\rho_{4}}\HZZ (G/G)
        = \upi_{2-2\sigma -2\lambda}\HZZ (G/G)\\
  &\qquad \mbox{with }y_{2,2}^{2}=y_{4,4}=\Res_{2}^{4}(x_{4,4})\\
  &\qquad \mbox{and }
x_{1,3}x_{3,1}=2x_{4,4}\\
x_{2,6} = x_{1,3}^{2}
  & \in   \upi_{2}\Sigma^{2\rho_{4}}\HZZ (G/G)
        = \upi_{-2\sigma -2\lambda}\HZZ (G/G).
\end{align*}

These elements and their restrictions generate
$\upi_{*}\Sigma^{m\rho_{4}}\HZZ$ for $m=1$ and 2.  For
$m>2$ the groups are generated by products of these elements. 

The element
\begin{displaymath}
z_{4,0} = \Res^{2}_{1}(y_{4,0}) =\Res^{2}_{1} (u_{\rho_{4}} )\in 
\upi_{4}\Sigma^{\rho_{4}}\HZZ (G/\ee)
\end{displaymath}

\noindent is invertible with $\gamma (y_{4,0})=-y_{4,0}$,
$z_{4,0}^{2}=z_{8,0}=\Res_{1}^{4}(x_{8,0})$ and
\begin{displaymath}
z_{-4m,0}:= z_{4,0}^{-m}=e_{m\rho_{4}}
       \in \upi_{-4m}\Sigma^{-m\rho_{4}}\HZZ (G/\ee)
\qquad \mbox{for }m>0, 
\end{displaymath}

\noindent where $e_{m\rho_{4}}$ is as in Definition
\ref{def-aeu}. These elements and their transfers generate the groups
in
\begin{displaymath}
\upi_{-4m}\Sigma^{-m\rho_{4}}\HZZ \qquad \mbox{for }m>0. 
\end{displaymath}

\begin{thm}\label{thm-div}
{\bf Divisibilities in the negative regular slices for $C_{4}$.}
There are the following infinite divisibilities in the third quadrant
of the {\SS} in Figure \ref{fig-sseq-3}.

\begin{enumerate}

\item [(i)] \label{div:i} 
$x_{-4,0}=\Tr_{1}^{4}(z_{-4,0})$ is
divisible by any monomial in $x_{1,3}$ and $x_{4,4}$, meaning that
\begin{displaymath}
x_{1,3}^{i}x_{4,4}^{j}x_{-4-4j-i,-4j-3k}=x_{-4,0}
\qquad \mbox{for }i,j\geq 0. 
\end{displaymath}

\noindent Moreover, no other basis element killed by $x_{3,1}$ and
$x_{4,4}$ has this property.
 
\item [(ii)] \label{div:ii} 
$x_{-4,0}$, and $x_{-7,-1}$ are divisible
by any monomial in $x_{4,4}$, $x_{6,2}$ and $x_{8,0}$, subject to the
relation $x_{6,2}^{2}=x_{8,0}x_{4,4}$.  Note here that
$x_{3,1}^{2}=2x_{6,2}$.

\noindent Moreover, no other basis element killed by $x_{4,4}$,
$x_{6,2}$ and $x_{8,0}$ has this property.

\item [(iii)] \label{div:iv} 
$y_{-7,-1}=\Res_{2}^{4}(x_{-7,-1})$ is
divisible by any monomial in $y_{2,2}$ and $y_{4,0}$, meaning that
\begin{displaymath}
y_{2,2}^{j}y_{4,0}^{k}y_{-7-2j-4k,-1-2j}=y_{-7,-1}
\qquad \mbox{for }j,k\geq 0. 
\end{displaymath}

\noindent Moreover, no other basis element killed by $y_{2,2}$,and
$y_{4,0}$ has this property.

\end{enumerate}
\end{thm}

We will prove Theorem \ref{thm-div} as a corollary of a more general
statement (Lemma \ref{lem-div} and Corollary \ref{cor-div}) in which
we consider all {\rep}s of the form $m\lambda +n\sigma $ for $m,n\geq
0$.  Let
\begin{displaymath}
\underline{R}=\bigoplus_{m,n\geq 0}\underline{H}_{*}S^{m\lambda +n\sigma }. 
\end{displaymath}

\noindent It is generated by the elements of (\ref{eq-au}) subject to
the relations of (\ref{eq-au-rel}).

In the larger ring 
\begin{displaymath}
\underline{\tilde{R}}
     =\bigoplus_{m,n\in \Z \atop mn\geq 0}
          \underline{H}_{*}S^{m\lambda +n\sigma }, 
\end{displaymath}

\noindent the elements $u_{\sigma }$, $\overline{u}_{\sigma } $ and
$\overline{\overline{u} }_{\lambda } $ are invertible with
\begin{align*}
e_{\sigma }
  & = u_{\sigma }^{-1} \in \underline{H}_{-1}S^{-\sigma } (G/G')   
       &e_{\lambda }
     = \overline{\overline{u} }_{\lambda }^{-1}
  \in \underline{H}_{-2}S^{-\lambda } (G/e).
\end{align*}

Define spectra $L_{m}$ and $K_{n}$ to be the cofibers of $a_{m\lambda
}$ and $a_{n\sigma }$.  Thus we have cofiber sequences
\begin{displaymath}
\xymatrix
@R=2mm
@C=15mm
{
{\Sigma^{-1}L_{m}}\ar[r]^(.5){c_{m\lambda }}
    &{S^{0}}\ar[r]^(.5){a_{m\lambda }}
        &{S^{m\lambda }}\ar[r]^(.5){b_{m\lambda }}
            &{L_{m}}\\
{\Sigma^{-1}K_{n}}\ar[r]^(.5){c_{n\sigma}}
    &{S^{0}}\ar[r]^(.5){a_{n\sigma}}
        &{S^{n\sigma}}\ar[r]^(.5){b_{n\sigma}}
            &{K_{n}}
}
\end{displaymath}

\noindent Dualizing gives
\begin{displaymath}
\xymatrix
@R=2mm
@C=15mm
{
{DL_{m}}\ar[r]^(.5){Db_{m\lambda }}
    &{S^{-m\lambda }}\ar[r]^(.5){Da_{m\lambda }}
        &{S^{0}}\ar[r]^(.5){Dc_{m\lambda }}
            &{\Sigma DL_{m}}\\
{DK_{n}}\ar[r]^(.5){Db_{n\sigma}}
    &{S^{-n\sigma }}\ar[r]^(.5){Da_{n\sigma}}
        &{S^{0}}\ar[r]^(.5){Dc_{n\sigma}}
            &{\Sigma  DK_{n}}
}
\end{displaymath}

\noindent The maps $Da_{m\lambda }$ and $Da_{n\sigma }$ are the same
as desuspensions of $a_{m\lambda }$ and $a_{n\sigma }$, which implies that 
\begin{displaymath}
DL_{m}=\Sigma^{-1-m\lambda }L_{m}
\qquad \aand 
DK_{n}=\Sigma^{-1-n\sigma }K_{n}.
\end{displaymath}

\noindent Inspection of the cellular chain complexes for $L_{m}$ and
$K_{n}$ and certain of their suspensions reveals that
\begin{displaymath}
\Sigma^{2-\lambda }L_{m}\wedge \HZZ
=L_{m}\wedge \HZZ=\Sigma^{2-2\sigma }L_{m}\wedge \HZZ
\end{displaymath}

\noindent and
\begin{displaymath}
\Sigma^{2-2\sigma }K_{n}\wedge \HZZ=K_{n}\wedge \HZZ,
\end{displaymath}

\noindent while $\Sigma^{1-\sigma }$ alters both $L_{m}\wedge \HZZ$
and $K_{n}\wedge \HZZ$.  We will denote $\Sigma^{k (1-\sigma
)}L_{m}\wedge \HZZ$ by $L_{m}^{(-1)^{k}}\wedge \HZZ$ and similarly for
$K_{n}$.
   
\begin{figure}
\begin{center}
\includegraphics[width=5cm]{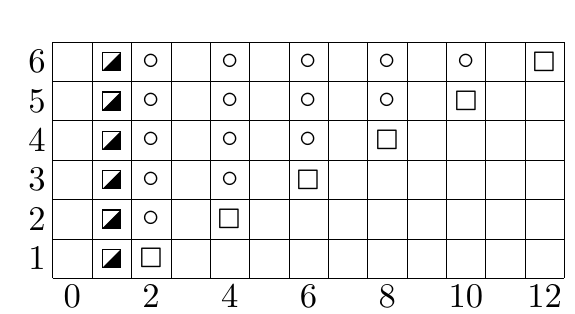}
\includegraphics[width=5cm]{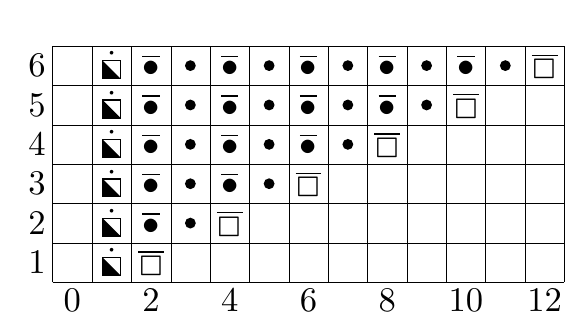}
\caption[Charts for {${\underline{\sl H}_{i}L_{m}^{\pm }}$}.]{Charts for
$\underline{H}_{i}L_{m}^{\pm }$.  The horizontal coordinate is $i$ and
the vertical one is $m$. $L_{m}$ is on the left and $L_{m}^{-}$ is on
the right.}  \label{fig-Lm}
\end{center}
\end{figure}

\begin{figure}
\begin{center}
\includegraphics[width=5cm]{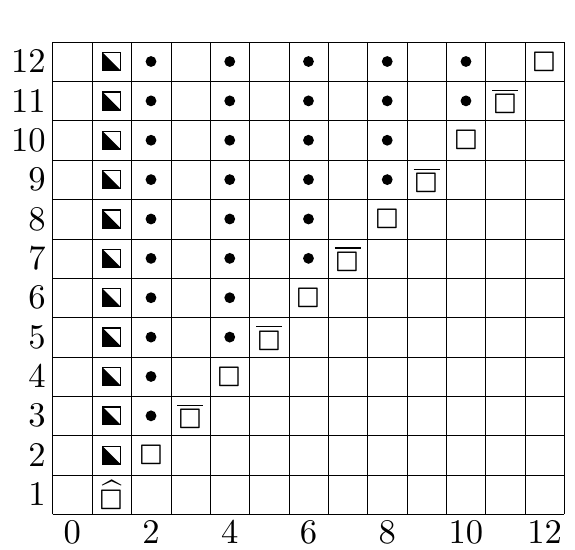}
\includegraphics[width=5cm]{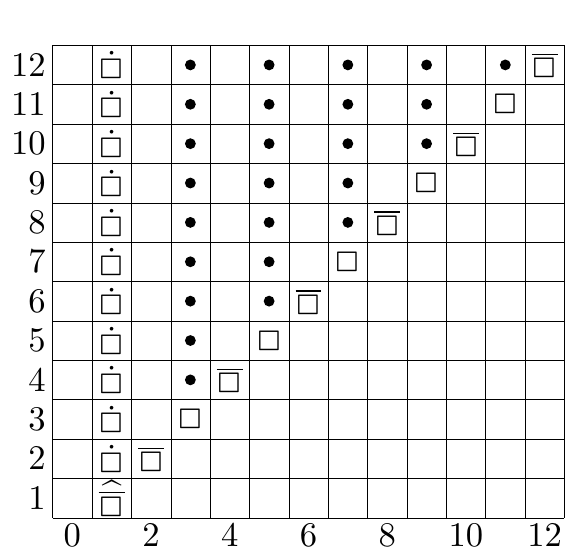} 
\caption[Charts for {${\underline{\sl H}_{i}K_{n}^{\pm }}$}.] {Charts for
$\underline{H}_{i}K_{n}^{\pm }$.  The horizontal coordinate is $i$ and
the vertical one is $n$. $K_{n}$ is on the left and $K_{n}^{-}$ is on
the right.}  \label{fig-Kn}
\end{center}
\end{figure}

\begin{figure}
\begin{center}
\includegraphics[height=4cm]{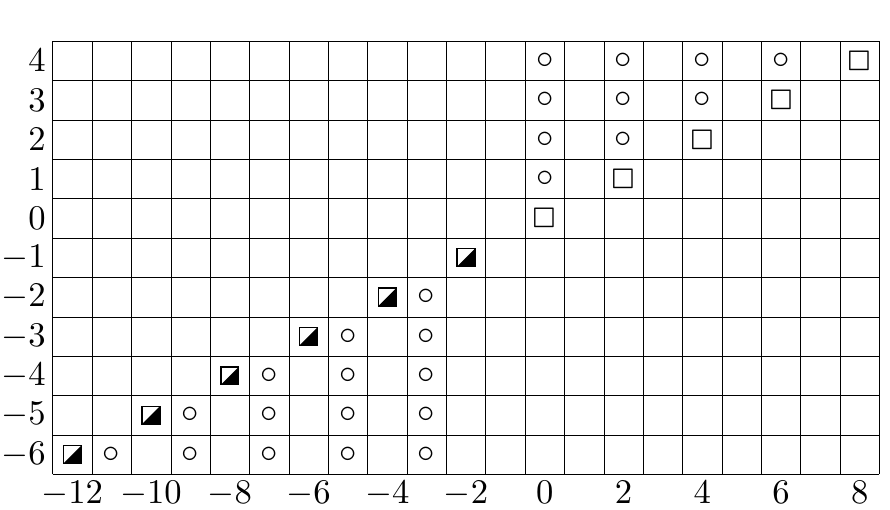}
\includegraphics[height=4cm]{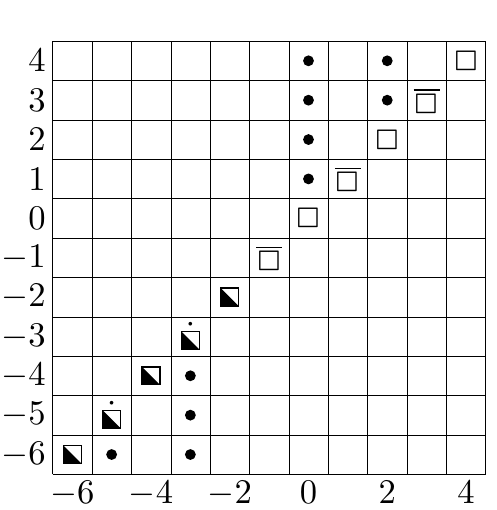} \caption[Charts for
${{\underline{\sl H}_{i}S^{m\lambda }}$} and {${\underline{\sl
H}_{i}S^{n\sigma }}$}.] {Charts for ${\underline{H}_{i}S^{m\lambda }}$
and ${\underline{H}_{i}S^{n\sigma }}$. The horizontal coordinates are
$i$ and the vertical ones are $m$ and $n$. $S^{m\lambda }$ is on the
left and $S^{n\sigma }$ is on the right.  }  \label{fig-spheres}
\end{center}
\end{figure}

The homology groups of $L_{m}^{\pm }$ and $K_{m}^{\pm }$ for $m,n>0$
are indicated in Figures \ref{fig-Lm} and \ref{fig-Kn}, and those for
$S^{m\lambda }$ and $S^{n\sigma }$ are shown in
Figure\ref{fig-spheres}.
    
In the following diagrams we will use the same notation for a map and
its smash product with any identity map.  Let $V=m\lambda +n\sigma $
with $m,n>0$, and let $R_{V}$ denote the fiber of $a_{V}$.  Since
$a_{V}$ is self-dual up to susupension, we have
$DR_{V}=\Sigma^{-1-V}R_{V}$.  In the following each row and column is
a cofiber sequence.
\begin{numequation}\label{eq-newfirst}
\begin{split}
\xymatrix
@R=5mm
@C=10mm
{
{}
    &{}
        &{\Sigma^{n\sigma -1}L_{m}}\ar@{=}[r]\ar[d]^(.5){c_{m\lambda }}
            &{\Sigma^{n\sigma -1}L_{m}}\ar[d]^(.5){}\\
{\Sigma^{-1}K_{n}}\ar[r]^(.5){c_{n\sigma }}\ar[d]^(.5){}
    &{S^{0}}\ar[r]^(.5){a_{n\sigma }}\ar@{=}[d]^(.5){}
        &{S^{n\sigma }}\ar[r]^(.5){b_{n\sigma }}\ar[d]^(.5){a_{m\lambda }}
            &{K_{n}}\ar[d]^(.5){}\\
{\Sigma^{-1}R_{V}}\ar[r]^(.5){c_{V}}
    &{S^{0}}\ar[r]^(.5){a_{V }}
        &{S^{V}}\ar[r]^(.5){b_{V}}\ar[d]^(.5){b_{m\lambda }}
            &{R_{V}}\ar[d]^(.5){}\\
{}
    &{}
        &{\Sigma^{n\sigma }L_{m}}\ar@{=}[r]
            &{\Sigma^{n\sigma }L_{m}}
}
\end{split}
\end{numequation}%

\noindent The homology sequence for the third column is the easiest
way to compute $\underline{H}_{*}S^{V}$.  That column is
\begin{numequation}\label{eq-third-column}
\begin{split}
\xymatrix
@R=5mm
@C=10mm
{
{\Sigma^{n\sigma -1}L_{m}}\ar[r]^(.6){c_{m\lambda }}
    &{S^{n\sigma }}\ar[r]^(.5){a_{m\lambda }}
        &{S^{V}}\ar[r]^(.4){b_{m\lambda }}
            &{\Sigma^{n\sigma }L_{m}},
}
\end{split}
\end{numequation}%

\noindent which dualizes to 
\begin{displaymath}
\xymatrix
@R=2mm
@C=10mm
{
{\Sigma^{1-n\sigma }DL_{m}}\ar@{=}[d]^(.5){}
    &{S^{-n\sigma }}\ar[l]_(.4){c_{m\lambda }}
        &{S^{-V}}\ar[l]_(.5){a_{m\lambda }}
            &{\Sigma^{-n\sigma }DL_{m}}\ar[l]_(.6){c_{m\lambda }}
                                       \ar@{=}[d]^(.5){}\\
{\Sigma^{-V}L_{m}}
    &   &   &{\Sigma^{-1-V}L_{m}}
}
\end{displaymath}

\noindent or 
\begin{numequation}\label{eq-dual-third-column}
\begin{split}
\xymatrix
@R=5mm
@C=10mm
{
{\Sigma^{ -1-V}L_{m}}\ar[r]^(.6){c_{m\lambda }}
    &{S^{-V}}\ar[r]^(.5){a_{m\lambda }}
        &{S^{-n\sigma }}\ar[r]^(.4){b_{m\lambda }}
            &{\Sigma^{-V}L_{m}}.
}
\end{split}
\end{numequation}%

For (\ref{eq-third-column}) the {\LES} in homology includes
\begin{displaymath}
\xymatrix
@R=5mm
@C=8mm
{
{\underline{H}_{i+1-n}L_{m}^{(-1)^{n}}}\ar[r]^(.6){c_{m\lambda }}
    &{\underline{H}_{i}S^{n\sigma }}\ar[r]^(.5){a_{m\lambda }}
        &{\underline{H}_{i}S^{V }}\ar[r]^(.4){b_{m\lambda }}
            &{\underline{H}_{i-n}L_{m}^{(-1)^{n}}}\ar[r]^(.5){c_{m\lambda }}
                &{\underline{H}_{i-1}S^{n\sigma }}
}
\end{displaymath}

{\bf Divisibility by $a_{\lambda }$.} Multiplication by $a_{\lambda }$ keads to 
\begin{displaymath}
\xymatrix
@R=5mm
@C=8mm
{
{\underline{H}_{i+1-n}L_{m}^{(-1)^{n}}}\ar[r]^(.6){c_{m\lambda }}
                                  \ar[d]^(.5){a'_{\lambda }}
    &{\underline{H}_{i}S^{n\sigma }}\ar[r]^(.5){a_{m\lambda }}\ar@{=}[d]^(.5){}
        &{\underline{H}_{i}S^{V }}\ar[r]^(.4){b_{m\lambda }}
                                  \ar[d]^(.5){a_{\lambda }}
            &{\underline{H}_{i-n}L_{m}^{(-1)^{n}}}\ar[r]^(.5){c_{m\lambda }}
                                  \ar[d]^(.5){a'_{\lambda }}
                &{\underline{H}_{i-1}S^{n\sigma }}\ar@{=}[d]^(.5){}\\
{\underline{H}_{i+1-n}L_{m'}^{(-1)^{n}}}\ar[r]^(.6){c_{m'\lambda }}
    &{\underline{H}_{i}S^{n\sigma }}\ar[r]^(.5){a_{m'\lambda }}
        &{\underline{H}_{i}S^{V+\lambda  }}\ar[r]^(.4){b_{m'\lambda }}
            &{\underline{H}_{i-n}L_{m'}^{(-1)^{n}}}\ar[r]^(.5){c_{m'\lambda }}
                &{\underline{H}_{i-1}S^{n\sigma },}
}
\end{displaymath}

\noindent where $m'=m+1$ and $a'_{\lambda }$ is induced by the
inclusion $L_{m}\to L_{m'}$.  

In the dual case we get 
\begin{numequation}\label{eq-a-lambda}
\begin{split}
\xymatrix
@R=5mm
@C=2mm
{
{\underline{H}_{i+1}S^{-n\sigma }}\ar[r]^(.4){b}\ar@{=}[d]^(.5){}
    &{\underline{H}_{i+1+|V|}L_{m}^{(-1)^{n}}}\ar[r]^(.6){c}
        &{\underline{H}_{i}S^{-V }}\ar[r]^(.5){a}
            &{\underline{H}_{i}S^{-n\sigma }}\ar[r]^(.4){b}
                                               \ar@{=}[d]^(.5){}
                &{\underline{H}_{i+|V|}L_{m}^{(-1)^{n}}}
\\
{\underline{H}_{i+1}S^{-n\sigma }}\ar[r]^(.4){b}
    &{\underline{H}_{i+3+|V|}L_{m'}^{(-1)^{n}}}\ar[r]^(.6){c}
                                                \ar[u]_(.5){Da'_{\lambda }}
        &{\underline{H}_{i}S^{-V-\lambda  }}\ar[r]^(.5){a}
                                            \ar[u]_(.5){a_{\lambda }}
            &{\underline{H}_{i}S^{-n\sigma }}\ar[r]^(.4){b}
                &{\underline{H}_{i+2+|V|}L_{m'}^{(-1)^{n}}}
                                            \ar[u]_(.5){Da'_{\lambda }}
}
\end{split}
\end{numequation}%

\noindent Here the subscripts on the horizontal maps ($m\lambda $ in
the top row and $m'\lambda $ in the bottom row) have been omitted
to save space. The five lemma implies that the middle
vertical map is onto when the left hand $Da'_{\lambda }$ is onto
and the right hand one is one to one.  
The left version of $Da'_{\lambda }$ is onto in every
case except $i=-|V|$ and the right version of it is one to one in all
cases except $i=-|V|$ and $i=-1-|V|$.  This is illustrated for small
$m$ in the following diagram in which trivial Mackey functors are
indicated by blank spaces.
\begin{displaymath}
\xymatrix
@R=-1mm
@C=1mm
{
{j}
    &{\underline{H}_{j}L_{1}}
        &{\underline{H}_{j}L_{2}}
            &{\underline{H}_{j}L_{3}}
                &{\underline{H}_{j}L_{4}}
                    &
    &{\underline{H}_{j}L_{1}^{-}}
        &{\underline{H}_{j}L_{2}^{-}}
            &{\underline{H}_{j}L_{3}^{-}}
                &{\underline{H}_{j}L_{4}^{-}}
\\
-1  &   &   &   &   &
    &   &   &   &\\
0   &   &   &   &   &
    &   &   &   &\\
1   
    &{\fourbox }
        &{\fourbox }\ar[luu]^(.5){}
            &{\fourbox }\ar[luu]^(.5){}
                &{\fourbox }\ar[luu]^(.5){}
                    &
    &{\dot{\twobox}}
        &{\dot{\twobox}}\ar[luu]^(.5){}
            &{\dot{\twobox}}\ar[luu]^(.5){}
                &{\dot{\twobox}}\ar[luu]^(.5){}\\
2   
    &{\Box}
        &{\circ} \ar[luu]^(.5){}
            &{\circ}\ar[luu]^(.5){}
                &{\circ} \ar[luu]^(.5){}
                    &
    &{\oBox}
        &{\obull}\ar[luu]^(.5){}
            &{\obull}\ar[luu]^(.5){}
                &{\obull}\ar[luu]^(.5){}\\
3   
    &   &   &   &   &
    &   &{\bullet}\ar[luu]^(.5){}
            &{\bullet}\ar[luu]^(.5){}
                &{\bullet}\ar[luu]^(.5){}\\
4   
    &   &{\Box}\ar@{=}[luu]^(.5){}
            &{\circ} \ar@{=}[luu]^(.5){}
                &{\circ}\ar@{=}[luu]^(.5){}
                    &
    &   &{\oBox}\ar@{=}[luu]^(.5){}
            &{\obull}\ar@{=}[luu]^(.5){}
                &{\obull}\ar@{=}[luu]^(.5){}\\
5   
    &   &   &   &   &
    &   &   &{\bullet}\ar@{=}[luu]^(.5){}
                &{\bullet}\ar@{=}[luu]^(.5){}\\
6   
    &   &   &{\Box} \ar@{=}[luu]^(.5){}
                &{\circ}\ar@{=}[luu]^(.5){}
                    &
    &   &   &{\oBox}\ar@{=}[luu]^(.5){}
                &{\obull}\ar@{=}[luu]^(.5){}\\ 
7   
    &   &   &   &   &
    &   &   &   &{\bullet}\ar@{=}[luu]^(.5){}\\
8   
    &   &   &   &{\Box} \ar@{=}[luu]^(.5){}
                    &
    &   &   &   &{\oBox}\ar@{=}[luu]^(.5){}\\
}
\end{displaymath}

\noindent It follows that the map $a_{\lambda }$ in
(\ref{eq-a-lambda}) is onto for all $i$ except $-|V|$.  {\em This is a
divisibility result.}  Note that $a_{\lambda }$ is trivial on
$\underline{H}_{*}X (G/e)$ for any $X$ since $\Res_{1}^{4}(a_{\lambda
})=0$.

\bigskip
{\bf Divisibility by $u_{\lambda }$.} For $u_{\lambda }$
multiplication we use the diagram
\begin{numequation}\label{eq-u-lambda}
\begin{split}
\xymatrix
@R=5mm
@C=3mm
{
{\underline{H}_{i+1}S^{-n\sigma }}\ar[r]^(.5){b}
    &{\underline{H}_{i+1}L_{m}^{(-1)^{n}}}\ar[r]^(.5){c}
        &{\underline{H}_{i}S^{-V}}\ar[r]^(.5){a}
            &{\underline{H}_{i}S^{-n\sigma }}\ar[r]^(.5){b}
                &{\underline{H}_{i}L_{m}^{(-1)^{n}}}
\\
{\underline{H}_{i-1}S^{-n\sigma-\lambda}}\ar[r]^(.5){b}
                                           \ar[u]^(.5){u_{\lambda }}
    &{\underline{H}_{i+1}L_{m}^{(-1)^{n}}}\ar[r]^(.5){c}
                                           \ar@{=}[u]^(.5){}
        &{\underline{H}_{i-2}S^{-V-\lambda }}\ar[r]^(.5){a}
                                           \ar[u]^(.5){u_{\lambda }}
            &{\underline{H}_{i-2}S^{-n\sigma-\lambda}}\ar[r]^(.5){b}
                                           \ar[u]^(.5){u_{\lambda }}
                &{\underline{H}_{i}L_{m}^{(-1)^{n}}}\ar@{=}[u]^(.5){}
}
\end{split}
\end{numequation}%

\noindent The rightmost $u_{\lambda }$ is onto in all cases except
$i=-n$ and $n$ even.  This is illustrated for $n=6$ and 7 in the
following diagram.
\begin{displaymath}
\xymatrix
@R=-1mm
@C=2mm
{
j   &-1 &-2 &-3 &-4 &-5 &-6 &-7 &-8 &-9\\
{\underline{H}_{j}S^{-6\sigma }}
    &   &   &{\bullet}
                &{}
                    &{\bullet}
                        &{\twobox }
\\
{\underline{H}_{j}S^{-6\sigma-\lambda  }}
    &   &   &{\bullet}\ar[llu]^(.5){}
                &{}
                    &{\bullet}\ar@{=}[llu]^(.5){}
                        &{} &{\circ}\ar[llu]^(.5){}
                                &{\fourbox }\ar[llu]^(.5){}
\\
{}\\
{\underline{H}_{j}S^{-7\sigma }}
    &   &   &{\bullet} 
                &{}
                    &{\bullet}
                        &   &{\dot{\twobox }}
\\
{\underline{H}_{j}S^{-7\sigma-\lambda  }}
    &   &   &{\bullet}\ar[llu]^(.5){}
                &{}
                    &{\bullet}\ar@{=}[llu]^(.5){}
                        &{} &{\bullet}\ar@{=}[llu]^(.5){}
                                &   &{\dot{\twobox }}\ar@{=}[llu]^(.5){}
\\
}
\end{displaymath}

\noindent

Thus the central $u_{\lambda }$ in (\ref{eq-u-lambda}) fails to be
onto only in when $i=-n$ and $n$ is even.

\bigskip {\bf Divisibility by $a_{\sigma }$.}  The corresponding
diagram is
\begin{displaymath}
\xymatrix
@R=5mm
@C=3mm
{
{\underline{H}_{i+1}S^{-n\sigma }}\ar[r]^(.4){b}
    &{\underline{H}_{i+1+|V|}L_{m}^{(-1)^{n}}}\ar[r]^(.6){c}
        &{\underline{H}_{i}S^{-V }}\ar[r]^(.5){a}
            &{\underline{H}_{i}S^{-n\sigma }}\ar[r]^(.4){b}
                &{\underline{H}_{i+|V|}L_{m}^{(-1)^{n}}}
\\
{\underline{H}_{i+1}S^{- n'\sigma }}\ar[r]^(.4){b}
                                            \ar[u]_(.5){a_{\sigma }}
    &{\underline{H}_{i+2+|V|}L_{m}^{(-1)^{n'}}}\ar[r]^(.6){c}
                                            \ar[u]_(.5){a_{\sigma}}
        &{\underline{H}_{i}S^{-V-\sigma   }}\ar[r]^(.5){a}
                                            \ar[u]_(.5){a_{\sigma }}
            &{\underline{H}_{i}S^{- n'\sigma }}\ar[r]^(.35){b}
                                            \ar[u]_(.5){a_{\sigma }}
                &{\underline{H}_{i+1+|V|}L_{m}^{(-1)^{n'}}}
                                            \ar[u]_(.5){a_{\sigma}}
}
\end{displaymath}


\noindent Here we have abbreviated ${n+1}$ by $n'$.  Since
$\Res_{2}^{4}(a_{\sigma })=0$, the map $a_{\sigma }$ must vanish on
$\underline{H}_{*}X (G/G')$ and $\underline{H}_{*}X (G/e)$. It can be
nontrivial only on $G/G$.  

By Lemma \ref{lem-hate}, the image of $a_{\sigma }$ is the kernel of
the restriction map $u_{\sigma }^{-1}\Res_{2}^{4}$ and the kernel of
$a_{\sigma }$ is the image of the transfer $\Tr_{2}^{4}$.  From Figure
\ref{fig-spheres} we see that $\Res_{2}^{4}$ kills
$\underline{H}_{i}S^{-n\sigma } (G/G)$ except the case $i=-n$ for even
$n$.  From Figure \ref{fig-Lm} we see that it kills
$\underline{H}_{j}L_{m}^{-} (G/G)$ for all $j$ and
$\underline{H}_{j}L_{m}(G/G)$ for odd $j>1$, but not the generators
for $j=1$ nor the ones for even values of $j$ from $2$ to $2m$.  The
transfer has nontrivial image in $\underline{H}_{j}L_{m}^{-}$ only for
$j=1$ and in $\underline{H}_{j}L_{m}$ only for $j=1$ and for even $j$
from 2 to $2m$.

{\em It follows that for odd $n$, each element of
$\underline{H}_{i}S^{-V} (G/G)$ is divisible by $a_{\sigma }$ except
when $i=-|V|=-2m-n$.  For even $n$ it is onto except when $i=-n$,
$i=-n-2m$, and $i$ odd from $1-n-2m$ to $-1-n$.  }

\bigskip

{\bf Divisibility by $u_{2\sigma }$.} For $u_{2\sigma }$
multiplication, the diagram is
\begin{displaymath}
\xymatrix
@R=5mm
@C=4mm
{
{\underline{H}_{i+1}S^{-n\sigma }}\ar[r]^(.5){b}
    &{\underline{H}_{i+1}L_{m}^{(-1)^{n}}}\ar[r]^(.5){c}
        &{\underline{H}_{i}S^{-V}}\ar[r]^(.5){a}
            &{\underline{H}_{i}S^{-n\sigma }}\ar[r]^(.5){b}
                &{\underline{H}_{i}L_{m}^{(-1)^{n}}}
\\
{\underline{H}_{i-1}S^{- (n+2)\sigma}}\ar[r]^(.5){b}
                                           \ar[u]^(.5){u_{2\sigma }}
    &{\underline{H}_{i+1}L_{m}^{(-1)^{n}}}\ar[r]^(.5){c}
                                           \ar@{=}[u]^(.5){}
        &{\underline{H}_{i-2}S^{-V-2\sigma  }}\ar[r]^(.5){a}
                                           \ar[u]^(.5){u_{2\sigma }}
            &{\underline{H}_{i-2}S^{- (n+2)\sigma}}\ar[r]^(.55){b}
                                           \ar[u]^(.5){u_{2\sigma }}
                &{\underline{H}_{i}L_{m}^{(-1)^{n}}}\ar@{=}[u]^(.5){}
}
\end{displaymath}

\noindent The rightmost $u_{2\sigma }$ is onto in all cases, so {\em every
element in $\underline{H}_{*}S^{-V}$ is divisible by $u_{2\sigma }$.}

The arguments above prove the following.

\begin{lem}\label{lem-div}
{\bf $RO (G)$-graded divisibility.}  Let $G=C_{4}$ and $V=m\lambda
+n\sigma $ for $m,n\geq 0$.  

\begin{enumerate}
\item [(i)] Each element in
$\underline{H}_{i}S^{-V} (G/G)$ or $\underline{H}_{i}S^{-V} (G/G')$ is
divisible by $a_{\lambda }$ or $\overline{a}_{\lambda }$ except when
$i=-|V|$.

\item [(ii)] Each element in
$\underline{H}_{i}S^{-V} (G/H)$  is
divisible by a suitable restriction of $u_{\lambda }$ except when
$i=-n$ for even $n$.

\item [(iii)] Each element in $\underline{H}_{i}S^{-V} (G/G)$ for odd
$n$ is divisible by $a_{\sigma }$ except when $i=-|V|$.  For even $n$ it is
divisible by $a_{\sigma }$ except when $i=-n$, $i=-|V|$ and $i$ is odd
from $i=1-|V|$ to $-1-n$.

\item [(iv)] Each element in $\underline{H}_{i}S^{-V} (G/H)$ is
divisible by a $u_{2\sigma }$, $u_{\sigma }$ or $\overline{u}_{\sigma } $.
\end{enumerate}
\end{lem}

In Theorem \ref{thm-div} we are looking for
divisibility by 
\begin{numequation}\label{eq-divisors}
\begin{split}\left\{ 
\begin{array}[]{rl}
x_{1,3}
  =  a_{\sigma }a_{\lambda }&\in 
      \underline{H}_{0}S^{\sigma +\lambda } (G/G)
    = \underline{H}_{1}S^{\rho } (G/G)\\ 
x_{4,4}
  = a_{\lambda}^{2}u_{2\sigma } &\in 
       \underline{H}_{2}S^{2\lambda +2\sigma } (G/G)
    = \underline{H}_{4}S^{2\rho } (G/G)\\ 
y_{2,2}
  = \overline{a} _{\lambda}u_{\sigma } &\in 
       \underline{H}_{1}S^{1\lambda +1\sigma } (G/G')
    = \underline{H}_{2}S^{\rho } (G/G)\\ 
x_{6,2}
  = a_{\lambda }u_{2\sigma }u_{\lambda }&\in 
  \underline{H}_{4}S^{2\lambda +2\sigma } (G/G)
   =\underline{H}_{6}S^{2\rho } (G/G) \\ 
x_{8,0}
  = u_{2\sigma }u_{\lambda }^{2}&\in 
  \underline{H}_{6}S^{2\lambda +2\sigma } (G/G)
   =\underline{H}_{8}S^{2\rho } (G/G) \\ 
y_{4,0}
  = u_{\sigma }\overline{u} _{\lambda }&\in 
  \underline{H}_{3}S^{\lambda +\sigma } (G/G')
   =\underline{H}_{4}S^{\rho } (G/G') 
\end{array}\right.
\end{split}
\end{numequation}%

\noindent In view of Lemma \ref{lem-div}(iv), we can ignore the
factors $u_{2\sigma }$ and $u_{\sigma }$ when analyzing such
divisibility.

\begin{cor}\label{cor-div}
{\bf Infinite divisibility by the divisors of (\ref{eq-divisors}).}
Let
\begin{displaymath}
V=m\lambda +n\sigma \qquad \mbox{for }m,n\geq 0. 
\end{displaymath}

\noindent Then
\begin{itemize}
\item [$\bullet$] Each element of $\underline{H}_{i}S^{-V} (G/G)$  is
infinitely divisible by ${x_{1,3}=a_{\sigma }a_{\lambda }}$ for
${i>-n}$ when $n$ is even and for $i\geq -n$ when $n$ is odd.

\item [$\bullet$] Each element of $\underline{H}_{i}S^{-V} (G/G)$ is
infinitely divisible by ${x_{4,4}=a_{\lambda }^{2}u_{2\sigma }}$ for
$i>-|V|$.

\item [$\bullet$] Each element of $\underline{H}_{i}S^{-V} (G/G')$ is
infinitely divisible by ${y_{2,2}=\overline{a} _{\lambda }u_{\sigma
}}$ for $i>-|V|$.

\item [$\bullet$] Each element of $\underline{H}_{i}S^{-V} (G/G)$ is
infinitely divisible by ${x_{6,2}=a_{\lambda }u_{2\sigma }u_{\lambda
}}$ for $i>-|V|$ when $n$ is odd and for $-|V|<i<-n$ when $n$ is even.

\item [$\bullet$] Each element of $\underline{H}_{i}S^{-V} (G/G)$ is
infinitely divisible by ${x_{8,0}=u_{2\sigma }u_{\lambda }^{2}}$ for
$i<-n$ when $n$ is even and for all $i$ when $n$ is odd.

\item [$\bullet$] Each element of $\underline{H}_{i}S^{-V} (G/G')$ is
infinitely divisible by ${y_{4,0}=u_{\sigma }\overline{u} _{\lambda
}}$ for $i<-n$ when $n$ is even and for all $i$ when $n$ is odd.
\end{itemize}
\end{cor}

This implies Theorem \ref{thm-div}.

\section{The spectra $\kR$ and $\kH$}\label{sec-kH}


Before defining our spectrum we need to recall some definitions and
formulas from \cite{HHR}.  Let $H\subset G$ be finite groups.  In
\cite[\S2.2.3]{HHR} we define a norm functor $N_{H}^{G}$ from the
category of $H$-spectra to that of $G$-spectra.  Roughly speaking, for
an $H$-spectrum $X$, $N_{H}^{G}X$ is the $G$-spectrum underlain by the
smash power $X^{(|G/H|)}$ with $G$ permuting the factors and $H$ leaving
each one invariant.  When $G$ is cyclic, we will denote the orders of
$G$ and $H$ by $g$ and $h$, and the norm functor by $N_{h}^{g}$.

There is a $C_{2}$-spectrum $MU_{\reals}$ underlain by the complex
cobordism spectrum $MU$ with group action given by complex
conjugation.  Its construction is spelled out in
\cite[\S B.12]{HHR}. For a finite cyclic 2-group $G$ we define
\begin{displaymath}
MU^{((G))}= N_{2}^{g}MU_{\reals}.
\end{displaymath}

\noindent Choose a generator $\gamma $ of $G$. In
\cite[(5.47)]{HHR} we defined generators
\begin{numequation}
\label{eq-rbar}
\orr_{k}=\orr^{G}_{k}\in
\upi^{C_{2}}_{k\rho_{2}}i_{C_{2}}^{*}MU^{((G))} (C_{2}/C_{2})
\cong \upi_{C_{2},k\rho_{2}}MU^{((G))} (G/G)
\end{numequation}%

\noindent (note that this group is a module over $G/C_{2}$) and
\begin{align*}
 r_{k}=\underline{r}_{1}^{2}(\orr_{k})
 & \in  \pi_{\ee, 2k}^{u}MU^{((G))} (G/G) 
 \cong  \upi^{\ee}_{2k}MU^{((G))} (\ee/\ee)  
  =\pi_{2k}^{u}MU^{((G))}.
\end{align*}

\noindent The Hurewicz images of the $\orr_{k}$ (for which we use the
same notation) are defined in terms of the coefficients (see
Definition \ref{def-graded})
\begin{displaymath}
\om_{k}\in
\upi^{C_{2}}_{k\rho_{2}}\HZZ_{(2)}\wedge MU^{((G))} (C_{2}/C_{2}) 
 = \upi_{C_{2},k\rho_{2}}\HZZ_{(2)}\wedge MU^{((G))} (G/G) 
\end{displaymath}

\noindent of the logarithm of the formal group law $\overline{F} $
associated with the left unit map from $MU$ to $MU^{((G))}$.  The
formula is
\begin{displaymath}
\sum _{k\geq 0}\overline{r}_{k}x^{k+1} 
=\left(x+\sum_{\ell>0}\gamma(\overline{m}_{2^{\ell}-1})x^{2^{\ell}}\right)^{-1}
           \circ \log_{\overline{F} }(x)
\end{displaymath}

\noindent where 
\begin{displaymath}
\log_{\overline{F} }(x)=x+\sum_{k>0}\overline{m}_{k}x^{k+1}. 
\end{displaymath}

For small $k$ we have 
\begin{align*}
\orr_{1}
 & =   (1-\gamma ) (\om_{1})\\                        
\orr_{2}
 & =   \om_{2}-2\gamma (\om_{1} )(1-\gamma ) (\om_{1})\\
\orr_{3}
 & =   (1-\gamma ) (\om_{3})-\gamma (\om_{1}) 
    (5\gamma (\om_{1})^{2}-6\gamma (\om_{1})\om_{1}+\om_{1}^{2} +2\om_{2})
\end{align*}

Now let $G=C_{2}$ or $C_{4}$ and, in the latter case
$G'=C_{2}\subseteq G$. The generators $\orr^{G}_{k}$ are the
$\orr_{k}$ defined above.  We also have elements $\orr^{G'}_{k}$
defined by similar formulas with $\gamma $ replaced by $\gamma^{2}$;
recall that $\gamma^{2} (\om_{k})= (-1)^{k}\om_{k}$.  
They are the images of similar generators of 
\begin{displaymath}
\upi^{C_{2}}_{k\rho_{2}}MU^{((G'))} (C_{2}/C_{2})
 \cong \upi_{C_{2},k\rho_{2}}MU^{((G'))} (G'/G')
\end{displaymath}

\noindent under the left unit map 
\begin{displaymath}
MU^{((G'))}\to MU^{((G'))}\wedge MU^{((G'))}\cong i^{*}_{G'} MU^{((G))}.
\end{displaymath}

\noindent
Thus we have
\begin{align*}
\orr_{1}^{G'}
 & =   2 \om_{1} \\
\orr_{2}^{G'}
 & =   \om_{2}+4\om_{1}^{2}\\
\orr_{3}^{G'}
 & =   2 \om_{3}+2\om_{1}\om_{2}+12\om_{1}^{3}
\end{align*}

\noindent If we set $\orr_{2}=0$ and $\orr_{3}=0$, we get 
\begin{numequation}
\label{eq-r3}
\begin{split}
\left\{\begin{array}[]{rl}
\orr_{1}^{G'}
 &\hspace{-2.5mm} =    \orr_{1,0}+\orr_{1,1}\\                        
\orr_{2}^{G'}
 &\hspace{-2.5mm} =   3\orr_{1,0} \orr_{1,1}+\orr_{1,1}^{2}\\
\orr_{3}^{G'}
 &\hspace{-2.5mm} =  5 \orr_{1,0}^{2}\orr_{1,1} +5\orr_{1,0} \orr_{1,1}^{2} 
        +\orr_{1,1}^{3}
            = \orr_{1,1}( 5 \orr_{1,0}^{2} +5\orr_{1,0} \orr_{1,1}
        +\orr_{1,1}^{2})\\
\gamma (\orr_{3}^{G'})
 &\hspace{-2.5mm} = -\orr_{1,0}( 5 \orr_{1,1}^{2} -5\orr_{1,0} \orr_{1,1}
        +\orr_{1,0}^{2})\\
\lefteqn{-\orr_{3}^{G'}\gamma (\orr_{3}^{G'})/\orr_{1,0} \orr_{1,1}}
          \qquad\qquad\\
 &\hspace{-2.5mm} =  \left(5 \orr_{1,1}^2-5
\orr_{1,0} \orr_{1,1}+\orr_{1,0}^2\right)  \left(\orr_{1,1}^2+5
\orr_{1,0} \orr_{1,1}+5 \orr_{1,0}^2\right)\\
&\hspace{-2.5mm} =  (5 \orr_{1,0}^{4}-20\orr_{1,0}^{3}\orr_{1,1}
                 +\orr_{1,0}^{2}\orr_{1,1}^{2}
          +20\orr_{1,0}\orr_{1,1}^{3}+5\orr_{1,1}^{4})\\ 
&\hspace{-2.5mm} =  \left(5 (\orr_{1,0}^{2}-\orr_{1,1}^{2})^{2}
        -20\orr_{1,0}\orr_{1,1}(\orr_{1,0}^{2}-\orr_{1,1}^{2})
             +11 (\orr_{1,0}\orr_{1,1})^{2} \right)
\end{array}  \right.\\
\end{split}
\end{numequation}%

\noindent where $\orr_{1,0}=\orr_{1}$ and $\orr_{1,1}=\gamma  (\orr_{1})$.

\begin{defin}\label{def-kH} 
{\bf $\kR$, $\KR$, $\kH$ and $\KH$.} 
The $C_{2}$-spectrum $k_{\reals}$ (connective real $K$-theory), is the
spectrum obtained from $MU_{\reals}$ by killing the $r_{n}$s for
$n\geq 2$. Its periodic counterpart $\KR$ is the telescope obtained
from $\kR$ by inverting\linebreak  ${\orr_{1}\in \upi_{\rho_{2}}\kR
(C_{2}/C_{2})}$.

The $C_{4}$-spectrum $\kH$ is obtained
from $MU^{((C_{4}))}$ by killing the $r_{n}$s and their conjugates for
$n\geq 2$.  Its periodic counterpart $\KH$ is the telescope obtained
from $\kH$ by inverting a certain element $D\in \upi_{4\rho_{4}}\kH
(C_{4}/C_{4})$ defined below in (\ref{eq-D}) and Table
\ref{tab-pi*}.
\end{defin}

The  image of $D$ in $\upi^{C_{2}}_{8\rho_{2}}\kH
(C_{2}/C_{2})\cong \upi_{C_{2},8\rho_{2}}\kH(C_{4}/C_{4})$ is
\begin{numequation}\label{eq-r24D}
\begin{split}
\left\{ \begin{array}{rl}
\underline{r}_{2}^{4} (D)
 &\hspace{-2.5mm} = \orr_{1,0}\orr_{1,1}\orr_{3}^{G'}\gamma (\orr_{3}^{G'})  \\
 &\hspace{-2.5mm} = \orr_{1,0}^{2}\orr_{1,1}^{2}\left( -5 \orr_{1,0}^{4}
+20 \orr_{1,0}^{3}\orr_{1,1}-\orr_{1,0}^{2}\orr_{1,1}^{2}
-20 \orr_{1,0} \orr_{1,1}^{3}-5 \orr_{1,1}^{4} \right)\\
 &\hspace{-2.5mm} = -\orr_{1,0}^{2}\orr_{1,1}^{2}
\left(5(\orr_{1,0}^{2}-\orr_{1,1}^{2})^{2}
-20\orr_{1,0}\orr_{1,1}(\orr_{1,0}^{2}-\orr_{1,1}^{2})\right.\\
 &\qquad \qquad \qquad \left. +11
(\orr_{1,0}\orr_{1,1})^{2} \right).
\end{array} \right.
\end{split}
\end{numequation}%

\noindent It is fixed by the action of $G/G'$, while its factors
$\orr_{1,0}\orr_{1,1}$ and $\orr_{3}^{G'}\gamma (\orr_{3}^{G'})$ are
each negated by the action of the generator $\gamma $.

We remark that while $MU^{((C_{4}))}$ is $MU_{\reals}\wedge
MU_{\reals}$ as a $C_{2}$-spectrum, $\kH$ is {\em not}
$k_{\reals}\wedge k_{\reals}$ as a $C_{2}$-spectrum.  The former has
torsion free underlying homotopy but the latter does not.
\bigskip

\section{The slice {\SS} for $\KR$}\label{sec-Dugger}

In this section we describe the slice {\SS} for $\KR$.  These results
are originally due to Dugger \cite{Dugger}, to which we refer for many
of the proofs.  This case is far simpler than that of $\KH$, but it is
very instructive.

\begin{thm}\label{thm-Dugger-slice}
{\bf The slice $\EE_{2}$-terms for $\KR$ and $\kR$. }
The slices of $\KR$ are
\begin{displaymath}
P_{t}^{t}\KR=\mycases{
\Sigma^{(t/2)\rho_{2}}\HZZ
        &\mbox{for $t$ even}\\
*       &\mbox{otherwise}
}
\end{displaymath}

\noindent For $\kR$ they are the same in nonnegative dimensions, and
contractible below dimension 0.
\end{thm}

Hence we know the integrally graded homotopy groups of these slices by
the results of \S \ref{sec-chain}, and they are shown in Figure
\ref{fig-sseq-1}.  It shows the $\EE_{2}$-term for the wedge of all of
the slices of $\KR$, and $\KR$ itself has the same $\EE_{2}$-term.  It
turns out that the differentials and Mackey functor extensions are
determined by the fact that $\upi_{*}\KR$ is 8-periodic, while the
$\EE_{2}$-term is far from it.  This explanation is admittedly
circular in that the proof of the Periodicity Theorem itself of
\cite[\S9]{HHR} relies on the existence of certain differentials
described below in (\ref{eq-slicediffs}).

\begin{figure}
\begin{center}
\includegraphics[width=12cm]{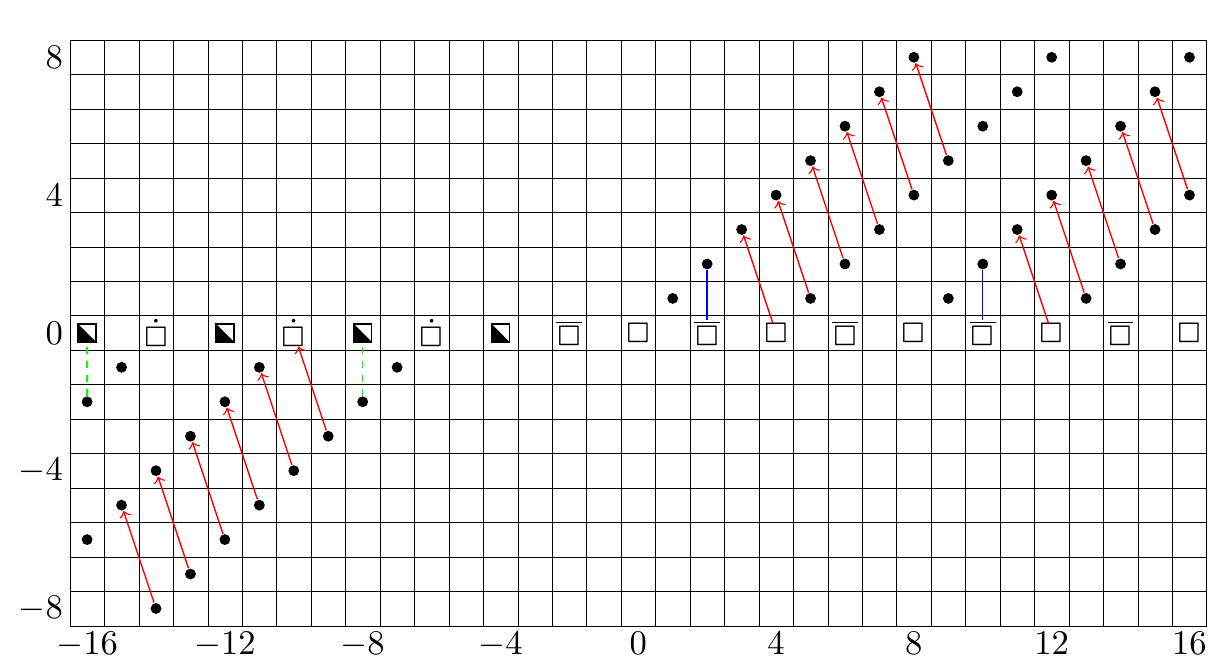} \caption[The slice {\SS} for
{$\KR$}.]{The slice {\SS} for $\KR$. Compare with Figure
\ref{fig-sseq-1}. Exotic transfers and restrictions are indicated
respectively by solid blue and dashed green lines.  Differentials are
in red.}  \label{fig-KR}
\end{center}
\end{figure}

\begin{thm}\label{thm-KR} 

{\bf The slice {\SS} for $\KR$.} The differentials and extensions in the
{\SS} are as indicated in Figure \ref{fig-KR}.

\end{thm}

\proof
There are four phenomena we need to establish:

\begin{enumerate}
\item [(i)] The differentials in the first quadrant, which are
indicated by red lines.

\item [(ii)] The differentials in the third quadrant.

\item [(iii)]  The exotic transfers in the first quadrant, which are
indicated by blue lines.

\item [(iv)]  The exotic restrictions in the third quadrant, which are
indicated by dashed green lines.
\end{enumerate}

For (i), note that there is a nontrivial element in
$\EE_{2}^{3,6} (G/G)$, which is part of the 3-stem, but
nothing in the $(-5)$-stem.  This means the former element must be
killed by a differential, and the only possiblilty is the one
indicated.  The other differentials in the first quadrant follow from
this one and the multiplicative structure.

For (ii), we know know that $\upi_{7}\KR=0$, so the same must be true
of $\upi_{-9}$.  Hence the element in $\EE_{2}^{-3,-12}$
cannot survive, leading to the indicated third quadrant differentials.

For (iii), note that $\upi_{2}$ and $\upi_{-6}$ must be the same as
Mackey functors.  This forces the indicated exotic transfers.  For
each $m\geq 0$ one has a nonsplit short exact sequence of $C_{2}$
Mackey functors
\begin{displaymath}
\xymatrix
@R=5mm
@C=10mm
{
0\ar[r]^(.5){}
    &{\EE_{2}^{2,8m+4}}\ar[r]^(.5){}\ar@{=}[d]^(.5){}
        &{\upi_{8m+2}\KR}\ar[r]^(.5){}\ar@{=}[d]^(.5){}
            &{\EE_{2}^{0,8m+2}}\ar[r]^(.5){}\ar@{=}[d]^(.5){}
                &0\\
    &\bullet
        &\dot{\Box}
            &\overline{ \Box}
}
\end{displaymath}

For (iv), note that $\upi_{-8}$ and $\upi_{0}$ must also agree.  This
forces the indicated exotic restrictions.  For each $m<0$ one has a
nonsplit short exact sequence
\begin{displaymath}
\xymatrix
@R=5mm
@C=10mm
{
0\ar[r]^(.5){}
    &{\EE_{2}^{0,8m}}\ar[r]^(.5){}\ar@{=}[d]^(.5){}
        &{\upi_{8m}\KR}\ar[r]^(.5){}\ar@{=}[d]^(.5){}
            &{\EE_{2}^{-2,8m-2}}\ar[r]^(.5){}\ar@{=}[d]^(.5){}
                &0\\
    &{\twobox}
        &\dot{\Box}
            &\bullet
}
\end{displaymath} 
\qed\bigskip 

In order to describe $\upi_{*}\KR$ as a graded Green functor, meaning
a graded Mackey functor with multiplication, we recall some notation
from \S\ref{sec-chain}(i) and Definition \ref{def-aeu}.  For $G=C_{2}$
we have elements
\begin{numequation}\label{eq-upi-HZ}
\begin{split}\left\{
\begin{array}{rlrl}
a=a_{\sigma }
 &\hspace{-3mm} \in \upi_{-\sigma }\HZZ (G/G)   \\
u=u_{2\sigma }
 &\hspace{-3mm} \in \upi_{2-2\sigma }\HZZ (G/G)   \\
x=u_{\sigma }
 &\hspace{-3mm} \in \upi_{1-\sigma }\HZZ (G/\ee) &&\mbox{with $x^{2}=\Res(u)$}  \\
z_{n} = e_{2n\rho_{2} }
 &\hspace{-3mm} \in   \upi_{2n (\sigma -1)}\HZZ (G/\ee)
                &&\mbox{for $n>0$}\\
a^{-i}\Tr (x^{-2n-1})
  &\hspace{-3mm} \in \upi_{(2n+1) (\sigma -1)+ i\sigma}\HZZ (G/G)
                &&\mbox{for $n>0$}
\end{array} \right.
\end{split}
\end{numequation}%

\noindent We will use the same symbols for the representatives of
these elements in the slice $E_{2}$-term.  The filtrations of $u$, $x$
and $z_{n}$ are zero while that of $a$ is one.  It follows that
$a^{-i}\Tr (x^{-2n-1})$ has filtration $-i$.  The element $x$ is invertible.

In $\underline{E}_{2}^{*,*}$ we have relations in 
\begin{numequation}\label{eq-upi-HZ-relns}
\begin{split}
\left\{
\begin{array}[]{rclrcl}
2a  &\hspace{-2.5mm} = & 0\qquad &                       
\Res (a) 
    &\hspace{-2.5mm} = & 0\\
z_{n} &\hspace{-2.5mm} = & x^{-2n}\qquad &
\Tr (x^{n}) &\hspace{-2.5mm} = &\mycases{
2u^{n/2}
        &\hspace{-1cm}\mbox{for $n$ even and $n\geq 0$}\\
\Tr_{}^{}(z_{-n/2})\neq 0\\
        &\hspace{-1cm}\mbox{for $n$ even and $n< 0$}\\
0       &\hspace{-1cm}\mbox{for $n$ odd and $n>-3$}\\
\neq 0  &\hspace{-1cm}\mbox{for $n$ odd and $n\leq -3$.}
} 
\end{array}
\right.
\end{split}
\end{numequation}%

We also have the element $\orr_{1}\in \upi_{1+\sigma }\kR (G/G)$, the
image of the element of the same name in $\orr_{1}\in
\upi_{1+\sigma}MU_{\reals} (G/G)$ of (\ref{eq-rbar}).  We use the same
symbol for its representative $\underline{E}_{2}^{0,1+\sigma }
(G/G)$. Then we have integrally graded elements
\begin{align*}
\eta 
  = a \orr_{1}
     &\in \underline{E}_{2}^{1,2} (G/G)  \\
v_{1}
  = x\cdot \Res(\orr_{1})
     &\in \underline{E}_{2}^{0,2} (G/\ee)
                   \qquad \mbox{with }\gamma (v_{1})=-v_{1}   \\
u\orr_{1}^{2}
     &\in \underline{E}_{2}^{0,4} (G/G) \\
w=2u\orr_{1}^{2}
     &\in \underline{E}_{2}^{0,4} (G/G) \\
b=u^{2}\orr_{1}^{4}
     &\in \underline{E}_{2}^{0,8} (G/G) 
                       \qquad \mbox{with }w^{2}=4b,
\end{align*}

\noindent where $\eta $ and $v_{1}$ are the images of the elements of
the same name in $\pi_{1}S^{0}$ and $\pi_{2}k$, and $w$ and $b$ are
permanent cycles.  The elements $x$, $v_{1}$ and $b$ are invertible.
Note that for $n<0$,
\begin{align*}
\underline{E}_{2}^{0,2n} (G/G)
 & = \mycases{
0      &\mbox{for }n=1\\
\Z \mbox{ generated by } \Tr_{}^{}(v_{1}^{-n})  
       &\mbox{for $n$ even}\\
\Z/2 \mbox{ generated by } \Tr_{}^{}(v_{1}^{-n})  
       &\mbox{for $n$ odd and }n<-1
}  
\end{align*}

\noindent so each group is killed by $\eta =a\orr_{1}$ by \ref{lem-hate}.

Then we have 
\begin{align*}
d_{3} (u)
 & =  a^{3}\orr_{1} 
     \qquad \mbox{by (\ref{eq-slicediffs2}) below,} \\
\mbox{so } 
d_{3} (u\orr_{1}^{2}) = d_{3} (u)\orr_{1}^{2}
 & =  a^{3}\orr_{1}^{3} =\eta^{3}\\
\Tr_{1}^{2}(x)
 & = a^{2}\orr_{1}  
    \qquad \mbox{by (\ref{eq-exotic-transfers}), raising filtration by 2,}  \\ 
\mbox{so } 
\Tr_{1}^{2}(v_{1})
 & = \eta^{2}.
\end{align*}

  Thus we get

\begin{thm}\label{thm-upi-KR}
{\bf The homotopy of $\KR$ as an integrally graded Green functor}.
With notation as above,
\begin{align*}
\upi_{*}\KR (G/\ee)
 & = \Z[v_{1}^{\pm 1}]  \\
\upi_{*}\KR (G/G)
 & = \Z[b^{\pm 1}, w, \eta ]/ (2\eta ,\eta^{3},w\eta ,w^{2}-4b)
\end{align*}

\noindent with
\begin{align*}
\Tr(v_{1}^{i})
 & = \mycases{
2b^{j}
       &\mbox{for }i=4j\\
\eta^{2}b^{j}
       &\mbox{for }i=4j+1\\
wb^{j}
       &\mbox{for }i=4j+2\\
0      &\mbox{for }i=4j+3
}  \\
\Res(b)
 & =  v_{1}^{4},\quad 
\Res(w)
  = 2 v,\mbox{ and } 
\Res(\eta )
  = 0.
\end{align*}

\noindent For each $j<0$, $b^{j}$ has filtration $-2$ and supports an
exotic restriction in the slice spectral sequence as indicated in
Figure \ref{fig-KR}.  Both $v_{1}\Res(b^{j})$ and
$\eta^{2}b^{j}$ have filtration zero, so the transfer relating them is
does not raise filtration.
\end{thm}

Now we will describe the $RO (G)$-graded slice spectral sequence and
homotopy of $\KR$.  The former is trigraded since $RO (G)$ itself is
bigraded, being isomorphic as an abelian group to $\Z \oplus \Z$.  For
each integer $k$, one can imagine a chart similar to Figure
\ref{fig-KR} converging to the graded Mackey functor $\upi_{k\sigma
+*}\KR$. Figure \ref{fig-KR} itself is the one for $k=0$.  The product
of elements in the $k$th and $\ell $th charts lies in the $(k+\ell )$th
chart.  We have elements as in (\ref{eq-upi-HZ})
\begin{align*}
a=a_{\sigma }
 & \in \underline{E}_{2}^{1,1-\sigma } (G/G)   \\
u=u_{2\sigma }
 & \in \underline{E}_{2}^{0,2-2\sigma } (G/G)   \\
x=u_{\sigma }
 & \in \underline{E}_{2}^{0,1-\sigma }  (G/\ee) 
                &&\mbox{with $\gamma (x)=-1$ and $x^{2}=\Res(u)$}  \\
z_{n} = x^{-2n}
 & \in  \underline{E}_{2}^{0,-2n+2n\sigma  } (G/\ee)
                &&\mbox{for $n>0$}\\
a^{-i}\Tr (x^{-2n-1})
  & \in \underline{E}_{2}^{-i,-i-2n+2n\sigma}  (G/G)
                &&\mbox{for $i\geq 0$ and $n>0$}\\
\orr_{1}
  & \in \underline{E}_{2}^{0,1+\sigma }  (G/G)
\end{align*}

\noindent where $a$, $x$, $z_{n}$ and $\orr_{1}$ are permanent cycles,
both $x$ and $\orr_{1}$ are invertible, and there are  relations as in
(\ref{eq-upi-HZ-relns}). We also know that
\begin{align*}
d_{3} (u)
 & =  a^{3}\orr_{1} \qquad \mbox{by (\ref{eq-slicediffs2}) below} \\
\aand 
\Tr_{1}^{2}(x)
 & = a^{2}\orr_{1}  \qquad \mbox{by (\ref{eq-exotic-transfers})} .
\end{align*}

\begin{thm}\label{thm-ROG-graded}
{\bf The $RO (G)$-graded slice spectral sequence for $\KR$} can be
obtained by tensoring that of Figure \ref{fig-KR} with
$\Z[\orr_{1}^{\pm 1}]$, that is for any integer $k$
\begin{align*}
\underline{E}_{2}^{s,t +k\sigma } (G/G)
 & \cong  \orr_{1}^{k} \underline{E}_{2}^{s,t -k }(G/G ) \\
\aand 
\underline{E}_{2}^{s,t +k\sigma } (G/\ee)
 & \cong  \Res(\orr_{1}^{k}) \underline{E}_{2}^{s,t -k }(G/\ee )
\end{align*}

\noindent and $\upi_{t+k\sigma}\KR$ has a similar description.
\end{thm}

\proof The element $\orr_{1}$ and its restriction are invertible
permanent cycles, so multiplication by either induces an isomorphism
in the spectral sequence.  \qed\bigskip

\begin{rem}\label{rem-a3}
In {\bf the $RO (G)$-graded slice {\SS } for $\kR$} one has $d_{3}
(u)=\orr_{1}a^{3}$, but $a^{3}$ itself, and indeed all higher powers
of $a$, survive to $\underline{E}_{4}=\underline{E}_{\infty }$.  Hence
the $\underline{E}_{\infty }$-term of this {\SS} does {\bf not} have
the horizontal vanishing line that we see in $\underline{E}_{4}$-term
of Figure \ref{fig-KR}. However when we pass from $\kR$ to $\KR$,
$\orr_{1}$ becomes invertible and we have
\begin{displaymath}
d_{3} (\orr_{1}^{-1}u)=a^{3}.
\end{displaymath}
\end{rem}

We can keep track of the groups in this trigraded {\SS } with the
help of four variable {\Ps} $g (\underline{E}_{r} (G/G))\in \Z[[x,y,z,t]]$
in which the rank of $\underline{E}_{r}^{s,i+j\sigma } (G/G)$ is the
coefficient in $\Z[[t]]$ of $x^{i-s}y^{j}z^{s}$.  The variable $t$
keeps track of powers of two. Thus  a copy of the integers
is represented by $1/ (1-t)$ or (when it is the kernel of a differential
of the form $\Z \to \Z/2$) $t/ (1-t)$.  Let 
\begin{numequation}\label{eq-aur}
\begin{split}
\wa =y^{-1}z,\qquad 
\uu =x^{2}y^{-1}\qquad \aand 
\rr =xy.
\end{split}
\end{numequation}

\noindent  Since 
\begin{displaymath}
\underline{E}_{2} (G/G) =\Z[a,u,\orr_{1}]/ (2a),
\end{displaymath}

\noindent we have 
\begin{align*}
g (\underline{E}_{2} (G/G))
 & = \left(\frac{1}{1-t}+ \frac{\wa}{1-\wa}\right)
  \frac{1}{(1-\uu) (1-\rr)}  \\
g (\underline{E}_{4} (G/G))
 & =  g (\underline{E}_{2} (G/G))
   -\frac{\uu+\rr\wa^{3}}
         {(1-\wa) (1-\uu^{2}) (1-\rr)}.
\end{align*}

\noindent We subtract the indicated expression from $g
(\underline{E}_{2} (G/G))$ because we have differerentials
\begin{displaymath}
d_{3} (a^{i}\orr_{1}^{j}u^{2k+1})= a^{i+3}\orr_{1}^{j+1}u^{2k}
\qquad \mbox{for all }i,j,k\geq 0. 
\end{displaymath}

\noindent Pursuing this further we get
\begin{align*}
g (\underline{E}_{4} (G/G))
 & =  \left(\frac{1}{1-t}+ \frac{\wa}{1-\wa}\right)
  \frac{1}{(1-\uu) (1-\rr)}
   -\frac{\uu}
         { (1-\uu^{2}) (1-\rr)}\\
 & \qquad 
   -\frac{a\uu+\rr\wa^{3}}
         {(1-\wa) (1-\uu^{2}) (1-\rr)} \\ 
 & = \frac{1+\uu-\uu (1-t)}
          {(1-t) (1-\uu^{2}) (1-\rr)} 
     +\frac{\wa (1+\uu)
            -\wa (\uu+\wa^{2}\rr)}
           {(1-\wa) (1-\uu^{2}) (1-\rr)}\\
 & = \frac{1+t\uu}
          {(1-t) (1-\uu^{2}) (1-\rr)} 
     +\frac{\wa-\wa^{3}+\wa^{3}-\wa^{3}\rr}
           {(1-\wa) (1-\uu^{2}) (1-\rr)}\\
 & = \frac{1+t\uu}
          {(1-t) (1-\uu^{2}) (1-\rr)} 
     +\frac{\wa+\wa^{2}}
           { (1-\uu^{2}) (1-\rr)} 
     +\frac{\wa^{3}}
           {(1-\wa) (1-\uu^{2})} .
\end{align*}

\noindent The third term of this expression represents the elements of
filtration above two (referred to in \ref{rem-a3}) which disappear
when we pass to $\KR$.  The first term  represents the elements of
filtration zero, which include
\begin{numequation}\label{eq-[2u]]}
\begin{split}
1,\qquad  
[2u]\in \langle 2,\,a ,\,a^{2}\orr_{1} \rangle
\qquad \aand 
[u^{2}]\in \langle \,a ,\,a^{2}\orr_{1}\,a ,\,a^{2}\orr_{1} \rangle
\end{split}
\end{numequation}

\noindent Here we use the notation $[2u]$ and $[u^{2}]$ to indicate
the images in $\underline{E}_{4}$ of the elements $2u$ and $u^{2}$ in
$\underline{E}_{2}$; see Remark \ref{rem-abuse} below.  The former
{\em not} divisible by 2 and the latter is not a square since $u$
itself is not present in $\underline{E}_{4}$, where the Massey
products are defined.  For an introduction to Massey products, we
refer the reader to \cite[A1.4]{Rav:MU}.

We now make a similar computation where we enlarge $\underline{E}_{2} (G/G)$
by adjoining $\orr_{1}^{-1}u$ and denote the resulting {\SS } terms by
$\underline{E}'_{2}$ and $\underline{E}'_{4}$.  

Let 
\begin{displaymath}
\www = \rr^{-1}\uu= xy^{-3}.
\end{displaymath}

\noindent  Then since 
\begin{displaymath}
\underline{E}'_{2} (G/G) =\Z[a,\orr_{1}^{-1}u,\orr_{1}]/ (2a),
\end{displaymath}

\noindent we have 
\begin{align*}
g (\underline{E}'_{2} (G/G))
 & = \left(\frac{1}{1-t}+ \frac{\wa}{1-\wa}\right)
  \frac{1}{(1-\www) (1-\rr)}  \\
g (\underline{E}_{4} (G/G))
 & =  g (\underline{E}_{2} (G/G))
   -\frac{\www+\wa^{3}}
         {(1-\wa) (1-\www^{2}) (1-\rr)} \\ 
 & =  \left(\frac{1}{1-t}+ \frac{\wa}{1-\wa}\right)
  \frac{1}{(1-\www) (1-\rr)}
   -\frac{\www}
         { (1-\www^{2}) (1-\rr)}\\
 & \qquad 
   -\frac{a\www+\wa^{3}}
         {(1-\wa) (1-\www^{2}) (1-\rr)} \\ 
 & = \frac{1+\www-\www (1-t)}
          {(1-t) (1-\www^{2}) (1-\rr)} 
     +\frac{\wa (1+\www)
            -\wa (\www+\wa^{2})}
           {(1-\wa) (1-\www^{2}) (1-\rr)}\\
 & = \frac{1+t\www}
          {(1-t) (1-\www^{2}) (1-\rr)} 
     +\frac{\wa+\wa^{2}}
           { (1-\www^{2}) (1-\rr)}  
\end{align*}

\noindent and there is nothing in $\underline{E}'_{4}$ with filtration
above two.  As far as we know there is no modification of the spectrum
$\kR$ corresponding to this modification of $\underline{E}_{r}$.
However the map $\underline{E}_{r}\kR \to \underline{E}_{r}\KR$
clearly factors through $\underline{E}'_{r}$

\section{Some elements in the homotopy groups of $\kH$ and
$\KH$}\label{sec-more}

For $G=C_{4}$ we will often use a (second) subscript $\epsilon $ on
elements such as $r_{n}$ to indicate the action of a generator $\gamma
$ of $G=C_{4}$, so $\gamma (x_{\epsilon })=x_{1+\epsilon }$ and
$x_{2+\epsilon }=\pm x_{\epsilon }$.  Then we have
\begin{numequation}
\label{eq-pikH}
\pi_{*}^{u}\kH=\upi_{*}\kH (G/\ee)=\upi_{\ee,*}\kH (G/G)
=\Z[r_{1},\,\gamma (r_{1})]=\Z[r_{1,0},\,r_{1,1}]
\end{numequation}%

\noindent where $\gamma^{2} (r_{1,\epsilon })=-r_{1,\epsilon }$. Here
we use $r_{1,\epsilon }$ and $\orr_{1,\epsilon }$ to denote the images
of elements of the same name in the homotopy of $MU^{((G))}$.

\begin{numequation}\label{eq-bigraded}
\begin{split}
\includegraphics[width=9cm]{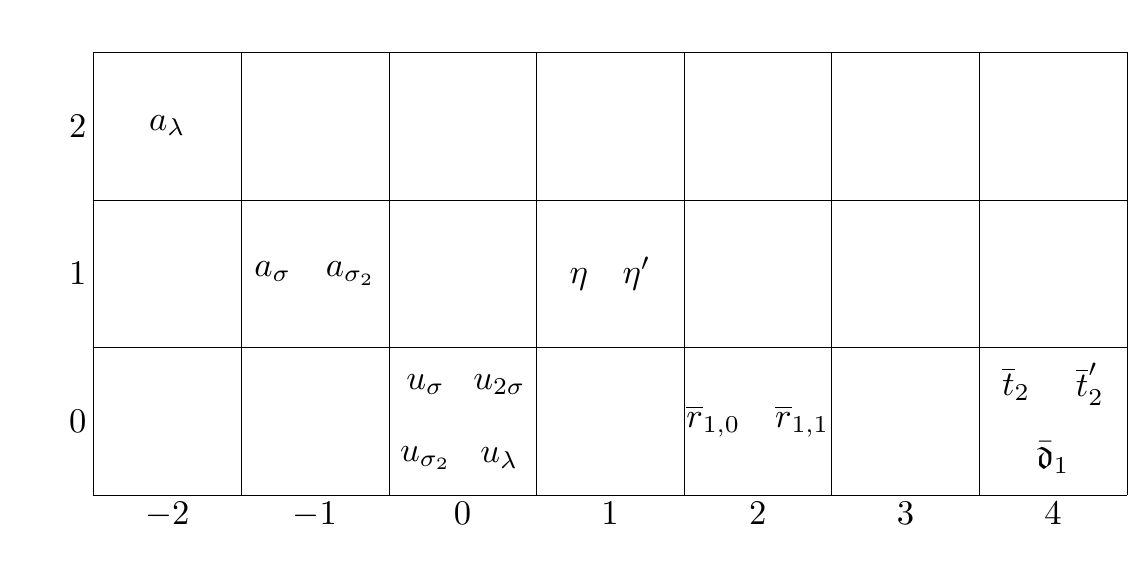}
\end{split}
\end{numequation}%

\noindent Here the vertical coordinate is $s$ and the horizontal
coordinate is $|t|-s$. More information about these elements can be
found in Table \ref{tab-pi*} below.

\noindent

We are using the following notational convention.  When
$x=\Tr_{2}^{4}(y)$ for some element ${y\in \upi_{\star}\kH (G/G')}$,
we will write $x'=\Tr_{2}^{4}(u_{\sigma}y)$.  Examples above include
the cases $x=\eta $ and $x=\ot_{2}$.  The primes could be iterated,
{\ie } we might write $x^{(k)}=\Tr_{2}^{4}(u_{\sigma}^{k}y)$, but
this turns out to be unnecessary.

The group action (by $G'$ on $\orr_{1,\epsilon }$, $a_{\sigma_{2}}$
and $u_{\sigma_{2}}$, and by $G$ on all the others) fixes each
generator but $u_{\sigma}$ and $u_{\sigma_{2}}$.
For them the action is given by
\begin{displaymath}
\xymatrix
@R=1mm 
@C=8mm
{
u_{\sigma}\ar@{<->}[r]^(.5){\gamma}
    &-u_{\sigma}
        &{\mbox{and} }
            &u_{\sigma_{2} }\ar@{<->}[r]^(.5){\gamma^{2}}
                &-u_{\sigma_{2} }
}
\end{displaymath}

\noindent by Theorem \ref{thm-module}.  This is compatible with the
following $G$-action:
\begin{displaymath}
\xymatrix
@R=5mm
@C=10mm
{
r_{1,0}\ar[r]^(.5){\gamma}
    &r_{1,1}\ar[d]^(.5){\gamma}\\
-r_{1,1}\ar[u]^(.5){\gamma} &-r_{1,0}\ar[l]^(.5){\gamma} }
\qquad \mbox{where } r_{1,\epsilon} =
\underline{r}_{1}^{2}(\orr_{1,\epsilon })\in \upi_{\ee,2}\kH (G/G) .
\end{displaymath}

We will see below (Theorem \ref{thm-d3ulambda}) that $d_{5}
(u_{2\sigma})=a_{\sigma}^{3}a_{\lambda}\normrbar_{1}$ and
$[u_{2\sigma}^{2}]$ is a permanent cycle.  Since all transfers are
killed by $a_{\sigma}$ multiplication (Lemma \ref{lem-hate}), this
implies that $[u_{2\sigma}x]$ is a permanent
cycle representing the Toda bracket
\begin{displaymath}
[u_{2\sigma}x]
 = [u_{2\sigma}\Tr_{2}^{4}(y)]
 = \langle x,\,a_{\sigma} ,\,a_{\sigma}^{2}a_{\lambda}\normrbar_{1} \rangle.
\end{displaymath}

\noindent This element is $x''$ since in $\EE_{2}$ we have
(using the Frobenius relation (\ref{eq-Frob}))
\begin{displaymath}
x''=\Tr_{2}^{4}(u_{\sigma}^{2}y)=\Tr_{2}^{4}(\Res_{2}^{4}(u_{2\sigma})y)
 = u_{2\sigma}\Tr_{2}^{4}(y)= u_{2\sigma}x.
\end{displaymath}

\noindent Similarly $x'''=u_{2\sigma}x'$. For $k\geq 4$,
$x^{(k)}=u_{2\sigma}^{2}x^{(k-4)}$ in $\upi_{\star}$ as well as
$\EE_{2}$.

\bigskip

The Periodicity Theorem \cite[Thm. 9.19]{HHR} states that inverting a class
in $\upi_{4\rho_{4}}\kH (G/G)$ whose image under
$\ur_{2}^{4}\Res_{2}^{4}$ is divisible by
$\orr_{3,0}^{G'}\orr_{3,1}^{G'}$ (see (\ref{eq-r3})) and
$\orr_{1,0}\orr_{1,1}=\orr_{1,0}^{G}\orr_{1,1}^{G}$ makes
$u_{8\rho_{4}}$ a permanent cycle.  One such class is

\noindent
\begin{numequation}\label{eq-D}
\begin{split}
\begin{array}[]{rl}
D&\hspace{-2.5mm} = N_{2}^{4} (\normrbar_{2}^{G'})\normrbar_{1}^{G} 
  = u_{2\sigma}^{-2}(\ur_{2}^{4}\Res_{2}^{4})^{-1}
            \left(\orr_{1,0}^{G}\orr_{1,1}^{G}
                  \orr_{3,0}^{G'}\orr_{3,1}^{G'} \right)\\
 &\hspace{-2.5mm} = \normrbar_{1}^{2} (-5\ot_{2}^{2}+20 \ot_{2}\normrbar_{1}
      +9\normrbar_{1}^{2}) \in \upi_{4\rho_{4}}\kH (G/G),\\
 &\qquad \mbox{where }\ot_{2}=\Tr_{2}^{4}(u_{\sigma }^{-1}[\orr_{1,0}^{2}])
          \mbox{ and $\normrbar_{1}$ is as in (\ref{eq-normrbar}) below,} 
\end{array}
\end{split}
\end{numequation}%

\noindent 
and $\KH = D^{-1}\kH$.  Then we know that $\Sigma^{32}\KH$ is
equivalent to $\KH$.

The Slice and Reduction Theorems \cite[Thms. 6.1 and 6.5]{HHR} imply that the
$2k$th slice of $\kH$ is the $2k$th wedge summand of
\begin{displaymath}
\HZZ \wedge N_{2}^{4}\left(\bigvee_{i\geq 0}S^{i\rho_{2}} \right).
\end{displaymath}

\noindent It follows that over $G'$ the $2k$th slice is a wedge of
$k+1$ copies of $\HZZ \wedge S^{k\rho_{2}}$.  Over $G$ we get the
wedge of the appropriate number of copies of $G_{+}\smashove{G'}\HZZ
\wedge S^{k\rho_{2}}$, wedged with a single copy of $\HZZ \wedge
S^{(k/2)\rho_{4}}$ for even $k$.  This is spelled out in Theorem
\ref{thm-sliceE2} below.

The group $\upi^{G'}_{\rho_{2}}\kH (G'/\ee)$ is {\em
not} in the image of the group action restriction
$\ur_{2}^{4}$ because $\rho_{2}$ is not the restriction of a
{\rep} of $G$.  However, $\pi_{2}^{u}\kH$ is refined (in the
sense of \cite[Def. 5.28]{HHR}) by a map from
\begin{numequation}
\label{eq-s1}
\xymatrix 
@R=5mm
@C=10mm
{
S_{\rho_{2}} : = G_{+}\smashove{G'}S^{\rho_{2}}
         \ar[r]^(.7){\os_{1} }
    &\kH.
}
\end{numequation}%

\noindent The Reduction Theorem implies that the 2-slice
$P_{2}^{2}\kH$ is $S_{\rho_{2}}\wedge \HZZ$.  We know that 
\begin{displaymath}
\upi_{2} ( S_{\rho_{2}}\wedge \HZZ) = \widehat{\oBox }.
\end{displaymath}

\noindent We use the symbols $r_{1}$ and $\gamma (r_{1})$ to denote
the generators of the underlying abelian group of
$\widehat{\oBox } (G/\ee)=\Z[G/G']_{-}$.  These elements have
trivial fixed point transfers and
\begin{displaymath}
\upi_{2} ( S_{\rho_{2}}\wedge \HZZ) (G/G')=0.
\end{displaymath}

Table \ref{tab-pi*} describes some elements in the slice {\SS} for $\kH$ in 
low dimensions, which we now discuss.

Given an element in $\pi_{\star}MU^{((G))}$, we will often use the
same symbol to denote its image in $\pi_{\star}\kH$. For example, in
\cite[\S9.1]{HHR}
\begin{numequation}
\label{eq-normrbar}
\normrbar_{n}\in \pi_{(2^{n}-1)\rho_{4}}^{G}MU^{((G))}
      =\upi_{(2^{n}-1)\rho_{4}}^{G}MU^{((G))} (G/G)
\end{numequation}%

\noindent was defined to be the composite
\begin{displaymath}
\xymatrix
@R=5mm
@C=10mm
{
S^{(2^{n}-1)\rho_{4}}\ar@{=}[r]^(.5){}
    &N_{2}^{4}S^{(2^{n}-1)\rho_{2}}\ar[rr]^(.5){N_{2}^{4}\orr_{2^{n}-1} }
        &   &N_{2}^{4}MU^{((G))}\ar[r]^(.5){}
                &MU^{((G))}.
}
\end{displaymath}

\noindent We will use the same symbol to denote its image in
$\upi_{(2^{n}-1)\rho_{4}}^{G}\kH (G/G)$. 

The element $\eta \in \pi_{1}S^{0}$ (coming from the Hopf map
$S^{3}\to S^{2}$) has image ${a_{\sigma}\orr_{1} \in
\upi^{G'}_{1}\kR (G'/G')}$.  There are two corresponding
elements
\begin{displaymath}
\eta_{\epsilon }\in \upi^{G'}_{1}\kH (G'/G')
\qquad \mbox{for } \epsilon =0,1. 
\end{displaymath}

\noindent We use the same symbol for their preimages under
$\ur_{2}^{4}$ in $\upi^{G}_{1}\kH (G/G')$, and there we have
\begin{displaymath}
\eta_{\epsilon}=a_{\sigma_{2}}\orr_{1,\epsilon}.
\end{displaymath}

\noindent We denote by $\eta $ again the image of either under the
transfer $\Tr_{2}^{4}$, so 
\begin{displaymath}
{\Res_{2}^{4}(\eta)=\eta_{0}+\eta_{1}}.
\end{displaymath}

Its cube is killed by a $d_{3}$ in the slice {\SS}, as is the sum of
any two monomials of degree 3 in the $\eta_{\epsilon }$.  It follows
that in $\EE_{4}$ each such monomial is equal to
$\eta_{0}^{3}$.  It has a nontrivial transfer, which we denote by $x_{3}$.

In \cite[Def. 5.51]{HHR} we defined 
\begin{numequation}
\label{eq-fk}
f_{k}=a_{\overline{\rho } }^{k}N_{2}^{g} (\orr_{k})
  \in \upi_{k}MU^{((G))} (G/G)
\end{numequation}%

\noindent for a finite cyclic 2-group $G$.  In particular,
$f_{2^{n-1}}=a_{\overline{\rho } }^{2^{n}-1}\normrbar_{n}$ for
$\normrbar_{n}$ as in (\ref{eq-normrbar}). The slice filtration of
$f_{k}$ is $k(g-1)$ and we will see below (Lemma \ref{lem-hate} and,
for $G=C_{4}$, Theorem \ref{thm-d3ulambda}) that
\begin{numequation}
\label{eq-Tr(usigma)}
\Tr_{G'}^{G}(u_{\sigma}) = a_{\sigma}f_{1}.
\end{numequation}%

Note that $u_{\sigma}\in \EE_{2}^{0,1-\sigma} (G/G')$
since the maximal subgroup for which the sign {\rep} $\sigma $ is
oriented is $G'$, on which it restricts to the trivial {\rep} of
degree 1.  This group depends only on the restriction of the $RO
(G)$-grading to $G'$, and the isomorphism extends to differentials as
well.  This means that $u_{\sigma}$ is a place holder corresponding
to the permanent cycle $1\in \EE_{2}^{0,0} (G/G')$. 

For $G=C_{4}$ (\ref{eq-Tr(usigma)}) implies
\begin{displaymath}
\Tr_{2}^{4}(u_{\sigma}) = a_{\sigma}f_{1}
= a_{\sigma}^{2}a_{\lambda}\normrbar_{1}. 
\end{displaymath}

\noindent For example 
\begin{align*}
\Tr_{2}^{4}(\eta_{0}\eta_{1})
 & =   \Tr_{2}^{4}(a_{\sigma_{2}}^{2}\orr_{1,0}\orr_{1,1})
   =   \Tr_{2}^{4}(u_{\sigma}\Res_{2}^{4}(a_{\lambda}\normrbar_{1}))\\
 & =   \Tr_{2}^{4}(u_{\sigma})a_{\lambda}\normrbar_{1}
   =    a_{\sigma}f_{1}a_{\lambda}\normrbar_{1}
   =    f_{1}^{2}
\end{align*}

The Hopf element $\nu \in \pi_{3}S^{0}$ has image 
\begin{displaymath}
a_{\sigma}u_{\lambda}\normrbar_{1}\in \upi_{3}\kH (G/G),
\end{displaymath}

\noindent so we also denote the latter by $\nu $.  (We will see below
in (\ref{eq-d-ul}) that $u_{\lambda}$ is not a permanent cycle, but
$\alphax := a_{\sigma}u_{\lambda}$ is (\ref{eq-alphax}).) It has an
exotic restriction $\eta_{0}^{3}$ (filtration jump two), which implies
that
\begin{displaymath}
2\nu =\Tr_{2}^{4}(\Res_{2}^{4}(\nu ))=\Tr_{2}^{4}(\eta_{0}^{3})=x_{3}.
\end{displaymath}

\noindent One way to see this is to use the Periodicity Theorem to
equate $\upi_{3}\kH$ with $\upi_{-29}\kH$,
which can be shown to be the Mackey functor $\circ $ in slice
filtration $-32$.  Another argument not relying on periodicity is given
below in Theorem \ref{thm-d3ulambda}.

The exotic restriction on $\nu $ implies
\begin{displaymath}
\Res_{2}^{4}(\nu^{2})=\eta_{0}^{6},
\end{displaymath}

\noindent with filtration jump 4.

\begin{thm}\label{thm-Hur}
{\bf The Hurewicz image.} The elements $\nu \in \upi_{3}\kH (G/G)$,
\linebreak ${\epsilon \in \upi_{8}\kH (G/G)}$, $\kappa \in
\upi_{14}\kH (G/G)$, and $\overline{\kappa} \in \upi_{20}\kH (G/G)$
are the images of elements of the same names in $\pi_{*}S^{0}$.  The
image of the Hopf map $\eta\in\pi_{1}S^{0}$ is either
$\eta=\Tr_{2}^{4}(\eta_{\epsilon})$ or its sum with $f_{1}$.
\end{thm}

We refer the reader to \cite[Table A3.3]{Rav:MU} for more information
about these elements.

\proof Suppose we know this for $\nu $ and $\overline{\kappa} $.
Then $\Delta_{1}^{-4}\nu $ is represented by an element of filtration
$-3$ whose product with $\nu^{2}$ is nontrivial.  This implies that
$\nu^{3}$ has nontrivial image in $\underline{\pi }_{9}\kH (G/G)$.
This is a nontrivial multiplicative extension in the first quadrant,
but not in the third.  The spectral sequence representative of
$\nu^{3}$ has filtration 11 instead of 3.  We will see later that
$\nu^{3}=2n$ where $n$ has filtration 1, and $\nu^{3}$ is the transfer
of an element in filtration 1.

Since $\nu^{3}=\eta \epsilon $ in $\pi_{*}S^{0}$, this implies that
$\eta $ and $\epsilon $ are both detected and have the images stated
in Table \ref{tab-pi*}.  It follows that $\epsilon \overline{\kappa
}$ has nontrivial image here. Since $\kappa^{2}=\epsilon
\overline{\kappa}$ in $\pi_{*}S^{0}$, $\kappa $ must also be
detected.  Its only possible image is the one indicated.


Both $\nu $ and $\overline{\kappa}$ have images of order $8$ in
$\pi_{*}TMF$ and its $K (2)$ localization.  The latter is the homotopy
fixed point set of an action of the binary tetrahedral group $G_{24}$
acting on $E_{2}$.  This in turn is a retract of the homotopy fixed
point set of the quaternion group $Q_{8}$. A restriction and transfer
argument shows that both elements have order at least 4 in the
homotopy fixed point set of $C_{4}\subset Q_{8}$.  

There is an orientation map $MU\to E_{2}$, which extends to a
$C_{2}$-{\eqvr} map $MU_{\reals}\to E_{2}$.  Norming up and
multiplying on the right gives us a $C_{4}$-{\eqvr} map
$N_{2}^{4}MU_{\reals}\to E_{2}$.  This $C_{4}$-action on the target is
compatible with the $G_{24}$-action leading to $L_{K (2)}TMF$.

The image of $\eta\in\pi_{1}S^{0}$ must restrict to
$\eta_{0}+\eta_{1}$, so modulo the kernel of $\Res_{2}^{4}$ it is the
element $\Tr_{2}^{4}(\eta_{\epsilon})$, which we are calling $\eta$.
The kernel of $\Res_{2}^{4}$ is generated by $f_{1}$.  \qed\bigskip

We now discuss the norm $N_{2}^{4}$, which is a functor from the
category of $C_{2}$-spectra to that of $C_{4}$ spectra.  As explained
above in connection with Corollary \ref{cor-normdiff}, for a
$C_{4}$-ring spectrum $X$ we have an internal norm  
\begin{displaymath}
\upi_{V}^{G'}i^{*}_{G'}X (G'/G')
\cong \upi_{G',V}^{G'}X (G/G)
 \to \upi_{\Ind_{2}^{4}V }^{G}X (G/G)
\end{displaymath}

\noindent and a similar functor on the slice spectral sequence for
$X$.  It preserves multiplication but not addition.  Its source is a
module over $G/G'$, which acts trivially on its target. Consider the
diagram
\begin{displaymath}
\xymatrix
@R=8mm
@C=8mm
{
{\upi_{G',V}X (G/G)}\ar[r]^(.45){\cong }
    &{\upi_{V}^{G'}i^{*}_{G'}X (G'/G')}\ar[r]^(.5){N_{2}^{4}}
        &{\upi_{\Ind_{2}^{4}V }^{G}X (G/G)}\ar[d]^(.5){\Res_{2}^{4}}\\
{\upi_{G',2V}X (G/G)}\ar[r]^(.45){\cong }
    &{\upi^{G'}_{2V}i^{*}_{G'}X (G'/G')}\ar[r]^(.5){\cong }
        &{\upi_{\Ind_{2}^{4}V }^{G}X (G/G')}
}
\end{displaymath}

\noindent

For $x\in \upi_{V}^{G'}i^{*}_{G'}X (G'/G')$ we have $x\gamma (x)\in
\upi_{2V}^{G'}i^{*}_{G'}X (G'/G')$ and $2V$ is the restriction of some
$W\in RO (G)$.  The group $\upi_{W}^{G}X (G/G')$ depends only on the
restriction of $W$ to $RO (G')$.  If $W'\in RO (G)$ is another virtual
{\rep} restricting to $2V$, then $W-W'=k (1-\sigma )$ for some integer
$k$.  The canonical isomorphism between $\upi_{W}^{G}X (G/G')$ and
$\upi_{W'}^{G}X (G/G')$ is given by multiplication by
$u_{\sigma}^{k}$.

\begin{defin}\label{def-bracket}
{\bf A second use of square bracket notation}. For $0\leq i\leq 2d$, let $f
(\orr_{1,0},\orr_{1,1})$ be a homogeneous polynomial of degree $2d-i$,
so
\begin{displaymath}
a_{\sigma_{2}}^{i}f (\orr_{1,0},\orr_{1,1})\in
\upi_{(2d-i)+ (2d-2i)\sigma_{2}}^{G'}i^{*}_{G'}\kH (G'/G')
\end{displaymath}

\noindent We will denote by
$[a_{\sigma_{2}}^{i}f(\orr_{1,0},\orr_{1,1})]$ its preimage in
$\upi_{2d-i+ (d-i)\lambda }\kH (G/G')$ under the isomorphism of
(\ref{eq-easy-iso}).
\end{defin}

The first use of square bracket notation is that of Remark
\ref{rem-abuse}. Note that ${\orr_{1,\epsilon }\in
\upi_{\rho_{2}}^{G'}i^{*}_{G'}\kH}$ is not the target of such an
isomorphism since ${\rho_{2}\in RO (G')}$ is not the restriction of
any element in $RO (G)$, hence the requirement that $f$ has even
degree.

We will denote $u_{\sigma }^{-1}[\orr_{1,\epsilon }^{2}]\in
\upi_{\rho_{4}}^{G}\kH (G/G')$ by $\overline{s}_{2,\epsilon }$.  Then
we have $\gamma(\overline{s}_{2,0})=-\overline{s}_{2,1}$ and
$\gamma(\overline{s}_{2,1})=-\overline{s}_{2,0}$.  We define
\begin{displaymath}
\ot_{2}:= (-1)^{\epsilon }\Tr_{2}^{4}(\overline{s}_{2,\epsilon }),
\end{displaymath}

\noindent which is independent of $\epsilon $, and we have 
\begin{displaymath}
\Res_{2}^{4}(\ot_{2})=\overline{s}_{2,0}-\overline{s}_{2,1}.
\end{displaymath}

\noindent

Then we have 
\begin{displaymath}
\Res_{2}^{4}(N_{2}^{4} (\orr_{1,0}))
 =\Res_{2}^{4}(\normrbar_{1}^{})
 =u_{\sigma}^{-1}[\orr_{1,0}\orr_{1,1}]
 \in \upi_{_{\rho_{4} }}\kH (G/G') .
\end{displaymath}

\noindent More generally for integers $m$ and $n$ 
\begin{align*}
\lefteqn{\Res_{2}^{4}(N_{2}^{4} (m\orr_{1,0}+n\orr_{1,1}))}\qquad\qquad\\
 & = u_{\sigma}^{-1}[(m\orr_{1,0}+n\orr_{1,1})
                                  (m\orr_{1,1}-n\orr_{1,0})]\\
 & = u_{\sigma}^{-1} ((m^{2}-n^{2}) [\orr_{1,0}\orr_{1,1}]
                   +mn ([\orr_{1,1}^{2}]-[\orr_{1,0}^{2}]))\\
 & = (m^{2}-n^{2})\Res_{2}^{4}(\normrbar_{1}^{})
                    -mn\,\Res_{2}^{4}(\ot_{2})
\end{align*}

\noindent so 
\begin{numequation}\label{eq-norm-orr}
\begin{split}
N_{2}^{4} (m\orr_{1,0}+n\orr_{1,1})
 = (m^{2}-n^{2})\normrbar_{1} -mn\ot_{2}.
\end{split}
\end{numequation}%

Similarly for integers $a$, $b$ and $c$,
\begin{align*}
\lefteqn{u_{\sigma }^{2}\Res_{2}^{4}(N_{2}^{4}
  (a\orr_{1,0}^{2}+b\orr_{1,0}\orr_{1,1}+c\orr_{1,1}^{2}))}\quad\\
 & = [(a\orr_{1,0}^{2}+b\orr_{1,0}\orr_{1,1}+c\orr_{1,1}^{2})
             (a\orr_{1,1}^{2}-b\orr_{1,0}\orr_{1,1}+c\orr_{1,0}^{2})]  \\
 & = [ac (\orr_{1,0}^{4}+\orr_{1,1}^{4})
              +b (c-a)\orr_{1,0}\orr_{1,1} (\orr_{1,0}^{2}-\orr_{1,1}^{2})
          + (a^{2}-b^{2}+c^{2})\orr_{1,0}^{2}\orr_{1,1}^{2}] \\
  & = [ac (\orr_{1,0}^{2}-\orr_{1,1}^{2})^{2}
              +b (c-a)\orr_{1,0}\orr_{1,1} (\orr_{1,0}^{2}-\orr_{1,1}^{2})
          + ((a+c)^{2}-b^{2})\orr_{1,0}^{2}\orr_{1,1}^{2}]  
\end{align*}

\noindent so 
\begin{numequation}\label{eq-norm-quad}
\begin{split}
N_{2}^{4}
  (a\orr_{1,0}^{2}+b\orr_{1,0}\orr_{1,1}+c\orr_{1,1}^{2})
= ac\,\ot_{2}^{2}+b (c-a)\normrbar_{1}\ot_{2}
                         +((a+c)^{2}-b^{2})\normrbar_{1}^{2}
\end{split}
\end{numequation}%

For future reference we need
\begin{align*}
\lefteqn{N_{2}^{4} ( 5 \orr_{1,0}^{2}\orr_{1,1} +5\orr_{1,0}\orr_{1,1}^{2}
            +\orr_{1,1}^{3}) }\qquad\qquad\\
 & = N_{2}^{4} (\orr_{1,1})N_{2}^{4} (5 \orr_{1,0}^{2} +5\orr_{1,0}\orr_{1,1}
            +\orr_{1,1}^{2}))  \\
 & = -\normrbar_{1}(5\ot_{2}^{2}-20\normrbar_{1}\ot_{2}
                         +11\normrbar_{1}^{2})
\end{align*}

\noindent Compare with (\ref{eq-r3}).

We also denote by
\begin{displaymath}
\eta_{\epsilon }
 =[a_{\sigma_{2}}\orr_{1,\epsilon }]
 \in \upi_{1}\kH (G/G')
\end{displaymath}

\noindent the preimage of $a_{\sigma_{2}}\orr_{1,\epsilon } \in
\upi_{1}^{G'}i^{*}_{G'}\kH (G'/G')$ and by $[a_{\sigma_{2}}^{2}]\in
\upi_{-\lambda }\kH (G/G')$ the preimage of $a_{\sigma_{2}}^{2}$.  The
latter is $\Res_{2}^{4}(a_{\lambda })$.

The values of $N_{2}^{4} (a_{\sigma_{2}})$ and $N_{2}^{4}
(u_{2\sigma_{2}})$ are given by Lemma \ref{lem-norm-au}, namely
\begin{align*}
N_{2}^{4} (a_{\sigma_{2}})
 & =  a_{\lambda } \\
\aand 
N_{2}^{4} (u_{2\sigma_{2}})
 & =  u_{2\lambda }/u_{2\sigma }.
\end{align*}

\noindent

\begin{center}
\begin{longtable}{|p{5cm}|p{6.8cm}|} 
    & \kill \caption[Some elements in the slice {\SS} and homotopy
groups of {$\kH$}.]
{Some elements in the slice {\SS}  and homotopy groups of $\kH$, listed in 
order of ascending filtration.}\\
\hline \multicolumn{1}{|c|}{Element}
    &\multicolumn{1}{|c|}{Description}\\
\hline\hline
 \multicolumn{2}{|c|}{Filtration 0}\\
\hline\endfirsthead  
\caption[Some elements in the slice {\SS} and homotopy groups of  {$\kH$}, continued.]
{Some elements in the slice {\SS} and homotopy groups of  $\kH$, continued.}\\
\hline   
\multicolumn{1}{|c|}{Element}    
    &\multicolumn{1}{|c|}{Description}\\
\hline\endhead 
\hline \multicolumn{2}{|r|}{Continued on next page} \\ \hline
\endfoot
\hline 
\endlastfoot 
\label{tab-pi*}
$\orr_{1,\epsilon}
      \in \upi_{\rho_{2}}^{G'}i_{G'}^{*}\kH (G'/G')$
\newline $\phantom{\orr_{1,\epsilon}\,\,} \cong  \upi_{G', \rho_{2}}\kH (G/G)$ 
\newline 
with $\orr_{1,2}=-\orr_{1,0}$
   &Images from (\ref{eq-rbar}) defined in \cite[(5.47)]{HHR}\\
\hline 
$r_{1,\epsilon}
 \in \upi_{\ee,2}\kH (G/G)$\newline 
$\phantom{\orr_{1,\epsilon}}\cong \upi_{G,2}\kH (G/\ee)\cong \pi_{2}^{u}\kH$
   &$\underline{r}_{1}^{2}(\orr_{1,\epsilon})$, generating   \newline 
    $\upi_{2}^{G}\kH/\mbox{torsion}=\widehat{\oBox } $ \\
\hline $u_{2\sigma}\in \EE_{2}^{0,2-2\sigma} (G/G)$ with
   &Element corresponding to
\newline  $\phantom{aaaaa} u_{2\sigma}\in \upi_{2-2\sigma}\HZZ (G/G)$\\
$d_{5} (u_{2\sigma}) = a_{\sigma}^{3}a_{\lambda}\normrbar_{1}^{}$
   &Slice differential of (\ref{eq-slicediffs2}) \\
$[2u_{2\sigma}]
         =\langle 2,\,a_{\sigma},\,a_{\sigma}^{2}a_{\lambda}\normrbar_{1}
                  \rangle$
\newline $\phantom{[2u_{2\sigma}]} \in \EE_{6}^{0,2-2\sigma} (G/G)$
   &Image of $2u_{2\sigma}$ in $\EE_{6}^{0,2-2\sigma} (G/G)$, 
which is a permanent cycle\\
$[u_{2\sigma}^{2}]
      =\langle a_{\sigma}^{3}a_{\lambda} ,\,\normrbar_{1} ,\, 
               a_{\sigma}^{3}a_{\lambda} ,\,\normrbar_{1} \rangle$\newline 
$\phantom{[u_{2\sigma}^{2}]} \in \EE_{6}^{0,4-4\sigma} (G/G)$
   &Image of $u_{2\sigma}^{2}$ in $\EE_{6}^{0,2-2\sigma} (G/G)$, 
which is a permanent cycle\\
\hline $u_{\sigma}\in \upi_{1-\sigma}\kH (G/G') $ 
 \newline  $\cong  \upi_{G',0}\kH (G/G)$ with
\newline    $\Res_{2}^{4}(u_{2\sigma})=u_{\sigma}^{2}$,
\newline    $\gamma (u_{\sigma})=-u_{\sigma}$ 
   & Isomorphic image of
     \newline  $\phantom{aaaaa}
1\in \upi_{0}\kH (G/G')\cong  \upi_{G',0}\kH (G/G)$\\
$\Tr_{2}^{4}(u_{\sigma}^{4k+1})=a_{\sigma}f_{1}u_{2\sigma}^{2k}$
\newline $\phantom{aaaaa}$(exotic transfer)
\newline $\Tr_{2}^{4}(u_{\sigma}^{2k})=2u_{2\sigma}^{k}$  
\newline $\Tr_{2}^{4}(u_{\sigma}^{4k+3})=0$ 
    &Follows from Theorem \ref{thm-exotic} and $d_{5} (u_{2\sigma})$ in
 (\ref{eq-slicediffs2}) \\
\hline $u_{\lambda}\in \EE_{2}^{0,2-\lambda} (G/G)$ with
    &Element corresponding to\newline
     $\phantom{aaaaa}u_{\lambda}\in \upi_{2-\lambda}\HZZ (G/G)$\\
$[2u_{\lambda}]\in \upi_{2-\lambda}\KH (G/G)$
    &$\langle 2,\,\eta  ,\,a_{\lambda} \rangle$\\
$a_{\sigma}^{3}u_{\lambda}=0$
    &Follows from the gold relation,\newline  Lemma \ref{lem-aeu}(vii)\\
$d_{3} (u_{\lambda})=\eta a_{\lambda}
         =\Tr_{2}^{4}([a_{\sigma_{2}}^{3}\orr_{1,0}])$ 
    &Slice differential of Theorem \ref{thm-d3ulambda}\\
$d_{5} ([u_{\lambda}^{2}]) =\alphax  a_{\lambda}^{2}\normrbar_{1}$ 
    &Slice differential of Theorem \ref{thm-d3ulambda}\\
$d_{7} ([2u_{\lambda}^{2}]) =\eta ' a_{\lambda}^{3}\normrbar_{1}$ 
    & $2\alphax  a_{\lambda}^{2}\normrbar_{1}$\\
$[4u_{\lambda}^{2}]\in \upi_{4-2\lambda}\KH (G/G)$
    &$\langle 2,\,\eta  ,\,a_{\lambda} \rangle^{2}
        =\langle 2,\,\eta ' ,\,a_{\lambda}^{3}\normrbar_{1} \rangle$\\
$[2a_{\sigma}u_{\lambda}^{2}]\in \upi_{4-\sigma -2\lambda}\KH (G/G)$
    &$\langle a_{\sigma},\,\eta ' ,\,a_{\lambda}^{3}\normrbar_{1} \rangle$\\
$d_{7} ([u_{\lambda}^{4}]) 
     =\langle \eta',\,\alphax  ,\,a_{\lambda}^{2}\normrbar_{1}\rangle
      a_{\lambda}^{3}\normrbar_{1}$ 
   &$[2u_{\lambda}^{2}d (u_{\lambda}^{2})]$\\
$[2u_{\lambda}^{4}]\in \upi_{8-4\lambda}\KH (G/G)$
   &$\Tr_{2}^{4}(\ou_{\lambda}^{4} )$\\
\hline 
$\ou _{\lambda}\in \EE_{2}^{0,2-\lambda}
(G/G')$ with 
    &$\Res_{2}^{4}(u_{\lambda})$\\
$d_{3} (\ou_{\lambda} )
   = [a_{\sigma_{2}}^{3} (\orr_{1,0}+\orr_{1,1})]$\newline 
$\phantom{d_{3} (\ou_{\lambda} )}
    =\Res_{2}^{4}(a_{\lambda })(\eta_{0}+\eta_{1})$
    &$\Res_{2}^{4}(d_{3} ( u_{\lambda}))$\\
$[2\ou_{\lambda}]\in \upi_{2-\lambda}\KH (G/G') $
    &$[\langle 2,\,a_{\sigma_{2}}^{3} ,\,\orr_{1,0}+\orr_{1,1} \rangle]
       =\langle 2,\,[a_{\sigma_{2}}^{2}] ,\,\eta_{0}+\eta_{1} \rangle$\\
$d_{7} ([\ou_{\lambda}^{2}] )
   = a_{\sigma_{2}}^{7}\orr_{1,0}^{3}$
    &$\Res_{2}^{4}(d_{5} ( u_{\lambda}^{2}))$\\
$[2\ou_{\lambda}^{2}]\in \upi_{4-2\lambda}\KH (G/G') $
    &$[\langle 2,\,a_{\sigma_{2}}^{7} ,\,\orr_{1,0}^{3} \rangle]
      =\langle 2,\,[a_{\sigma_{2}}^{2}]^{2} ,\,\eta_{0}^{3} \rangle$\\
$[\ou_{\lambda}^{4}]\in \upi_{8-4\lambda}\KH (G/G') $
    &$[\langle a_{\sigma_{2}}^{7},\,\orr_{1,0}^{3} ,\, 
         a_{\sigma_{2}}^{7},\,\orr_{1,0}^{3}\rangle]$\newline 
   $\phantom{a_{\sigma_{2}}^{7},\,\orr_{1,0}^{3}}
       =\langle [a_{\sigma_{2}}^{2}]^{2} ,\,\eta_{0}^{3},\,
                    [a_{\sigma_{2}}^{2}]^{2} ,\,\eta_{0}^{3} \rangle$\\
\hline
$u_{\sigma_{2}}\in \upi_{(G',1-\sigma_{2})}\kH (G/e)$ with\newline 
   $\Res_{1}^{2}(\ou_{\lambda} )=u_{\sigma_{2}}^{2}$,\newline 
   $\gamma^{2} (u_{\sigma_{2}})=-u_{\sigma_{2}}$ and \newline 
$\Tr_{1}^{2}(u_{\sigma_{2}})=a_{\sigma_{2}}^{2}(\orr_{1,0}+\orr_{1,1})$\newline 
    $\phantom{aaaaa}$(exotic transfer)
   & Isomorphic image of $1\in \upi_{0}\kH (G/e)$\\
\hline 
$\os_{2,\epsilon }
      \in \upi_{\rho_{4}}^{G}\kH (G/G')$
   &$u_{\sigma}^{-1}[\orr_{1,\epsilon }^{2}]$\\
\hline 
$\normrbar_{1}
      \in \upi_{\rho_{4}}^{G}\kH (G/G)$ with \newline 
$\Res_{2}^{4}(\normrbar_{1})
=u_{\sigma}^{-1}[\orr_{1,0}\orr_{1,1}] $
   &Image from (\ref{eq-normrbar}) defined in \cite[\S9.1]{HHR}\\
\hline
$\ot_{2}
 \in \upi_{\rho_{4}}^{G}\kH (G/G) $ with\newline 
$\Res_{2}^{4}(\ot_{2})= \os_{2,0}-\os_{2,1}$
   &$(-1)^{\epsilon }\Tr_{2}^{4}(\os_{2,\epsilon })$
            for either value of $\epsilon $\\
\hline $\ot_{2}' \in \upi_{2+\lambda}^{G}\kH (G/G) $ with\newline
$\Res_{2}^{4}(\ot_{2}')= [\orr_{1,0}^{2}]+[\orr_{1,1 }^{2}]$
&$\Tr_{2}^{4}([\orr_{1,\epsilon }^{2}])$
        for either value of $\epsilon $\\
\hline $D\in \upi_{4\rho_{4}}\kH (G/G)$, \newline 
   the  periodicity element
   &$-\normrbar_{1}^{2} (5\ot_{2}^{2}-20 \ot_{2}\normrbar_{1}
                         +11\normrbar_{1}^{2})$\\
\hline 
$\Sigma_{2, \epsilon}
      \in \EE_{2}^{0,4}\kH (G/G')$ with \newline 
     $\Sigma_{2,2}=\Sigma_{2,0}$ and\newline 
 $d_{3} (\Sigma_{2, \epsilon})  = \eta_{\epsilon } ^{2} (\eta_{0}+\eta_{1} )$
   &$(-1)^{\epsilon } u_{\rho_{4}}\os_{2, \epsilon}
    = (-1)^{\epsilon }\ou_{\lambda} [\orr_{1,\epsilon }^{2}]$\\
\hline 
$T_{2}\in \EE_{2}^{0,4}\kH (G/G)$ with \newline 
      ${\Res_{2}^{4}(T_{2})= \Sigma_{2,0}-\Sigma_{2,1}}$ and\newline 
       $d_{3} (T_{2}) =\eta^{3}$
   &$\Tr_{2}^{4}(\Sigma_{2,\epsilon})
= (-1)^{\epsilon }u_{\lambda}
    \Tr_{2}^{4}([\orr_{1,\epsilon }^{2}])$\newline 
           for either value of $\epsilon $\\
\hline 
$T_{4}\in \EE_{2}^{0,8}\kH (G/G)$ with \newline 
    $T_{4}^{2}=\Delta_{1} (T_{2}^{2}-4\Delta_{1})$,\newline
$\Res_{2}^{4}(T_{4})= (\Sigma_{2,0}-\Sigma_{2,1})\delta_{1}$  and \newline
$d_{3} (T_{4})=0$
   &$(-1)^{\epsilon }\Tr_{2}^{4}(\Sigma_{2,\epsilon}\delta_{1})
= u_{2\sigma}u_{\lambda}^{2}\ot_{2}\normrbar_{1}^{}$\newline 
           for either value of $\epsilon $  \\
\hline 
$\delta_{1}\in \EE_{2}^{0,4}\kH (G/G')$ with\newline 
$\gamma (\delta_{1})=-\delta_{1}$, $\Tr_{2}^{4}(\delta_{1})=0$ \newline 
 and 
 $d_{3} (\delta_{1})=\eta_{0}\eta_{1} (\eta_{0}+\eta_{1})$
   &$u_{\rho_{4}}\Res_{2}^{4}(\normrbar_{1})
        =\ou_{\lambda}[\orr_{1,0}\orr_{1,1}]$\\
\hline 
$\Delta_{1}\in \EE_{2}^{0,8}\kH (G/G)$ with \newline 
$\Res_{2}^{4}(\Delta_{1})=\delta_{1}^{2}$, \newline 
$\Res_{1}^{4}(\Delta_{1})=r_{1,0}^{2}r_{1,1}^{2}$ and \newline 
$d_{5} (\Delta_{1})=\nu x_{4}$
   &$u_{2\rho_{4}}\normrbar_{1}^{2}
         = u_{2\sigma}u_{\lambda}^{2}\normrbar_{1}^{2}$\\
\hline \hline 
\multicolumn{2}{|c|}{Filtration 1}\\
\hline 
$a_{\sigma_{2}}
\in \upi_{G', -\sigma_{2}}\kH (G/G)$\newline 
$\phantom{a_{\sigma_{2}}}\cong \upi_{-\sigma_{2}}^{G'}\kH (G'/G')$\newline   
with $2a_{\sigma_{2}} =0$
   & See Definition \ref{def-aeu} \\
\hline 
$\eta_{\epsilon }
\in \upi_{1}\kH (G/G')$\newline 
$\phantom{\eta_{\epsilon }}\cong \upi^{G'}_{1}\kH (G'/G')$
 with $2\eta_{\epsilon } =0$
   & $[a_{\sigma_{2}}\orr_{1,\epsilon} ]$\\
\hline 
$\eta \in \upi^{G}_{1}\kH (G/G)$ with\newline  
$\Res_{2}^{4}(\eta )= \eta_{0}+\eta_{1}$ \newline 
$\phantom{\Res_{2}^{4}(\eta )}\in \upi^{G}_{1}\kH (G/G')$ 
   & $\Tr_{2}^{4}(\eta_{\epsilon })
    =\Tr_{2}^{4}([a_{\sigma_{2}}\orr_{1,0}])
    =\Tr_{2}^{4}([a_{\sigma_{2}}\orr_{1,1}])$ \\
\hline 
$\eta'\in \underline{\pi }_{2-\sigma}\kH (G/G)$ with\newline 
$\Res_{2}^{4}(\eta')
=u_{\sigma}(\eta_{0}+\eta_{1})$
                  \newline 
$\phantom{\Res_{2}^{4}(\eta')}\in \upi^{G}_{2-\sigma}\kH (G/G')$
   &$\Tr_{2}^{4}(\eta_{0}u_{\sigma})
    =\Tr_{2}^{4}([a_{\sigma_{2}}\orr_{1,0}]u_{\sigma})$\newline 
    $\phantom{\Tr_{2}^{4}(\eta_{0}u_{\sigma})}
    =\Tr_{2}^{4}([a_{\sigma_{2}}\orr_{1,1}]u_{\sigma})$\\
\hline 
$\alphax \in \upi_{2-\sigma -\lambda}\kH (G/G)$ with
    &$[a_{\sigma}u_{\lambda} ]
     =\langle a_{\sigma},\,\eta ,\,a_{\lambda}\rangle$\\
$\Res_{2}^{4}(\alphax  ) = u_{\sigma}[a_{\sigma_{2} }^{3}\orr_{1,0}]$
\newline $\phantom{aaaaa}$ (exotic restriction)
    &Follows from Theorem \ref{thm-exotic} and $d_{3} (u_{\lambda})$ in
\newline  Theorem \ref{thm-d3ulambda}  \\
$2\alphax  =\eta 'a_{\lambda}$
    & Transfer of the above\\
$\eta \alphax  =a_{\lambda}\langle 2,\,a_{\sigma} ,\,f_{1}^{2} \rangle$\newline 
$\phantom{\eta \alphax }
    =a_{\lambda}\langle 2,\,a_{\sigma} ,\,\Tr_{2}^{4}(\eta_{0}\eta_{1}) \rangle$
    &\\
$\eta' \alphax  =0$
    &\\
$a_{\sigma } \alphax  =a_{\lambda }\Tr_{2}^{4}(u_{\sigma }^{2})$
    &$[a_{\sigma }^{2}u_{\lambda }]
          = a_{\lambda }[2u_{2\sigma }]$ 
        by the gold relation, Lemma \ref{lem-aeu}(vii)\\
$a_{\sigma }^{2} \alphax  =0$
    &$[a_{\sigma }^{3}u_{\lambda }]
          =a_{\lambda }a_{\sigma }\Tr_{2}^{4}(u_{\sigma }^{2})=0$\\
\hline 
$\betax \in \upi_{4-3\sigma -\lambda}\kH (G/G)$ with
    &$[\alphax  u_{2\sigma}]
           =\langle\alphax  ,\,a_{\sigma}^{2} ,\,f_{1} \rangle$\\
$\Res_{2}^{4}(\betax )=a_{\sigma_{2}}^{3}u_{\sigma}^{3}\orr_{1,0}$
    &Follows from value of $\Res_{2}^{4}(\alphax  )$\\
$2\betax=a_{\lambda}\langle \eta ',\,a_{\sigma}^{2} ,\,f_{1} \rangle$
    & Transfer of the above\\
$d_{5} (u_{2\sigma}u_{\lambda}^{2}) = \betax a_{\lambda}^{2}\normrbar_{1}$
    &\\
$\eta \betax = 2a_{\lambda}u_{2\sigma}^{2}\normrbar_{1}^{}$
\newline $\phantom{aaaaa}$ (exotic multiplication)
    &\\
$\eta' \betax = a_{\sigma}^{2}a_{\lambda}^{3}u_{2\sigma}^{2}\normrbar_{1}^{2}$
\newline $\phantom{aaaaa}$ (exotic multiplication)
    &\\
\hline $\nu \in \upi_{3}\kH (G/G)$ with
    &$a_{\sigma}u_{\lambda}\normrbar_{1}=\alphax \normrbar_{1}$, generating
     $\circ=\upi_{3}\kH$\\
$\Res_{2}^{4}(\nu)=\eta_{0}^{3}$
and $2\nu =x_{3}$ 
\newline (exotic restriction
\newline and group extension) 
    &
     \newline Follows from those on $\alphax $\\
\hline 
 \hline 
\multicolumn{2}{|c|}{Filtration 2}\\
\hline $[a_{\sigma_{2}}^{2}]\in \upi_{-\lambda }\kH (G/G')$
&Preimage of $a_{\sigma_{2}}^{2}
    \in \upi_{-2\sigma_{2}}i_{G'}^{*}\kH (G'/G')$\\
\hline $a_{\lambda }\in \upi_{-\lambda }\kH (G/G)$ with\newline 
$4a_{\lambda }=0$ and 
$\Res_{2}^{4}(a_{\lambda })=[a_{\sigma_{2}}^{2}]$
& See Definition \ref{def-aeu}\\
\hline $ \eta_{\epsilon }^{2} ,\, \eta_{0}\eta_{1}
\in \underline{\pi }^{G}_{2}\kH (G/G')$ with\newline
$\Tr_{2}^{4}(\eta_{\epsilon }^{2})= (-1)^{\epsilon }a_{\lambda}\ot_{2}'$
 and \newline
$\Tr_{2}^{4}(\eta_{0}\eta_{1})=f_{1}^{2}$ (exotic transfer)
&$u_{\sigma}[a_{\sigma_{2}}^{2}]\os_{2,\epsilon}$ and
$u_{\sigma}[a_{\sigma_{2}}^{2}]\Res_{2}^{4}(\normrbar_{1})$, \newline
generating the torsion $\widehat{\bullet}\oplus \JJ $ in
$\upi^{G}_{2}\kH$\\
\hline 
$\eta^{2}
=a_{\lambda} (\ot_{2}'+a_{\sigma }^{2}a_{\lambda }\normrbar_{1}^{2})
=a_{\lambda} \ot_{2}'+f_{1}^{2} $
 & $a_{\lambda}\ot_{2}'$ has order 2 by Lemma \ref{lem-hate}\\
$\eta \eta '=a_{\lambda}[u_{2\sigma} \ot_{2}]$\newline  
$(\eta ')^{2}=a_{\lambda}[u_{2\sigma}\ot_{2}']$
&See (\ref{eq-order2}) for the definition of
$[u_{2\sigma} \ot_{2}]$ and $[u_{2\sigma} \ot_{2}']$ 
\\
\hline $\nu^{2}\in \upi_{6}\kH (G/G)$
    &$2a_{\lambda}u_{\lambda}u_{2\sigma}\normrbar_{1}^{2}
       =\langle 2,\,\eta  ,\,f_{1} ,\,f_{1}^{2} \rangle$\\
\hline 
$\kappa \in \upi_{14}\kH (G/G)$
   &$2a_{\lambda}u_{2\sigma}^{2}u_{\lambda}^{3}\normrbar_{1}^{4}$\\
\hline  \hline 
\multicolumn{2}{|c|}{Filtration 3}\\
\hline
$f_{1} \in \upi_{1}\kH (G/G)$ 
   &$a_{\sigma}a_{\lambda}\normrbar_{1}^{}$,
     generating  the summand $\bullet$ of $\upi_{1}\kH$\\
\hline 
$\eta_{0}^{3}=\eta_{0}^{2}\eta_{1}=\eta_{0}\eta_{1}^{2}=\eta_{1}^{3}$
\newline $\phantom{ab}\in  \upi^{G}_{3}\kH (G/G')$
   & $\eta_{\epsilon }u_{\sigma}[a_{\sigma_{2}}^{2}]\Res_{2}^{4}(\normrbar_{1})
 = \eta_{\epsilon }u_{\sigma}[a_{\sigma_{2}}^{2}]\os_{2,\epsilon }$\\
\hline $x_{3}\in \upi_{3}\kH (G/G)$\newline 
 with $\Res_{2}^{4}(x_{3})=0$
   & $\Tr_{2}^{4}(\eta_{0}^{2}\eta_{1})
       =a_{\lambda} \eta' \normrbar_{1}$\\
\hline \hline 
\multicolumn{2}{|c|}{Filtration 4}\\
\hline $x_{4}\in \EE_{2}^{4,8} (G/G)$\newline  
with $d_{5} (x_{4})=f_{1}^{3}$,\newline 
$\Res_{2}^{4}(x_{4} )=(\eta_{0}\eta_{1})^{2}=\eta_{0}^{4}$ \newline 
and $2x_{4}=f_{1}\nu $
    &$a_{\lambda}^{2}u_{2\sigma
}\normrbar_{1}^{2}$\\
\hline 
$\overline{\kappa} \in \upi_{20}\kH (G/G)$
   &$a_{\lambda}^{2}u_{2\sigma}^{3}u_{\lambda}^{4}\normrbar_{1}^{6}$\\
$2\overline{\kappa }
=\Tr_{2}^{4}(u_{\sigma }
     \Res_{2}^{4}(u_{2\sigma }^{2}u_{\lambda }^{5}\normrbar_{1}^{5}))$
\newline 
$\phantom{aaaaa}$(exotic transfer)
   &\\  
\hline  \hline 
\multicolumn{2}{|c|}{Filtration 8}\\
\hline $\epsilon \in \upi_{8}\kH (G/G)$
    &
     $x_{4}^{2}=\langle f_{1},\, f_{1}^{2} ,\, f_{1} ,\,f_{1}^{2} \rangle
        \in \EE_{6}^{8,16} (G/G)$\\
\hline  \hline 
\multicolumn{2}{|c|}{Filtration 11}\\
\hline $\nu^{3}=\eta \epsilon \in \upi_{9}\kH (G/G)$
    &Represents $f_{1}x_{4}^{2}\in \EE_{2}^{11,20} (G/G)$\\
\hline
\end{longtable}
\end{center}

\section{Slices for $\kH$ and $\KH$}\label{sec-slices}
 
In this section we will identify the slices for $\kH$ and $\KH$ and
the generators of their integrally graded homotopy groups.  For the
latter we will use the notation of Table \ref{tab-pi*}.  Let
\begin{numequation}\label{eq-Xmn}
\begin{split}
X_{m,n}=\mycases{
\Sigma^{m\rho_{4}}\HZZ
       &\mbox{for }m=n\\
G_{+}\smashove{G'}\Sigma^{(m+n)\rho_{2}}\HZZ
       &\mbox{for }m<n.
}
\end{split}
\end{numequation}%

\noindent The slices of $\kH$ are certain finite wedges of these, and
those of $\KH$ are a certain infinite wedges.  Fortunately we can
analyze these slices by considering just one value of $m$ at a time,
this index being preserved by the first differential $d_{3}$.  These
are illustrated below in Figures \ref{fig-sseq-7a}--\ref{fig-sseq-8b}.
They show both $\EE_{2}$ and $\EE_{4}$ in four
cases depending on the sign and parity of $m$.

\begin{thm}\label{thm-sliceE2}
{\bf The slice $\EE_{2}$-term for $\kH$.}
The slices of $\kH$ are
\begin{displaymath}
P_{t}^{t}\kH=\mycases{
\bigvee_{0\leq m\leq t/4}X_{m,t/2-m}
        &\mbox{for $t$ even and $t\geq 0$}\\
*       &\mbox{otherwise}
}
\end{displaymath}

\noindent where $X_{m,n}$ is as in (\ref{eq-Xmn}).

The structure of $\pi^{u}_{*}\kH$ as a $\Z[G]$-module (see
(\ref{eq-pikH})) leads to four types of orbits and slice summands:

\begin{enumerate}
\item [(1)] $\left\{(r_{1,0}r_{1,1})^{2\ell } \right\}$ leading to
$X_{2\ell ,2\ell }$ for $\ell \geq 0$; see the leftmost diagonal in
Figure \ref{fig-sseq-7a}.  On the 0-line we have a copy of $\Box$ (defined in
Table \ref{tab-C4Mackey}) generated under restrictions by
\begin{displaymath}
\Delta_{1}^{\ell }=u_{2\ell \rho_{4}}\normrbar_{1}^{2\ell}= 
u_{2\sigma}^{\ell }u_{\lambda}^{2\ell }\normrbar_{1}^{2\ell}
\in  \EE_{2}^{0,8\ell } (G/G).
\end{displaymath}

\noindent In positive filtrations we have
\begin{align*}
\circ &\subseteq  \EE_{2}^{2j,8\ell }
\qquad \mbox{generated by } \\
a_{\lambda}^{j}u_{2\sigma}^{\ell }
          u_{\lambda}^{2\ell -j}\normrbar_{1}^{2\ell }
 & \in   \EE_{2}^{2j,8\ell } (G/G)
           \qquad \mbox{for }0<j\leq 2\ell \mbox{ and}  \\
\bullet &\subseteq  \EE_{2}^{2k+4\ell ,8\ell }
\qquad \mbox{generated by } \\
a_{\sigma}^{2k}a_{\lambda}^{2\ell }u_{2\sigma}^{\ell -k}
          \normrbar_{1}^{2\ell }
 & \in   \EE_{2}^{2k+4\ell,8\ell } (G/G)
           \qquad \mbox{for }0<k\leq \ell.  
\end{align*}

\item [(2)] $\left\{(r_{1,0}r_{1,1})^{2\ell+1 } \right\}$ leading
to $X_{2\ell +1,2\ell +1}$ for $\ell \geq 0$; see the leftmost
diagonal in Figure \ref{fig-sseq-7b}.  On the 0-line we have a copy of
$\oBox$ generated under restrictions by
\begin{align*}
\delta_{1}^{2\ell +1 }
 = u_{\sigma}^{2\ell +1}
          \Res_{2}^{4}(u_{\lambda}\normrbar_{1})^{2\ell +1}
 &\in  \EE_{2}^{0,8\ell+4 } (G/G').
\end{align*}

\noindent In positive filtrations we have
\begin{align*}
\obull  &\subseteq  \EE_{2}^{2j,8\ell+4 }
\qquad \mbox{generated by } \\
u_{\sigma}^{2\ell +1}
     \Res_{2}^{4}(a_{\lambda}^{j}u_{\lambda}^{2\ell +1-j}
                  \normrbar_{1}^{2\ell +1})
 & \in   \EE_{2}^{2j,8\ell+4 } (G/G')
           \qquad \mbox{for }0<j\leq 2\ell+1 , \\
\bullet &\subseteq  \EE_{2}^{2j+1 ,8\ell+4 }
\qquad \mbox{generated by } \\
 a_{\sigma}a_{\lambda}^{j}
       u_{2\sigma}^{\ell}u_{\lambda}^{2\ell +1-j}
                  \normrbar_{1}^{2\ell +1}
 & \in   \EE_{2}^{2j+1,8\ell+4 } (G/G)
           \qquad \mbox{for }0\leq j\leq 2\ell+1\mbox{ and}  \\
\bullet &\subseteq  \EE_{2}^{2k+4\ell +3 ,8\ell+4 }
\qquad \mbox{generated by } \\
a_{\sigma}^{2k+1}a_{\lambda}^{2\ell +1}
       u_{2\sigma}^{\ell-k}
                  \normrbar_{1}^{2\ell +1}
 & \in   \EE_{2}^{2k+4\ell +3,8\ell+4 } (G/G)
           \qquad \mbox{for }0< k\leq \ell.  
\end{align*}

\item [(3)] $\left\{r_{1,0}^{i}r_{1,1}^{2\ell -i}, r_{1,0}^{2\ell
-i}r_{1,1}^{i} \right\}$ leading to $X_{i,2\ell -i}$ for $0\leq i<\ell
$; see other diagonals in Figure \ref{fig-sseq-7a}.  On the 0-line we
have a copy of $\widehat{\Box}$ generated (under $\Tr_{2}^{4}$,
$\Res_{1}^{2}$ and the group action) by

\begin{displaymath}
u_{\sigma  }^{\ell } \os_{2}^{\ell -i}
\Res_{2}^{4}(u_{\lambda}^{\ell }\normrbar_{1}^{i}) 
\in \EE_{2}^{0,4\ell } (G/G') 
\end{displaymath}

In positive filtrations we have
\begin{align*}
\widehat{\bullet}  &\subseteq  \EE_{2}^{2j,4\ell }
\qquad \mbox{generated by } \\
\lefteqn{u_{\sigma  }^{\ell } \os_{2}^{\ell -i}
\Res_{2}^{4}(a_{\lambda}^{j}u_{\lambda}^{\ell-j }\normrbar_{1}^{i}) 
}\qquad\qquad\\
  &\in \EE_{2}^{2j,4\ell } (G/G')\qquad \mbox{for }0<j\leq \ell  \\
  & =  \eta_{\epsilon }^{2j}u_{\sigma}^{\ell -j}\os_{2}^{\ell -i-j}
\Res_{2}^{4}(u_{\lambda}^{\ell -j}\normrbar_{1}^{i}) 
       \qquad \mbox{for }0 < j < \ell -i.
\end{align*}

\item [(4)] $\left\{r_{1,0}^{i}r_{1,1}^{2\ell+1 -i},
r_{1,0}^{2\ell+1 -i}r_{1,1}^{i} \right\}$ leading to
$X_{i,2\ell+1 -i}$ for $0\leq i\leq \ell $; see other diagonals in
Figure \ref{fig-sseq-7b}.  On the 0-line we have a copy of
$\widehat{\oBox}$ generated (under transfers and the group
action) by
\begin{displaymath}
r_{1,0}\Res_{1}^{2}(u_{\sigma  }^{\ell } \os_{2}^{\ell -i})
\Res_{1}^{4}(u_{\lambda}^{\ell }\normrbar_{1}^{i}) 
\in \EE_{2}^{0,4\ell +2} (G/\ee) 
\end{displaymath}

\noindent In positive filtrations we have
\begin{align*}
\widehat{\bullet}  &\subseteq  \EE_{2}^{2j+1,4\ell+2 }
\qquad \mbox{generated by } \\
\lefteqn{\eta_{\epsilon } u_{\sigma  }^{\ell } \os_{2}^{\ell -i}
\Res_{2}^{4}(a_{\lambda}^{j}u_{\lambda}^{\ell-j }\normrbar_{1}^{i}) 
}\qquad\qquad\\
  &\in \EE_{2}^{2j+1,4\ell+2 } (G/G')
      \qquad \mbox{for }0\leq j\leq \ell \\
  & =  \eta_{\epsilon }^{2j+1}u_{\sigma}^{\ell -j}\os_{2}^{\ell -i-j}
\Res_{2}^{4}(u_{\lambda}^{\ell -j}\normrbar_{1}^{i}) 
       \qquad \mbox{for }0 \leq  j \leq  \ell -i . 
\end{align*}
\end{enumerate}
\end{thm}

\begin{cor}\label{cor-subring}
{\bf A subring of the slice $E_{2}$-term.}  The ring
$\EE_{2}\kH (G/G')$ contains 
\begin{displaymath}
\Z[\delta_{1},\Sigma_{2,\epsilon},\eta_{\epsilon }
      \colon \epsilon =0,\,1]/
\left(2\eta _{\epsilon },\delta_{1}^{2}-\Sigma_{2,0}\Sigma_{2,1},
\eta_{\epsilon } \Sigma_{2,\epsilon +1}+\eta _{1+\epsilon }\delta_{1} \right);
\end{displaymath}

\noindent see Table \ref{tab-pi*} for the
definitions of its generators.  In particular the elements $\eta_{0} $
and $\eta_{1}$ are algebraically independent mod 2 with
\begin{displaymath}
\gamma^{\epsilon } (\eta_{0}^{m}\eta_{1} ^{n})
    \in \upi_{m+n}X_{m,n} (G/G')\qquad \mbox{for }m\leq n. 
\end{displaymath}

\noindent
The element $(\eta_{0}\eta_{1})^{2}$ is the fixed point restriction of 
\begin{displaymath}
u_{2\sigma}a_{\lambda}^{2}\normrbar_{1}^{2}
\in \EE_{2}^{4,8}\kH (G/G),
\end{displaymath}

\noindent which has order 4, and the transfer of the former is twice
the latter.  The element $\eta_{0}\eta_{1}$ is not in the image of
$\Res_{2}^{4}$ and has trivial transfer in $\EE_{2}$.
\end{cor}

\proof
We detect this subring with the monomorphism
\begin{displaymath}
\xymatrix
@R=1mm
@C=10mm
{
{\EE_{2}\kH (G/G')}\ar[r]^(.5){\ur_{2}^{4}}
   &{\EE_{2}\kH (G'/G')}\\
\eta_{\epsilon}\ar@{|->}[r]
   &{\,a_{\sigma}\orr_{1,\epsilon} }\\
\Sigma_{2,\epsilon}\ar@{|->}[r]
   &{\,u_{2\sigma}\orr_{1,\epsilon} ^{2}}\\
\delta_{1}\ar@{|->}[r]
   &{\,u_{2\sigma}\orr_{1,0}\orr_{1,1} },
}
\end{displaymath}

\noindent in which all the relations are transparent. 
\qed\bigskip 

\begin{cor}\label{cor-KH}
{\bf Slices for $\KH$.}  The slices of $\KH$ are
\begin{displaymath}
P_{t}^{t}\KH=\mycases{
\bigvee_{m\leq t/4}X_{m,t/2-m}
        &\mbox{for $t$ even}\\
*       &\mbox{otherwise}
}
\end{displaymath}

\noindent where $X_{m,n}$ is as in Theorem \ref{thm-sliceE2}.  Here
$m$ can be any integer, and we still require that $m\leq n$.
\end{cor}

\proof Recall that $\KH$ is obtained from $\kH$ by inverting a certain
element 
\begin{displaymath}
{D\in \upi_{4\rho_{4}}\kH (G/G)}
\end{displaymath}

\noindent described in Table \ref{tab-pi*}.
Thus $\KH$ is the homotopy colimit of the diagram 
\begin{displaymath}
\xymatrix
@R=5mm
@C=10mm
{
k_{[2]}\ar[r]^(.4){D}
   &\Sigma^{-4\rho_{4}}\kH\ar[r]^(.5){D}
      &\Sigma^{-8\rho_{4}}\kH\ar[r]^(.6){D}
         &\dotsb 
}
\end{displaymath}

\noindent Desuspending by $4\rho_{4}$ converts slices to slices, so
for even $t$ we have
\begin{align*}
P_{t}^{t}\KH
 & =   \lim_{k\to \infty }\Sigma^{-4k\rho_{4}}P^{t+16k}_{t+16k}\kH \\
 & =   \lim_{k\to \infty }\Sigma^{-4k\rho_{4}}
            \bigvee_{0\leq m\leq t/4+4k}X_{m,t/2+8k-m} \\
 & =   \lim_{k\to \infty }
            \bigvee_{0\leq m\leq t/4+4k}X_{m-4k,t/2+4k-m} \\
 & =   \lim_{k\to \infty }
            \bigvee_{-4k\leq m\leq t/4}X_{m,t/2-m} \\
 & =   \bigvee_{m\leq t/4}X_{m,t/2-m}.\hfill \qed 
\end{align*}
\bigskip
 
\begin{cor}\label{cor-filtration}
{\bf A filtration of $\kH$. }
Consider the diagram
\begin{displaymath}
\xymatrix
@R=5mm
@C=10mm
{
{\kH}\ar[d]^(.5){}
    &{\Sigma^{\rho_{4}}\kH}\ar[d]^(.5){} \ar[l]_(.5){\normrbar_{1}}
        &{\Sigma^{2\rho_{4}}\kH}\ar[d]^(.5){}\ar[l]_(.45){\normrbar_{1}}
            &\dotsb \ar[l]_(.4){\normrbar_{1}}\\
y_{0}
    &y_{1}=\Sigma^{\rho_{4}}y_{0}
        &y_{2}=\Sigma^{2\rho_{4}}y_{0}
} 
\end{displaymath}

\noindent where $y_{0}$ is the cofiber of the map induced by
$\normrbar_{1}$.  Then the slices of $y_{m}$ are
\begin{displaymath}
P_{t}^{t}y_{m}=\mycases{
X_{m,t/2-m}
        &\mbox{for $t$ even and $t\geq 4m$}\\
*       &\mbox{otherwise.}
}
\end{displaymath}
\end{cor}


\begin{cor}\label{cor-Filtration}
{\bf A filtration of $\KH$.}  Let $R=\Z_{(2)}[x]/ (11x^{2}-20x+5)$.
After tensoring with $R$ (by smashing with a suitable Moore spectrum
$M$) there is a diagram
\begin{displaymath}
\xymatrix
@R=5mm
@C=10mm
{\dotsb \ar[r]^(.5){}
    &\Sigma^{2\rho_{4}}\kH\ar[r]^(.5){f_{2}}\ar[d]^(.5){}
        &\Sigma^{\rho_{4}}\kH\ar[r]^(.5){f_{1}}\ar[d]^(.5){}
            &k_{[2]}\ar[r]^(.4){f_{0}}\ar[d]^(.5){}
                &\Sigma^{-\rho_{4}}\kH\ar[r]^(.6){f_{-1}}\ar[d]^(.5){}
                    &\dotsb 
\\
    &Y_{2}
        &Y_{1}
            &Y_{0}
                &Y_{-1}
                    &
}
\end{displaymath}

\noindent where the homotopy colimit of the top row is $\KH$ and each
$Y_{m}$ has slices similar to those of $y_{m}$ as in Corollary
\ref{cor-filtration}.
\end{cor}

\proof The periodicity element $D=-\normrbar_{1}^{2} (5\ot_{2}^{2}-20
\ot_{2}\normrbar_{1} +11\normrbar_{1}^{2})$ can be factored as 
\begin{displaymath}
D=D_{0}D_{1}D_{2}D_{3}
\end{displaymath}

\noindent where $D_{i}=a_{i}\normrbar_{1}+b_{i}\ot_{2}$ with $a_{i}\in
\Z_{(2)}^{\times }$ and $b_{i}\in R$.  Then
let $f_{4n+i}$ be multiplication by $D_{i}$.  It follows that the
composite of any four successive $f_{m}$s is $D$, making the colimit
$\KH$ as desired.  The fact that $a_{i}$ is a unit means that the
$Y$'s here have the same slices as the $y$'s in Corollary
\ref{cor-filtration}.  \qed\bigskip

\begin{rem}\label{rem-Witt}
The 2-adic completion of $R$ is
the Witt ring $W (\F{4})$ used in Morava $E_{2}$-theory.  This follows
from the fact that the roots of the quadratic polynomial involve
$\,\sqrt[]{5}$, which is in $W (\F{4})$ but is not a 2-adic integer.

Moreover if we assume that $D_{0}D_{1}= 5\ot_{2}^{2}-20
\ot_{2}\normrbar_{1} +11\normrbar_{1}^{2}$, then the composite maps
$f_{4n}f_{4n+1}$, as well as $f_{4n+2}$ and $f_{4n+3}$, can be
constructed without adjoining $\,\sqrt[]{5}$.
\end{rem}



It turns out that $y_{m}\wedge M$ and $Y_{m}$ for $m\geq 0$ not only have the
same slices, but the same slice spectral sequence, which is shown in
Figures \ref{fig-sseq-7a}--\ref{fig-sseq-8b}. See Remark \ref{rem-Ym}
below.  We do not know if they have the same homotopy type.

\section{Some differentials in the slice {\SS} for $\kH$}\label{sec-diffs}

Now we turn to differentials.  The only
generators in (\ref{eq-bigraded}) that are not permanent cycles are
the $u$'s.  We will see that it is easy to account for the elements in
$\EE_{2}^{0,|V|-V} (G/H)$ for proper subgroups $H$ of
$G=C_{4}$.  From (\ref{eq-bigraded}) we see that
\begin{numequation}\label{eq-sparse}
\EE_{2}^{s,t}= 0 \qquad \mbox{for $|t|$ odd.} 
\end{numequation}%

\noindent This sparseness condition implies that $d_{r}$ can be
nontrivial only for odd values of $r$.

Our starting point is the Slice Differentials Theorem of
\cite[Thm. 9.9]{HHR}, which is derived from the fact that the geometric
fixed point spectrum of $MU^{((G))}$ is $MO$.  It says that in the
slice {\SS} for $MU^{((G))}$ for an arbitrary finite cyclic 2-group
$G$ of order $g$, the first nontrivial differential on various powers
of $u_{2\sigma}$ is
\begin{numequation}
\label{eq-slicediffs}
d_{r} (u_{2\sigma}^{2^{k-1}})
 = a_{\sigma}^{2^{k}}a_{\orho}^{2^{k}-1}
                 N_{2}^{g} (\orr_{2^{k}-1}^{G})
\in \EE_{r}^{r,r+2^{k} (1-\sigma )-1}MU^{((G))}(G/G),
\end{numequation}%

\noindent where $r=1+(2^{k}-1)g$ and $\orho$ is the
reduced regular {\rep} of $G$.

In particular
\begin{numequation}
\label{eq-slicediffs2}
\begin{split}
\left\{
\begin{array}[]{rlll}
d_{5} (u_{2\sigma}) 
 &\hspace{-2.5mm}
   =  a_{\sigma}^{3}a_{\lambda}\normrbar_{1}
     &\hspace{-2.5mm}\in \EE_{5}^{5,6-2\sigma}MU^{((G))} (G/G)
           &\mbox{for }G=C_{4}\\  
d_{13} ([u_{2\sigma}^{2}]) 
 &\hspace{-2.5mm} 
    =  a_{\sigma}^{7}a_{\lambda}^{3}\normrbar_{2}
     &\hspace{-2.5mm}\in \EE_{13}^{13,16-4\sigma}MU^{((G))} (G/G)
           &\mbox{for }G=C_{4}\\  
d_{3} (u_{2\sigma}) 
 &\hspace{-2.5mm}
   =  a_{\sigma}^{3}\orr_{1} 
     &\hspace{-2.5mm}\in \EE_{3}^{3,4-2\sigma}MU_{\reals} (G/G)
           &\mbox{for }G=C_{2}\\  
d_{7} ([u_{2\sigma}^{2}]) 
 &\hspace{-2.5mm} 
   =  a_{\sigma}^{7}\orr_{3} 
     &\hspace{-2.5mm}\in \EE_{3}^{7,10-4\sigma}MU_{\reals} (G/G)
           &\mbox{for }G=C_{2}.
\end{array}
\right.
\end{split}
\end{numequation}%

The first of these leads directly to a similar differential in the
slice {\SS } for $\kH$.  The target of the second one has
trivial image in $\kH$ and we shall see that $[u_{2\sigma }^{2}]$
turns out to be a permanent cycle.

There are two ways to leverage the third and fourth differentials of
(\ref{eq-slicediffs2}) into information about $\kH$. 
\begin{enumerate}
\item [(i)] They both lead to differentials in the slice {\SS} for the
$C_{2}$ spectrum $i^{*}_{G'}\kH$.  They are spelled out ind 
(\ref{eq-d3-and-d7}) and will be studied in detail below in
\S\ref{sec-C2diffs}.  They completely determine the slice {\SS }
$\underline{E}_{*}^{*,\star} (G/G')$ for both $\kH$ and $\KH$.

Since $u_{\lambda }$ restricts to $\ou_{\lambda }$, which is
isomorphic to $u_{2\sigma_{2}}$, we get some information about
differentials on powers of $u_{\lambda }$.  The $d_{3}$ on
$u_{2\sigma_{2}}$ forces a $d_{3} (u_{\lambda })=\eta a_{\lambda }$.
The target of $d_{7} ([u_{2\sigma_{2}}^{2}])$ turns out to be the
exotic restriction of an an element in filtration 5, leading to $d_{5}
([u_{\lambda }]^{2})=\nu a_{\lambda }^{2}$.  We will also see that
even though $[u_{2\sigma_{2}}^{4}]$ is a permanent cycle, $[u_{\lambda
}^{4}]$ (its preimage under the restrtction map $\Res_{2}^{4}$) is
not.

\item [(ii)] One can norm up the differentials on $u_{2\sigma_{2}}$
and its square using Corollary \ref{cor-normdiff}, converting the
$d_{3}$ and $d_{7}$ to a $d_{5}$ and a $d_{13}$.  The source of the
latter is $[a_{\sigma }u_{\lambda }^{4}]$, which implies that
$[u_{\lambda }^{4}]$ is not a permanent cycle.

\end{enumerate}

The differentials of (\ref{eq-slicediffs2}) lead to Massey products
which are permanent cycles,
\begin{align*}
\langle 2,\,a_{\sigma}^{2} ,\,f_{1} \rangle 
 & =  [2u_{2\sigma}]
  =  \Tr_{G'}^{G}(u_{\sigma}^{2})
 \in \mycases{
\EE_{6 }^{0,2-2\sigma}MU^{((G))} (G/G)
           &\mbox{for }G=C_{4}\\
\EE_{4}^{0,2-2\sigma}MU_{\reals} (G/G)
           &\mbox{for }G=C_{2}
}\\
\langle 2,\,a_{\sigma}^{4} ,\,f_{3} \rangle 
 & =  [2u_{2\sigma}^{2}]
  =  \Tr_{G'}^{G}(u_{\sigma}^{4})
 \in \mycases{
\EE_{14 }^{0,4-4\sigma}MU^{((G))} (G/G)
           &\mbox{for }G=C_{4}\\ 
\EE_{8}^{0,4-4\sigma}MU_{\reals} (G/G)
           &\mbox{for }G=C_{4}
}
\end{align*}

\noindent and (by Theorem \ref{thm-exotic}) to exotic transfers
\begin{numequation}\label{eq-exotic-transfers}
\begin{split}\left\{
\begin{array}[]{rl}
a_{\sigma}f_{1}
 &\hspace{-2.5mm} = \left\{
\begin{array}{lll}
\Tr_{2}^{4} (u_{\sigma}) 
    &\in\EE_{\infty }^{4,5-\sigma}MU^{((G))} (G/G)
         &\mbox{(filtration jump 4)}\\
    &     \qquad \mbox{for }G=C_{4}\\
\Tr_{1}^{2} (u_{\sigma}) 
    &\in \EE_{\infty }^{2,3-\sigma}MU_{\reals} (G/G)
        &\mbox{(filtration jump 2)}\\
    &     \qquad \mbox{for }G=C_{2}  
\end{array} \right.  \\
a_{\sigma}^{3}f_{3}
 &\hspace{-2.5mm} = \left\{
\begin{array}{lll}
\Tr_{2}^{4} (u_{\sigma}^{3}) 
    &\in \EE_{\infty }^{12,15-3\sigma}MU_{\reals} (G/G)
        &\mbox{(filtration jump 12)}\\
    &     \qquad \mbox{for }G=C_{4}\\
\Tr_{1}^{2} (u_{\sigma}^{3}) 
    &\in \EE_{\infty }^{6,9-3\sigma}MU_{\reals} (G/G)
        &\mbox{(filtration jump 6)}\\
    &     \qquad \mbox{for }G=C_{2}  
\end{array} \right. 
\end{array} \right.
\end{split}
\end{numequation}%
 
Since $a_{\sigma }$ and $2a_{\lambda }$ kill transfers by Lemma
\ref{lem-hate}, we have Massey products,
\begin{numequation}\label{eq-order2}
\begin{split}
[u_{2\sigma }\Tr_{2}^{4}(x)] = \Tr_{2}^{4}(u_{\sigma }^{2}x)
      =\langle a_{\sigma }f_{1},\,a_{\sigma } ,\, \Tr_{2}^{4}(x)\rangle
\quad \mbox{with }2a_{\lambda } [u_{2\sigma }\Tr_{2}^{4}(x)]=0.
\end{split}
\end{numequation}%

Now, as before, let $G=C_{4}$ and $G'=C_{2}\subseteq G$.  We need to
translate the $d_{3}$ above in the slice {\SS} for $MU_{\reals}$ into
a statement about the one for $\kH$ as a $G'$-spectrum.  We have an
{\eqvr} multiplication map $m$ of $G'$-spectra
\begin{displaymath}
\xymatrix
@R=2mm
@C=10mm
{
   &MU^{((G))}\ar@{=}[d]^(.5){}\\
MU_{\reals}\ar[r]^(.4){\eta_{L}}
   & MU_{\reals}\wedge MU_{\reals}\ar[r]^(.6){m}
      & MU_{\reals}\\
{\orr_{1}^{G'}}\ar@{|->}[r]^(.5){}
   &\,{\orr_{1,0}^{G}+\orr_{1,1}^{G}}\ar@{|->}[r]^(.5){}
      &{\,\orr_{1}^{G'}}\\
   &a_{\sigma}^{3} (\orr_{1,0}^{G}+\orr_{1,1}^{G})\ar@{|->}[r]^(.5){}
      &\,a_{\sigma}^{3}\orr_{1}^{G'}\\
{\orr_{3}^{G'}}\ar@{|->}[r]^(.5){}
   &{\left(\begin{array}[]{r}
5\orr_{1,0}^{G}\orr_{1,1}^{G} (\orr_{1,0}^{G}+\orr_{1,1}^{G}) 
        +(\orr_{1,1}^{G})^{3}\\ \bmod (\orr_{2}^{G},\orr_{3}^{G})
     \end{array} \right)}
           \ar@{|->}[r]^(.5){}
      &\,\orr_{3}^{G'}\\
}
\end{displaymath}

\noindent where the elements lie in
$\upi^{G'}_{\rho_{2}}(\cdot)(G'/G')$ and
$\upi^{G'}_{3\rho_{2}}(\cdot)(G'/G')$.  In the slice {\SS} for
$MU^{((G))}$ as a $G'$-spectrum, $d_{3} (u_{2\sigma})$ and
$d_{7}(u_{2\sigma}^{2})$ must be $G$-invariant since $u_{2\sigma}$ is,
and they must map respectively to $a_{\sigma}^{3}\orr_{1}^{G'}$ and
$a_{\sigma}^{7}\orr_{3}^{G'}$, so we have
\begin{numequation}\label{eq-d3-and-d7}
\begin{split}
\left\{\begin{array}{rl}
d_{3} (u_{2\sigma_{2} }) = d_{3} (\ou_{\lambda})
 &\hspace{-2.5mm} =   a_{\sigma_{2} }^{3} (\orr_{1,0}^{G}+\orr_{1,1}^{G})
  =  a_{\sigma_{2} }^{2} (\eta_{0}+\eta_{1})\\
d_{7} ([u_{2\sigma_{2} }^{2}]) = d_{7} ([\ou_{\lambda}^{2}])
 &\hspace{-2.5mm} =   a_{\sigma_{2} }^{7}\left(5\orr_{1,0}^{G}\orr_{1,1}^{G}
           (\orr_{1,0}^{G}+\orr_{1,1}^{G})
           + (\orr_{1,1}^{G})^{3}+\dotsb  \right)\\
 &\hspace{-2.5mm} =    a_{\sigma_{2} }^{7} (\orr_{1,0}^{G})^{3}+\dotsb\\
 &
     \qquad \mbox{since $a_{\sigma_{2}}^{3}(\orr_{1,0}^{G}+\orr_{1,1}^{G})=0$
            in $\EE_{4}$} 
\end{array} \right.
\end{split}
\end{numequation}%

\noindent We get similar differentials in the slice spectral sequence
for $\kH$ as a $C_{2}$-spectrum in which the missing terms in $d_{7}
(\ou_{\lambda}^{2})$ vanish.

Pulling back along the isomorphism $\ur_{2}^{4}$ gives 
\begin{numequation}\label{eq-d-ul}
\begin{split}\left\{
\begin{array}[]{rl}
d_{3} (\Res_{2}^{4}(u_{\lambda}))
    =   d_{3} (\overline{u}_{\lambda} )
 &\hspace{-2.5mm} =   [a_{\sigma_{2}}^{2}] (\eta_{0}+\eta_{1}) 
  = \Res_{2}^{4}(a_{\lambda}\eta)\\
d_{7} (\Res_{2}^{4}(u_{\lambda}^{2}))
 =    d_{7} (\overline{u}_{\lambda}^{2} )
 &\hspace{-3mm}
   =  \Res_{2}^{4}(a_{\lambda }^{2})\eta_{0}^{3}
   =  \Res_{2}^{4}(a_{\lambda }^{2}\nu )
\end{array} 
\right. 
\end{split} 
\end{numequation}%

\noindent  These imply that 
\begin{displaymath}
d_{3} (u_{\lambda})=a_{\lambda}\eta
\qquad \aand 
d_{5} (u_{\lambda }^{2}) = a_{\lambda }^{2}\nu .
\end{displaymath}

The differential on $u_{\lambda}$ leads to the following Massey
products, the second two of which are permanent cycles.
\begin{numequation}\label{eq-alphax}
\begin{split}\left\{
\begin{array}[]{rl}
\left[u_{\lambda}^{2} \right]
  =  \langle a_{\lambda },\,\eta  ,\,  a_{\lambda },\,\eta\rangle 
         &\hspace{-3mm} \in \EE_{4}^{0,4-2\lambda} (G/G)\\
\left[2u_{\lambda} \right]
  =   \langle 2,\,\eta ,\, a_{\lambda} \rangle 
         &\hspace{-3mm} \in \EE_{4}^{0,2-\lambda} (G/G) \\
\alphax  :=  [a_{\sigma}u_{\lambda}]
  =   \langle a_{\sigma},\,\eta ,\,a_{\lambda} \rangle 
         &\hspace{-3mm} \in \EE_{4}^{1,3-\sigma-\lambda} (G/G)
\end{array} 
\right.
\end{split} 
\end{numequation}

\noindent where $\alphax  $ satisfies
\begin{align*}
a_{\sigma}^{2}\alphax  
 & =   \langle a_{\sigma}^{3},\,\eta  ,\, a_{\lambda}\rangle
 =  a_{\sigma }[a_{\sigma}^{2}u_{\lambda}]
 =   a_{\sigma} [2a_{\lambda}u_{2\sigma}] 
 =   [2a_{\sigma }a_{\lambda}u_{2\sigma}]
 =  0\\
\Res_{2}^{4}(\alphax )
 & =   [a_{\sigma_{2}}^{3}\orr_{1,\epsilon}]u_{\sigma}
  \in \EE_{4}^{3,5-\sigma -\lambda} (G/G') \\
&   \qquad \mbox{(exotic restriction with filtration jump 2 by Theorem 
                  \ref{thm-exotic}(i))} \\
2\alphax & = \Tr_{2}^{4}(\Res_{2}^{4}(\alphax ))
 = \Tr_{2}^{4}(u_{\sigma}
           [a_{\sigma_{2}}^{3}\orr_{1,\epsilon}])\\
 & =   \eta' a_{\lambda}
 \in  \EE_{4}^{3,5-\sigma -\lambda} (G/G)\\
&
\qquad \mbox{(exotic group extension with jump 2)} \\
\Tr_{2}^{4}(x)\alphax  
 & =   \Tr_{2}^{4}(x \cdot \Res_{2}^{4}(\alphax  ))  
   =   \Tr_{2}^{4}(x [a_{\sigma_{2}}^{3}\orr_{1,0}]
                 u_{\sigma})\\
\eta \alphax  
 & =   \Tr_{2}^{4}([a_{\sigma_{2}}\orr_{1,1}])\alphax  
   =   \Tr_{2}^{4}([a_{\sigma_{2}}^{4}\orr_{1,0}\orr_{1,1}] 
                          u_{\sigma})
 =  a_{\lambda}^{2}\normrbar_{1}^{}\Tr_{2}^{4}(u_{\sigma}^{2})\\
 & =    a_{\lambda}\normrbar_{1}^{}
              \langle 2,\,a_{\sigma} ,\,a_{\sigma}f_{1} \rangle
 =  \langle 2,\,a_{\sigma} ,\,f_{1}^{2} \rangle\\
\eta' \alphax  
 & =   a_{\lambda}^{2}\normrbar_{1}^{}\Tr_{2}^{4}(u_{\sigma}^{3})
 =  0\\
d_{7} ([\ou_{\lambda}^{2}])
 & = [a_{\sigma_{2}}^{7}\orr_{1,0}^{3}]
         \quad \mbox{in }\EE_{4}\\
 & = \Res_{2}^{4}(\alphax )\Res_{2}^{4}(a_{\lambda}^{2}\normrbar_{1})
   = \Res_{2}^{4}(\alphax a_{\lambda}^{2}\normrbar_{1})
   = \Res_{2}^{4}(d_{5} (u_{\lambda}^{2} ))\\ 
d_{5} ([u_{\lambda}^{2}])
 & =   \alphax  a_{\lambda}^{2}\normrbar_{1}=a_{\lambda }^{2}\nu \\
d_{7} ([2u_{\lambda}^{2}])
 & =  ( 2\alphax  ) a_{\lambda}^{2}\normrbar_{1}
 =   a_{\lambda}^{3}\eta' \normrbar_{1}.
\end{align*}
 
\noindent Note that $\nu =\alphax  \normrbar_{1}^{}$, with the exotic
restriction and group extension on $\alphax  $ being consistent with
those on $\nu $.

\bigskip

The differential on $[u_{\lambda}^{2}]$ yields Massey products
\begin{numequation}\label{eq-brackets3}
\begin{split}
\left\{ 
\begin{array}[]{rl}
[a_{\sigma}^{2}u_{\lambda}^{2}]
 &\hspace{-2.5mm} =   \langle a_{\sigma}^{2},\,\alphax  
                               ,\,a_{\lambda}^{2}\normrbar_{1} \rangle\\
\left[\eta' u_{\lambda}^{2} \right]
 &\hspace{-2.5mm} =   \langle \eta ',\,\alphax  
                               ,\,a_{\lambda}^{2}\normrbar_{1} \rangle.
\end{array} \right.
\end{split}
\end{numequation}%

\begin{thm}\label{thm-normed-slice-diffs}
{\bf Normed up slice differentials for $\kH$ and $\KH$.} In the slice
{\SS}s for $\kH$ and $\KH$,
\begin{align*}
d_{5} ( [a_{\sigma }u_{\lambda }^{2}])
 & = 0  \\
\aand 
d_{13} ([a_{\sigma }u_{\lambda }^{4}])
 & = a_{\lambda }^{7}[u_{2\sigma }^{2}]\normrbar_{1}^{3} .
\end{align*}
\end{thm}

\proof
The two slice differentials over $G'$ are 
\begin{align*}
d_{3} (u_{2\sigma_{2}})
 & = a_{\sigma_{2}}^{3}\orr_{1}^{G'}
   = a_{\sigma_{2}}^{3} (\orr_{1,0}+\orr_{1,1})  \\
\aand 
d_{7} ([u_{2\sigma_{2}}^{2}])
 & = a_{\sigma_{2}}^{7}\orr_{3}^{G'}
   = a_{\sigma_{2}}^{7} ( 5 \orr_{1,0}^{2}\orr_{1,1} +5\orr_{1,0}
\orr_{1,1}^{2} +\orr_{1,1}^{3}) 
\end{align*}

\noindent We need to find the norms of both sources and targets.
Lemma \ref{lem-norm-au} tells us that
\begin{align*}
N_{2}^{4} (a_{\sigma_{2}}^{k})
 & =  a_{\lambda }^{k}\\
\aand N_{2}^{4} (u_{2\sigma_{2}}^{k}) & = u_{\lambda }^{2k}/u_{2\sigma
}^{k} \qquad \mbox{in }\underline{E}_{2}(G/G) .
\end{align*}

\noindent Previous calculations give
\begin{align*}
N_{2}^{4} (\orr_{1,0}+\orr_{1,1})
 & =  -\overline{t}_{2}\qquad \mbox{by (\ref{eq-norm-orr})} \\
\aand 
N_{2}^{4} ( 5 \orr_{1,0}^{2}\orr_{1,1} +5\orr_{1,0}\orr_{1,1}^{2}
            +\orr_{1,1}^{3}) 
 & =   -\normrbar_{1} (5\ot_{2}^{2}-20 \ot_{2}\normrbar_{1}
                        +11\normrbar_{1}^{2})\\
 & \phantom{ =  -\overline{t}_{2}}\, \qquad \mbox{by (\ref{eq-norm-quad})} ,
\end{align*}

For the first differential, Corollary \ref{cor-normdiff} tells us that
\begin{align*}
a_{\lambda }^{3}\overline{t}_{2}
 & = d_{5} (a_{\sigma }u_{\lambda }^{2}/u_{2\sigma }) \\
 & = d_{5} (a_{\sigma } u_{\lambda }^{2})/u_{2\sigma }
     -a_{\sigma }u_{\lambda }^{2}d_{5} (u_{2\sigma })/[u_{2\sigma }^{2}]\\
 & = d_{5} (a_{\sigma } u_{\lambda }^{2})/u_{2\sigma }
         -a_{\sigma }u_{\lambda }^{2}a_{\sigma }^{3}
                    a_{\lambda }\normrbar_{1}/[u_{2\sigma }^{2}]
\end{align*}

\noindent Multiplying both sides by the permanent cycle
$[u_{2\sigma}^{2}]$ gives
\begin{align*}
[u_{2\sigma }d_{5} (a_{\sigma } u_{\lambda }^{2})]
 & = a_{\lambda }^{3}[u_{2\sigma }^{2}] \ot_{2}
 +  a_{\sigma }u_{\lambda }^{2}a_{\sigma }^{3}
                    a_{\lambda }\normrbar_{1} \\
 & = a_{\lambda }^{3}[u_{2\sigma }^{2}] \ot_{2}
      + 4a_{\lambda }^{3}[u_{2\sigma }^{2}]\normrbar_{1}\\
 & = a_{\lambda }^{3}[u_{2\sigma }^{2}] \ot_{2}\\
d_{5} (a_{\sigma } u_{\lambda }^{2})
 & =  a_{\lambda }^{3}[u_{2\sigma } \ot_{2}].
\end{align*}

\noindent We have seen that 
\begin{displaymath}
\eta \eta ' = a_{\lambda }[u_{2\sigma } \overline{t}_{2}].
\end{displaymath}

\noindent This implies that $a_{\lambda }^{2}[u_{2\sigma }
\overline{t}_{2}]$ vanishes in $\underline{E}_{5}$ since $a_{\lambda
}\eta $ is killed by $d_{3}$.  It follows that $d_{5} (a_{\sigma }
u_{\lambda }^{2}) = a_{\lambda }^{3}[u_{2\sigma } \ot_{2}]=0$ as claimed.

For the second differential we have 
\begin{align*}
d_{13} ([a_{\sigma }u_{\lambda }^{4}/u_{2\sigma }^{2}]) 
 & =  a_{\lambda }^{7} \normrbar_{1} (-5\ot_{2}^{2}+20 \ot_{2}\normrbar_{1}
      +9\normrbar_{1}^{2})\\
d_{13} ([a_{\sigma }u_{\lambda }^{4}])
 & =  a_{\lambda }^{7} [u_{2\sigma }^{2}]\normrbar_{1} 
        (-5\ot_{2}^{2}+20 \ot_{2}\normrbar_{1}
      +9\normrbar_{1}^{2})\\
 & = a_{\lambda }^{7}[u_{2\sigma}^{2}]\normrbar_{1}
                 (-\ot_{2}^{2}+\normrbar_{1}^{2})
\end{align*}

\noindent since $a_{\lambda }$ has order 4.  As we saw above, 
$a_{\lambda }^{2}[u_{2\sigma } \overline{t}_{2}]$ vanishes in
$\underline{E}_{5}$, so $d_{13} ([a_{\sigma }u_{\lambda }^{4}])$ is as
claimed.  \noindent\qed\bigskip

We can use this to find the differential on $[u_{\lambda}^{4}]$.
We have
\begin{numequation}\label{eq-d7}
\begin{split}\left\{
\begin{array}[]{rl}
d ([u_{\lambda}^{4}])
 &\hspace{-2.5mm} =   [2u_{\lambda}^{2}]d ([u_{\lambda}^{2}])
   =   [2u_{\lambda}^{2}]\alphax a_{\lambda}^{2}\normrbar_{1}^{}
   =   (2\alphax )a_{\lambda}^{2}[u_{\lambda}^{2}]\normrbar_{1}\\
 &\hspace{-2.5mm} =  \eta 'a_{\lambda}^{3}[u_{\lambda}^{2}]\normrbar_{1}^{}
   =   [\eta 'u_{\lambda}^{2}]a_{\lambda}^{3}\normrbar_{1}^{}
   =  \langle \eta',\,\alphax  ,\,a_{\lambda}^{2}\normrbar_{1}\rangle 
      a_{\lambda}^{3}\normrbar_{1}.
\end{array}
\right.
\end{split}
\end{numequation}%

The differential on $u_{2\sigma}$ yields
\begin{displaymath}
[x u_{2\sigma}]
=\langle x,\, a_{\sigma}^{2} 
                       ,\, f_{1}\rangle
\end{displaymath}

\noindent for any permanent cycle $x$ killed by $a_{\sigma}^{2}$.
Possible values of $x$ include 2, $\eta $, $\eta '$ (each of which is
killed by $a_{\sigma}$ as well) and $\alphax $.  For the last of these
we write
\begin{numequation}\label{eq-brackets2}
\betax :=  [\alphax  u_{2\sigma}]
  =   \langle \alphax ,\,a_{\sigma}^{2} 
                ,\, f_{1}\rangle
  =   \langle[ a_{\sigma} u_{\lambda }] ,\,a_{\sigma}^{2} 
                ,\, f_{1}\rangle
          \in \EE_{6}^{1,5-3\sigma-\lambda} (G/G),
\end{numequation}%

\noindent which satisfies
\begin{align*}
\Res_{2}^{4}(\betax )
 & =   a_{\sigma_{2}}^{3}u_{\sigma}^{3}\orr_{1,\epsilon}
  \in  \EE_{4}^{3,7-3\sigma -\lambda} (G/G')\\
&
   \qquad \mbox{(exotic restriction with jump 2)}\\
2\betax 
 & =   \Tr_{2}^{4}(\Res_{2}^{4}(\betax ))
 =  \eta'     a_{\lambda}u_{2\sigma}
  \in  \EE_{4}^{3,7-3\sigma -\lambda} (G/G)\\
&
\qquad \mbox{(exotic group extension with jump 2)}\\
d_{5} ([u_{2\sigma}u_{\lambda}^{2}])
 & =   a_{\sigma}^{3}a_{\lambda}u_{\lambda}^{2}\normrbar_{1}
        +\alphax a_{\lambda}^{2}u_{2\sigma}\normrbar_{1}
 =   (a_{\sigma}^{3}u_{\lambda}^{2}+\alphax u_{2\sigma})
          a_{\lambda}^{2}\normrbar_{1}\\
 & =  (2a_{\sigma}a_{\lambda}u_{\lambda}u_{2\sigma}+\betax )
          a_{\lambda}^{2}\normrbar_{1}
 =    \betax a_{\lambda}^{2}\normrbar_{1}\\
d_{7} ([2u_{2\sigma }u_{\lambda }^{2}])
 & =   2\betax\cdot  a_{\lambda}^{2}\normrbar_{1}
 =   \eta' a_{\lambda}^{3}u_{2\sigma}\normrbar_{1}\\
\Res_{2}^{4}(d_{5} ([u_{2\sigma}u_{\lambda}^{2}]))
 & =   u_{\sigma}^{3}a_{\sigma_{2}}^{3}\orr_{1,\epsilon}
             \Res_{2}^{4}(a_{\lambda}^{2}\normrbar_{1})
  =    u_{\sigma}^{2}a_{\sigma_{2}}^{7}\orr_{1,0}^{3}
  =   u_{\sigma}^{2}d_{7} (\overline{u}_{\lambda}^{2} ).
\end{align*}

\begin{thm}\label{thm-d3ulambda}
{\bf The differentials on powers of $u_{\lambda}$ and $u_{2\sigma}$.}
The following differentials occur in the slice {\SS} for $\kH$.  Here
$\ou_{\lambda}$ denotes $\Res_{2}^{4}(u_{\lambda})$.
\begin{align*}
d_{3} (u_{\lambda})
 & =   a_{\lambda}\eta  
   = \Tr_{2}^{4}(a_{\sigma_{2}}^{3}\orr_{1,\epsilon })\\
d_{3} (\ou_{\lambda})
 & =   \Res_{2}^{4}(a_{\lambda}) (\eta_{0}+\eta_{1})
   = \left[a_{\sigma_{2}}^{3} (\orr_{1,0}+\orr_{1,1}) \right]\\
d_{5} (u_{2\sigma})
 & =   a_{\sigma}^{3}a_{\lambda}\normrbar_{1}^{} \\
d_{5} ([u_{\lambda}^{2}])
 & =   a_{\lambda}^{2}a_{\sigma}u_{\lambda}\normrbar_{1}
  = a_{\lambda}^{2}\alphax \normrbar_{1}
  =  a_{\lambda}^{2}\nu\qquad 
           \mbox{for $\alphax$ 
                  as in (\ref{eq-alphax})} \\
d_{5} ([u_{2\sigma}u_{\lambda}^{2}])
 & =   a_{\sigma}^{3}a_{\lambda}u_{\lambda}^{2}\normrbar_{1}
        +\alphax a_{\lambda}^{2}u_{2\sigma}\normrbar_{1}
 =   (a_{\sigma}^{3}u_{\lambda}^{2}+\alphax u_{2\sigma})
          a_{\lambda}^{2}\normrbar_{1}\\
 & =  \betax a_{\lambda}^{2}\normrbar_{1}
       \qquad \mbox{for $\betax$ as in (\ref{eq-brackets2})} \\
d_{7} ([2u_{2\sigma}u_{\lambda}^{2}])
 & =   \eta' a_{\lambda}^{3}u_{2\sigma}\normrbar_{1}\\ 
d_{7} ([2u_{\lambda}^{2}])
 & =   2 a_{\lambda}^{2}\alphax \normrbar_{1}
  =  a_{\lambda}^{3}\eta '\normrbar_{1}^{}\\
d_{7} ([\ou_{\lambda}^{2}])
 & =   \Res_{2}^{4}(a_{\lambda}^{2})\eta_{0}^{3}
   = a_{\sigma_{2}}^{7}\orr_{1,0}^{3}\\
d_{7} ([u_{\lambda}^{4}])
 & =   [\eta 'u_{\lambda}^{2}]a_{\lambda}^{3}\normrbar_{1}
 = \langle \eta',\,\alphax  ,\,a_{\lambda}^{2}\normrbar_{1}\rangle 
      a_{\lambda}^{3}\normrbar_{1}.
\end{align*}

\noindent The elements
\begin{align*}
u_{\sigma},  
    &   &[2u_{\lambda}]
            & = \langle 2,\,\eta ,\,a_{\lambda} \rangle
                  ,\\
[2u_{2\sigma}]
    & =  \langle 2,\,a_{\sigma}^{2} ,\,f_{1} \rangle
            =\Tr_{2}^{4}(u_{\sigma}^{2}), 
        &[4u_{\lambda}^{2}]
            & = \langle 2,\,\eta ' 
                              ,\,a_{\lambda}^{3}\normrbar_{1} \rangle
            =\Tr_{1}^{4}(u_{\sigma_{2} }^{4}),\\
[u_{2\sigma}^{2} ]
    &=\langle a_{\sigma}^{2},\,f_{1} ,\,a_{\sigma}^{2},\,f_{1} \rangle
        & [2\ou_{\lambda}^{2}]
            & = \langle 2,\,a_{\sigma_{2} }^{6} 
                              ,\,a_{\sigma_{2}}\orr_{1,0}^{3} \rangle
            =\Tr_{1}^{2}(u_{\sigma_{2} }^{4}), \\
[2u_{2\sigma}u_{\lambda}]
    &=\langle [2u_{2\sigma}],\,\eta  ,\,a_{\lambda} \rangle,
        & [2u_{\lambda}^{4}]
            & = \langle 2,\,\eta ',\alphax 
                              ,\,a_{\lambda}^{5}\normrbar_{1}^{2} \rangle
            =\Tr_{2}^{4}(\ou_{\lambda}^{4}), \\
[\ou_{\lambda}^{4}]
            & = \langle a_{\sigma_{2}}^{7},\,\orr_{1,0}^{3} , \,
                             a_{\sigma_{2}}^{7},\,\orr_{1,0}^{3}\rangle 
        &\aand [u_{\lambda}^{8}]
            & =  \langle [\eta 'u_{\lambda }^{2}]
                          ,\,a_{\lambda}^{3}\normrbar_{1} 
                      ,\, [\eta 'u_{\lambda }^{2}]
                      ,\,a_{\lambda}^{3}\normrbar_{1}  \rangle
\end{align*}

\noindent are permanent cycles.


We also have the following exotic restriction and transfers.
\begin{align*}
\Res_{2}^{4}(a_{\sigma}u_{\lambda})
 & =   u_{\sigma}\Res_{2}^{4}(a_{\lambda})\eta_{\epsilon }
    =  u_{\sigma}a_{\sigma_{2}}^{3}\orr_{1,\epsilon }
      \qquad \mbox{(filtration jump 2)}  \\
\Tr_{2}^{4}(u_{\sigma}^{k})
 & =   \mycases{
a_{\sigma}^{2}a_{\lambda}\normrbar_{1}u_{2\sigma}^{(k-1)/2}
  = a_{\sigma}f_{1}u_{2\sigma}^{(k-1)/2}\\
\qquad \mbox{(filtration jump 4)}
       &\mbox{for $k\equiv 1$ mod 4}\hspace{4cm}\\
2u_{2\sigma}^{k/2}
       &\mbox{for $k$ even}\\
0      &\mbox{for $k\equiv 3$ mod 4}
}\\
\Tr_{1}^{2}(u_{\sigma_{2}}^{k})
 & = \mycases{
a_{\sigma_{2}}^{2} (\orr_{1,0}+\orr_{1,1})\ou_{\lambda}^{(k-1) /2}
 =   a_{\sigma_{2}} (\eta_{0}+\eta_{1})\ou_{\lambda}^{(k-1) /2}
              \\
\qquad \mbox{(filtration jump 2)}
       &\hspace{-2.4cm} \mbox{for $k\equiv 1$ mod 4 }\\
2\ou_{\lambda}^{k/2}
       &\hspace{-2.4cm}\mbox{for $k$ even}\\
a_{\sigma_{2}}^{6} \orr_{1,0}^{3}\ou_{\lambda}^{(k-3) /2}\\
\qquad \mbox{(filtration jump 6)}
       &\hspace{-2.4cm}\mbox{for $k\equiv 3$ mod 4 }
}
\end{align*}
\end{thm}

\proof All differentials were established above.  

The differential on $u_{\lambda}^{2}$ does {\em not} lead to an exotic
transfer because neither $[\ou_{\lambda}^{2}]$ nor
$[u_{\lambda} a_{\lambda}^{2}\normrbar_{1}]$ is a permanent cycle as
required by Theorem \ref{thm-exotic}.

We need to discuss the element $[2u_{2\sigma}u_{\lambda}] =\langle
[2u_{2\sigma}],\,\eta ,\,a_{\lambda} \rangle$.  To see that this Toda
bracket is defined, we need to verify that $[2u_{2\sigma}]\eta =0$.
For this we have
\begin{displaymath}
[2u_{2\sigma}] \eta
 = [2u_{2\sigma}]\Tr_{2}^{4}(\eta_{0})  
 = \Tr_{2}^{4}(2u_{\sigma}^{2}\eta_{0}) 
 = \Tr_{2}^{4}(0)=0.
\end{displaymath}

The exotic restriction and transfers are applications of Theorem
\ref{thm-exotic} to the differentials on $u_{\lambda}$ and on
$\left[u_{2\sigma}^{(k+1)/2} \right]$ and
$\left[\ou_{\lambda}^{(k+1)/2} \right]$ for odd $k$.  For even $k$ we
have
\begin{displaymath}
\Tr_{2}^{4}(u_{\sigma}^{k})
  = \Tr_{2}^{4}\left(\Res_{2}^{4}
                     \left(\left[u_{2\sigma}^{k/2}\right]\right) \right)
  = \left[2u_{2\sigma}^{k/2}  \right]
\quad \mbox{since } \Tr_{2}^{4}(\Res_{2}^{4}(x))= (1+\gamma )x,
\end{displaymath}

\noindent and similarly for even powers of $u_{\sigma_{2}}$.

As remarked above, we lose no information by inverting the class $D$,
which is divisible by $\normrbar_{1}$.  It is shown in
\cite[Thm. 9.16]{HHR} that inverting the latter makes
$u_{2\sigma}^{2}$ a permanent cycle.  One can also see this from
(\ref{eq-slicediffs2}).  Since $d_{5}
(u_{2\sigma})=a_{\sigma}^{3}a_{\lambda}\normrbar_{1}$, $d_{5}
(u_{2\sigma}\normrbar_{1}^{-1})=a_{\sigma}^{3}a_{\lambda}$.  This
means that $d_{13}
([u_{2\sigma}^{2}])=a_{\sigma}^{7}a_{\lambda}^{3}\normrbar_{3}$ is
trivial in $\EE_{6} (G/G)$.  It turns out that there is no possible
target for a higher differential.  \qed\bigskip \bigskip

\section{$\kH$ as a $C_{2}$-spectrum}
\label{sec-C2diffs}


Before studying the slice {\SS} for the $C_{4}$-spectrum $\kH$
further, it is helpful to explore its restriction to $G'=C_{2}$, for
which the $\Z$-bigraded portion
\begin{displaymath}
\EE_{2}^{*,*}i_{G'}^{*}\kH (G'/G')
\cong \EE_{2}^{*,(G',*)}\kH (G/G)
\cong \EE_{2}^{*,*}\kH (G/G')
\end{displaymath}

\noindent (see Theorem \ref{thm-module} for these isomorphisms) is the
isomorphic image of the subring of Corollary \ref{cor-subring}.  In
the following we identify $\Sigma_{2,\epsilon }$, $\delta_{1}$ and
$\orr_{1,\epsilon }$ (see Table \ref{tab-pi*}) with their images under
$\ur_{2}^{4}$. From the differentials of (\ref{eq-d3-and-d7}) we get
\begin{numequation}\label{eq-k2-C2-diffs}
\begin{split}
\left\{\begin{array}{rl}
d_{3} (\Sigma_{2,\epsilon })
 &\hspace{-2.5mm} =   \eta_{\epsilon }^{2} (\eta_{0}+\eta_{1})
   = a_{\sigma }^{3} \orr_{1,\epsilon }^{2} (\orr_{1,0}+\orr_{1,1}) \\
d_{3} (\delta_{1})
 &\hspace{-2.5mm} =   \eta_{0}^{2}\eta_{1}+\eta_{0}\eta_{1}^{2}
   = a_{\sigma }^{3} \orr_{1,0 }  \orr_{1,1 } (\orr_{1,0}+\orr_{1,1}) \\
d_{7} ([\delta_{1}^{2}])
 &\hspace{-2.5mm} =   d_{7} (u_{2\sigma}^{2}) \orr_{1,0}^{2}\orr_{1,1}^{2}
 = a_{\sigma}^{7}\orr_{3}^{G'} 
               \orr_{1,0}^{2}\orr_{1,1}^{2}\\
 &\hspace{-2.5mm} =  a_{\sigma}^{7} (5\orr_{1,0}^{2}\orr_{1,1} 
                       +5\orr_{1,0}\orr_{1,1}^{2}+\orr_{1,1}^{3})
               \orr_{1,0}^{2}\orr_{1,1}^{2}.
\end{array} \right.
\end{split}
\end{numequation}%

\noindent The $d_{3}$s above make all monomials in $\eta_{0}$ and
$\eta_{1}$ of any given degree $\geq 3$ the same in $\EE_{4}
(G/G')$ and $\EE_{4} (G'/G')$, so $d_{7}
(\delta_{1}^{2})=\eta_{0}^{7}$.  Similar calculations show that
\begin{displaymath}
d_{7} ([\Sigma_{2,\epsilon }^{2}])=\eta_{0}^{7}
 =  a_{\sigma}^{7}\orr_{1,0}^{7}.
\end{displaymath}

\noindent The image of the periodicity element $D$ here is as in
(\ref{eq-r24D}).

We have the following values of the transfer on powers of $u_{\sigma }$.
\begin{displaymath}
\Tr_{1}^{2}(u_{\sigma }^{i})=\mycases{
[2u_{2\sigma }^{i/2}]
       &\mbox{for $i$ even} \\
{}
[a_{\sigma }^{2} u_{2\sigma }^{(i-1) /2}] (\orr_{1,0}+\orr_{1,1})
       &\mbox{for $i\equiv 1$ mod 4} \\
{}
[u_{2\sigma }^{4}]^{(i-3) /8}a_{\sigma }^{6}\orr_{1,0}^{3} 
    = [ u_{2\sigma }^{4}]^{(i-3) /8}a_{\sigma }^{6}\orr_{1,1}^{3}
       &\mbox{for $i\equiv 3$ mod 8}\\
0      &\mbox{for $i\equiv 7$ mod 8}
}
\end{displaymath}

\noindent

\begin{figure}
\begin{center}
\includegraphics[width=10.9cm]{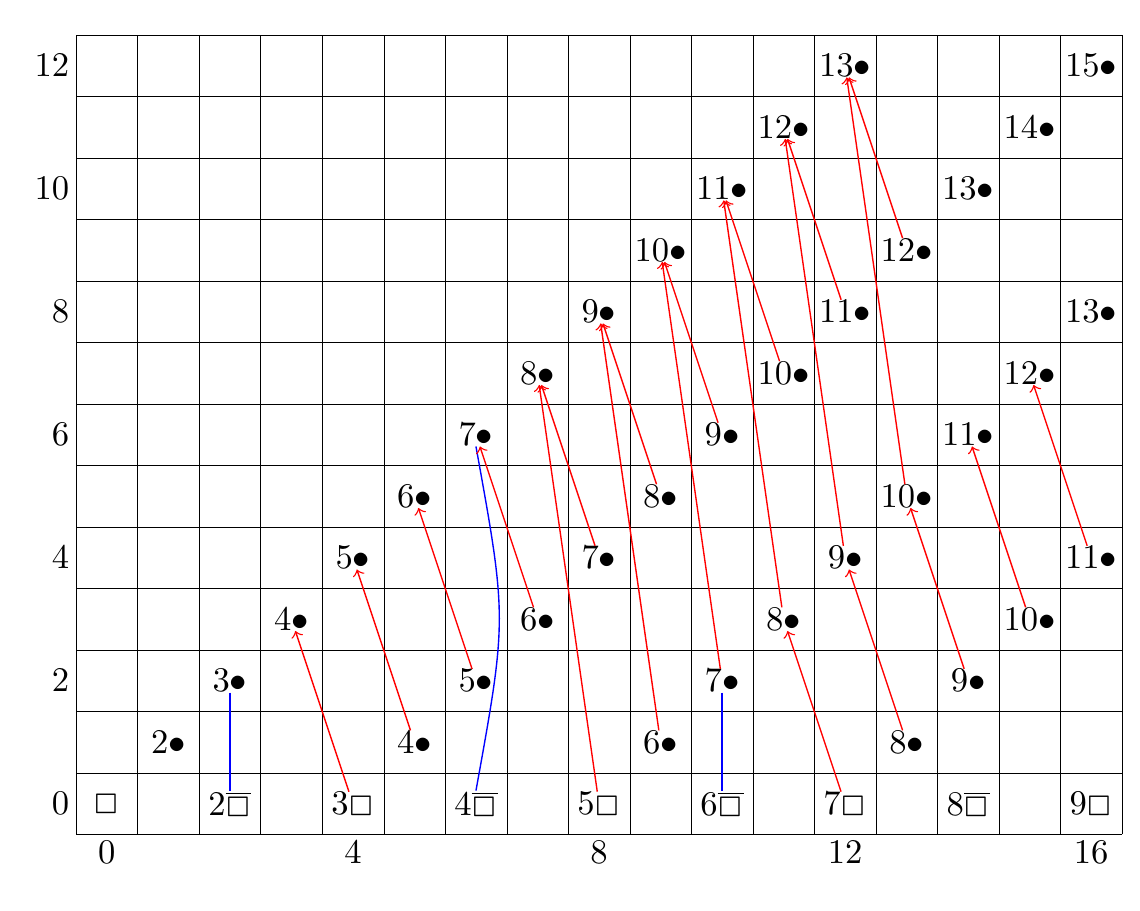} 

\caption[The slice {\SS} for $\kH$ as a $C_{2}$-spectrum. ]{The slice
{\SS} for $\kH$ as a $C_{2}$-spectrum.  The Mackey functor symbols are
defined in Table \ref{tab-C2Mackey}. The $C_{4}$-structure of the
Mackey functors is not indicated here.  In each bidegree we have a
direct sum of the indicated number of copies of the indicated Mackey
functor. Each $d_{3}$ has maximal rank, leaving a cokernel of rank 1,
and each $d_{7}$ has rank 1. Blue lines indicate exotic transfers. The
ones raising filtration by 2 have maximal rank while the ones raising
it by 6 have rank 1. The resulting
$\EE_{8}=\EE_{\infty }$-term is shown below.}

\includegraphics[width=11cm]{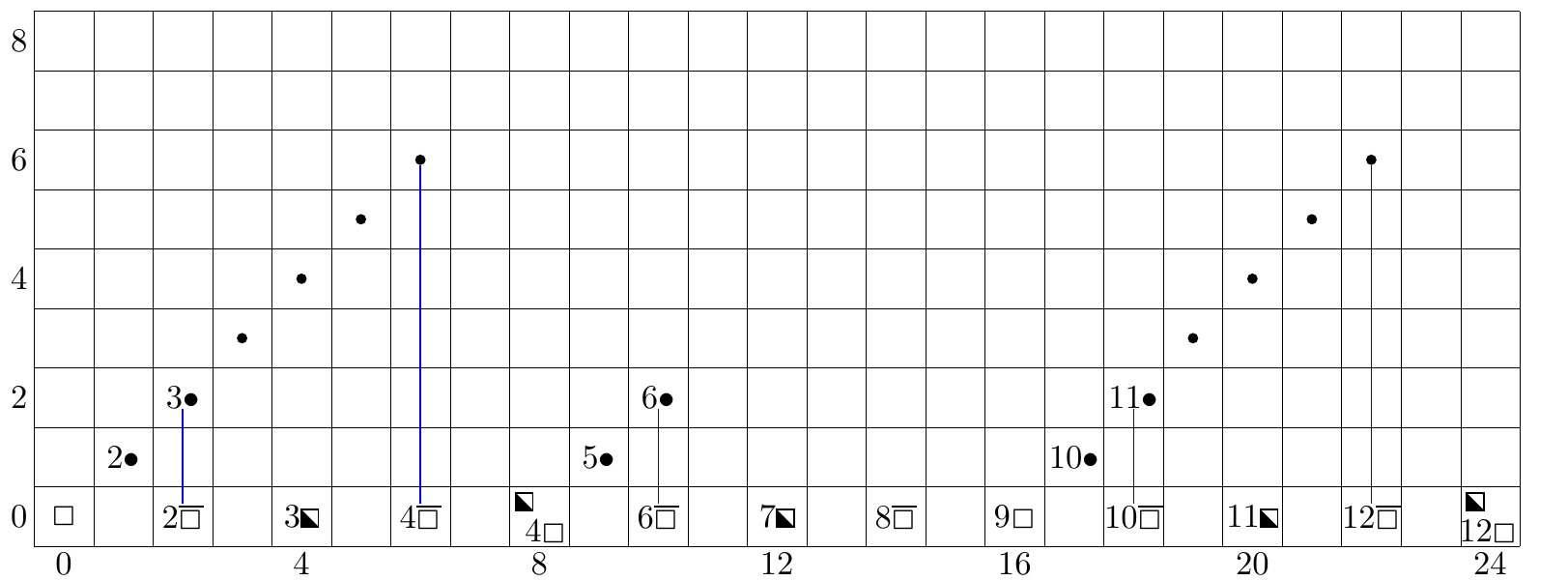} 
\vspace{3cm}
\label{fig-G'SS}
\end{center}

\end{figure}

This leads to the following, for which Figure
\ref{fig-G'SS} is a visual aid.

\begin{thm}\label{thm-G'SS}
{\bf The slice {\SS} for $\kH$ as a $C_{2}$-spectrum.}  Using the
notation of Table \ref{tab-C2Mackey} and Definition
\ref{def-enriched}, we have
\begin{align*}
\EE_{2}^{*,*} (G'/\ee)
 & =   \Z[r_{1,0}, r_{1,1}]\qquad 
   \mbox{with }r_{1,\epsilon }\in \EE_{2}^{0,2} (G'/\ee) \\
\EE_{2}^{*,*} (G'/G')
 & =   \Z[\delta_{1},\Sigma_{2,\epsilon},\eta_{\epsilon }
      \colon \epsilon =0,\,1]/\\
&\qquad 
\left(2\eta _{\epsilon },\delta_{1}^{2}-\Sigma_{2,0}\Sigma_{2,1},
\eta_{\epsilon } \Sigma_{2,\epsilon +1}+\eta _{1+\epsilon }\delta_{1} \right),
\end{align*}

\noindent so 
\begin{displaymath}
\EE_{2}^{s,t} =\mycases{
\Box \oplus \bigoplus_{\ell}\widetilde{\Box}
       &\mbox{for $(s,t)= (0,4\ell )$ with $\ell \geq 0$}\vspace{2mm}\\
\bigoplus_{\ell +1}\widetilde{\oBox}
       &\mbox{for $(s,t)= (0,4\ell+2 )$ with $\ell \geq 0$}\vspace{2mm}\\
\bullet\oplus \bigoplus_{u+\ell}\widetilde{\bullet}
       &\mbox{for $(s,t)= (2u,4\ell+4u )$ 
            with $\ell \geq 0$ and $u>0$}\vspace{2mm}\\
\bigoplus_{u+\ell}\widetilde{\bullet}
       &\mbox{for $(s,t)= (2u-1,4\ell+4u -2)$ 
            with $\ell \geq 0$ and $u> 0$}\vspace{2mm}\\
0       &\mbox{otherwise.}
}
\end{displaymath}

The first set of differentials and determined by 
\begin{displaymath}
d_{3} (\Sigma_{2,\epsilon })=\eta_{\epsilon }^{2} (\eta_{0}+\eta_{1})
\qquad \aand 
d_{3} (\delta_{1}) = \eta_{0}\eta_{1} (\eta_{0}+\eta_{1})
\end{displaymath}


\noindent and there is a second set of differentials determined by 
\begin{displaymath}
d_{7} (\Sigma_{2,\epsilon }^{2})=d_{7} (\delta_{1}^{2})=\eta_{0}^{7}
\end{displaymath}

\end{thm}

\begin{cor}\label{cor-G'SS}
{\bf Some nontrivial permanent cycles.} The elements  listed below in
$\EE_{2}^{s,8i+2s}\kH (G/G')$ are nontrivial permanent
cycles.  Their tranfers in $\EE_{2}^{s,8i+2s}\kH (G/G)$ are also
permanent cycles.
\begin{itemize}
\item [$\bullet$] $\Sigma_{2,\epsilon }^{2i-j}\delta_{1}^{j}$ for
$0\leq j\leq 2i$ ($4i+1$ elements of infinite order including
$\delta_{1}^{2i}$), $i$ even and $s=0$.

\item [$\bullet$] $\eta_{\epsilon }\Sigma_{2,\epsilon
}^{2i-j}\delta_{1}^{j}$ for $0\leq j<2i$ and $\eta_{\epsilon
}\delta_{1}^{2i}$ ($4i+2$ elements or order 2) for $i$ even and $s=1$.

\item [$\bullet$] $\eta_{\epsilon }^{2}\Sigma_{2,\epsilon
}^{2i-j}\delta_{1}^{j}$ for $0\leq j<2i$ and
$\delta_{1}^{2i}\left\{\eta_{0}^{2},\,\eta_{0}\eta_{1},\,\eta_{1}^{2}
\right\}$ ($4i+3$ elements or order 2) for $i$ even and $s=2$.

\item [$\bullet$] $\eta_{0}^{s}\delta_{1}^{2i}$ for $3\leq s\leq 6$ (4
elements or order 2) and $i$ even.

\item [$\bullet$]
$\Sigma_{2,\epsilon}^{2i-j}\delta_{1}^{j}+\delta_{1}^{2i}$ for $0\leq
j\leq 2i$ ($4i+1$ elements of infinite order including
$2\delta_{1}^{2i}$), $i$ odd and $s=0$.

\item [$\bullet$]
$\eta_{\epsilon}\Sigma_{2,\epsilon}^{2i-j}\delta_{1}^{j}+\delta_{1}^{2i}$
for $0\leq j\leq 2i-1$ and $\eta_{0}\delta_{1}^{2i-1}
(\Sigma_{2,1}+\delta_{1}) = \eta_{1}\delta_{1}^{2i-1}
(\Sigma_{2,0}+\delta_{1})$ ($4i+1$ elements of order 2), $i$ odd and $s=1$.

\item [$\bullet$]
$\eta_{\epsilon}^{2}\Sigma_{2,\epsilon}^{2i-j}\delta_{1}^{j}+\delta_{1}^{2i}$
for $0\leq j\leq 2i-1$,
$\eta_{0}^{2}\delta_{1}^{2i-1}(\Sigma_{2,1}+\delta_{1}) =
\eta_{0}\eta_{1}\delta_{1}^{2i-1}(\Sigma_{2,0}+\delta_{1})$ and
$\eta_{0}\eta_{1}\delta_{1}^{2i-1}(\Sigma_{2,1}+\delta_{1}) =
\eta_{1}^{2}\delta_{1}^{2i-1}(\Sigma_{2,0}+\delta_{1})$ ($4i+2$
elements of order 2) for $i$ odd and $s=2$.
\end{itemize}

In $\EE_{2}^{0,8i+4}\kH (G/G')$ we have $2\Sigma_{2,\epsilon
}^{2i+1-j}\delta_{1}^{j}$ for $0\leq j\leq 2i$ and $2\delta_{1}^{j}$,
$4i+3$ elements of infinite order, each in the image of the transfer
$\Tr_{1}^{2}$.
\end{cor}

\begin{rem}\label{rem-a7}
In {\bf the $RO (G)$-graded slice {\SS } for $k_{[2]}$} one has
\begin{displaymath}
d_{3}(u_{2\sigma })=a_{\sigma }^{3} (\orr_{1,0}+\orr_{1,1})
\qquad \aand 
d_{7}([u_{2\sigma }^{2}])=a_{\sigma }^{7} \orr_{3}^{G'}
 =a_{\sigma }^{7} \orr_{1,0}^{3},
\end{displaymath}

\noindent but $a^{7}$ itself, and indeed all higher powers of $a$,
survive to $\underline{E}_{8}=\underline{E}_{\infty }$.  Hence the
$\underline{E}_{\infty }$-term of this {\SS} does {\bf not} have the
horizontal vanishing line that we see in $\underline{E}_{8}$-term of
Figure \ref{fig-KR}. However when we pass from $k_{[2]}$ to $K_{[2]}$,
$\orr_{3}^{G'}=5\orr_{1,0}^{2}\orr_{1,1}
+5\orr_{1,0}\orr_{1,1}^{2}+\orr_{1,1}^{3}$ becomes invertible and we
have
\begin{displaymath}
d_{7} ((\orr_{3}^{G'})^{-1}[u_{2\sigma }^{2}])
 = d_{7} (\orr_{1,0}^{-3}[u_{2\sigma }^{2}])
 = a^{7}.
\end{displaymath}

\noindent On the other hand, $\orr_{1,0}+\orr_{1,1}$ is not invertible, so we cannot divide $u_{2\sigma }$ by it.
\end{rem}

We now give the {\Ps} computation analogous to the one following
Remark \ref{rem-a3}, using the notation of (\ref{eq-aur}).  In $RO
(G')$-graded slice {\SS } for $k_{[2]}$ we have
\begin{displaymath}
\underline{E}_{2} (G'/G') 
 =\Z[a_{\sigma },u_{2\sigma },\orr_{1,0},\orr_{1,2}]/ (2a_{\sigma }),
\end{displaymath}

\noindent so 
\begin{align*}
g (\underline{E}_{2} (G'/G'))
 & = \left(\frac{1}{1-t}+ \frac{\wa}{1-\wa}\right)
  \frac{1}{(1-\uu) (1-\rr)^{2}}  \\
g (\underline{E}_{4} (G'/G'))
 & =  g (\underline{E}_{2} (G'/G'))
   -\frac{\uu+\rr\wa^{3}}
         {(1-\wa) (1-\uu^{2}) (1-\rr)^{2}} \\ 
 & = \frac{1+t\uu}
          {(1-t) (1-\uu^{2}) (1-\rr)^{2}} 
     +\frac{\wa+\wa^{2}}
           { (1-\uu^{2}) (1-\rr)^{2}} 
     +\frac{\wa^{3}}
           {(1-\wa) (1-\uu^{2})(1-\rr)} 
\end{align*}

\noindent as before.  The next differential leads to 
\begin{align*}
g (\underline{E}_{8} (G'/G'))
 & =  g (\underline{E}_{4} (G'/G')) -\frac{\uu^{2}+\rr^{3}\wa^{7}}
         {(1-\wa) (1-\uu^{4}) (1-\rr)} \\
 & =  g (\underline{E}_{4} (G'/G')) -\frac{\uu^{2}}
         { (1-\uu^{4}) (1-\rr)}\\
 & \qquad  -\frac{\uu^{2}\wa}
         {(1-\wa) (1-\uu^{4}) (1-\rr)} -\frac{\rr^{3}\wa^{7}}
         {(1-\wa) (1-\uu^{4}) (1-\rr)} \\
 & =  \frac{1+t\uu}{(1-t) (1-\uu^{2}) (1-\rr)^{2}} 
     +\frac{\wa+\wa^{2}}{ (1-\uu^{2}) (1-\rr)^{2}}\\
 &\qquad  
     +\frac{\wa^{3}}{(1-\wa) (1-\uu^{2})(1-\rr)} 
    -\frac{\uu^{2}}{ (1-\uu^{4}) (1-\rr)}  \\
 & \qquad
          -\frac{\uu^{2} (\wa+\wa^{2})}{ (1-\uu^{4}) (1-\rr)} 
          -\frac{\uu^{2}\wa^{3}}{(1-\wa) (1-\uu^{4}) (1-\rr)} 
          -\frac{\rr^{3}\wa^{7}}{(1-\wa) (1-\uu^{4}) (1-\rr)} \\
 & =   \frac{(1+t\uu) (1+\uu^{2})- (1-t) (1-\rr)\uu^{2}}
            {(1-t) (1-\uu^{4}) (1-\rr)^{2}} 
       + \frac{(\wa+\wa^{2}) ( 1+\uu^{2} -\uu^{2} (1-\rr))}
              { (1-\uu^{4}) (1-\rr)^{2}}  \\
 &\qquad + \frac{\wa^{3} (1+\uu^{2})-\uu^{2}\wa^{3}-\rr^{3}\wa^{7}}
                {(1-\wa) (1-\uu^{4})(1-\rr)}\\
 & =   \frac{1+t\uu+ (t+\rr-t\rr)\uu^{2}+t\uu^{3}}
            {(1-t) (1-\uu^{4}) (1-\rr)^{2}}  
        + \frac{(\wa+\wa^{2}) (1+\uu^{2}\rr)}
              {(1-\uu^{4}) (1-\rr)^{2}}  \\
 &\qquad +\frac{\wa-\wa^{7}+\wa^{7}-\rr^{3}\wa^{7}}
               {(1-\wa) (1-\uu^{4})(1-\rr)}\\
 & =  \frac{1+t\uu+ (t+\rr-t\rr)\uu^{2}+t\uu^{3}}
            {(1-t) (1-\uu^{4}) (1-\rr)^{2}} 
      + \frac{(\wa+\wa^{2}) (1+\uu^{2}\rr)}
              {(1-\uu^{4}) (1-\rr)^{2}}  \\
 &\qquad +\frac{\wa^{3}+\wa^{4}+\wa^{6}+\wa^{6}}
               { (1-\uu^{4})(1-\rr)} 
         +\frac{\wa^{7} (1+\rr+\rr^{2})}
               {(1-\wa) (1-\uu^{4})}
\end{align*}

\noindent The fourth term of this expression represents the elements
with filtration above six, and the first term represents the elements of
filtration 0.  The latter include
\begin{align*}
[2u_{2\sigma }]
 & \in \langle 2,\,a_{\sigma }^{2} 
                ,\,a_{\sigma }(\orr_{1,0}+\orr_{1,1}) \rangle ,  \\
[2u_{2\sigma }^{2}]
 & \in \langle 2,\,a_{\sigma } ,\,a_{\sigma }^{6}\orr_{1,0}^{3} \rangle ,  \\
[(\orr_{1,0}+\orr_{1,1}) u_{2\sigma }^{2}]
 & \in   \langle a_{\sigma }^{4} ,\,a_{\sigma }^{}\orr_{1,0}^{3}
                ,\, \orr_{1,0}+\orr_{1,1}\rangle   \\
 & \qquad \mbox{with } 
      (\orr_{1,0}+\orr_{1,1}) [2u_{2\sigma }^{2}]
               =2 [(\orr_{1,0}+\orr_{1,1}) u_{2\sigma }^{2}], \\ 
[2u_{2\sigma }^{3}]
 & \in \langle 2,\, a_{\sigma }^{2}(\orr_{1,0}+\orr_{1,1}) 
                ,\, a_{\sigma}^{2},\,a_{\sigma}^{6}\orr_{1,0}^{3}\rangle \\ 
\aand 
[u_{2\sigma }^{4}]
 & \in  \langle a_{\sigma }^{4} ,\,a_{\sigma }^{3}\orr_{1,0}^{3} 
                \,a_{\sigma }^{4} ,\,a_{\sigma }^{3}\orr_{1,0}^{3} \rangle 
\end{align*}

\noindent with notation as in \ref{rem-abuse}.

As indicated in \ref{rem-a7}, we can get rid of them by formally
adjoining $w:=(\orr_{3}^{G'})^{-1}u_{2\sigma }^{2}$ to $\underline{E}_{2}
(G'/G')$.  As before we denote the enlarged {\SS } terms by
$\underline{E}'_{r} (G'/G')$ This time let
\begin{displaymath}
\www= \rr^{-3}\uu^{2}=xy^{-7}.
\end{displaymath}

\noindent  Then we have 
\begin{displaymath}
\underline{E}'_{r} (G'/G')
  =\left(\frac{1-\uu^{2}}{1-\www} \right)\underline{E}_{r} (G'/G')
\qquad \mbox{for $r=2$ and $r=4$} 
\end{displaymath}

\noindent and 
\begin{align*}
g (\underline{E}'_{8} (G'/G'))
 & =  g (\underline{E}'_{4} (G'/G'))
           -\frac{\www+\wa^{7}}{(1-\wa) (1-\www^{2}) (1-\rr)} \\
 & =  g (\underline{E}'_{4} (G'/G'))
           -\frac{\www}{ (1-\www^{2}) (1-\rr)}  
           -\frac{(\wa+\wa^{2})\www}{ (1-\www^{2}) (1-\rr)}  \\
 &\qquad    -\frac{\wa^{3}\www+\wa^{7}}{(1-\wa) (1-\www^{2}) (1-\rr)} \\
 & =   \frac{1+t\uu}{(1-t) (1-\www) (1-\rr)^{2}} 
     +\frac{\wa+\wa^{2}}{ (1-\www) (1-\rr)^{2}}\\
 &\qquad  
     +\frac{\wa^{3}}{(1-\wa) (1-\www)(1-\rr)}
      -\frac{\www}{ (1-\www^{2}) (1-\rr)} \\
 &\qquad   -\frac{(\wa+\wa^{2})\www}{ (1-\www^{2}) (1-\rr)} 
        -\frac{\wa^{3}\www+\wa^{7}}{(1-\wa) (1-\www^{2}) (1-\rr)}  \\
 & =   \frac{(1+t\uu) (1+\www) - (1-t) (1-\rr)\www}
            {(1-t) (1-\www^{2}) (1-\rr)^{2}}
      +\frac{(\wa+\wa^{2}) (1- (1-\rr)\www)}{ (1-\www^{2}) (1-\rr)^{2}}\\
 &\qquad  +\frac{\wa^{3} (1+\www) -\wa^{3}\www-\wa^{7}}
                {(1-\wa) (1-\www^{2})(1-\rr)}  \\ 
 & =   \frac{1+t\uu+ (t+\rr-t\rr)\www+t\www\uu}
            {(1-t) (1-\www^{2}) (1-\rr)^{2}}
      +\frac{(\wa+\wa^{2}) (1+\rr\www)}{ (1-\www^{2}) (1-\rr)^{2}}\\
 &\qquad  +\frac{\wa^{3}+\wa^{4}+\wa^{5}+\wa^{6}}{ (1-\www^{2})(1-\rr)} .
\end{align*}

\noindent Again the first term represents the elements of
filtration 0.  These include
\begin{align*}
[2u_{2\sigma }]
 & \in \langle 2,\,a_{\sigma }^{2} 
                ,\,a_{\sigma }(\orr_{1,0}+\orr_{1,1}) \rangle ,  \\
[2w]
 & \in \langle 2,\,a_{\sigma } ,\,a_{\sigma }^{6}\rangle ,  \\
[(\orr_{1,0}+\orr_{1,1})w]
 & \in   \langle a_{\sigma }^{4} ,\,a_{\sigma }^{3}
                ,\, \orr_{1,0}+\orr_{1,1}\rangle   \\
 & \qquad \mbox{with } 
      (\orr_{1,0}+\orr_{1,1}) [2w]
               =2 [(\orr_{1,0}+\orr_{1,1}) w], \\ 
[2u_{2\sigma }w]
 & \in \langle 2,\, a_{\sigma }^{2}(\orr_{1,0}+\orr_{1,1}) 
                ,\, a_{\sigma}^{2},\,a_{\sigma}^{6}\rangle \\ 
\aand 
[w^{2}]
 & \in  \langle a_{\sigma }^{4} ,\,a_{\sigma }^{3}
                \,a_{\sigma }^{4} ,\,a_{\sigma }^{3} \rangle 
\end{align*}

\noindent where, as indicated above,
$w=(\orr_{3}^{G'})^{-1}u_{2\sigma}^{2}$.

From these

\section{The effect of the first differentials over $C_{4}$}
\label{sec-C4diffs}


Theorem \ref{thm-sliceE2} lists elements in the slice {\SS} for
$\kH$ over $C_{4}$ in terms of
\begin{displaymath}
r_{1}, \,  \os_{2} , \,\normrbar_{1};\quad 
\eta , \,a_{\sigma}, \,a_{\lambda};\quad 
u_{\lambda},\,u_{\sigma},\,\mbox{ and }u_{2\sigma}.
\end{displaymath}

\noindent All but the $u$'s are permanent cycles, and the action of
$d_{3}$ on $u_{\lambda}$, $u_{\sigma}$ and $u_{2\sigma}$ is
described above in Theorem \ref{thm-d3ulambda}.

\begin{prop}\label{prop-d3}
{\bf $d_{3}$ on elements in Theorem \ref{thm-sliceE2}.}
We have the following $d_{3}$s, subject to the conditions on $i$, $j$,
$k$ and $\ell $ of Theorem \ref{thm-sliceE2}:

\begin{itemize}
\item [$\bullet$] On $X_{2\ell ,2\ell }$:
\begin{align*}
d_{3} (a_{\lambda}^{j}u_{2\sigma}^{\ell }
          u_{\lambda}^{2\ell -j}\normrbar_{1}^{2\ell })
     & = \mycases{
a_{\lambda}^{j+1}\eta u_{2\sigma}^{\ell}u_{\lambda}^{2\ell -j-1}
                                \normrbar_{1}^{2\ell}\hspace{-1.7cm}\\
\qquad 
\in \upi_{*}X_{2\ell,2\ell +1} (G/G)
        &\mbox{for $j$ odd}\\
0\hspace{2cm}&\mbox{for $j$ even}
}\\
   d_{3} (a_{\sigma}^{2k}a_{\lambda}^{2\ell }u_{2\sigma}^{\ell-k }
          \normrbar_{1}^{2\ell })
     & = 0
\end{align*}

\item [$\bullet$] On $X_{2\ell+1 ,2\ell+1 }$:
\begin{align*}
d_{3} (\delta_{1}^{2\ell +1})
   & = \eta u_{\sigma}^{2\ell+1}
       \Res_{2}^{4}(a_{\lambda}u_{\lambda}^{2\ell }\normrbar_{1}^{2\ell+1})\\
   & \in \upi_{*}X_{2\ell+1,2\ell +2} (G/G')\\
\lefteqn{ d_{3} (u_{\sigma}^{2\ell +1}
     \Res_{2}^{4}(a_{\lambda}^{j}u_{\lambda}^{2\ell +1-j}
   \normrbar_{1}^{2\ell +1}))}\qquad\qquad\\
   & = \mycases{
\eta u_{\sigma}^{2\ell +1}
     \Res_{2}^{4}(a_{\lambda}^{j+1}u_{\lambda}^{2\ell -j}
                  \normrbar_{1}^{2\ell +1})\\
     \qquad 
         \in \upi_{*}X_{2\ell+1,2\ell +2} (G/G')
        &\mbox{for $j$ even}\\
0\hspace{2cm}&\mbox{for $j$ odd}
}\\
\lefteqn{   d_{3} (a_{\sigma}a_{\lambda}^{j}
       u_{\sigma}^{2\ell}u_{\lambda}^{2\ell +1-j}
                  \normrbar_{1}^{2\ell +1})}\qquad\qquad\\
     & = \mycases{
\eta a_{\sigma}a_{\lambda}^{j+1}
       u_{\sigma}^{2\ell}u_{\lambda}^{2\ell -j}
                  \normrbar_{1}^{2\ell +1}\\ 
     \qquad 
         \in \upi_{*}X_{2\ell+1,2\ell +2} (G/G)
        &\mbox{for $j$ even}\\
0\hspace{2cm}&\mbox{for $j$ odd}
}\\
\lefteqn{d_{3} (a_{\sigma}^{2k+1}a_{\lambda}^{2\ell +1}
       u_{2\sigma}^{\ell-k}
                  \normrbar_{1}^{2\ell +1})}\qquad\qquad\\ 
     & =  0
\end{align*}

\item [$\bullet$] On $X_{i ,2\ell-i }$:
\begin{align*}
d_{3} (u_{\sigma  }^{\ell } \os_{2}^{\ell -i}
\Res_{2}^{4}(u_{\lambda}^{\ell }\normrbar_{1}^{i}) )
     & = \mycases{
\eta^{3} u_{\sigma  }^{\ell-1 } \os_{2}^{\ell -i-1}
\Res_{2}^{4}(u_{\lambda}^{\ell-1 }\normrbar_{1}^{i})\\
    \qquad 
         \in \upi_{*}X_{i,2\ell +1-i} (G/G')\\
       &\hspace{-2cm}\mbox{for  $\ell $ odd}\\
0      &\hspace{-2cm}\mbox{for $\ell $ even}
}\\
 d_{3} (\eta^{2j}u_{\sigma}^{\ell -j}\os_{2}^{\ell -i-j}
\Res_{2}^{4}(u_{\lambda}^{\ell -j}\normrbar_{1}^{i}) )
     & = \mycases{
\eta^{2j+1}u_{\sigma}^{\ell -j}\os_{2}^{\ell -i-j}
\Res_{2}^{4}(a_{\lambda}u_{\lambda}^{\ell -j-1}\normrbar_{1}^{i})
            \\
   \qquad 
         \in \upi_{*}X_{i,2\ell +1-i} (G/G')\\
       &\hspace{-3.0cm} \mbox{for $\ell -j$ odd}\\
0
       &\hspace{-3.0cm} \mbox{for $\ell -j$ even}
}
\end{align*}

\item [$\bullet$] On $X_{i ,2\ell+1-i }$:
\begin{align*}
 d_{3} (r_{1}\Res_{1}^{2}(u_{\sigma  }^{\ell } \os_{2}^{\ell -i})
\Res_{1}^{4}(u_{\lambda}^{\ell }\normrbar_{1}^{i}) )
     & = 0\\
 d_{3} ( \eta^{2j+1}u_{\sigma}^{\ell -j}\os_{2}^{\ell -i-j}
\Res_{2}^{4}(u_{\lambda}^{\ell -j}\normrbar_{1}^{i}) ) 
     & = \mycases{
\eta^{2j+2}u_{\sigma}^{\ell -j}\os_{2}^{\ell -i-j}
\Res_{2}^{4}(a_{\lambda}u_{\lambda}^{\ell -j-1}\normrbar_{1}^{i}) 
      \\
   \qquad 
         \in \upi_{*}X_{i,2\ell +2-i} (G/G')\hspace{-3cm}\\
       &\hspace{-3cm} \mbox{for $\ell -j$ odd}\\
0      &\hspace{-3cm} \mbox{for $\ell -j$ even}
}
\end{align*}

\end{itemize}
\end{prop}

Note that in each case the first index of $X$ is unchanged by the
differential, and the second one is increased by one.  Since $X_{m,n}$
is a summand of the $2 (m+n)$th slice, each $d_{3}$ raises the slice
degree by 2 as expected.

\begin{figure} 
\begin{center} 
\includegraphics[width=11cm]{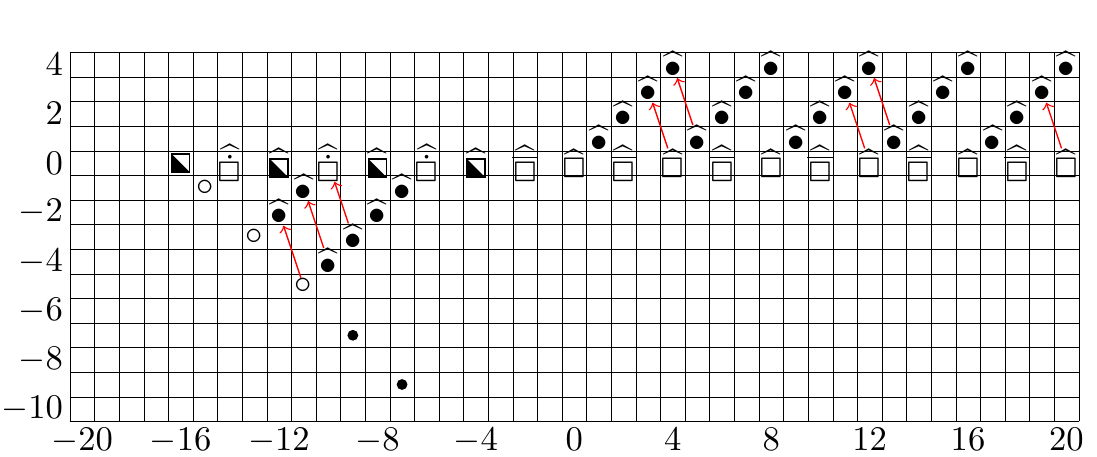}
\includegraphics[width=11cm]{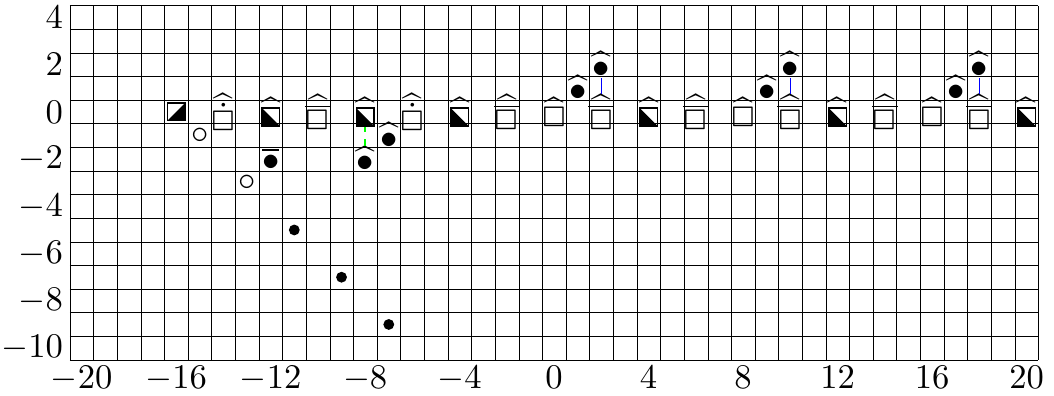} \caption[{The slice
$\EE_{2}$ and $\EE_{4}$-terms for $X_{-4,n}$ for $n\geq -4$.}]  {The
subquotient of the slice $\EE_{2}$ and $\EE_{4}$-terms for $\kH$ for
the slice summands $X_{-4,n}$ for $n\geq -4$.  This is also the slice
{\SS } for $Y_{-4}$ (after tensoring with $R$) as in Corollary
\ref{cor-Filtration}.}  \label{fig-sseq-8a}
\end{center}    
\end{figure}

\begin{figure} 
\begin{center}
\includegraphics[width=11cm]{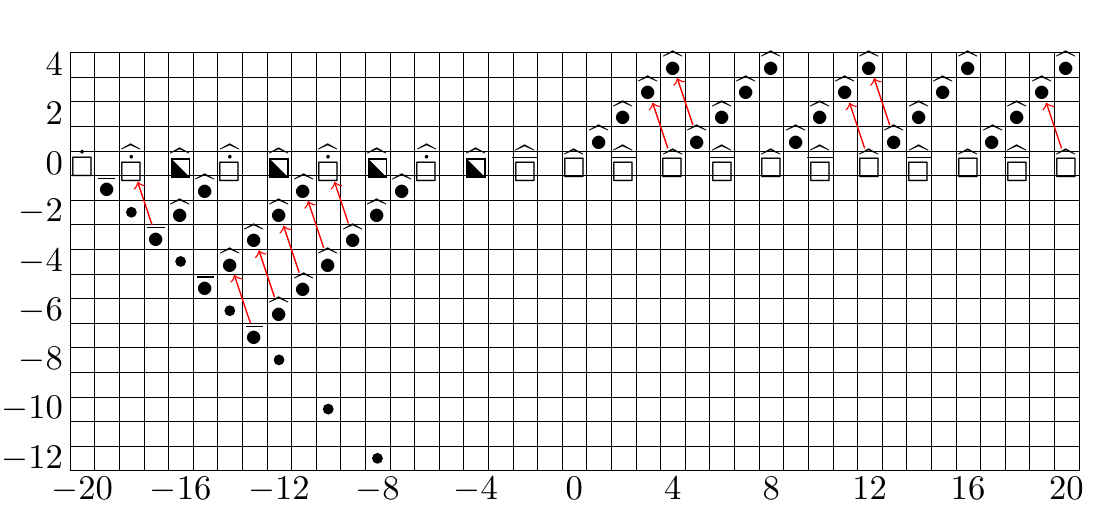}
\includegraphics[width=11cm]{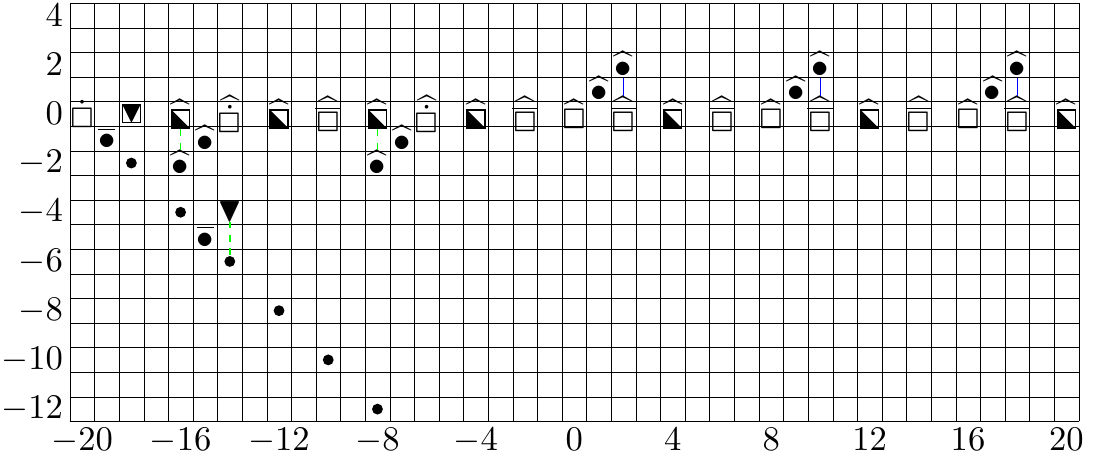} \caption[{The slice
$\EE_{2}$ and $\EE_{4}$-terms for $X_{-5,n}$ for $n\geq -5$.}]  {The
subquotient of the slice $\EE_{2}$ and $\EE_{4}$-terms for $\kH$ for
the slice summands $X_{-5,n}$ for $n\geq -5$. This is also the slice
{\SS } for $Y_{-5}$ (after tensoring with $R$) as in Corollary
\ref{cor-Filtration}. } \label{fig-sseq-8b}
\end{center} 
\end{figure}

\begin{rem}\label{rem-Ym}
{\bf The spectra $y_{m}$ and $Y_{m}$ of Corollaries \ref{cor-filtration} and
\ref{cor-Filtration}}.  Similar statements can be proved for the case
$\ell <0$.  We leave the details to the reader, but illustrate the
results in Figures \ref{fig-sseq-8a} and \ref{fig-sseq-8b}.

The source of each differential in \ref{prop-d3} is the product of
some element in $\upi_{\star}\HZZ$ with a power of $\normrbar_{1}^{}$
or $\delta_{1}$.  The target is the product of a different element in
$\upi_{\star}\HZZ$ with the same power.  This means they are
differentials in the slice {\SS } for the spectra $y_{m}$ of
\ref{cor-filtration}.

Similar differentials occur when we replace $\normrbar_{1}^{i}$ by any
homogeneous polynomial of degree $i$ in $\normrbar_{1}^{}$ and
$\ot_{2}$ in which the coefficient of $\normrbar_{1}^{i}$ is odd.
This means they are also differentials in the slice {\SS } for the
spectra $Y_{m}$ of \ref{cor-Filtration}.
\end{rem}

\begin{figure}
\begin{center}
\includegraphics[width=11.11cm]{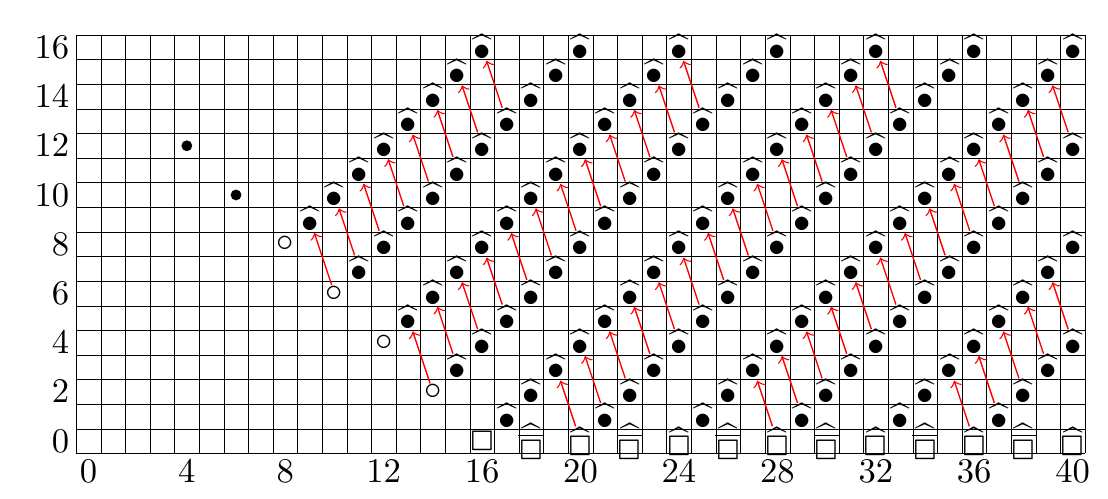}\\
\includegraphics[width=11cm]{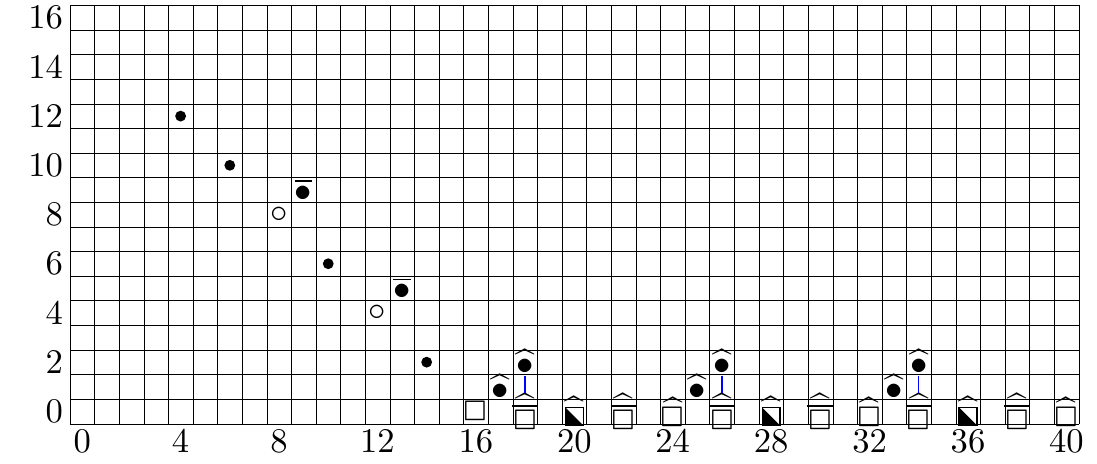} \caption[The slice
{$\EE_{2}$}- and $\EE_{4}$-terms for $X_{4,n}$ for $n\geq 4$.]{The
subquotient of the slice $\EE_{2}$- and $\EE_{4}$-terms for $\kH$ for
the slice summands $X_{4,n}$ for $n\geq 4$. Exotic transfers are shown
in blue and differentials are in red. The symbols are defined in Table
\ref{tab-C4Mackey}.  This is also the slice {\SS } for $y_{4}$ as in
Corollary \ref{cor-filtration} and $Y_{4}$ (after tensoring with $R$) as in
Corollary \ref{cor-Filtration}.}  \label{fig-sseq-7a}
\end{center}
\end{figure}

\begin{figure} 
\begin{center}
\includegraphics[width=11.11cm]{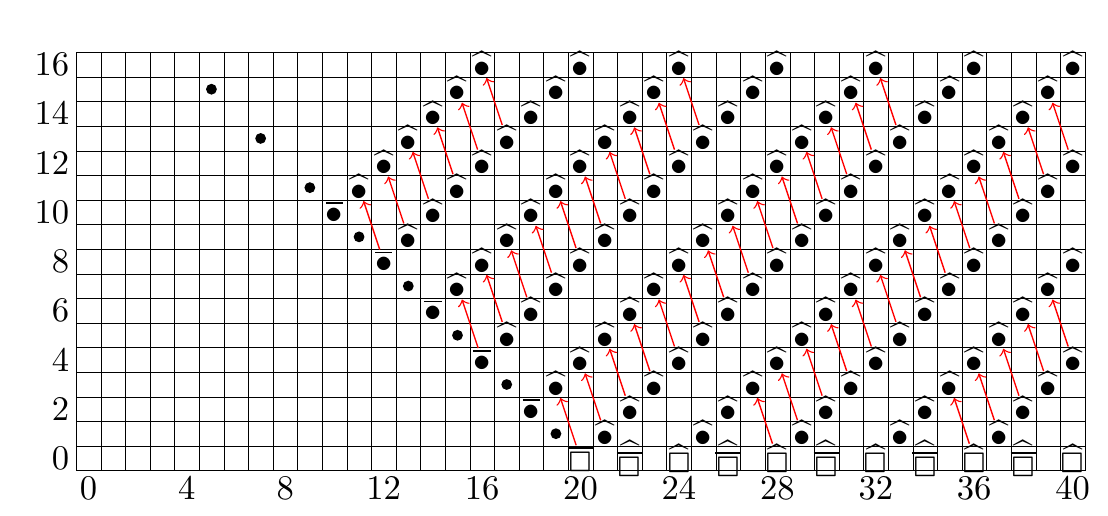}\\ 
\includegraphics[width=11cm]{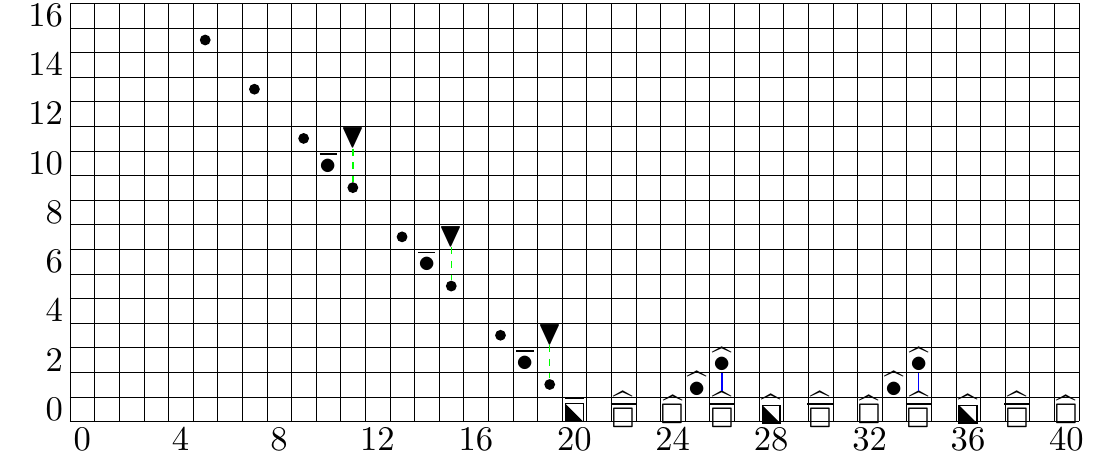} \caption[The slice
$\EE_{2}$ and $\EE_{4}$-terms for for $X_{5,n}$ for $n\geq
5$.]{The subquotient of the slice $\EE_{2}$ and $\EE_{4}$-terms for
$\kH$ for the slice summands $X_{5,n}$ for $n\geq 5$.  Exotic
restrictions and transfers are shown in dashed green and solid blue
lines respectively.  This is also the slice {\SS } for $y_{5}$ as in
Corollary \ref{cor-filtration} and for $Y_{5}$ (after tensoring with
$R$) as in Corollary \ref{cor-Filtration}.}  \label{fig-sseq-7b}
\end{center}
\end{figure}

These differentials are illustrated in the upper charts in  Figures
\ref{fig-sseq-7a}--\ref{fig-sseq-8b}.  In order to pass to $\EE_{4}$
we need the following exact sequences of Mackey functors.
\begin{displaymath}
\xymatrix
@R=5mm
@C=10mm
{
0 \ar[r]^(.5){}
  &\bullet\ar[r]^(.5){} 
     &\circ\ar[r]^(.5){d_{3}} 
        &{\widehat{\bullet}}\ar[r]^(.5){} 
           &{\obull }\ar[r]^(.5){} 
              &0\\
0\ar[r]^(.5){} 
  &{\widehat{\twobox}}\ar[r]^(.5){} 
     &{\widehat{\Box}}\ar[r]^(.5){d_{3}} 
        &{\widehat{\bullet} }\ar[r]^(.5){} 
           &0\\
  &0\ar[r]^(.5){} 
     &\obull\ar[r]^(.5){d_{3}} 
        &{\widehat{\bullet}}\ar[r]^(.5){} 
           &{\JJ }\ar[r]^(.5){} 
              &0\\
0 \ar[r]^(.5){}
  &\odbox\ar[r]^(.5){} 
     &\oBox\ar[r]^(.5){d_{3}} 
        &{\widehat{\bullet}}\ar[r]^(.5){} 
           &\JJ\ar[r]^(.5){} 
              &0
}
\end{displaymath}

The resulting subquotients of $\EE_{4}$ are shown in the lower charts
of Figures \ref{fig-sseq-7a}--\ref{fig-sseq-8b} and described below in
Theorem \ref{thm-sliceE4}.  In the latter the slice summands are
organized as shown in the Figures rather than by orbit type as in
Theorem \ref{thm-sliceE2}.

\begin{thm}\label{thm-sliceE4}
{\bf The slice $\EE_{4}$-term for $\kH$.}  The elements of
Theorem \ref{thm-sliceE2} surviving to $\EE_{4}$, which live
in the appropriate subquotients of $\upi_{*}X_{m,n}$, are
as follows.

\begin{enumerate}
\item 
In $\upi_{*}X_{2\ell ,2\ell }$ (see the
leftmost diagonal in Figure \ref{fig-sseq-7a}), on the 0-line we still
have a copy of $\Box$ generated under fixed point restrictions by
$\Delta_{1}^{\ell }\in \EE_{4}^{0,8\ell }$.  In positive
filtrations we have

\begin{align*}
\circ &\subseteq  \EE_{4}^{2j,8\ell }
\qquad \mbox{generated by } \\
&\mycases{
a_{\lambda}^{j}u_{2\sigma}^{\ell }
          u_{\lambda}^{2\ell -j}\normrbar_{1}^{2\ell }
  \in   \EE_{4}^{2j,8\ell } (G/G)\\
       &\hspace{-5cm} \mbox{for $j$ even and }0<j\leq 2\ell,  \\
2a_{\lambda}^{j}u_{2\sigma}^{\ell }
          u_{\lambda}^{2\ell -j}\normrbar_{1}^{2\ell }
=a_{\sigma}^{2}a_{\lambda}^{j-1}u_{2\sigma}^{\ell +1}
          u_{\lambda}^{2\ell -j-1}\normrbar_{1}^{2\ell }
  \in   \EE_{4}^{2j,8\ell } (G/G)\\
       &\hspace{-5cm} \mbox{for  $j$ odd and }0<j\leq 2\ell \mbox{ and}}\\
\bullet &\subseteq  \EE_{4}^{2k+2\ell ,8\ell }
\qquad \mbox{generated by } \\
&a_{\sigma}^{2k}a_{\lambda}^{2\ell }u_{2\sigma}^{\ell -k}
          \normrbar_{1}^{2\ell }
  \in   \EE_{4}^{2j+2k,8\ell } (G/G)
           \qquad \mbox{for }0<k\leq \ell.  
\end{align*}

\item 
In $\upi_{*}X_{2\ell ,2\ell +1}$ (see the
second leftmost diagonal in Figure \ref{fig-sseq-7a}), in filtration 0
we have $\widehat{\oBox }$, generated (under transfers and
the group action) by
\begin{displaymath}
r_{1}\Res_{1}^{2}(u_{\sigma  }^{2\ell } 
\Res_{1}^{4}(u_{\lambda}^{2\ell }\normrbar_{1}^{2\ell }) 
\in \EE_{2}^{0,8\ell+2} (G/\ee) .
\end{displaymath}

\noindent In positive filtrations we have
\begin{displaymath}
\begin{array}[]{rll}
\widehat{\bullet}
 & \subseteq
      \EE_{4}^{1,8\ell +2}
         &\mbox{generated (under transfers and the group action) by} \\
 &       & \eta u_{\sigma}^{2\ell }\Res_{2}^{4}(u_{\lambda}
                \normrbar_{1})^{2\ell }  =  \EE_{4}^{1,8\ell +2} (G/G')\\
\obull 
 & \subseteq   
      \EE_{4}^{4k+1,8\ell +2}
         &\mbox{for $0<k\leq \ell $ generated by} \\                        
 &       &  x=\eta^{4k+1} u_{\sigma}^{2\ell-2k}
      \Res_{2}^{4}(u_{\lambda}\normrbar_{1})^{2\ell -2k}
             \in  \EE_{4}^{4k+1,8\ell +2} (G/G')\\
 &       & \mbox{with $(1-\gamma )x=\Tr_{2}^{4}(x)=0$.}
\end{array}
\end{displaymath}

\item 
In $\upi_{*}X_{2\ell+1 ,2\ell +1}$ (see the leftmost diagonal
in Figure \ref{fig-sseq-7b}), on the 0-line we have a copy of $\odbox$
generated under fixed point $\Delta_{1}^{(2\ell +1)/2 }\in
\EE_{4}^{0,8\ell+4}$.  In positive filtrations we have
\begin{displaymath}
\begin{array}[]{rll}
\obull  
    &\subseteq 
         \EE_{2}^{2j,8\ell+4 }
            & \mbox{generated by }\\
    &      & u_{\sigma}^{2\ell +1}
     \Res_{2}^{4}(a_{\lambda}^{j}u_{\lambda}^{2\ell +1-j}
                  \normrbar_{1}^{2\ell +1})\in  \EE_{2}^{2j,8\ell+4 } (G/G')\\
    &      &\mbox{for }0<j\leq 2\ell+1 , \\
\bullet 
    &\subseteq 
         \EE_{2}^{2j+1 ,8\ell+4 }
            & \mbox{generated by }\\
    &      & a_{\sigma}a_{\lambda}^{j}
              u_{\sigma}^{2\ell}u_{\lambda}^{2\ell +1-j}
                  \normrbar_{1}^{2\ell +1}\in \EE_{2}^{2j+2k,8\ell+4 } (G/G)\\
    &      & \mbox{for }0\leq j\leq 2\ell+1\mbox{ and} \\
\bullet 
    &\subseteq 
         \EE_{2}^{2k+4\ell +3 ,8\ell+4 }
            &\mbox{generated by }\\
    &      &a_{\sigma}^{2k+1}a_{\lambda}^{2\ell +1}
       u_{2\sigma}^{\ell-k}
                  \normrbar_{1}^{2\ell +1}
              \in   \EE_{2}^{2k+4\ell +2,8\ell+4 } (G/G)\\
    &      &\mbox{for }0< k\leq 2\ell+1.  
\end{array}
\end{displaymath}

\item 
In $\upi_{*}X_{2\ell+1 ,2\ell +2}$ (see the
second leftmost diagonal in Figure \ref{fig-sseq-7b}), in filtration 0
we have $\widehat{\oBox }$, generated (under transfers and
the group action) by
\begin{displaymath}
r_{1}\Res_{1}^{2}(u_{\sigma  }^{2\ell+1 } 
\Res_{1}^{4}(u_{\lambda}^{2\ell+1 }\normrbar_{1}^{2\ell+1 }) 
\in \EE_{4}^{0,8\ell+6} (G/\ee) .
\end{displaymath}

\noindent In positive filtrations we have
\begin{displaymath}
\begin{array}[]{rll}
\JJ & \subseteq  
        \EE_{4}^{4k+3,8\ell+6}
           &\mbox{for $0\leq k\leq \ell $ generated under transfer by } \\  
    &      &x=\eta^{4k+3}\Delta_{1}^{\ell -k}
                 \in   \EE_{4}^{4k+3,8\ell+6} (G/G')\\  
    &      &\mbox{with $(1-\gamma )x=0$.}
\end{array}
\end{displaymath}

\noindent The generator of $\EE_{4}^{4k+3,8\ell+6} (G/G')$ is the
exotic restriction of the one in $\EE_{4}^{4k+1,8\ell+4} (G/G)$.

\item 
In $\upi_{*}X_{m,m+i}$ for $i\geq 2$ (see the
rest of Figures \ref{fig-sseq-7a} and \ref{fig-sseq-7b}), in
filtration 0 we have
\begin{displaymath}
\begin{array}[]{rll}
\widehat{\oBox }
    & \subseteq   
          \EE_{4}^{0,4m+4j+2}
           &\mbox{generated under transfers and group action by }\\  
    &      &r_{1}\Res_{1}^{2}(u_{\sigma  }^{m +j} \os_{2}^{j})
                \Res_{1}^{4}(u_{\lambda}^{m +j}\normrbar_{1}^{m })\\  
    &      &\quad\in \EE_{4}^{0,4m+4j+2} (G/\ee)\mbox{ for $j\geq 0$,}\\
\widehat{\twobox}
    & \subseteq   
          \EE_{4}^{0,8\ell +4}
           &\mbox{generated under transfers and group action by }\\  
    &      &r_{1}\Res_{1}^{2}(u_{\sigma  }^{m +j} \os_{2}^{j})
             \Res_{1}^{4}(u_{\lambda}^{m +j}\normrbar_{1}^{m })\\  
    &      &\quad\in     \EE_{4}^{0,8\ell  +4} (G/\ee)
             \mbox{ for $\ell \geq m/2$ and} \\
\widehat{\Box}
    & \subseteq   
          \EE_{4}^{0,8\ell }
           &\mbox{generated under transfers, restriction and group action by}\\
    &      &x_{8\ell ,m}
             =\Sigma_{2,0}^{2\ell -m}\delta_{1}^{m}+\ell \delta_{1}^{2\ell }
              \mbox{ where } \\
    &      &\Sigma_{2,\epsilon }=u_{\rho_{2}}\os_{2,\epsilon } \\
    &      &\mbox{and } \delta_{1}=u_{\rho_{2}}\Res_{2}^{4}(\normrbar_{1})\\ 
    &      &\quad\in \EE_{4}^{0,8\ell  } (G/G')\mbox{ for $0\leq m\leq 2\ell -1$}.
\end{array}
\end{displaymath}

\noindent In positive filtrations we have 
\begin{displaymath}
\begin{array}[]{rll}
\widehat{\bullet}
    &  \subseteq 
         \EE_{4}^{2,8\ell  +4}
           &\mbox{generated under transfers and group action by }\\ 
    &      &\eta_{0}^{2}\Res_{2}^{4}(\Delta_{1}^{\ell})
               =   \eta_{0}^{2}\delta_{1}^{2\ell } = 
\eta_{0}^{2}u_{\sigma}^{2\ell}\Res_{2}^{4}(u_{\lambda}\normrbar_{1})^{2\ell}\\ 
    &      &\quad\in   \EE_{4}^{2,8\ell +4} (G/G')
\quad \mbox{and}\\                       
\widehat{\bullet}
    &  \subseteq 
         \EE_{4}^{s,8\ell +2s}
           &\mbox{generated under transfers and group action by }\\ 
    &      &\eta_{\epsilon }^{s}x_{8\ell ,m}\in\EE_{4}^{s,8\ell +2s} (G/G')\\ 
    &      &  \mbox{for $s=1,2$ and $0\leq m\leq 2\ell-1 $}. 
\end{array}
\end{displaymath}

\noindent Each generator of $\EE_{4}^{2,8\ell +4} (G/G')$ is an exotic
transfer of one in\linebreak  $\EE_{4}^{0,8\ell +2} (G/e)$.
\end{enumerate}
\end{thm}

\begin{prop}\label{prop-perm}
{\bf Some nontrivial permanent cycles.} The elements listed in Theorem
\ref{thm-sliceE4}(v) other than $\eta_{\epsilon }^{2}\delta_{1}^{2\ell
}$ are all nontrivial permanent cycles.
\end{prop}

\proof Each such element is either in the image of
$\EE_{4}^{0,*} (G/\ee)$ under the transfer and
therefore a nontrivial permanent cycle, or it is one of the ones listed
in Corollary \ref{cor-G'SS}.  \qed\bigskip

{\em In subsequent discussions and charts, starting with Figure
\ref{fig-E4}, we will omit the elements in Proposition
\ref{prop-perm}.  These elements all occur in $\EE_{4}^{s,t}$ for
$0\leq s\leq 2$.}

Analogous statements can be made about the slice {\SS} for $\KH$. Each
of its slices is a certain infinite wedge spelled out in Corollary
\ref{cor-KH}.  Their homotopy groups are determined by the chain
complex calculations of Section \ref{sec-chain} and illustrated in
Figures \ref{fig-sseq-1} (with Mackey functor induction
$\uparrow_{2}^{4}$ applied) and \ref{fig-sseq-3}.  Analogs of Figures
\ref{fig-sseq-7a}--\ref{fig-sseq-7b} are shown in Figures
\ref{fig-sseq-8a}--\ref{fig-sseq-8b}.  In each figure, exotic
transfers and restrictions are indicated by blue and dashed green lines
respectively.  As in the $\kH$ case, most of the elements shown in
this chart can be ignored for the purpose of calculating higher
differentials.  {\em In the third quadrant the elements we are
ignoring all occur in $\EE_{4}^{s,t}$ for $-2\leq s\leq 0$.}

The resulting reduced $\EE_{4}$ for $\KH$ is shown in Figure
\ref{fig-KH}.  The information shown there is very useful for
computing differentials and extensions.  The periodicity theorem tells
us that $\upi_{n}\KH$ and $\upi_{n-32}\KH$ are isomorphic.  For $0\leq
n<32$ these groups appear in the first and third quadrants
respectively, and the information visible in the {\SS} can be quite
different.

For example, we see that $\upi_{0}\KH$ has summand of the form $\Box$,
while $\upi_{-32}\KH$ has a subgroup isomorphic to $\fourbox $.  The
quotient $\Box/\fourbox $ is isomorphic to $\circ $.  This leads to
the exotic restrictions and transfer in dimension $-32$ shown in
Figure \ref{fig-KH}.  Information that is transparent in dimension 0
implies subtle information in dimension $-32$.  Conversely, we see
easily that $\upi_{-4}\KH=\dot{\twobox}$ while $\upi_{28}\KH$ has a
quotient isomorphic to $\odbox$.  This leads to the ``long transfer''
(which raises filtration by 12) in dimension 28.


\section{Higher differentials and exotic Mackey functor extensions}
\label{sec-higher}

\begin{figure}
\begin{center}
\includegraphics[width=11cm]{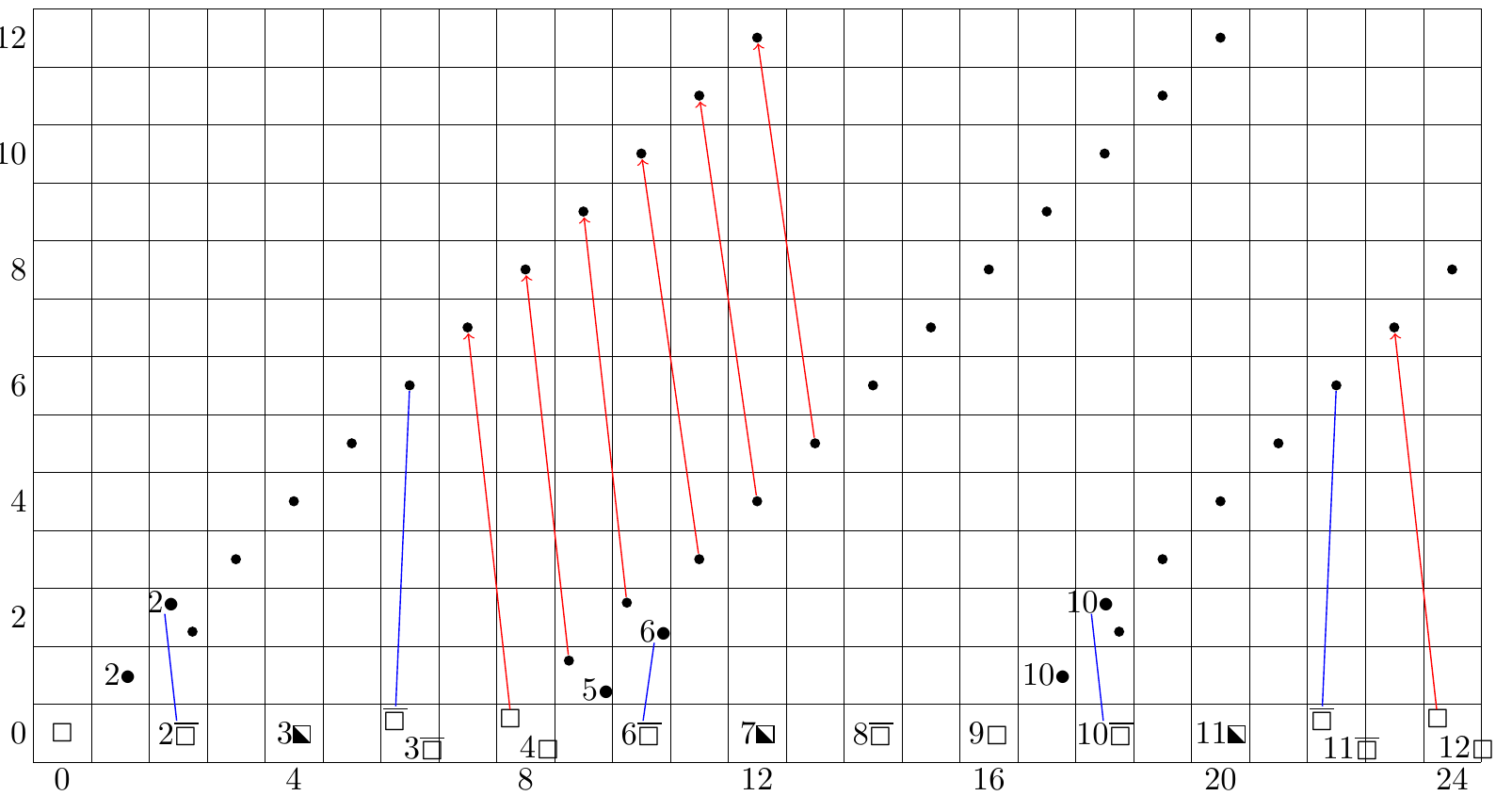} \caption[The slice
$\EE_{7}$-term for the $C_{2}$-spectrum $i_{G'}^{*}\kH$.]{The slice
$\EE_{7}$-term for the $C_{2}$-spectrum $i_{G'}^{*}\kH$.  The Mackey
functor symbols are defined in Table \ref{tab-C2Mackey}. A number
$n$ in front of a symbol indicates an $n$-fold direct sum.  Blue lines
indicate exotic transfers and red lines indicate differentials.}  
\label{fig-G'SS-3}
\end{center}
\end{figure}

We can use the results of the \S\ref{sec-C2diffs} to study higher
differentials and extensions.  The $\EE_{7}$-term implied by them is
illustrated in Figure \ref{fig-G'SS-3}. For each $\ell ,s\geq 0$ there is
a generator
\begin{displaymath}
y_{8\ell +s,s}:=\eta_{0}^{s} \delta_{1}^{2\ell }
\in \EE_{7}^{s,8\ell +2s} (G/G')
\end{displaymath}

\noindent with
\begin{displaymath}
d_{7} (y_{16k+s+8,s}) = y_{16k+s+7,s+7}.
\end{displaymath}

\noindent 
Recall that 
\begin{displaymath}
\delta_{1}=\ou_{\lambda }\orr_{1,0}\orr_{1,1}
\in \underline{E}_{2}^{0,4}\kH (G/G')
\cong \underline{E}_{2}^{0,(G',4)}\kH (G/G),
\end{displaymath}

\noindent and in the latter group we denote $\ou_{\lambda }$ by
$u_{2\sigma }$.  We have 
\begin{displaymath}
d_{3} (\delta_{1})= d_{3} (\ou_{\lambda })\orr_{1,0}\orr_{1,1}
\cong d_{3} (u_{2\sigma  })\orr_{1,0}\orr_{1,1}
=a_{\sigma }^{3} (\orr_{1,0}+\orr_{1,1})\orr_{1,0}\orr_{1,1}.
\end{displaymath}

\noindent

If the source has the form $\Res_{2}^{4}(x_{16k+s+8,s})$, then such
an $x$ must support a nontrivial $d_{r}$ for $r\leq 7$.  If it has a
nontrivial transfer $x'_{16k+s+8,s}$, then such an $x'$ cannot support an
earlier differential, and we must have
\begin{displaymath}
d_{r} (x'_{16k+s+8,s})
  =\Tr_{2}^{4}(d_{7} (y_{16k+s+8,s}))
  =\Tr_{2}^{4}(y_{16k+s+7,s+7})   \qquad \mbox{for some $r\geq 7$}.
\end{displaymath}

\noindent We could get a higher differential (meaning $r>7$) if
$y_{16k+s+7,s+7}$ supports an exotic transfer.

\begin{figure}   
\begin{center} 
\includegraphics[width=11cm]{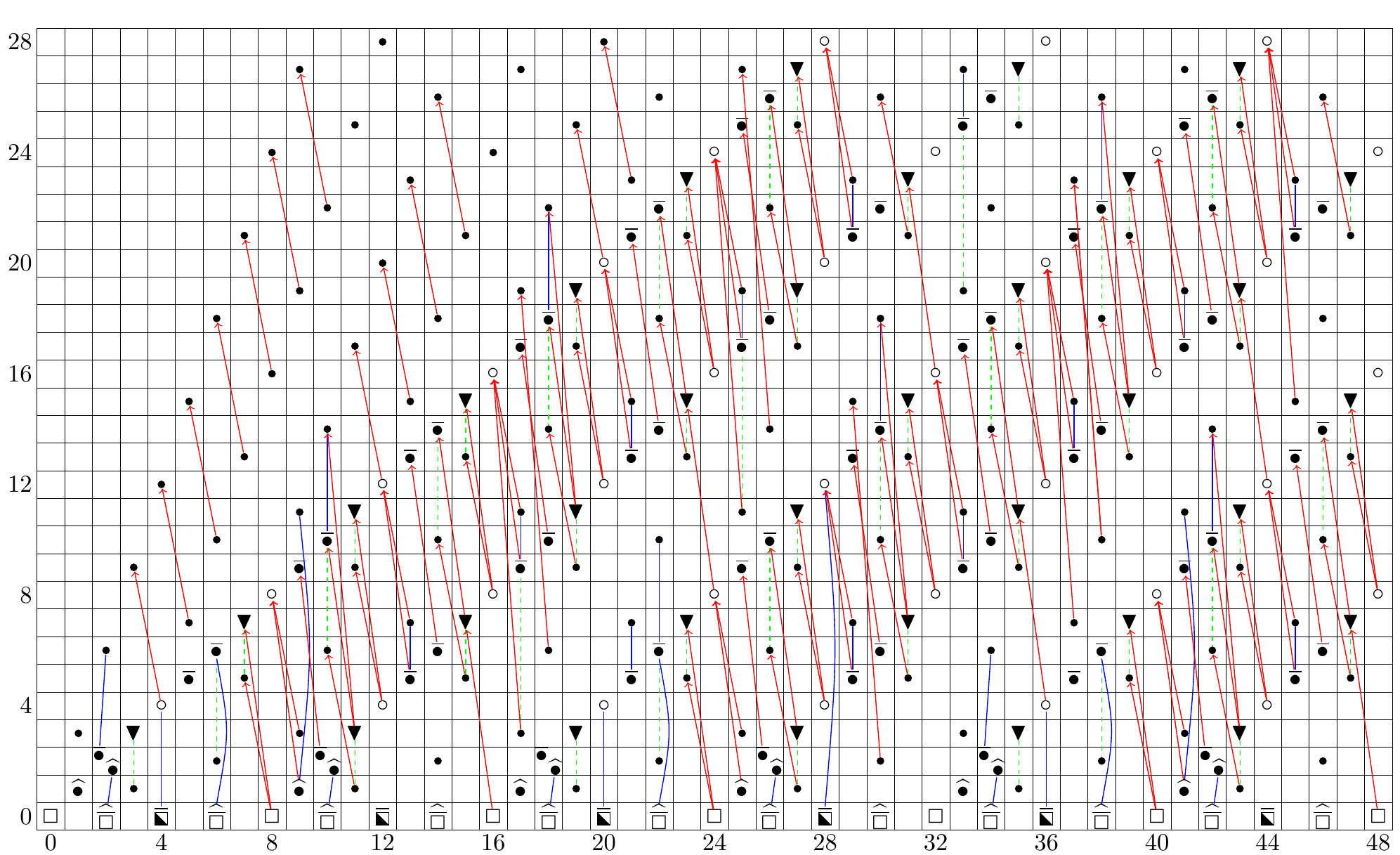} \caption[The reduced
$\EE_{4}$-term of the slice {\SS} for $\kH$.]{The $\EE_{4}$-term of
the slice {\SS} for $\kH$ with elements of Proposition \ref{prop-perm}
removed. Differentials are shown in red.  Exotic transfers and
restrictions are shown as solid blue and dashed green lines
respectively. The Mackey functor symbols are as in Table
\ref{tab-C4Mackey}.}  \label{fig-E4}
\end{center}
\end{figure}

We have seen (Figure \ref{fig-E4} and Theorem \ref{thm-sliceE4}) that
for $s\geq 3$ and $k\geq 0$,
\begin{numequation}\label{eq-eighth-diagonls}
\begin{split}
\EE_{5}^{s,16k+8+2s}
 = \mycases{
\circ
       &\mbox{for $s\equiv 0$ mod 4}\\
\overline{\bullet} 
       &\mbox{for $s\equiv 1,2$ mod 4}\\
\JJ
       &\mbox{for $s\equiv 3$ mod 4.}
}
\end{split}
\end{numequation}%

\noindent For $s=1,2$, $\EE_{5}^{s,16k+8+2s}$ has
$\overline{\bullet}$ as a direct summand. For $s=0$ it has $\Box$ as a
summand, and the differentials on it factor through its quotient
$\circ$; see (\ref{eq-C4-SES}).

\begin{landscape}
\begin{figure}
\begin{center}
\includegraphics[width=18cm]{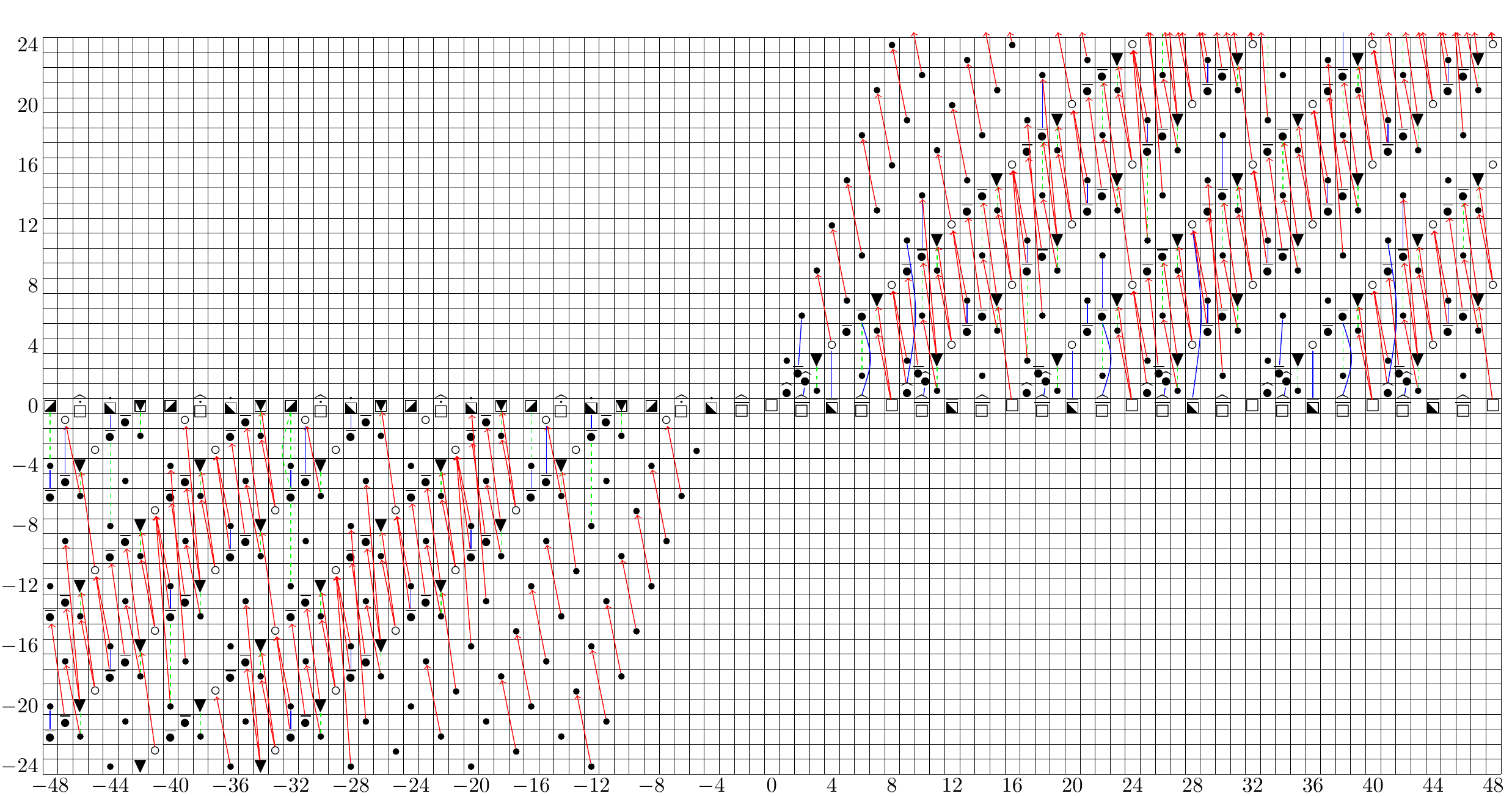} \caption[The reduced
$\EE_{4}$-term of the slice {\SS} for the periodic spectrum
$\KH$.]{The reduced $\EE_{4}$-term of the slice {\SS} for the periodic
spectrum $\KH$.  Differentials are shown in red.  Exotic transfers and
restrictions are shown in solid blue and dashed green vertical lines
respectively. The Mackey functor symbols are indicated in the table
below Figure \ref{fig-E14a}.  }
  
\label{fig-KH}
\end{center}
\end{figure}
\end{landscape}

The corresponding statement in the third quadrant is 
\begin{displaymath}
\EE_{5}^{-s,-16k-2s-24}
 = \mycases{
\circ
       &\mbox{for $s\equiv 3$ mod 4}\\
\overline{\bullet} 
       &\mbox{for $s\equiv 1,2$ mod 4}\\
\JJ
       &\mbox{for $s\equiv 0$ mod 4.}
}
\end{displaymath}

\noindent for $s\geq 3$ and $k\geq 0$.  For $s=1,2$ the groups have
similar summands, and for $s=0$ there is a summand of the form
$\JJbox$, which has $\JJ$ as a subgroup; again see (\ref{eq-C4-SES}).
This is illustrated in Figure \ref{fig-KH}.

\begin{thm}\label{thm-d7}
{\bf Differentials for $C_{4}$ related to the $d_{7}$s for $C_{2}$.}
The differential
\begin{displaymath}
d_{7} (y_{16k+s+8,s})=y_{16k+s+7,s+7}\qquad \mbox{with }s\geq 3 
\end{displaymath}

\noindent has the following implications for the congruence
classes of $s$ modulo 4.

\begin{enumerate}
\item [(i)] For $s\equiv 0$, $\EE_{7}^{s,16k+8+2s}=\circ$ and
$\EE_{7}^{s+7,16k+14+2s}=\JJ$.  Hence $y_{16k+s+8,s}$ is a
restriction with a nontrivial transfer, and
\begin{align*}
d_{5} (x_{16k+s+8,s})
 & = x_{16k+s+7,s+5}  \\
\aand 
d_{7} (2x_{16k+s+8,s})&=d_{7} (\Tr_{2}^{4}(y_{16k+s+8,s}))\\
 & = \Tr_{2}^{4}(y_{16k+s+7,s+7})= x_{16k+s+7,s+7}.
\end{align*}

\item [(ii)] For $s\equiv 1$,
\begin{align*}
d_{7} (y_{16k+s+8,s})
 & =  y_{16k+s+7,s+7}\\
\aand 
d_{5} (x_{16k+s+8,s+2})
 & =  \Tr_{2}^{4}(y_{16k+s+7,s+7}) = 2x_{16k+s+7,s+7}
\end{align*}

\noindent  This leaves the fate of $x_{16k+s+7,s+7}$ undecided; see below.


\item [(iii)] For $s\equiv 2$, 
$\EE_{7}^{s,16k+8+2s}=\overline{\bullet} $ and
$\EE_{7}^{s+7,16k+14+2s}=\overline{\bullet} $.  Neither the
source nor target is a restriction or has a nontrivial transfer, so no
additional differentials are implied.

\item [(iv)] For $s\equiv 3$,  $\EE_{7}^{s,16k+8+2s}=\JJ $ and
$\EE_{7}^{s+7,16k+14+2s}=\overline{\bullet} $.  In this case
the source is an exotic restriction; again see Figure
\ref{fig-sseq-7b}). Thus we have 
\begin{align*}
d_{7} (y_{16k+s+8,s})
 & =  y_{16k+s+7,s+7}\\
\aand 
d_{5} (x_{16k+s+8,s-2})
 & =  x_{16k+s++7,s+3}\\
 \mbox{with }
\Res_{2}^{4}(x_{16k+s++7,s+3})
 &  =y_{16k+s+7,s+7}.
\end{align*}
\noindent Moreover, $\Tr_{2}^{4}(y_{16k+8+s,s})$ is nontrivial and it
supports a nontrivial $d_{11}$ when $4k+s\equiv 3$ mod 8.  The other
case, $4k+s\equiv 7$, will be discussed below.
\end{enumerate}
\end{thm}

\proof (i) The target Mackey functor is $\JJ$ and $y_{16k+s+7,s+7}$ is
the exotic restriction of $x_{16k+s+7,s+5}$; see Figure
\ref{fig-sseq-7b} and Theorem \ref{thm-sliceE4}.  The indicated
$d_{5}$ and $d_{7}$ follow.

(ii)   The differential is nontrivial on the $G/G'$ component of
\begin{displaymath}
\xymatrix
@R=5mm
@C=10mm
{
{\obull=\EE_{7}^{s,16k+8+2s}}\ar[r]^(.5){d_{7}}
    &{\EE_{7}^{s+7,16k+14+2s}=\circ}
}
\end{displaymath}

\noindent Thus the target has a nontrivial transfer, so the source
must have an exotic transfer.  The only option is $x_{16k+s+8,s+2}$,
and the result follows.

(iv) We prove the statement about $d_{11}$ by showing that
\begin{displaymath}
{y_{16k+s+7,s+7}=\eta_{0}^{s+7}\delta_{1}^{4k}}
\end{displaymath}

\noindent supports an exotic
transfer that raises filtration by 4.  First note that 
\begin{align*}
\Tr_{2}^{4}(\eta_{0}\eta_{1}) 
 & = \Tr_{2}^{4}(a_{\sigma_{2}}^{2}\overline{r}_{1,0}\overline{r}_{1,0} )
   = \Tr_{2}^{4}(u_{\sigma }\Res_{2}^{4}(a_{\lambda }\normrbar_{1})) 
   = \Tr_{2}^{4}(u_{\sigma })a_{\lambda }\normrbar_{1}  \\
 & = a_{\sigma }a_{\lambda }\normrbar_{1}a_{\lambda }\normrbar_{1} 
        \qquad \mbox{by (\ref{eq-exotic-transfers})}.
\end{align*}

\noindent Next note that the three elements
\begin{displaymath}
y_{8,8}=\eta_{0}^{8}=\Res_{2}^{4}(\epsilon ),\quad 
y_{20,4}= \eta_{0}^{4}\delta_{1}^{4}=\Res_{2}^{4}(\overline{\kappa } )
\quad \aand 
y_{32,0}= \delta_{1}^{8}=\Res_{2}^{4} (\Delta^{4})
\end{displaymath}

\noindent are all permanent cycles, so the same
is true of all
\begin{displaymath}
y_{16m +4\ell ,4\ell }=\eta_{0}^{4\ell }\delta_{1}^{4m }
\qquad \mbox{for $m,\ell \geq 0$ and $m+\ell $ even.} 
\end{displaymath}

\noindent  It follows that for such $\ell $ and $m$,
\begin{align*}
\eta_{0}\eta_{1}y_{16m+4\ell ,4\ell }
 & = \eta_{0}\eta_{1}\eta_{0}^{4\ell }\delta_{1}^{4m}
=\eta_{0}^{4\ell+2 }\delta_{1}^{4m}
=y_{16m+4\ell+2 ,4\ell+2 }  \\
 & =   \eta_{0}\eta_{1}\Res_{2}^{4}(x_{16m+4\ell ,4\ell }),
\end{align*}

\noindent so
\begin{displaymath}
\Tr_{2}^{4}(y_{16m+4\ell+2 ,4\ell+2 })
 = \Tr_{2}^{4}(\eta_{0}\eta_{1})x_{16m+4\ell ,4\ell }
 = f_{1}^{2}x_{16m+4\ell ,4\ell }.
\end{displaymath}

\noindent  This is the desired exotic transfer.
\qed\bigskip

We now turn to the unsettled part of \ref{thm-d7}(iv).

\begin{thm}\label{thm-d7a}
{\bf The fate of $x_{16k+s+8,s}$ for $4k+s\equiv 7$ mod 8 and $s\geq
7$.}  Each of these elements is the target of a $d_{7}$ and hence a
permanent cycle.
\end{thm}

\proof Consider the element $\Delta_{1}^{2}\in
\underline{E}_{2}^{0,16} (G/G)$.  We will show that
\begin{displaymath}
d_{7} (\Delta_{1}^{2}) = x_{15,7}=\Tr_{2}^{4}(y_{15,7}).
\end{displaymath}
 
\noindent This is the case $k=0$ and $s=7$.  The remaining cases will
follow via repeated multiplication by $\epsilon $, $\overline{\kappa } $ and
$\Delta_{1}^{4}$.  

We begin by looking at 
\begin{displaymath}
\Delta_{1}=u_{2\sigma }u_{\lambda }^{2}\normrbar_{1}^{2}.
\end{displaymath}

\noindent  From Theorem \ref{thm-d3ulambda} we have 
\begin{displaymath}
d_{5} (u_{2\sigma }) = a_{\sigma}^{3}a_{\lambda}\normrbar_{1}
\qquad \aand 
d_{5} (u_{\lambda}^{2}) = a_{\sigma }a_{\lambda}^{2}u_{\lambda }\normrbar_{1}
\end{displaymath}

\noindent Using the gold relation $a_{\sigma }^{2}u_{\lambda }=2a_{\lambda
}u_{2\sigma }$, we have
\begin{align*}
d_{5} (\Delta_{1})
 = d_{5} (u_{2\sigma }u_{\lambda }^{2})\normrbar_{1}^{}
 & = (a_{\sigma}^{3}a_{\lambda}u_{\lambda}^{2}\normrbar_{1}^{}
         + a_{\sigma }a_{\lambda}^{2}u_{\lambda }u_{2\sigma }\normrbar_{1})
       \normrbar_{1}^{}  \\
 & = a_{\sigma }a_{\lambda }u_{\lambda }(a_{\sigma}^{2}u_{\lambda}
 + a_{\lambda}u_{2\sigma })\normrbar_{1}^{2}  \\
 & = a_{\sigma }a_{\lambda }u_{\lambda }(2a_{\lambda }u_{2\sigma }
 + a_{\lambda}u_{2\sigma })\normrbar_{1}^{2}  \\
 & = a_{\sigma }a_{\lambda }^{2}u_{\lambda }u_{2\sigma }\normrbar_{1}^{2} 
 \qquad \mbox{since }2a_{\sigma }=0  \\
  & = \nu  x_{4}.
\end{align*}

\noindent Since $\nu $ supports an exotic group extension, $2\nu
=x_{3}$, we have
\begin{displaymath}
2d_{5} (\Delta_{1})= d_{7} (2\Delta_{1})= x_{3}x_{4}.
\end{displaymath}

\noindent From this it follows that 
\begin{displaymath}
d_{7} (\Delta_{1}^{2})=\Delta_{1}d_{7} (2\Delta_{1})=x_{15,7}
\end{displaymath}

\noindent as claimed.
\qed\bigskip

\begin{figure}    
\begin{center} 
\includegraphics[width=11cm]{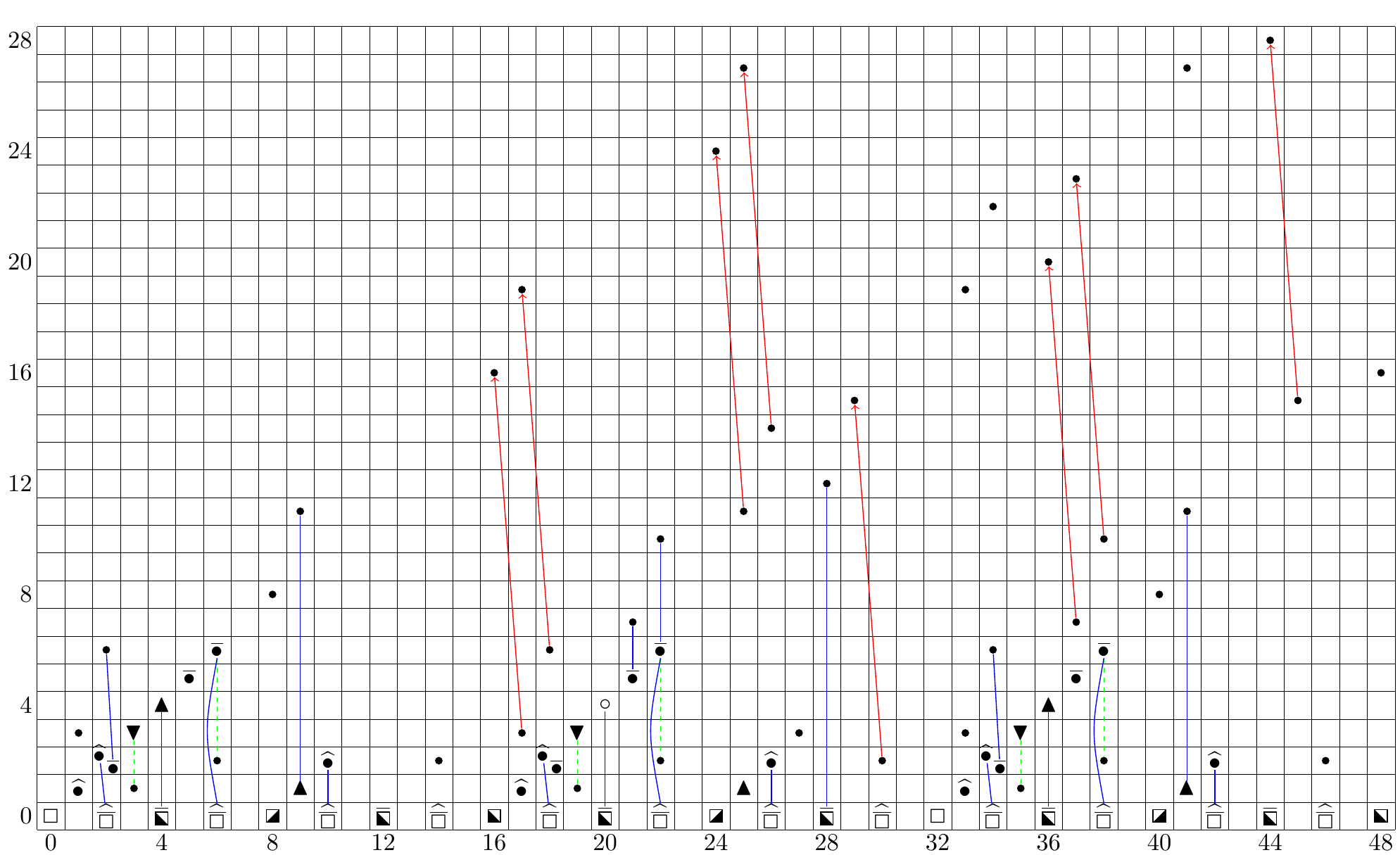} \caption[The reduced
$\EE_{4}$-term of the slice {\SS} for $\kH$.]{The $\EE_{12}$-term of
the slice {\SS} for $\kH$ with elements of Proposition \ref{prop-perm}
removed. Differentials are shown in red.  Exotic transfers and
restrictions are shown as solid blue and dashed green lines
respectively. The Mackey functor symbols are as in Table
\ref{tab-C4Mackey}.}  \label{fig-E12}
\end{center} 
\end{figure}

The resulting reduced $\underline{E}_{12}$-term is shown in Figure
\ref{fig-E12}.  It is sparse enough that the only possible remaining
differentials are the indicated $d_{13}$s.  In order to establish them
we need the following.



 
The surviving class in $\underline{E}_{12}^{20,3} (G/G)$ is 
\begin{displaymath}
x_{17,3}=f_{1}\Delta_{1}^{2}
 = a_{\sigma }a_{\lambda }\normrbar_{1}\cdot
[u_{2\sigma }^{2}]u_{\lambda }^{4}\normrbar_{1}^{4} = (a_{\sigma
}u_{\lambda }^{4}) (a_{\lambda }[u_{2\sigma }^{2}]\normrbar_{1}^{5}).
\end{displaymath}

\noindent The second factor is a permanent cycle, so Theorem
\ref{thm-normed-slice-diffs} gives 
\begin{displaymath}
d_{13} (f_{1}\Delta_{1}^{2})
 = (a_{\lambda }^{7}[u_{2\sigma }^{2}]\normrbar_{1}^{3})
     (a_{\lambda }[u_{2\sigma }^{2}]\normrbar_{1}^{5})
=a_{\lambda }^{8}[u_{2\sigma }^{2}]\normrbar_{1}^{8}
 = \epsilon^{2} = x_{4}^{4}.
\end{displaymath}

\noindent The surviving class in  $\underline{E}_{12}^{32, 2} (G/G)$ is 
\begin{displaymath}
x_{30,2}
  =a_{\sigma }^{2}u_{2\sigma }^{3}u_{\lambda }^{8}\normrbar_{1}^{8}
  \in \underline{E}_{12}^{32,2} (G/G)
\end{displaymath}

\noindent satisfies
\begin{displaymath}
\epsilon x_{30,2}
 =f_{1}\overline{\kappa }x_{17,3}
= f_{1}^{2}x_{4}\Delta_{1}^{2}, 
\end{displaymath}

\noindent  so we have proved the following.

\begin{thm}\label{thm-d13}
{\bf $d_{13}$s in the slice {\SS } for $\kH$.}  There are
differentials
\begin{displaymath}
d_{13} (f_{1}^{\epsilon }x_{4}^{m}\Delta_{1}^{2n})
 = f_{1}^{\epsilon -1}x_{4}^{m+4}\Delta_{1}^{2 (n-1)}
\end{displaymath}

\noindent for $\epsilon =1,2$, $m+n$ odd, $n\geq 1$ and $m
\geq 1-\epsilon $. The {\SS } collapses from $\underline{E}_{14}$.

\end{thm}

\bigskip

To finish the calculation we have

\begin{landscape}
\begin{figure}
\begin{center}
\includegraphics[width=18.2cm]{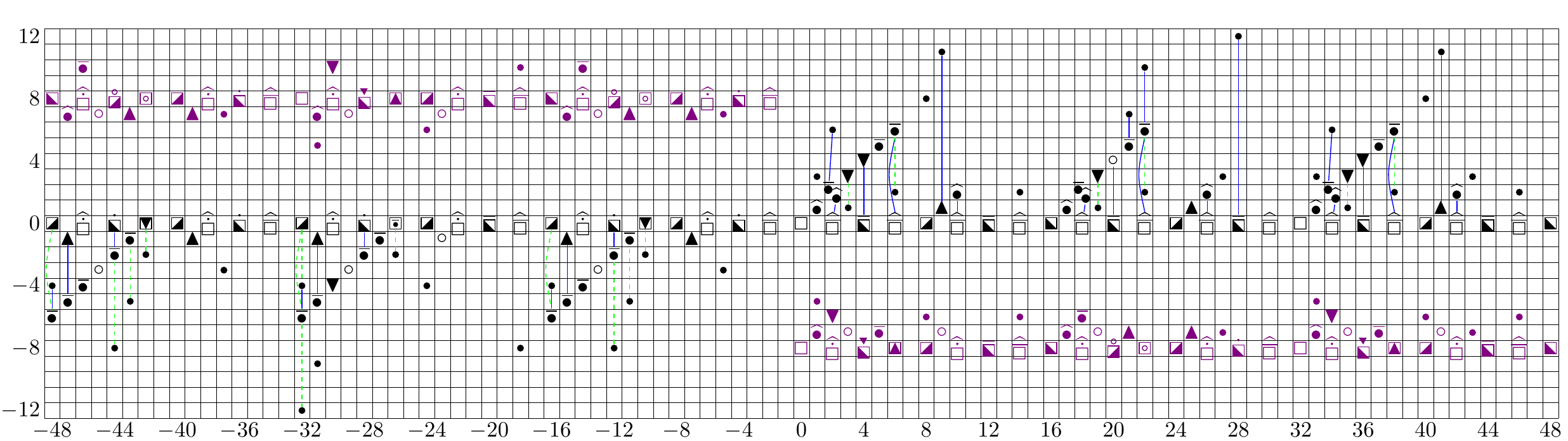} 
\caption[The reduced $\EE_{\infty }$-term of the slice {\SS} for
$\KH$.]{The reduced $\EE_{14}=\EE_{\infty }$-term of the slice {\SS} for
$\KH$. The exotic Mackey functor extensions lead to the Mackey
functors shown in violet in the second and fourth quadrants. The Mackey
functor symbols are indicated below. 
}   
\includegraphics[width=18.2cm]{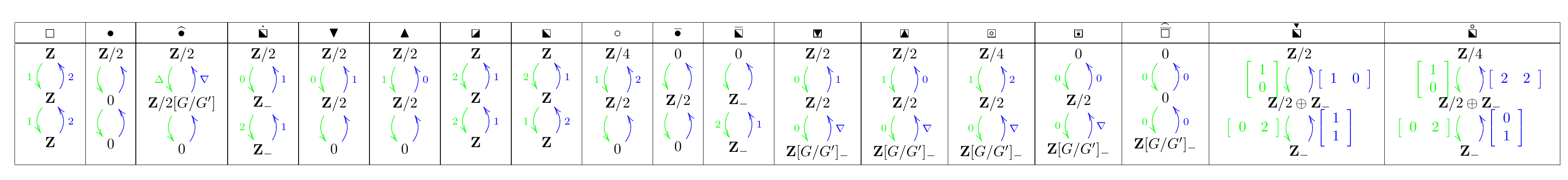} \label{fig-E14a}
\end{center}
\end{figure}
\end{landscape}

\begin{thm}\label{thm-double}  
{\bf Exotic transfers from and restrictions to the 0-line.} In
$\upi_{*}\kH$, for $i\geq 0$ we have
\begin{align*}   
\Tr_{1}^{2} (r_{1,\epsilon }r_{1,0}^{4i}r_{1,1}^{4i})
 & =  \eta _{\epsilon }^{2}\orr_{1,0}^{4i}\orr_{1,1}^{4i}  
      &&\in \upi_{8i+2}
          &\mbox{(filtration jump 2)}\\ 
\Tr_{1}^{4} (r_{1,0}^{8i+1}r_{1,1}^{8i+1})
 & =   2x_{4}\Delta_{1}^{4i}  
      &&\in \upi_{32i+4}
          &\mbox{(filtration jump 4)}\\ 
\Tr_{1}^{2} ((r_{1,0}^{3}+r_{1,1}^{3}) r_{1,0}^{8i}r_{1,1}^{8i})
 & =  \eta _{0}^{3}\eta _{1}^{3}\delta_{1}^{8i}  
      &&\in \upi_{32i+6}
          &\mbox{(filtration jump 6)}\\ 
\Tr_{1}^{4}(r_{1,0}^{8i+5}r_{1,1}^{8i+5})
 & =   2x_{4}\Delta_{1}^{4i+2}
      &&\in \upi_{32i+20}
          &\mbox{(filtration jump 4)}\\
\Tr_{1}^{2} ((r_{1,0}^{3}+r_{1,1}^{3}) r_{1,0}^{8i+4}r_{1,1}^{8i+4})
 & =  \eta _{0}^{3}\eta _{1}^{3}\delta_{1}^{8i+4}  
      &&\in \upi_{32i+22}
          &\mbox{(filtration jump 6)}\\ 
\aand  
\Tr_{2}^{4}(2\delta_{1}^{8i+7})
 & =   x_{4}^{3}\Delta_{1}^{4i+2}
      &&\in \upi_{32i+28}
          &\mbox{(filtration jump 12,}\\
 &    &&  &\mbox{the long transfer)}. 
\end{align*}
 
Let $\underline{M}_{k}$ denote the reduced value of $\upi_{k}\kH$,
meaning the one obtained by removing the elements of Proposition
\ref{prop-perm}.  Its values are shown in purple in Figure
\ref{fig-E14a}, and each has at most two summands. For even $k$ one
of them contains torsion free elements, and we denote it by
$\underline{M}'_{k}$.  Its values depend on $k$ mod 32 and are as
follows, with symbols as in Table \ref{tab-C4Mackey}.
 
\begin{center}
\begin{tabular}[]{|c||c|c|c|c|c|c|c|c|c|c|c|c|c|c|c|c|}
\hline 
$k$
   &0 &2 &4 &6 &8 &10&12&14 
   &16&18&20&22&24&26&28&30\\
\hline 
$\underline{M}'_{k}$
   &$\Box$
      &$\widehat{\dot{\Box}}$
         &$\JJoldiagbox$
            &$\JJJbox $
               &$\fourbox$
                  &$\widehat{\dot{\Box}}$ 
                     &$\overline{\twobox} $
                        &$\widehat{\oBox}$
   &$\twobox $
      &$\widehat{\dot{\Box}}$ 
         &$\circoldiagbox$
            &$\circbox$
               &$\fourbox $
                  &$\widehat{\dot{\Box}}$ 
                     &$\dot{\twobox}$
                        &$\widehat{\oBox}$\\
\hline 
\end{tabular}
\end{center}
\end{thm}


\proof  We have two tools at our disposal: the
periodicity theorem and Theorem \ref{thm-exotic}, which relates exotic
transfers to differentials.  

Figure \ref{fig-KH} shows that $\underline{M}'_{k}$ has the
indicated value for $-8\leq k\leq 0$ because the same is true of
$\EE_{4}^{0,k}$ and there is no room for any exotic extensions. On the
other hand $\EE_{4}^{0,k+32}$ does not have the same value for $k=-8$, $k=-6$
and $k=-4$.  This comparison via periodicity forces 
\begin{itemize}
\item [$\bullet$] the indicated $d_{5}$ and $d_{7}$ in dimension 24,
which together convert $\Box$ to $\fourbox$. These were also established in
Theorem \ref{thm-d7}.

\item [$\bullet$] the short transfer in dimension 26, which converts
$\widehat{\oBox}$ to $\widehat{\dot{\Box}}$. It also follows from the
the results of Section \ref{sec-C2diffs}.

\item [$\bullet$] the long transfer in dimension 28, which converts
$\overline{\twobox}$ to $\dot{\twobox}$.

\end{itemize}

The differential corresponding to the long transfer is 
\begin{align*}
d_{13} ([2u_{\lambda }^{7}])
 & =  a_{\sigma }a_{\lambda }^{6}u_{2\sigma }u_{\lambda }^{4}
              \normrbar_{1}^{3},  \\
\mbox{so}\qquad  
d_{13} (a_{\sigma } [2u_{\lambda }^{7}])
 & =  a_{\sigma }^{2}a_{\lambda }^{6}u_{2\sigma }u_{\lambda }^{4}
              \normrbar_{1}^{3}  
   =  2a_{\lambda }^{7}[u_{2\sigma }^{2}]u_{\lambda }^{3}
                \normrbar_{1}^{3}  
\end{align*}

\noindent This compares well with the $d_{13}$ of Theorem 
\ref{thm-normed-slice-diffs}, namely 
\begin{displaymath}
d_{13} (a_{\sigma }[u_{\lambda }^{4}])
    =a_{\lambda }^{7}[u_{2\sigma }^{2}]\normrbar_{1}^{3}.
\end{displaymath}

The statements in dimensions 4 and 20 have similar proofs, and we will
only give the details for the former.  It is based on comparing the
$\EE_{4}$-term for $\KH$ in dimensions $-28$ and $4$.  They must
converge to the same thing by periodicity.  From the slice
$\EE_{4}$-term in dimension 4 we see there is a short exact sequence
\begin{numequation}\label{eq-4-stem}
\begin{split}
\xymatrix
@R=4mm
@C=18mm
{
0\ar[r]
    &\JJ \ar[r]
        &{\underline{M}'_{4}}\ar[r]
            &{\odbox}\ar[r]
                &0\\
    &\Z/2\ar@/_/[dd]_(.5){0}\ar@{=}[r]
        &\Z/2\ar@/_/[dd]_(.5){\left[\begin{array}{c}1\\0\end{array} \right]}
             \ar[r]
            &0\ar@/_/[dd]\\
{}\\
    &\Z/2\ar@/_/[dd]\ar@/_/[uu]_(.5){1}
         \ar[r]^(.4){\left[\begin{array}{c}1\\0\end{array} \right]}
        &\Z/2\oplus \Z_{-}
            \ar@/_/[dd]_(.5){\left[\begin{array}{cc}0&2\end{array} \right]}
            \ar@/_/[uu]_(.5){\left[\begin{array}{cc}1&a\end{array} \right]}
            \ar[r]^(.6){\left[\begin{array}{cc}0&1\end{array} \right]}
            &\Z_{-}\ar@/_/[dd]_(.5){2}\ar@/_/[uu]\\
{}\\
    &0\ar@/_/[uu]\ar[r]
        &\Z_{-}\ar@/_/[uu]_(.5){\left[\begin{array}{c}b\\1\end{array} \right]}
               \ar@{=}[r]
            &\Z_{-}\ar@/_/[uu]_(.5){1},
}
\end{split}
\end{numequation}%

\noindent while the $(-28)$-stem gives 
\begin{displaymath}
\xymatrix
@R=4mm
@C=17mm
{
0\ar[r]
    &\dot{\twobox}\ar[r]
        &{\underline{M}'_{4}}\ar[r]
            &\obull\ar[r]
                &0 \\
    &\Z/2\ar@/_/[dd]_(.5){0}\ar@{=}[r]^(.5){}
        &\Z/2\ar@/_/[dd]_(.5){\left[\begin{array}{c}1\\0\end{array} \right]}
             \ar[r]^(.5){}   
            &0\ar@/_/[dd]_(.5){} \\
{}\\
    &\Z_{-}\ar@/_/[dd]_(.5){2}\ar@/_/[uu]_(.5){1}
         \ar[r]^(.4){\left[\begin{array}{c}c\\1\end{array} \right]}
        &\Z/2\oplus\Z_{-}
            \ar@/_/[dd]_(.5){\left[\begin{array}{cc}0&2\end{array} \right]}
            \ar@/_/[uu]_(.5){\left[\begin{array}{cc}1&a\end{array} \right]}
            \ar[r]^(.6){\left[\begin{array}{cc}1&d\end{array} \right]}
            &\Z/2\ar@/_/[uu]_(.5){}\ar@/_/[dd]_(.5){}\\
{}\\
    &\Z_{-}\ar@/_/[uu]_(.5){1}\ar@{=}[r]
        &\Z_{-}\ar@/_/[uu]_(.5){\left[\begin{array}{c}b\\1\end{array} \right]}
               \ar[r]^(.5){}
            &0\ar@/_/[uu]_(.5){}
}
\end{displaymath}

\noindent The commutativity of the second diagram requires that
\begin{displaymath}
a+b=c=1  \qquad \aand b+d=c+d=0,
\end{displaymath}

\noindent giving $(a,b,c,d)= (0,1,1,1)$.  The diagram
for $M_{4}$ is that of $\JJoldiagbox$ in Table \ref{tab-C4Mackey}. 


In dimension 20 the short exact sequence of (\ref{eq-4-stem}) is replaced by 
\begin{displaymath}
\xymatrix
@R=5mm
@C=10mm
{
0\ar[r] 
    &\circ \ar[r]
        &{\underline{M}'_{20}}\ar[r]
            &{\odbox}\ar[r]
                &0
}
\end{displaymath}

\noindent and the resulting diagram for $\underline{M}'_{20}$ is that
of $\circoldiagbox$.



Similar arguments can be made in dimensions 6 and 22.
\qed\bigskip

We could prove a similar statement about exotic restrictions hitting
the 0-line in the third quadrant in dimensions congruent to 0, 4, 6,
14, 16, 20 (where there is an exotic transfer) and 22.  The problem is
naming the elements involved.

\begin{table}[h]
\caption
[Infinite Mackey functors in the slice {\SS }.]
{Infinite Mackey functors in the reduced
$\EE_{\infty}$-term for $\KH$. In each even degree there is an
infinite Mackey functor on the 0-line that is related to a summand of
$\upi_{2k}\KH$ in the manor indicated. The rows in each diagram are 
short or 4-term exact sequences with the summand appearing in both
rows.}  \label{tab-0-line}
 
\begin{center} 
\begin{tabular}{|c|c||c|c|} 
\hline 
Dimension
    &Third quadrant
        &Dimension
            &Third quadrant\\
mod 32
    &First quadrant
        &mod 32
            &First quadrant\\
\hline 
0   &$
\xymatrix
@R=1mm
@C=7mm
{
{\fourbox} \ar[r]^(.5){\phantom{d_{7}}}
    &\Box \ar[r]^(.5){}\ar@{=}[d]^(.5){}
        &\circ \\
0\ar[r]^(.5){}
    &\Box \ar@{=}[r]^(.5){}
        &\Box
}$
        &16 &$
\xymatrix
@R=1mm
@C=7mm
{
{\fourbox} \ar[r]^(.5){}
    &{\twobox} \ar[r]^(.5){}\ar@{=}[d]^(.5){}
        &{\JJ}\\
    &{\twobox} \ar[r]^(.5){}
        &\Box\ar[r]^(.5){d_{7}}
            &\bullet
}$\\
\hline 
2, 10   &$
\xymatrix
@R=1mm
@C=7mm
{
{\widehat{\dot{\Box}}}
               \ar[r]^(.5){}
    &{\widehat{\dot{\Box}}}
               \ar[r]^(.5){}\ar@{=}[d]^(.5){}
        &0 \\
{\widehat{\bullet}}
               \ar[r]^(.5){}
    &{\widehat{\dot{\Box}}}
               \ar[r]^(.5){}
        &{\widehat{\oBox}}
}$
        &18, 26 &$
\xymatrix
@R=1mm
@C=7mm
{
{\widehat{\dot{\Box}}}
               \ar[r]^(.5){}
    &{\widehat{\dot{\Box}}}
               \ar[r]^(.5){}\ar@{=}[d]^(.5){}
        &0 \\
{\widehat{\bullet}}
               \ar[r]^(.5){}
    &{\widehat{\dot{\Box}}}
               \ar[r]^(.5){}
        &{\widehat{\oBox}}
}$\\
\hline 
4   &$
\xymatrix
@R=1mm
@C=7mm
{
{\dot{\twobox}}
               \ar[r]^(.5){}
    &{\JJoldiagbox }
               \ar[r]^(.5){}\ar@{=}[d]^(.5){}
        &{\obull}\\
{\JJ}
               \ar[r]^(.5){}
    &{\JJoldiagbox }
               \ar[r]^(.5){}
        &{\overline{\twobox} }
}$
        &20 &$
\xymatrix
@R=1mm
@C=7mm 
{
{\dot{\twobox}}
               \ar[r]^(.5){}
    &{\circoldiagbox }
               \ar[r]^(.5){}\ar@{=}[d]^(.5){}
        &{\circ} \\
{\circ}
               \ar[r]^(.5){}
    &{\circoldiagbox }
               \ar[r]^(.5){}
        &{\overline{\twobox} }
}$\\
\hline 
6   &$
\xymatrix
@R=1mm
@C=7mm
{
{\bullet}\ar[r]^(.5){d_{7}}
    &{\JJbox } 
               \ar[r]^(.5){}
        &{\JJJbox } 
               \ar[r]^(.5){}\ar@{=}[d]^(.5){}
            &{\bullet}\\
    &{\widehat{\oBox}}
               \ar[r]^(.5){}
        &{\JJJbox }
               \ar[r]^(.5){}
            &{\JJJ}
}$
        &22 &$
\xymatrix
@R=1mm
@C=7mm 
{
{\JJbox }
               \ar[r]^(.5){}
    &{\circbox}
               \ar[r]^(.5){}\ar@{=}[d]^(.5){}
        &{\bullet} \\
{\circ}
               \ar[r]^(.5){}
    &{\circbox}
               \ar[r]^(.5){}
        &{\widehat{\oBox}}
}$\\
\hline 
8   &$
\xymatrix
@R=1mm
@C=7mm
{
{\fourbox}
         \ar[r]^(.5){\phantom{d_{7}}}
    &{\fourbox}
           \ar[r]^(.5){}\ar@{=}[d]^(.5){}
        &0\\
    &{\fourbox}
         \ar[r]^(.5){}
        &{\Box}
               \ar[r]^(.5){d_{5}, d_{7}}
            &{\circ}
}$
        &24 &$
\xymatrix
@R=1mm
@C=7mm
{
{\fourbox}
         \ar[r]^(.5){}
    &{\fourbox}
           \ar[r]^(.5){}\ar@{=}[d]^(.5){}
        &0\\
    &{\fourbox}
         \ar[r]^(.5){}
        &{\Box}
               \ar[r]^(.5){d_{5}, d_{7}}
            &{\circ}
}$\\
 \hline 
12  &$
\xymatrix
@R=1mm
@C=7mm
{
{\bullet}
               \ar[r]^(.5){d_{13}}
    &{\dot{\twobox}}
               \ar[r]^(.5){}
        &{\overline{\twobox} }\ar@{=}[d]^(.5){} \\
    &0
               \ar[r]^(.5){}
        &{\overline{\twobox} }
               \ar[r]^(.5){}
            &{\overline{\twobox} }
}$
        &28 &$
\xymatrix
@R=1mm
@C=7mm
{
{\dot{\twobox}}
               \ar[r]^(.5){}
    &{\dot{\twobox}}
               \ar[r]^(.5){}\ar@{=}[d]^(.5){}
        &0 \\
{\bullet}
               \ar[r]^(.5){}
    &{\dot{\twobox}}
               \ar[r]^(.5){}
        &{\overline{\twobox} }
}$\\
\hline 
14  &$
\xymatrix
@R=1mm
@C=7mm
{
{\JJ} \ar[r]^(.5){d_{7}}
    &{\JJbox }
               \ar[r]^(.5){}
        &{\widehat{\oBox}}
               \ar@{=}[d]^(.5){}
            & \\
    &0
               \ar[r]^(.5){}
        &{\widehat{\oBox}}
               \ar[r]^(.5){}
            &{\widehat{\oBox}}
}$
        &30 &$
\xymatrix
@R=1mm
@C=7mm
{
{\widehat{\oBox}}
               \ar[r]^(.5){}
    &{\widehat{\oBox}}
               \ar[r]^(.5){}\ar@{=}[d]^(.5){}
        &0 \\
0
               \ar[r]^(.5){}
    &{\widehat{\oBox}}
               \ar[r]^(.5){}
        &{\widehat{\oBox}}
}$\\
\hline 
\end{tabular}
\end{center}
 
\end{table}

In Table \ref{tab-0-line} we show short or 4-term exact sequences in the 16 even
dimensional congruence classes.  In each case the value of
$\underline{M}'_{k}$ is the symbol appearing in both rows of the diagram.
For even $k$ with $0\leq k<32$, we typically have short exact
sequences
\begin{displaymath}
\xymatrix
@R=2mm
@C=10mm
{
0\ar[r]^(.5){}
    &{\EE_{4}^{0,k-32}}\ar[r]^(.5){}
        &{\underline{M}'_{k}}\ar[r]^(.5){}\ar@{=}[d]^(.5){}
            &{\mbox{quotient} }\ar[r]^(.5){}
                &{0}\\
0\ar[r]^(.5){}
    &{\mbox{subgroup} }\ar[r]^(.5){}
        &{\underline{M}'_{k}}\ar[r]^(.5){}
            &{\EE_{4}^{0,k}}\ar[r]^(.5){}
                &{0,} 
}
\end{displaymath}

\noindent where the quotient or subgroup is finite and may be spread
over several filtrations.  This happens for the quotient in dimensions
$-32$, $-16$ and $-12$, and for the subgroup in dimensions 6 and 22.

This is the situation in dimensions where no differential hits
[originates on] the 0-line in the third [first] quadrant.  When such a
differential occurs, we may need a 4-term sequence, such as the one in
dimension -22.

In dimensions 8 and 24 there is more than one such differential,
the targets being a quotient and subgroup of the Mackey functor
$\circ=\Box/\fourbox$. 

In dimension $-18$ we have a $d_{7}$ hitting the 0-line.  Its source
is written as ${\circ \subseteq \EE_{4}^{-7,-24}}$ in Figure
\ref{fig-KH}.  Its generator supports a $d_{5}$, leaving a copy of
$\JJ$ in $\EE_{7}^{-7,-24}$.

There is no case in which we have such
differentials in both the first and third quadrants.

\begin{cor}\label{cor-higher}
{\bf The $\EE_{\infty }$-term of the slice {\SS} for $\KH$.}
The surviving elements in the spectral sequence for $\KH$ are shown in
Figure \ref{fig-E14a}.
\end{cor}

\begin{thebibliography}{Dug05}

\bibitem[Ada84]{Ad:Prereq}
J.~F. Adams.
\newblock Prerequisites (on equivariant stable homotopy) for {C}arlsson's
  lecture.
\newblock In {\em Algebraic topology, {A}arhus 1982 ({A}arhus, 1982)}, volume
  1051 of {\em Lecture Notes in Math.}, pages 483--532. Springer, Berlin, 1984.

\bibitem[Ati66]{Atiyah:KR}
M.~F. Atiyah.
\newblock {$K$}-theory and reality.
\newblock {\em Quart. J. Math. Oxford Ser. (2)}, 17:367--386, 1966.

\bibitem[Dug05]{Dugger}
Daniel Dugger.
\newblock An {A}tiyah-{H}irzebruch spectral sequence for {$KR$}-theory.
\newblock {\em $K$-Theory}, 35(3-4):213--256 (2006), 2005.

\bibitem[GM95]{Greenlees-May}
J.~P.~C. Greenlees and J.~P. May.
\newblock Equivariant stable homotopy theory.
\newblock In {\em Handbook of algebraic topology}, pages 277--323.
  North-Holland, Amsterdam, 1995.

\bibitem[HHRa]{HHR}
Michael~A. Hill, Michael~J. Hopkins, and Douglas~C. Ravenel.
\newblock The non-existence of elements of {Kervaire} invariant one.
\newblock Online at http://arxiv.org/abs/arXiv:0908.3724v2 and on the third
  author's home page.

\bibitem[HHRb]{HHR:RO(G)}
Michael~A. Hill, Michael~J. Hopkins, and Douglas~C. Ravenel.
\newblock The slice spectral sequence for {$RO(C_{p^{n}})$}-graded suspensions
  of {$H\underline{\bf Z}$ I}.
\newblock To appear.


\bibitem[Lew88]{Lewis:ROG}
L.~Gaunce Lewis, Jr.
\newblock The {$R{\rm O}(G)$}-graded equivariant ordinary cohomology of complex
  projective spaces with linear {${\bf Z}/p$} actions.
\newblock In {\em Algebraic topology and transformation groups ({G}\"ottingen,
  1987)}, volume 1361 of {\em Lecture Notes in Math.}, pages 53--122. Springer,
  Berlin, 1988.


\bibitem[Rav86]{Rav:MU}
Douglas~C. Ravenel.
\newblock {\em Complex cobordism and stable homotopy groups of spheres}, volume
  121 of {\em Pure and Applied Mathematics}.
\newblock Academic Press Inc., Orlando, FL, 1986.
\newblock Errata and second edition available online at author's home page.

\bibitem[TW95]{Thevenaz-Webb}
Jacques Th{\'e}venaz and Peter Webb.
\newblock The structure of {M}ackey functors.
\newblock {\em Trans. Amer. Math. Soc.}, 347(6):1865--1961, 1995.

\end{thebibliography}

\end{document}